\documentclass[final]{siamltex}

\include{macros_edge}

\usepackage{color}
\definecolor{blue}{rgb}{0,0,1}
\definecolor{red}{rgb}{1,0,0}
\definecolor{green}{rgb}{0,1,0}
\definecolor{orange}{rgb}{1,.5,0}
\definecolor{plum}{rgb}{0.62, 0.0, 0.77}

\newcommand\SGvar[1]{\textcolor{magenta}{{ #1}}}
\newcommand\sg[1]{\textcolor{black}{{ #1}}}

\newcommand\sbullet[1][.4]{\mathbin{\vcenter{\hbox{\scalebox{#1}{$\bullet$}}}}}
\newcommand\dotexp{\, _{\sbullet} \, \widehat{}\ \,}

\newcommand{\btrue}{b_{\text{\scriptsize{true}}}}
\newcommand{\inv}{^{-1}}
\newcommand{\lj}{^{(\ell),j}}
\newcommand{\R}{\mathbb{R}}
\newcommand{\rhp}{\mbox{reshape}}
\newcommand{\vvec}{\mbox{vec}}
\newcommand{\whA}{\widehat{A}}
\newcommand{\whb}{\widehat{b}}
\newcommand{\whd}{\widehat{d}}
\newcommand{\whw}{\widehat{w}}
\newcommand{\xreg}{x_{\text{\scriptsize{reg}}}}
\newcommand{\xtrue}{x_{\text{\scriptsize{true}}}}

\title{An Inner-Outer Iterative Method for Edge Preservation in Image Restoration and Reconstruction}

\author{
Silvia Gazzola\thanks{Department of Mathematical Sciences, University of Bath, Bath, UK ({\tt S.Gazzola@bath.ac.uk}). SG's effort for this paper is supported in part by the EPSRC under grants EP/P005985/1 and EP/T001593/1, and the NSF under grant DMS-1906664.} \and
Misha E. Kilmer\thanks{Department of Mathematics, Tufts University, Medford, MA 02155 ({\tt misha.kilmer@tufts.edu}).} \and
         James~G.~Nagy\thanks{Department of Mathematics, Emory University, Atlanta, GA ({\tt jnagy@emory.edu}). 
         JGN's effort for this paper is supported in part by the U.S. National Science Foundation under grant DMS-1819042 and the NIH under grant 1R13EB028700-01.} \and
         Oguz Semerici\thanks{Spotify ({\tt oguz.semerci@gmail.com })} \and
         Eric L. Miller\thanks{Department of Electrical and Computer Engineering, Tufts University, Medford, MA 02155  ({\tt elmiller@tufts.edu}).
         ELM's effort for this paper is based upon work supported by the U.S. Department of Homeland Security, Science and Technology Directorate, Office of University Programs, under Grant Award 2013-ST-061-ED0001. The views and conclusions contained in this document are those of the authors and should not be interpreted as necessarily representing the official policies, either expressed or implied, of the U.S. Department of Homeland Security.}
	}

\begin{document}

\maketitle
%%%%%%%%%%%%%%%ABSTRACT %%%%%%%%%%%%%%%%%
\begin{abstract}
We present a new inner-outer iterative algorithm for edge enhancement in \sg{imaging problems}. 
%reconstruction and restoration 
%problems.    
At  each outer iteration, we formulate a Tikhonov-regularized problem where the \sg{penalization} is \sg{expressed in the} 2-norm \sg{and involves a regularization operator} designed to improve edge resolution as the outer iterations progress, through an adaptive 
%sequential refinement 
process. An efficient hybrid regularization method is used to \sg{project the Tikhonov-regularized problem onto approximation subspaces of increasing dimensions (inner iterations), while conveniently choosing the regularization parameter (by %are chosen very efficiently on 
applying well-known techniques, such as the discrepancy principle or the ${\mathcal L}$-curve criterion, to the projected problem)}. \sg{This procedure }   
results in an automated algorithm for edge recovery that does not involve regularization parameter tuning by the user, \sg{nor} repeated calls to sophisticated optimization algorithms, \sg{and is therefore 
particularly attractive from a computational point of view.} A key \sg{to the success of the new algorithm} is the design of the regularization operator through the use of an \sg{adaptive} diagonal weighting matrix that effectively enforces smoothness only where \sg{needed}. We demonstrate
the value of our approach on applications in X-ray CT image reconstruction and in \sg{image} deblurring, and show that it
can be computationally much more attractive than other well-known strategies for edge preservation, while providing \sg{solutions} of \sg{greater or equal} quality.
\end{abstract}

%%%%%%%%%%%%%%%%%%%% KEYWORDS %x%%%%%%%%%%%
\begin{keywords}
edge \sg{enhancement}, hybrid regularization methods, Krylov subspace methods, parameter choice strategies, iterative reweighting, tomography, image deblurring
\end{keywords}
%%%%%%%%%%%%%%%% AMS NUMBERS %%%%%%%%%%%%
\begin{AMS}
\sg{65F08, 65F10, 65F22}
\end{AMS}

%%%%%%%%%%%%%%% PAGE STYLE AND LABELING %%%%%%%%%%%%%
\thispagestyle{plain}
\markboth{O. Semerici, M. E. Kilmer, E. L. Miller, S. Gazzola}{Iterative Edge Preservation}

%%%%%%%%%%%%%%%%%%%%%%%%%%%%%%%%%%%%%%%
\section{Introduction} \label{sec:Intro}
In this paper, we consider a new inner-outer iterative algorithm for edge \sg{enhancement} in image restoration and reconstruction problems.  Specifically, we consider generating regularized solutions $\xreg$ to 
\begin{equation} \label{eq:model}
Ax  = \btrue + \eta =b , \quad \mbox{ with }\quad \btrue := A\xtrue ,\end{equation}
where $A \in \mathbb{R}^{m \times n}$ is a known \sg{ill-posed} forward operator,  
$b$ represents the known measured data \sg{corrupted by some unknown noise} $\eta$, and $\xtrue$ is the unknown, vectorized version of the true image to which we would like to generate an approximation.  \sg{Here and in the following we assume that $\eta$ is Gaussian white noise.} 
     
One of the most well-known regularization methods is Tikhonov
regularization, which consists in computing
\begin{equation} \label{eq:tik}   \sg{\xreg = \arg} \min_{x} \| A x - b \|_2^2 + \lambda^2 R(x),  \end{equation}
where $R(x)$ is referred to as the regularization term. The regularization
parameter $\lambda$ determines the trade-off between the fidelity to the model, and \sg{the} damping of the noise through
enforcement of a priori information via $R$.   Indeed, the role of $\lambda$ cannot be overemphasized, since the quality of the solution estimate is highly dependent on its value.   Unfortunately, in practice, good values of $\lambda$ are rarely known a priori.  When $R(x)$ is a non-quadratic term, the most widely used approach is to solve (\ref{eq:tik}) with a suitable optimization method for many values of $\lambda$, and use a heuristic approach, based on the output, to decide which provides the best quality solution.  For some choices of $R(x)$ that try to preserve edges in $x$, this can be a very expensive strategy, as we discuss below. 

%One common choice for regularization 
%operator is  $R(x) := \|L x \|_2^2$
%where $L$ is usually taken to be either the identity or a discrete gradient or discrete Laplacian operator.   
%Due to the use of the 2-norm on the constraint term $\| L x \|_2^2$, solutions tend to be 
%smooth, assuming an appropriate value of $\lambda$ is known.   In some applications, a smooth regularized
%solution is indeed desirable.   
%Solving this regularized least squares problem with a suitable Krylov method is straightforward if $\lambda$ is known, as $A$ and $L$ are usually structured.  
To motivate the work in this paper, let us first consider the case $R(x) = \| Lx \|_2^2$, where $L$ is taken to be the identity (standard form Tikhonov) or a discrete derivative-type operator (general form Tikhonov), such as the discrete gradient or Laplacian.    
%Assume that the 
%optimal value of $\lambda$ is known.  
Assuming that an appropriate value of $\lambda$ is known, the corresponding solution $\xreg$ will necessarily be smooth given these choices of $L$\sg{, as the magnitude of $x$ or derivatives thereof is penalized}.  For $\lambda$ fixed, the solution to the optimization problem is straightforward to compute.   For large scale problems, this is usually done via a Krylov subspace-based iterative \sg{projection} method, such as CGLS, LSQR \cite{PaSa82a,PaSa82b} or LSMR \cite{FoSa11}.  At each iteration, these iterative methods project the original problem (\ref{eq:tik}) onto Krylov subspaces of increasing dimension, they \sg{only} require \sg{one matrix-vector} product with \sg{each of the} $A, A^T, L, L^T$, and they have short term recurrences (so storage is minimal).   Provided \sg{that} the number of iterations is small relative to the dimension of the problem, and that fast matrix-vector products are possible, these are relatively efficient solvers.  

%However, this approach has the disadvantage of a new solve for each new $\lambda$.  
%
%If $\lambda$ is not known a priori, several alternatives are possible.  
%One option is to solve (\ref{eq:tik}) for several fixed values of $\lambda$, then based on a suitable metric or parameter selection heuristic \cite{ }, select the regularized solution that is most desirable.   For large scale problems, the solves are accomplished 

\sg{The} so-called hybrid methods \cite{Bjo88,BjGrDo94,ChPa15,GaNa14,GaNoRu15,KiOLe01,OlSi81} \sg{exploit} the efficiency of iterative solvers such as those mentioned in the previous paragraph \sg{to select a good regularization} parameter when $R(x)$ has the form $\| L x \|_2^2$.  Iterative solvers produce a small-dimensional projected problem that retains certain features of the original large-scale problem. The order of the projected problem equals the number of performed iterations. Additional regularization (for instance, Tikhonov regularization) is then \sg{efficiently} applied to the \sg{small-dimensional} projected problem. 
%that typically is much smaller in dimension.  
From there, a regularized solution to the original large-scale problem can be obtained.  
%% OLD bits
%The \sg{upside} of the hybrid approach is that the regularization parameter is selected for the small-scale problem, and if a technique such as the discrepancy principle, the generalized cross validation, or the ${\mathcal L}$-curve criterion is used, the \sg{quantities} for the heuristic \sg{are} available essentially for free with almost no extra floating point operations or storage.  \sg{(These strategies can be applied with a negligible computational cost, by manipulating the projected quantities of much smaller size)} \SG{if we include the discrepancy principle, should we still refer to it as ``heuristic''?} 
%% NEW bits
The \sg{upside} of the hybrid approach is that the regularization parameter is selected \sg{at each iteration} for the small-\sg{dimensional} problem \sg{only}: in this setting, some well-known techniques such as the discrepancy principle or the ${\mathcal L}$-curve criterion can be \sg{applied with a negligible computational cost (both in terms of floating point operations and storage), by manipulating the much smaller projected quantities.} 
However, with the usual choices of $L$ \sg{mentioned above}, the solutions will be smoother than desired.    
% However, for the problems of interest to us in this paper, we desire regularized solutions with sharp edges consistent
%with those in $x_{true}$, and the aforementioned approaches which rely on Krylov subspaces generated by $A$ typically do not possess those feature.  
          
Two well-known alternative choices for $R(x)$ when one desires regularized solutions with sharp edges are $TV(x)$ (\sg{isotropic} total variation) and $\| Lx \|_p^p$, 
 with $L$ the discrete gradient operator and $p$ close to 1 (see \cite{Vogelbook} and references therein). 
 %% OLDEST bits
 % Iteratively reweighted least norm approaches or fixed point approaches are among the algorithms designed to solve \cite{eq:reg} when $R(x) = TV(x)$ or $R(x) = \|Lx \|_p^p$ \cite{list of citations}.  
 For a single fixed $\lambda$, the algorithms employed to solve the optimization problem (\ref{eq:tik}) usually require much more computational effort than if the 2-norm \sg{regularization} is used.   \sg{Moreover}, 
the biggest pitfall of using these regularization operators is that $\lambda$ is not known a priori.  Thus, many optimization problems need to be solved for a discrete set of $\lambda$ in order to \sg{compute a meaningful solution}. %arrive at the solution.     

\sg{A natural and well-established way of handling problem (\ref{eq:tik}) with $R(x)=\|Lx\|_p^p$, $p\neq 2$, is to reformulate it as a sequence of quadratic problems: the $\ell$th quadratic problem involves a regularization term of the form $R(x)=\|M^{(\ell)}x\|_2^2$, where the $\ell$th regularization matrix $M^{(\ell)}$ is defined with respect to $L$ and the approximate solution of the $(\ell-1)$th quadratic problem in the sequence. Here and in the following we refer to this approach as 
\emph{iteratively reweighted norm} (IRN) method: this was first proposed in \cite{IRNekki}, and then extended in \cite{renaut2017hybrid,wohlberg2008lp,IRNtv}.} Recently, an Arnoldi-Tikhonov (hybrid) method for $R(x): = \| x \|_1$ was proposed \cite{GaNa14}, which can be used when the operator $A$ is square. This approach combines a reweighting strategy (to approximate $R(x)$) and a hybrid strategy applied with the discrepancy principle to choose the regularization parameter. Since the weights are updated at each iteration of the hybrid method, and since the weights are formally regarded as variable preconditioners after transforming the current quadratic Tikhonov problem into standard form, the reweighted problem is efficiently projected using the flexible Arnoldi algorithm: because of this, only one cycle of iterations should be performed, while, for all the other IRN methods mentioned so far, inner-outer iteration cycles are needed when dealing with large-scale problems. This method based on flexible Krylov subspaces has been extended to work with $R(x)=TV(x)$ in \cite{gazzolasabate2019}, and with rectangular matrices $A$ as well as a sparsity-under-transform regularization term in \cite{chunggazzola2018}. 
      
     In this paper we propose a technique that is similar to the IRN methods, in the sense that we, too, generate a sequence of quadratic problems, with a different weight matrix at each outer-iteration, and employ a hybrid approach on each regularized problem in the sequence.  However, our method differs from each of the methods listed above in that: (a) our regularization matrix is different -- indeed it is rectangular and possibly not full column rank; (b) the operator $A$ may be rectangular; (c) due to (a) and (b), we must use a different hybrid projection algorithm.   
%           
% 
%Another edge-based alternative which capitalizes on solving two-norm problems is proposed in \cite{ChenKilmerHansen14}.  Their idea is to recover the smooth part of the solution and the edge part of the solution separately by expanding $x = Wy = W_1 y_1 + W_2 y_2$ where $R(W_1) \perp R(W_2)$, the columns of $W$ are orthonormal and resemble the singular vectors of the operator and leverage structure and fast transforms.   As this approach does not strictly fit into the class of regulariztaion problems of type (\ref{eq:tik}), 
%we will not consider it further here.  
%% For deblurring in particular, $A$ often has special structure, and thus $W$ with the desired properties can be generated and applied quickly (see \cite{ChungKilmerOleary15} and references therein).    A least squares problem is solved for $y_1$, and then $\| L x \|_p^p$ is minimized to find $y_2$.  This has the advantage that the regularizaiton parameter is determined on the `cheaper' problem of recovering the smooth 
%% solution, but solving the secondary problem can be time consuming.  
%However, what our present approach has in common with the approach in \cite{ChenKilmerHansen14} is that  
%the aim is to invoke the use of numerical linear algebraic tools to get a straightforward algorithm that is easier to implement and interpret perhaps than, say, TV.   
%       
%           
More precisely, we propose an automated, \sg{inner-outer iterative} approach, consisting in solving %in which we solve 
a sequence of regularized least squares 
 problems with \sg{2}-norm \sg{regularization of the form}
 \begin{equation} 
 \label{eq:MK} \min_{x \in \Gamma^{(\ell)}} \| A x - b \|_2^2 + \lambda^2 \| M^{(\ell)} x \|_2^2, \quad\sg{\ell=1,2,3,\dots}
 \end{equation}
 where a near-optimal value of $\lambda$ for the $\ell$th problem (\ref{eq:MK}), which we denote as $\lambda_{*,\ell}$, is determined on-the-fly for each 
 regularization operator $M^{(\ell)}$.   Here, $\Gamma^{(\ell)}$ denotes a specific projection space in which we will look for a \sg{projected} solution.   Because this solution space \sg{depends in part on} $M^{(\ell)}$ itself (but not on the regularization parameter) we use a superscript on it as well. %\SG{make sure that the concept of projection is well-explained.}
To compute $\Gamma^{(\ell)}$ and 
 $\lambda_{*,\ell}$ we employ the hybrid regularization approach \sg{(inner iterations)} of \cite{KiHaEs07}. We give a
 method for designing $M^{(\ell)}$ adaptively (outer iterations) in such a way that edges are enhanced as $\ell$ increases.  
 % For fixed $k$, we 
 %will solve the regularized problem, 
 We discuss expected behavior on a class of images, and 
present results \sg{and comparisons} on applications in X-ray CT and image deblurring to demonstrate that our algorithm 
 has the potential to produce high-quality images in a computationally efficient manner.    

This paper is organized as follows. In Section \ref{sec:background}, we give definitions, notations, and motivations
for studying our new inner-outer \sg{iterative} approach.  We describe hybrid iterative regularization algorithms in Section \ref{sec:hybrid}, highlighting the particular method in \cite{KiHaEs07} that we leverage in our new approach.  In Section \ref{sec:algorithm}, we develop our \sg{inner-outer iterative} algorithm, and provide some analysis. Section \ref{sec:results} is devoted to numerical results.  We give conclusions and \sg{outline} future work in Section \ref{sec:conclusions}.

%%%%%%%%%%%%%%%%%%%%%%%%%%%%%%%%%%%%%%%
\section{Background and Motivation} \label{sec:background}
We begin by motivating the need to take $M^{(\ell)}$ in (\ref{eq:MK}) as something other than the identity or the discrete gradient operator. 
Then, we discuss other algorithms proposed in the literature that have been used for edge-preservation. 

\subsection{Filtering Methods}
%\MEK{Some info on plain vanilla Tikhonov, why edges are wiped out with a 2-norm constraint. } 

Consider the model (\ref{eq:model}) and assume that $A$ has rank $\rho$. Consider the singular value decomposition (SVD) of $A = U \Sigma V^T = \sum_{i=1}^{\rho} \sigma_i u_i v_i^T$. 
\sg{Then,} the minimum norm least squares solution to $A x = b$ is easily seen to be given by
\begin{equation}\label{eq:xSVD}%\[ 
  x = \sum_{i=1}^{\rho} \frac{u_i^T b}{\sigma_i} v_i = 
  \sum_{i=1}^{\rho}  \left(\frac{{u_i^T\btrue} }{ \sigma_i}  v_i + 
                \frac{u_i^T \sg{\eta}}{\sigma_i} v_i\right) , 
\end{equation}%\]
and $\sigma_1 \gg \sigma_\rho > 0$. 
\sg{Under the assumption of} white noise \sg{$\eta$}, the values of $|u_i^T \sg{\eta}|$ are approximately constant, whereas the $|u_i^T \btrue|$ are large and dominant for small indices $i$ but, eventually, decay toward zero and the noise becomes dominant (this property in commonly known as discrete Picard condition, see \cite{Hansen97}). 

Regularization methods like TSVD, and Tikhonov with $R(x) = \| x \|_2^2$ in (\ref{eq:tik}), 
%(or, equivalently, $M_k=I$ in (\ref{eq:MK})) 
work by damping the contribution from the terms corresponding to large indices in (\ref{eq:xSVD}), i.e., they are essentially filtering methods and compute a regularized solution of the form
\begin{equation}\label{xregfilt}
\xreg = \sum \phi_i \frac{u_i^Tb}{\sigma_i} v_i\,,
\end{equation}
where $\phi_i$ denotes a \sg{scalar} that decreases with \sg{increasing} $i$; see \cite{Hansen97,HaNaOL06}. However, the edge content from the image is encoded in the \sg{high-frequency} terms that are damped or discarded \sg{in (\ref{xregfilt})}, and hence such regularized
solutions tend to be smooth.  
%\SG{perhaps we should dwell a bit more on the behavior of the $\sigma_i$'s and the $v_i$'s.}

\begin{figure}
\begin{center}
\begin{tabular}{cc}
\includegraphics[width=5cm]{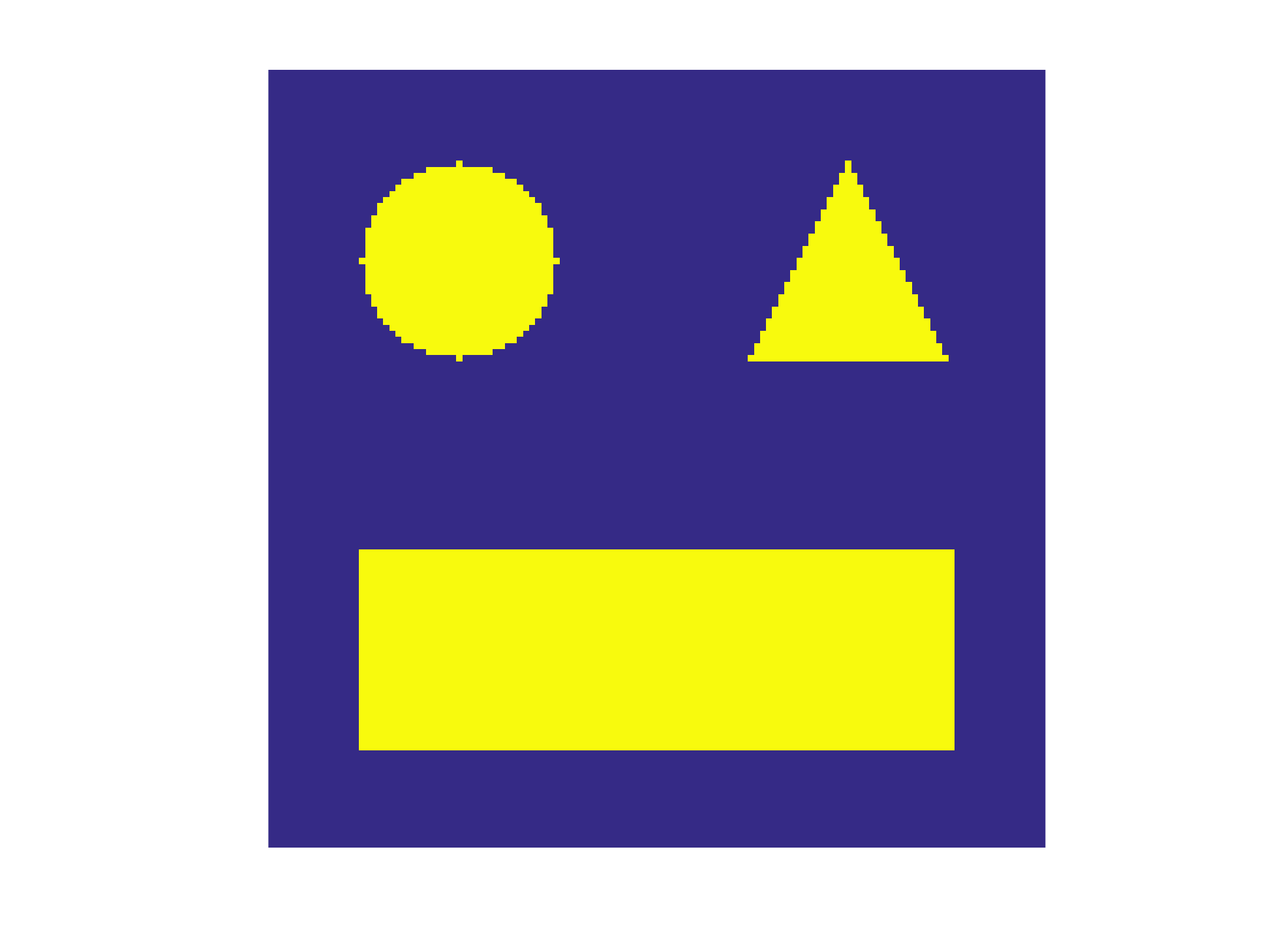} & 
\includegraphics[width=5cm]{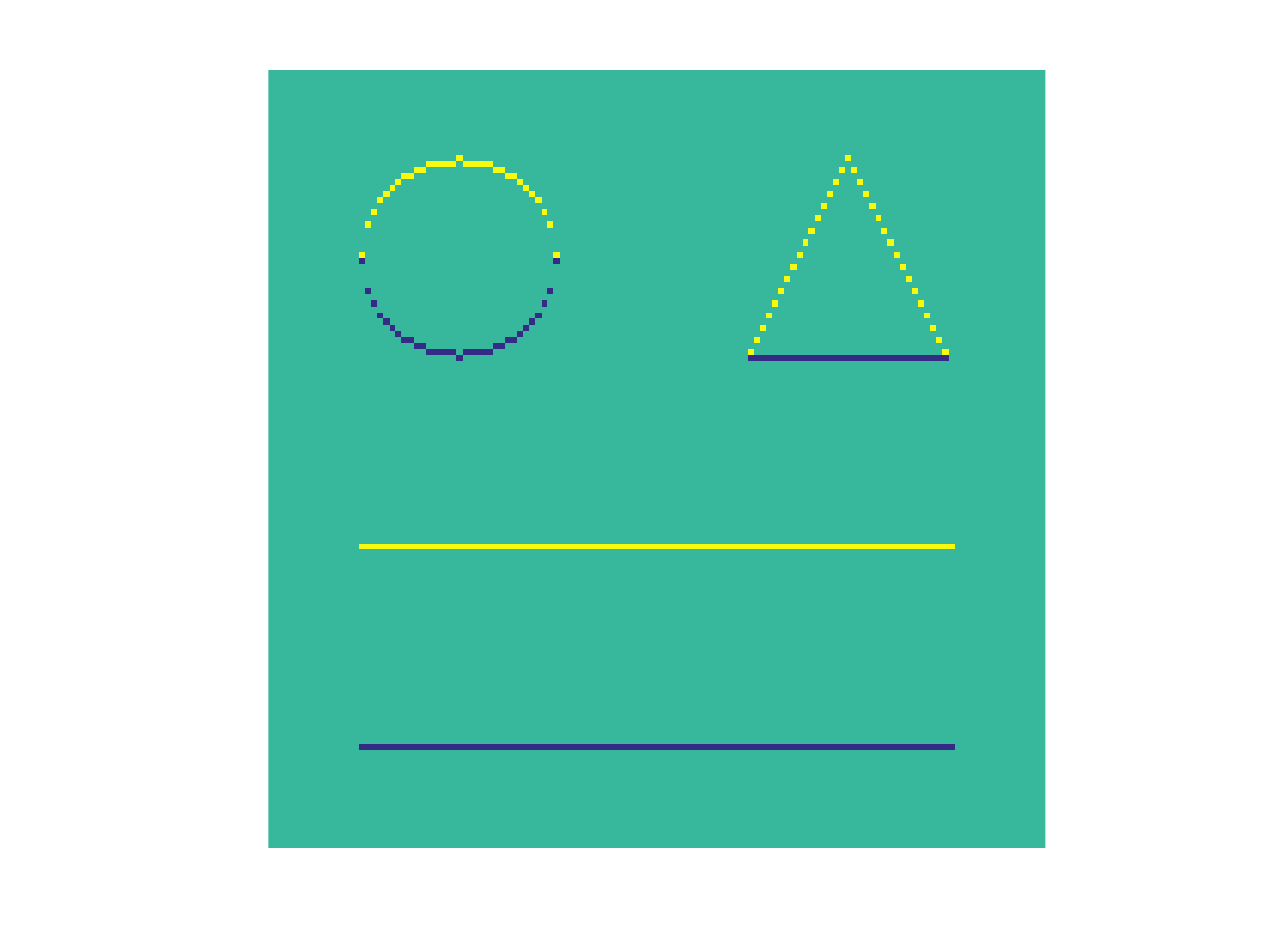} \\ 
\includegraphics[width=5cm]{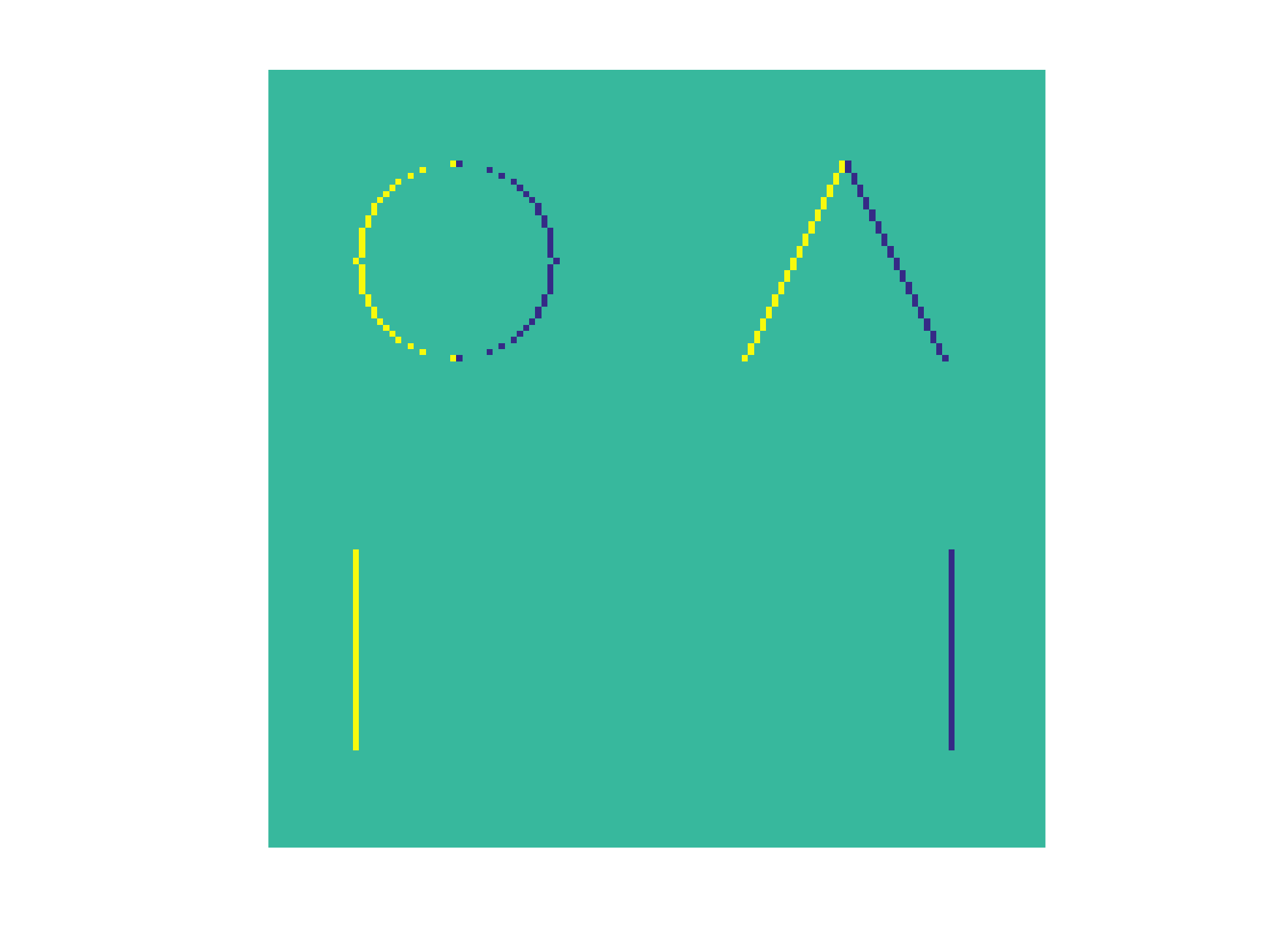} & 
\includegraphics[width=5cm]{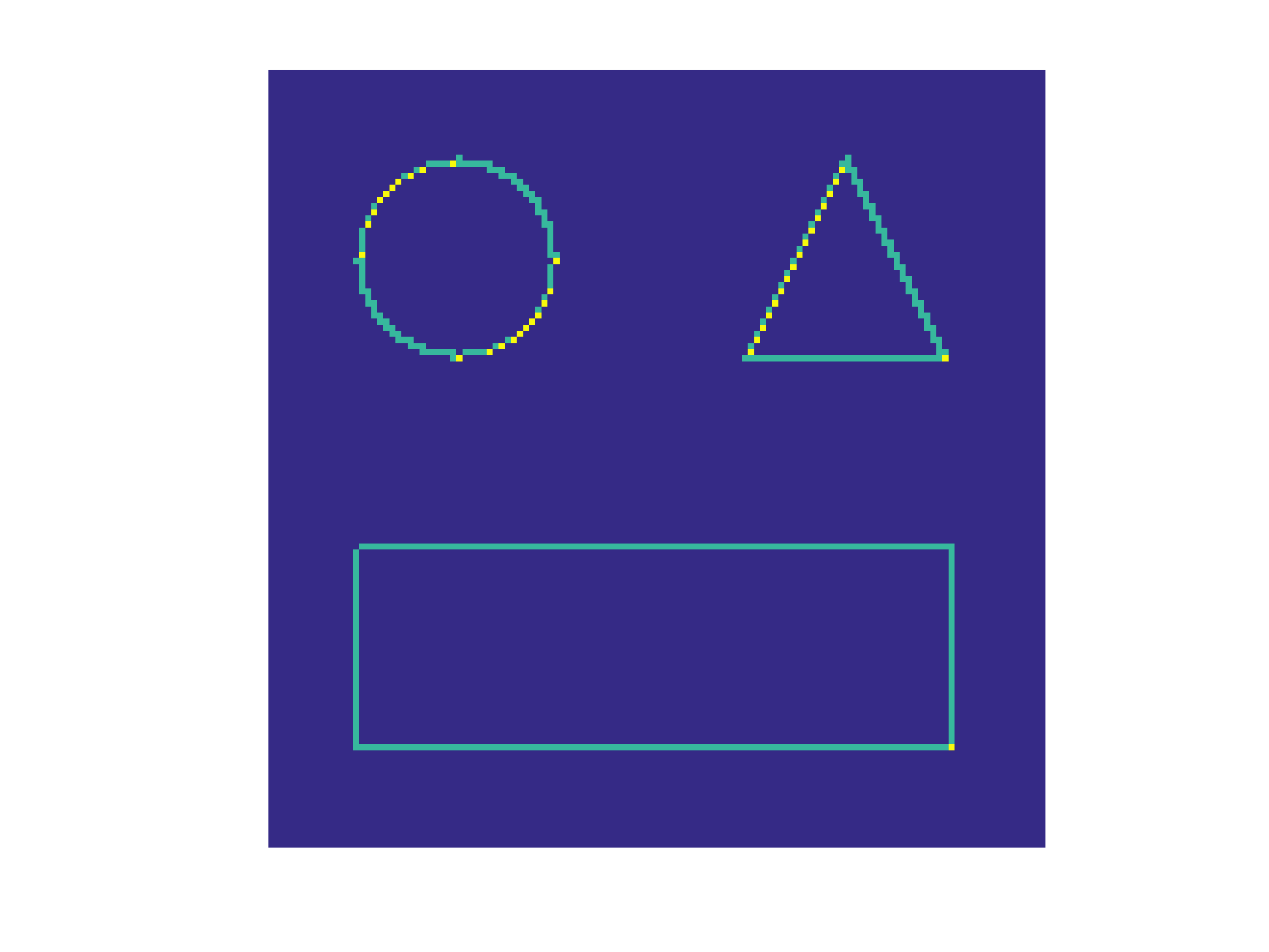}
\end{tabular}
\caption{\label{fig:onenorm} Top left:  exact image $X=\rhp(x)$. Top right:  vertical derivatives $L_v X$ displayed as an image. 
Bottom left:  horizontal derivatives $X L_h^T$ displayed as an image. Bottom right: \sg{pixel-wise magnitude of the gradient $[L_v X]_{i,j}^2 + [XL_h^T]_{i,j}^2$}. Extra rows and columns
of zeros were added as needed to keep the images appearing square for ease of comparison.}
\end{center}
\end{figure}

Next, let us consider what happens when $R(x) = \| L x \|_2^2$, where $L$ is the discrete gradient operator.  To fix notation, 
we assume the following.  
Let $X$ be an $N_v \times N_h$ image\sg{, and let $x=\vvec(X)$ be its vectorized version, obtained by stacking its columns.} Define 
\[ L_v = \bea{ccccc} -1 & 1 & 0 & 0 & \cdots \\ 0 & -1 & 1 & 0 & \cdots \\
                                  0 & 0 & -1 & 1 &  \cdots \\
                                  \vdots & \ddots & \ddots & \ddots &\ddots \\
                                  0 & 0 & \cdots & 0 & - 1 \eea  \]
to be the $N_v \!-\!1 \times N_v$ discrete first derivative operator.  We define $L_h\in \mathbb{R}^{(N_h-1)\times N_h}$ similarly. Using $\otimes$ to represent Kronecker product between two matrices, note that
\begin{itemize}
 \item $L_v X = \mbox{reshape}(( I \otimes L_v) \vvec(X)) \in \mathbb{R}^{(N_v -1) \times N_h}$ is the image that approximates all vertical derivatives;
 \item $X L_h^T = \mbox{reshape}((L_h \otimes I)\vvec(X)) \in \mathbb{R}^{N_v \times (N_h -1)}$ is the image that approximates horizontal derivatives.
 \end{itemize}
\sg{Here %$\vvec(\cdot)$ is the operator that takes a 2D image and reshapes it as a vector by stacking its columns; 
$\rhp(\cdot)$ can be regarded as the inverse of the $\vvec(\cdot)$ operator.} Now let %$L = \bea{cc} I \otimes L_v \\ L_h \otimes I \eea$,
\begin{equation}\label{eq:gradient}
L = \bea{cc} I \otimes L_v \\ L_h \otimes I \eea\,.
\end{equation}
The matrix $L$ is called the \emph{discrete gradient operator} and, clearly, $Lx \in \mathbb{R}^{N_h(N_v-1)+(N_h-1)N_v}$ contains the edge information of $x$. 
%For perspective, 
Figure \ref{fig:onenorm} illustrates the vertical derivative image, horizontal derivative image, and \sg{magnitude of the} gradient image of a  \sg{standard} test image of size $128\times 128$ pixels from the \sg{MATLAB toolbox} ``IR Tools'' \cite{Gazzola2019}.  
%\SG{maybe replace this with something from IR Tools.}

Consider (\ref{eq:tik}) with $R(x) = \| L x \|_2^2$.   Since $L x$ contains the edge information, adding a multiple
of $R(x)$ to the data misfit component of the cost functional has the effect of producing regularized solutions 
where jumps are {penalized}.  For images that are supposed to be smooth, this is not a problem. But assume for a moment that $X$ is a  \sg{piecewise constant image with jump discontinuities (as the one displayed in the top-left corner of Figure \ref{fig:onenorm})}. Additive noise
has the effect of creating jumps, so $R(x)$ smooths out the discontinuities in a region of otherwise constant value, which is desired.  Unfortunately, because the values of $Lx$ are large near \sg{true} jump discontinuities, the price one pays is that the
jump discontinuity is also smoothed out.    

\subsection{Edge Preservation Algorithms} 
As mentioned in Section \ref{sec:Intro}, two popular choices for $R(x)$ in (\ref{eq:tik}) are Total Variation \cite{RuOsFa92}  and $\| L x \|_p^p$ for
 $1 \le p < 2$ \cite{Hansen97}. 
%One straightforward alternative regularization method that is used to damp noise while retaining edge information
%is to set $R(x) = \| L x \|_p^p$, taking $1 < p \ll 2$.    For $p$ very close to one, convergence can be very slow, particularly when closing in on the solution, since the system that must be solved to find the search direction starts to become very ill-conditioned as edges are detected.   Also, as mentioned previously, one typically would have to solve (\ref{eq:tik}) for multiple values of $\lambda$.  
The Total Variation functional \sg{can be expressed, in a discrete setting, as}
\begin{equation}\label{eq:TV}
R(x) = \mbox{TV}(x)=\sum_{i=1}^{N_h -1} \sum_{j=1}^{N_v -1} \sqrt{ |L_v X|_{ij}^2 + | X L_h^T|_{ij}^2}.
\end{equation}
Clearly, either choice for $R$ necessitates using an iterative optimization technique and,  at each step of this routine, %optimization, 
a linear system or least squares problem \sg{may} be solved to generate the next search direction.  
The lagged diffusivity fixed point iteration (see \cite{VogelOman} for algorithm details) is a popular choice, but other options such as steepest descent, Newton's method, and primal-dual approaches are also possible (see \cite{Vogelbook} for details and references).  
For $R(x) = \| L x \|_1$, one can employ iteratively reweighted methods (from the Newton family) or interior point methods \cite{LiSantosa}.   
 Other approximations to 
TV may also be employed (see, for example, \sg{\cite{MLSQR}} and the discussion in \cite{Vogelbook}).   

However, \sg{a} key issue \sg{that is somewhat disregarded within the edge-preserving solvers listed so far} is the choice of regularization parameter $\lambda$.   \sg{Indeed,} the optimization codes written to solve (\ref{eq:tik}) with the edge-preserving $R(x)$  outlined in the previous paragraph {\it all assume that the value of $\lambda$ is fixed}.   Since the {\it optimal value of $\lambda$ is not known a priori}, the procedure \sg{usually advocated} is to solve multiple optimization problems for different values of $\lambda$ \sg{picked from a discrete set}, and then choose the value of $\lambda$ that gives a ``best'' image \sg{according to} some metrics. This \sg{involves} many calls to the optimization routine; moreover, the convergence behavior of a particular run of the optimization \sg{method} will vary with $\lambda$. Thus, \sg{on one side} it would appear that recovery of edge information via any of the aforementioned methods would be very computationally expensive, because we have to solve multiple optimization problems for multiple $\lambda$ \sg{(one for each value of $\lambda$)}; \sg{moreover, as hinted above,} each one of the optimization problems may require multiple linear system solves to determine the outcome.   

On the other \sg{side}, if we take $R(x) = \| x \|_2$ or $R(x) = \| L x \|_2 $, the cost functional \sg{in} (\ref{eq:tik}) is
quadratic and 
%(at least relative to the algorithms referenced above) 
it is relatively straightforward to compute the solution to the optimization problem  (i.e., in this case the cost of one \sg{full} solve is roughly the cost of a single optimization step for the edge-preserving $R$'s). Even so, \sg{especially when it comes to the choice of the regularization parameter}, further savings in the case of a quadratic problem (\ref{eq:tik}) can be obtained through hybrid iterative methods. Hybrid regularization methods operate by projecting a quadratic problem onto smaller-dimensional subspaces where certain properties of noise and signal separation can be preserved.  \sg{When employing a hybrid approach}, (additional) regularization is applied to the small projected problem instead, and the regularization parameter is chosen by applying a suitable heuristic to the projected problem. This is {significantly} cheaper than solving a quadratic problem (\ref{eq:tik}) for each new value of $\lambda$. This is discussed in detail  in the next section. However, \sg{adopting} this family of \sg{quadratic} $R(x)$'s,  we trade edge content for computational efficiency.                          

In this paper we propose to merge the best of these two worlds, defining an \sg{adaptive} sequence of quadratic problems
\[  \min_x \| A x - b \|_2^2 + \lambda^2 \| \underbrace{D^{(\ell)} L}_{=:M^{(\ell)}} x \|_2^2,\quad\ell=1,2,3,\dots \]
  whereby a hybrid approach can be employed to automatically choose the regularization parameter in an efficient way for each quadratic problem in the sequence.    Our choice of $M^{(\ell)}$ has features in common with the methods proposed in \sg{\cite{Bardsley, IRNtv}}, in that it involves different weightings of the gradient operator to emphasize edges in current estimated solutions.  However, it differs from those methods in that the effects of our $M^{(\ell)}$ are cumulative, i.e., \sg{information from all the previous steps is retained}.    
\sg{To run  our algorithm we only need to be able} to perform matrix vector products with the forward operator $A$ in (\ref{eq:model}) and with $L$ in (\ref{eq:gradient}), and their respective transposes. Moreover, our approach uses hybrid methods to automatically determine regularization parameters.

\section{Hybrid Iterative Methods}\label{sec:hybrid}

In this section, we give the general outline for two \sg{(hybrid)} approaches that could ultimately be used within our inner-outer \sg{iterative strategy}. \sg{The first one is based on LSQR, while the second one is based on a joint bidiagonalization algorithm.} In our problems, $D^{(\ell)}L$ will not have full column rank, which makes the second option a more elegant choice and easier to implement (\sg{for this reason, the second option is used to perform the numerical experiments displayed in Section \ref{sec:results}}). However, for completeness and ease of exposition, we present the LSQR-based hybrid approach first.  For more specific details, and for additional methods from the hybrid regularization family, we refer 
%the reader 
to \cite{GaNa14,GaNoRu15}.

Consider the problem
\begin{equation}\label{eq:tikGen}%\[ 
\min_{x} \|A x - b \|_2^2 + \lambda^2 \| M x \|_2^2, 
\end{equation}%\]
where,  
%$\lambda$ is known or TBD and 
for the remainder of this section, we assume \sg{that} $M$ remains fixed.   This problem is mathematically equivalent to 
the linear least squares problem
\begin{equation}\label{eq:tikGenLS} %\[ 
\min_{x} \left\| \bea{c} A \\ \lambda M \eea x - \bea{c} b \\ 0 \eea \right\|_2^2.  
\end{equation}%\]
The hybrid methods described in the following subsections \sg{compute approximate solutions to (\ref{eq:tikGen}) and (\ref{eq:tikGenLS})} as functions of $\lambda$, %\sg{and of the dimension of a projection subspace}, 
and select good choices of $\lambda$ in efficient ways. %, as \sg{detailed} below.

\subsection{Hybrid-LSQR Method for Standard Form Tikhonov} \label{ssec:hybridstd}
 Let us first address solving (\ref{eq:tikGen}) %the first version 
 iteratively by LSQR \cite{PaSa82a,PaSa82b}. Implicitly, this approach amounts to transforming problem (\ref{eq:tikGen}) into 
standard form and \sg{constraining the solution to belong to a $k$-dimensional subspace $\mathcal{K}_k$ at the $k$th iteration, i.e., solving}
\begin{equation} \label{eq:stdform} 
\min_{y \in \mathcal{K}_k} \| AM_A^{\dagger}y - b \|_2^2 + \lambda^2 \| y \|_2^2, 
\end{equation}
 where $M_A^{\dagger}$ is the $A$-weighted psedo-inverse \sg{of $M$}, \sg{defined as}
 \begin{equation}\label{eq:MAdagger}
   M_A^\dagger = \left(I - \left(A\left(A-M^\dagger\right)\right)^\dagger A \right) M^\dagger\,.
\end{equation}
In the above definition, the use of the dagger symbol without an accompanying subscript (e.g., $M^\dagger$) is meant
to denote the standard Moore-Penrose pseudo-inverse (see, for instance, \cite{GvL}). Note that, if $M$ is invertible,
then $M_A^\dagger = M^{-1}$, and, if $M$ has full column rank, then $M_A^\dagger = M^\dagger$. If 
$M$ is an underdetermined matrix, then $M_A^\dagger$ may not be the same as $M^\dagger$.
For further details on $A$-weighted pseudo-inverses in the context of Tikhonov regularization,
we refer %the reader 
to \cite{Hansen97}. %\newline

%\SG{Previous version: \\
%Implicitly, this approach amounts to transforming the equation to 
%standard form, and then putting a restriction on where the solution to the standard form equation
%can live:
%\[% \begin{equation} \label{eq:stdform} 
%\min_{y \in \mathcal{K}_k} \| AM^{\dagger}y - b \|_2^2 + \lambda^2 \| y \|_2^2, 
%\]%\end{equation}
% where $\mathcal{K}_k$ is a $k$-dimensional subspace and
% $M^{\dagger}$ should be understood as follows:
% \begin{itemize}
%   \item If $M$ is invertible, $M^{\dagger} = M^{-1}$.
%   \item As the (A-weighted) pseudo-inverse.  For our application, we must use the A-weighted psuedoinverse 
%   \MEK{add explanation and conditions}.   
%   \end{itemize}}
%   \JGN{I think this is a bit complicated to define without using a subscript, $M_A^\dagger$.} \SG{I agree with JN's proposed edits -- they are already incorporated above.} \newline

\sg{When adopting} LSQR, $\mathcal{K}_k$ represents the $k$-dimensional Krylov subspace
%\SG{Previous version:
%\[ 
%\mathcal{K}_k = \mbox{span} \{C^T b, (C^T C) A^Tb, (C^TC)^2 A^T b, \ldots, (C^TC)^{k-1} A^T b \}, 
%\]
%generated by the operator and the right-hand side in (\ref{eq:stdform}), where $C = A M_A^{\dagger}$. Also, I believe that we should always have $C^Tb$ (some instances of $A^Tb$ appear).}\\
%\SG{Proposed change (heavier notation, but more direct -- $\mathcal{K}_k$ does not appear elsewhere in the manuscript):}
\[ 
\mathcal{K}_k = \mbox{span} \{(A M_A^{\dagger})^T b, ((A M_A^{\dagger})^T (A M_A^{\dagger})) (A M_A^{\dagger})^Tb, %((A M_A^{\dagger})^TC)^2 A^T b, 
\ldots, ((A M_A^{\dagger})^T(A M_A^{\dagger}))^{k-1} (A M_A^{\dagger})^T b \}, 
\]
generated by $(A M_A^{\dagger})^T (A M_A^{\dagger})$ and $(A M_A^{\dagger})^Tb$. %the \sg{coefficient matrix} and the right-hand side in (\ref{eq:stdform}). 
The $k$th iteration of the LSQR algorithm updates the matrix recurrence
\[ (A M_A^{\dagger}) V_k = U_{k+1} B_k \]
where $B_k$ is $(k+1) \times k$ \sg{lower} bidiagonal, %in exact arithmetic 
$V_k$ has $k$ orthonormal columns that span $\mathcal{K}_k$, 
and $U_{k+1}$ has $k+1$ orthonormal columns with $U_{k+1}(\beta e_1) = b$.    
\sg{Taking $y=V_kz$ (i.e., imposing $y\in\mathcal{K}_k$), and exploiting} unitary invariance of the 2-norm allows us to rewrite (\ref{eq:stdform}) as the  $O(k)$-sized least squares problem
%\begin{equation}\label{eq:stdformProj} 
%\min_{z}  \left\| \bea{c} B_k \\ \lambda I \eea  z -  \bea{c} \beta e_1 \\ 0 \eea \right\|_2^2 ,
%\end{equation}
\begin{equation}\label{eq:stdformProj} 
\sg{\min_{z\in\R^k}  \| B_kz - \beta e_1\|_2^2 + \lambda^2 \|z \|_2^2 ,}
\end{equation}
%assuming $y^{(\lambda)} = V_k z^{(\lambda)}$, 
which is a Tikhonov-regularized projected problem in standard form. \sg{Denoting by $z^{(\lambda,k)}$ the minimizer of the above problem, $y^{(\lambda,k)}=V_kz^{(\lambda,k)}$ is the \sg{minimizer} of (\ref{eq:stdform})}. 
If $k$ is large \sg{enough, so} that the 
        spectral behavior of the bidiagonal matrix $\sg{B_k}$ mimics that of the operator $A M_A^{\dagger}$, then the impact of the choice of $\lambda$ \sg{on the quality of the solution} can be observed. \sg{However, as hinted in Section \ref{sec:Intro}, $k$ should be assumed small enough, so that problem (\ref{eq:stdformProj}) can be efficiently generated and solved.} 
%\SG{It seems to me that we might write slightly contrasting statements in the following, depending on the parameter choice strategy we would like to employ. If we consider the discrepancy principle, we care about $k$ being small enough so that a zero finder can be efficiently applied to the projected discrepancy, and we do not need to solve each projected problem for a set of discrete values of $\lambda$ that is fixed in advance (at least with a classical implementation of the discrepancy principle). On the contrary, if we consider the ${\mathcal L}$-curve we would like to update the relevant quantities as $k$ increases and for a discrete set of values of $\lambda$ simultaneously (as suggested in \cite{KiHaEs07}), and potentially we do not care so much about $k$ increasing.}        \\
% 
 For fixed $k$, we can use any suitable \sg{parameter choice strategy (such as the discrepancy principle or the ${\mathcal L}$-curve criterion)} on the projected problem (\ref{eq:stdformProj}) to choose $\lambda$. 
% \sg{For instance, if a good estimate of the noise norm $\|\eta\|_2$ is available, we can apply the discrepancy principle, i.e., apply a zero-finder to compute the value of $\lambda>0$ that satisfies the nonlinear equation $\|B_k z^{(\lambda,k)} - \beta e_1\|=\tau\|\eta\|_2$, where $\tau>1$ is a safety constant: if $k$ is small enough, the computational cost of these operations is negligible. Alternatively, }
%%to choose the $\lambda$ (for $k$ fixed)
%using the ${\mathcal L}$-curve would require the pairs $(\|B_k z^{(\lambda,k)} - \beta e_1\|, \|z^{(\lambda,k)} \|)$: these norms can be computed from the ($k-1$)th ones 
%%previous estimate (i.e., for $k-1$) 
%with just a few floating point operations for a finite set of $\lambda$ simultaneously.

\sg{Still assuming $k$ fixed, and} regardless of the parameter selection method, let $\lambda_{k}$ denote the estimated optimal parameter. \sg{Assuming $k$ is large enough, }
%(again, we are assuming $k$ is fixed here).  
% For example, we can  take $\lambda$ as the value, $\lambda_*$, at the corner of the ${\mathcal L}$-curve \cite{Lcurve}.  
we \sg{then} find $z^{(\lambda_{k},k)}$ by solving (\ref{eq:stdformProj}), \sg{we} set $y^{(\lambda_{k},k)} = V_k z^{(\lambda_{k},k)}$, and \sg{we} transform from standard form (\ref{eq:stdform}) back to general form  (\ref{eq:tikGen}), to get \sg{an approximation} $x^{(\lambda_{k},k)}$ of an optimal %the desired 
regularized solution of (\ref{eq:tikGen}). 
%, which is the desired regularized image estimate.   
    Computing $y^{(\lambda_{k},k)}$ may be done with short term recurrences on the iteration $k$ (i.e., in $k$, for fixed $\lambda$), %(i.e., for $\lambda$ fixed, and $k$ varying),  
    without storage of $V_k$.  
        See, for example, \cite{KiOLe01} for details.  %\MEK{add one line ref also Julianne's work}\\
%        \sg{In Section \ref{sec:algorithm} we will use the lighter notations $z^{(\lambda)}=z^{(\lambda,k)}$, $y^{(\lambda)}=y^{(\lambda,k)}$, and $x^{(\lambda)}=x^{(\lambda,k)}$ %, depend also on $k$, we suppress this notationally, 
%to avoid confusion when also the regularization operators will change within the outer iterations. }
        
%\SG{Previously: Although $z^{(\lambda)}$ and $y^{(\lambda)}$ depend also on $k$, we suppress this notationally, to avoid confusion when we begin to consider the case of changing regularization operators in Section \ref{sec:algorithm}.}  
%        

%        As noted in \cite{KilmerOleary}, the ${\mathcal L}$-curve for the projected problem is the ${\mathcal L}$-curve for (\ref{eq:stdform}) (assuming exact arithmetic) when solved this way, since $\|B_k z^{(\lambda)} - \beta e_1\| = \| A M^{\dagger} y^{(\lambda)} - b \|$ and $\|z^{(\lambda)} \| = \|y^{(\lambda)} \|$.   
        
\sg{To evaluate the computational cost of this class of hybrid algorithms we should consider that} each iteration requires a matrix-vector product with $A M_A^{\dagger}$ and one with its transpose. Note that $A M_A^{\dagger} v$ is computed as the two subsequent ``products'' $A(M_A^{\dagger} v)$.   The product $M_A^{\dagger} v$ is actually computed implicitly:  for instance, if $M_A^{\dagger}= M^{-1}$, then
         $M^{-1}v = g$ \sg{is the same as} $v = Mg$, so in practice the system $Mg =v$ is solved for $g$, and $M^{-1}$ is never explicitly computed. Thus, one usually refers to the ``product'' of $M_A^{\dagger} v$ as ``applying $M_A^{\dagger}$''.    Unfortunately, applying $M_A^{\dagger}$ in the case of the $A$-weighted pseudo-inverse \sg{defined as in (\ref{eq:MAdagger})} requires additional solver calls. \sg{Moreover, in general, each solver call} can involve invoking iterative solvers. Therefore, analyzing the computational cost of the products $M_A^{\dagger}v$, and of this hybrid algorithm as a whole, may be difficult. 
         % and it is therefore difficult to analyze in terms of computational cost.   
         In any case, the overall cost in producing a regularized solution depends predominantly on $k$, and the cost of an application of $M_A^{\dagger}$ and matrix-vector products with $A$ (likewise, their transposes); we must highlight again that the computational cost is \sg{essentially} 
independent of the number of discrete values of $\lambda$ tested to determine the optimal 
         regularization parameter.
%\SG{are we sure that this is the case, even if, say, $n^3/k^3$ values of $\lambda$ should be tested? (e.g., each iteration of the discrepancy principle costs $O(k^3)$; if we need  $n^3/k^3$ iterations of the zero finder (which will never be the case!), then the computational cost of the $k$th hybrid iteration is comparable to the computational cost of computing the SVD of $A$.}\\
%%
%\SG{Previously: However, it is for all 
%intents and purposes {\it independent of the number of discrete values of $\lambda$ tested} to determine the optimal regularization parameter.}
         
\subsection{Hybrid Method for General Form Tikhonov}  \label{ssec:hybridgen}
         
	%However, 
	Applying $M_A^{\dagger}$ is often non-trivial \cite{Hansen97}, especially when $M$ is defined as a product of matrices (see also \cite{gazzolasabate2019}). %For this case, 
The authors of \cite{KiHaEs07} develop an \sg{alternative} hybrid method based on a joint bidiagonalization algorithm, which we now describe.
	Here, $M$ can be rectangular and does not need to be full rank. 
Avoiding transformation into standard form (and the need of applying $M_A^{\dagger}$), and starting instead from (\ref{eq:tikGenLS}), the formulation
 %\be \label{eq:one}
 \begin{equation}\label{eq:genTikLSproj}
   \min_{x \in \mathcal{Z}_k}
  \left\| \left[ \begin{array}{c} A \\ \lambda M \end{array} \right] x - 
                                 \left[ \begin{array}{c}  b \\ 0 \end{array} \right] \right\|_2 , \qquad
   \mathcal{Z}_k = \mbox{span}\{z_1,\ldots,z_k\} ,
 \end{equation}
 is considered, where a partial joint bidiagonalization for $A$ and $M$ is updated at each iteration as follows
%% OLD bits
%  \begin{itemize}
%  \item $A Z_k = U_{k+1} B_k$ and
%  \item $M Z_k = \hat{U}_k \hat{B}_k$, with $Z_k = [z_1,\ldots,z_k]$. 
%  \item $U_{k+1}(\beta_1 e_1) = b, \qquad \beta_1 = \| b \|_2 $
%  \end{itemize}
%% NEW bits
\begin{equation}\label{eq:jointbd}
A Z_k = U_{k+1} B_k, \qquad M Z_k = \hat{U}_k \hat{B}_k,
\end{equation}
where $Z_k = [z_1,\ldots,z_k]$, $U_{k+1}(\beta e_1) = b$, $\beta = \| b \|_2$. %\\
Taking $x = Z_k w$ in (\ref{eq:genTikLSproj}) \sg{results in the following equivalent problems}
% \begin{eqnarray} %\label{eq:proj}
%  \min_{x \in \mathcal{Z}_k} \left\| \left[ \begin{array}{c} A \\ \lambda M \end{array} \right] x - 
%                                                         \left[ \begin{array}{c} b \\ 0 \end{array} \right] 
%   \right\|_2 & = & \min_{g^{(\lambda)} } \left\| \left[ \begin{array}{c} A Z_k \\ \lambda M Z_k \end{array} \right] g^{(\lambda)} 
%   - \left[ \begin{array}{c} b \\ 0 \end{array} \right] \right\|_2 \nonumber \\
%  & = & \min_{g^{(\lambda)}} \left\| \left[ \begin{array}{c} B_k \\ \lambda \hat{B}_k \end{array} \right] g^{(\lambda)} -
%   \left[ \begin{array}{c} \beta_1 e_1 \\ 0 \end{array} \right] \right\|_2.  \label{eq:equiv}
% \end{eqnarray}
\[
  \min_{x \in \mathcal{Z}_k} \left\| \left[ \begin{array}{c} A \\ \lambda M \end{array} \right] x - 
                                                         \left[ \begin{array}{c} b \\ 0 \end{array} \right] 
   \right\|_2 = \min_{w\in\R^k} \left\| \left[ \begin{array}{c} A Z_k \\ \lambda M Z_k \end{array} \right] w 
   - \left[ \begin{array}{c} b \\ 0 \end{array} \right] \right\|_2 \nonumber \\
  = \min_{w\in\R^k} \left\| \left[ \begin{array}{c} B_k \\ \lambda \hat{B}_k \end{array} \right] w -
   \left[ \begin{array}{c} \beta e_1 \\ 0 \end{array} \right] \right\|_2\,.
\]
The last problem in the above equations is a Tikhonov-regularized problem of dimension $O(k)$ involving only bidiagonal matrices and, for $k$ small relative to the original problem dimension, it is significantly cheaper to solve. 
\sg{Let us {denote by $w^{(\lambda,k)}$ the minimizer of this $O(k)$-dimensional problem: then $x^{(\lambda,k)}=Z_kw^{(\lambda,k)}$ is the \sg{minimizer} of (\ref{eq:genTikLSproj}). }} The bidiagonal structure leads to recurrence updates (in $k$, for fixed $\lambda$) for $x^{(\lambda,k)}$, for the residual norm $\|A x^{(\lambda,k)} - b \|_2$ and for the regularization norm term $\|M x^{(\lambda,k)} \|_2$, as functions of $\lambda$.  So, for each fixed $k$, the regularization parameter can be computed efficiently using strategies such as the discrepancy principle or the ${\mathcal L}$-curve criterion, which only require those quantities.

Thus, as long as $k$ is large enough that the bidiagonal pair inherits certain spectral properties of the original operators, the impact of $\lambda$ as a regularization parameter for the $k$th system can be observed, \sg{a suitable parameter $\lambda_{k}$ can be selected, and the solution 
%Since the $x^{(\lambda_*)}$ can be generated by short term recurrences, 
$x^{(\lambda_{k},k)}=Z_kw^{(\lambda_{k},k)}$ of (\ref{eq:genTikLSproj}) can be readily obtained by employing short term recurrences: this solution is an estimate of an optimal regularized solution of (\ref{eq:tikGenLS})}. The reader is referred to \cite{KiHaEs07} for further details.   

The algorithm requires only that products with $A, M$ and their transposes can be computed.   The cost of generating the bidiagonal pair is somewhat difficult to quantify. \sg{The $k$th step requires one call to LSQR to produce orthogonal projections that are needed to update the partial joint bidiagonalization (\ref{eq:jointbd}), and each LSQR iteration requires 4 matrix-vector products (with each of $A$, $A^T$, $M$, $M^T$). }
%,  but it also requires 2 calls to LSQR to produce orthogonal projections. 
If LSQR needs $m_k$ steps, then the cost of a call to LSQR is proportional to $m_k$ times the cost of those 4 matrix-vector products. The value of $m_k$ is typically small relative to problem dimension but, depending on $k$, it may affect the overall reconstruction time. 
%the overall impact can affect the reconstruction time. 
Promising methods for keeping $m_k$ small are currently under investigation. 
%\SG{missing reference here, with the acronym ``KiEdS19''.} %\cite{KiEdS19}.   
Nevertheless, bounds on the behavior \sg{of $k$ and $m_k$} are possible, %(\SG{should we put a reference?}), 
and for some classes of operators (e.g., see the CT image examples in the Section \ref{sec:results}) both $k$ and $m_k$ are very small. Even in cases where $k$ and/or $m_k$ are larger, when this approach is used as part of the inner-outer edge-preserving \sg{iterative} scheme below, the overall performance in producing quality images without any parameter tuning favors our approach versus multiple calls per $\lambda$ to sophisticated optimization routines, whose behavior and performance is at least as difficult %, if not more difficult, 
to analyze.     

%%%%%%%%%%%%%%%%%%%%%%%%%%%%%%%%%%%%%%%
\section{New Inner-Outer Iterative Algorithm for Edge Preservation} \label{sec:algorithm} %\sg{(previous title: Inner-Outer Iterative Algorithm).}

So far, we have seen that hybrid iterative methods for 2-norm Tikhonov-regularized problems are computationally appealing for two reasons: {(a) their cost is associated with the number of iterations, $k$, and the ability to compute
matrix-vector products, thus requiring little storage \sg{(although a precise computational cost estimate can be difficult; see Sections \ref{ssec:hybridstd} and \ref{ssec:hybridgen})}}; 
%\SG{this sounds a bit ambiguous, I propose changing to}\\
%\sg{(1) each iteration only requires matrix-vector products with the coefficient and regularization matrices, so that their computational cost grows linearly with the number of iterations $k$, which is typically small, and little storage is necessary}; 
(b) a good value for $\lambda$ can be computed on-the-fly essentially for free, with \sg{strategies} like the discrepancy principle or the ${\mathcal L}$-curve criterion. On the other hand, choosing the 2-norm of the gradient as a regularizer has a smoothing effect; changing \sg{to} the $p$-norm of the gradient, or changing to TV regularization, \sg{would enhance edges but} %means we need to solve 
requires the solution of a sequence of difficult optimization problems for many values of $\lambda$. 

We therefore look to build an approach that allows us to leverage the \sg{capabilities} of hybrid algorithms, with their computational efficiency and ability to choose \sg{regularization} parameters on-the-fly, but has the \sg{adaptivity} to capture edge information in the reconstruction process. Clearly, based on the discussion in \sg{Section \ref{sec:background}}, this means we cannot use a regularization term of the form $R(x)= \|Lx \|_2^2$.     
 
We propose to \sg{solve} the sequence of problems
\begin{equation} \label{eq:seqtrue}   
x^{(*,k_{\ell})} = \arg \min_{x \in \Gamma_{k_\ell}^{(\ell)}} \| A x - b \|_2^2 +              \lambda^{2}_{*,\ell} \| D^{(\ell)} L x \|_2^2 , \quad\sg{\ell=1,2,3,\dots}
\end{equation}
where $\lambda_{*,\ell}$ denotes the ``optimal'' regularization parameter for the $\ell$th problem, chosen according to some criterion. The $\ell$th problem computes a solution over a subspace of dimension $k_{\ell}$ (i.e., subspace dimension can vary with $\ell$, hence the subscript on $\Gamma$). In this section, we will answer the following \sg{questions}:
 \begin{itemize}
    \item How can $D^{(\ell)}$ be expressed, and how should it evolve to produce \sg{edge-enhanced} images?
    \item How can $\lambda_{*,\ell}$ be estimated efficiently?
    \item How can $x^{(*,k_{\ell})}$ be obtained efficiently?
 \end{itemize}

\subsection{Diagonal Weighting}\label{ssec:weights}
  
  Suppose \sg{$x^{(\ast,0)}$} is an estimate of $x$. 
  %\sg{ obtained, say, from running (\ref{eq:seqtrue}) with $D^{(0)}=I$}. 
  Assume further that it is a good enough estimate that at least one edge is visible.   If we compute the gradient image (as a matrix-vector product) and normalize it, i.e., 
\begin{equation}\label{eq:NormalizeGrad}
   g^{(0)} = \SGvar{|L x^{(\ast,0)}|} / \| L x^{(\ast,0)} \|_{\infty},
\end{equation}
where $| \cdot |$ in the numerator is used to denote element-wise absolute value, then
this vector will have the largest values equal to $1$ where \sg{the dominant edges} are located \sg{(these are precisely the values we {\it do not want to penalize})} and smaller nonnegative values where pixels are still corrupted by noise \sg{(these are precisely the values we {\it do want to penalize})}. If we now consider the ``image''
\begin{equation}\label{eq:GradientMap}
   d^{(1)} : = {\bf 1} - (g^{(0)})\dotexp p,  \qquad p > 0,
\end{equation}
where ${\bf 1}$ is the vector of all ones and $\dotexp$ denotes element-wise exponentiation, then the values we do not want to penalize have been mapped to \sg{the smallest values $\geq 0$}, and the noisy parts \sg{of the image} we want to penalize %with noise 
have been mapped to the largest values \sg{$\leq 1$}.\\  
%\MEK{Eric, I seem to remember this was motivated by other literature?  Can you elaborate. Should cite that, plus Oguz' thesis.}\\  
%\SG{the results must depend on the value of $p > 0$ -- I am commenting on this.}\\ 
\sg{Choosing $p\gg 1$ in (\ref{eq:GradientMap}) results in more penalization of the supposedly smooth regions (as more entries in the vector $d^{(1)}$ are forced to be close to 1), and therefore in an overall smoother reconstruction. Choosing $p\ll 1$ in (\ref{eq:GradientMap}) results in less penalization of the supposedly smooth regions (as less entries in the vector $d^{(1)}$ are forced to be close to 1), and therefore in an overall less smooth reconstruction.} 
Thus if, for $\ell = 1$, we use % (i.e., at the second outer iteration), we use
   \[ D^{(1)} := \mbox{diag}(d^{(1)}) \]
  in (\ref{eq:seqtrue}), this has the effect of enforcing smoothness on parts of the image where we still expect to be able to see, and improve, smoothness; \sg{at the same time,} it does not wash out the ``true'' edges, because we have (almost) zeroed out those components of the gradient image, \sg{as prescribed by the diagonal weighting in (\ref{eq:seqtrue})}. % (\ref{eq:GradientMap}). % by applying the diagonal weighting matrix.    
    
Suppose for the moment that $\lambda_{*,1}$ is known, and that $x^{(*,1)}$ has been determined (as described in the next subsections). If $x^{(*,1)}$ is an improvement over \sg{$x^{(\ast,0)}$}, then the normalized gradient image \linebreak[4]$g^{(1)} = |L x^{(*,1)}|/\| L x^{(*,1)} \|_\infty$ should reveal even more edge information.   
We set  $d^{(2)} = {\bf 1} - (g^{(1)})\dotexp p$, so that again
edge information we do not want to penalize has been mapped to \sg{the smallest values $\geq 0$}, and noise-contaminated components have been mapped to the largest values \sg{$\leq 1$}. 
%We form $\Lambda = \mbox{diag}(d^{(2)})$, and then 
%\[ D^{(2)} = \Lambda D^{(1)} . \]
%\SG{can we avoid introducing $\Lambda$? i.e., here and in the following just}
\sg{We then take \[ D^{(2)} = \mbox{diag}(d^{(2)}) D^{(1)} . \]}
Notice that $D^{(2)}$ encodes information about both the previous solution estimates $x^{(\ast,0)}$ and $x^{(\ast,1)}$, and its entries still remain between 0 and 1. %\SG{the following needs to be made more precise.} 
\sg{The null space of $D^{(2)}L$ is spanned by the constant vector $\bf 1$, and by the piecewise-constant vectors having edges corresponding to the dominant edges of $x^{(*,1)}$ and $x^{(\ast,0)}$.} %We then solve (\ref{eq:seqtrue}) for $\ell = 2$, as described in the next subsection. 
This process is iterated. %, and is summarized in Section \ref{ssec:mainalgo}. 
More insight into this choice of diagonal weighting matrices is given in \cite{PhD_thesis}. More algorithmic details are unfolded in the following Sections \ref{ssec:choice} and \ref{ssec:mainalgo}.
    
\subsection{Choosing the Regularization Parameter and Computing the Solution}\label{ssec:choice} 
%\sg{(previous title: Parameter and Solution Computation)} 

We now further elaborate on techniques for the computation of the regularization parameters and the corresponding solutions for the $\ell$th problem in (\ref{eq:seqtrue}). 

We can use a hybrid method to determine $\lambda_{*,\ell}$ and $x^{(*,k_{\ell})}$. In our numerical examples, we 
use the method described in Section \ref{ssec:hybridgen}, together with \sg{either the discrepancy principle or }the ${\mathcal L}$-curve criterion \sg{for} the projected problem (\ref{eq:genTikLSproj}). The following statements can be easily reformulated to work in connection with the hybrid method described in Section \ref{ssec:hybridstd}.

\sg{The discrepancy principle can be applied if a good estimate of the norm $\|\eta\|_2$ of the noise affecting the data is available. When used at the $k$th iteration of the hybrid solver for (\ref{eq:seqtrue}), it consists in determining the parameter $\lambda_k$ such that 
\begin{equation}\label{eq:dp}
\left\|B_kw^{(\lambda_k,k)}-\beta e_1\right\|_2=\tau\|\eta\|_2=\left\|Ax^{(\lambda_k,k)} - b\right\|_2\,,
\end{equation}
where $\tau>1$, $\tau\simeq 1$, is a safety threshold, and where the factorizations (\ref{eq:jointbd}) and the properties of the matrices appearing therein have been exploited. 
%Here the notation $x^{(k,\lambda_k)}=Z_kw^{(k,\lambda_k)}$, analogous to the one used in Section \ref{ssec:hybridgen}, has been temporarily resumed. 
Equation (\ref{eq:dp}) is nonlinear with respect to $\lambda_k$, and an appropriate zero-finder should be employed (see, for instance, \cite{ZeroFind}). Problem (\ref{eq:dp}) can be solved at a negligible additional computational cost, provided that $k$ is relatively much smaller than the original problem size $\min\{n,m\}$, and by updating some relevant factorizations for $B_k$ as $k$ increases. The inner iterations (in $k$) are stopped when no significative variations in $\lambda_k$ are detected for two consecutive values of $k$. Alternatively, one can employ the so-called ``secant method'' \cite{GaNoRu15}, which updates $\lambda_k$ in such a way that stopping by the discrepancy principle is ensured. When the stopping criterion is satisfied, we take $k_\ell=k$, and the optimal parameter $\lambda_{\ast,\ell}$ for the $\ell$th problem in (\ref{eq:seqtrue}) is set to $\lambda_{k_\ell}$. %, i.e., the $\lambda_{k}$ satisfying the stopping criterion for the discrepancy principle and the corresponding approximate solution 
Also, $\Gamma_{k_\ell}^{(\ell)}=\mathcal{Z}_{k_\ell}$, and  $x^{(\ast,k_{\ell})}=Z_{k_\ell}w^{(k_\ell,\lambda_{k_\ell})}$.}
%
%$x^{(\lambda_k,k)}$ can be easily computed if the matrix $Z_k$ is stored. %\SG{double-check this last statement, and provide a reference.} 
%The ; 

\sg{The ${\mathcal L}$-curve criterion can be employed at the $k$th iteration of the hybrid solver for (\ref{eq:seqtrue}), even when an estimate for $\|\eta\|_2$ is not available: in this case, the horizontal axis is determined by $\| {B}_k w^{(\lambda,k)} - \beta e_1 \|_2 = \| A x^{(\lambda,k)} - b \|_2$, and the vertical axis is determined by $\| \hat{B}_k w^{(\lambda,k)} \|_2 = \| M^{(\ell)} x^{(\lambda,k)} \|_2$, where $M^{(\ell)} = D^{(\ell)} L$; the above quantities are evaluated, at each iteration $k$, for a given fixed discrete set of values of $\lambda$, \sg{which correspond to points on the ${\mathcal L}$-curve}. The inner iterations (in $k$) are stopped when the corner of an ${\mathcal L}$ is visible for a few iterations, and the estimated corner value is constant for a few iterations. When this stopping criterion is satisfied, we take $k_\ell=k$, $\Gamma_{k_\ell}^{(\ell)}=\mathcal{Z}_{k_\ell}$, and $\lambda_{\ast,\ell}=\lambda_{k_{\ell}}$ as the optimal parameter  for the $\ell$th problem (\ref{eq:seqtrue}). 
%is then $\lambda_{k_\ell}$.  
The solution $x^{(\ast,\ell)}$ is generated by a short term recurrence, and it is therefore very storage efficient. We refer to \cite{KiHaEs07} for further computational details.}

 \sg{We stress that somewhat different computational approaches should be adopted when considering the discrepancy principle or the ${\mathcal L}$-curve criterion within the inner hybrid solver \sg{for the $\ell$th problem in (\ref{eq:seqtrue})}. This is because, at the $k$th hybrid iteration, the former computes approximate solutions that correspond to regularization parameters that are sequentially selected by the zero finder applied to (\ref{eq:dp}), while the latter simultaneously computes approximate solutions that correspond to a predetermined set of regularization parameters. However, in both cases, the computational overhead of selecting $\lambda_{\ast,\ell}$ is negligible if ${k_\ell}\ll\min\{m,n\}$.} 
 
\sg{Similarly to the discrepancy principle and the ${\mathcal L}$-curve criterion, every parameter choice rule 
%This fact makes methods like the discrepancy principle and the ${\mathcal L}$-curve 
that relies on the current residual or solution (semi)norms is an attractive option for hybrid methods.} %for selection methods.  
 
Clearly, this process of \sg{adaptively} computing optimal solutions and updating the regularization operator to account for newly acquired edge information can be repeated for the $(\ell+1)$th problem in (\ref{eq:seqtrue}). A sketch of the resulting algorithm is presented in the next section.

  \subsection{Summary of the Main Algorithm}\label{ssec:mainalgo}

As already highlighted, our new algorithm iteratively improves available solution estimates by recovering and exploiting edge information on the go. This algorithm is inherently based on an inner-outer iteration scheme, where inner iterations are handled with a hybrid algorithm, and the outer iterations update the available regularization term. The main operations performed by the new method are summarized in Algorithm \ref{alg:main}. 
 
%We consider iteratively improving our solution estimates, picking up edge information as we go, 

%\begin{algorithm}[h] \label{alg:main}
%\caption{Inner-Outer iterations for recovering $N_v \times N_h$ images}
%Set $D^{(0)} = I$ of size $(N_h (N_v-1) + N_v(N_h-1))$. Choose $p>0$ for (\ref{eq:GradientMap}).  %\\
%%Set $L$ to multiple of discrete gradient operator \\
%
%For $\ell=1,\ldots,$ until a stopping criterion is satisfied\\
%\ \hspace*{2mm}  $M = D^{(\ell-1)} L$ \\
%\ \hspace*{2mm} Call the hybrid algorithm to choose $\lambda_*$ and return $x_{*,\ell}$. \\
%\ \hspace*{2mm} Compute $g = |L x^{(*,\ell)}|/\|L x^{(*,\ell)}\|_{\infty}$, the rescaled gradient image of the current estimate. \\
%\ \hspace*{2mm}  Set $d := 1 - g \dotexp p$ (produces map to [0,1])  \\
%%\ \hspace*{2mm} Set $\Lambda = \mbox{diag}(d)$ \\
%\ \hspace*{2mm} Set $D^{(\ell)} \leftarrow \mbox{diag}(d) D^{(\ell-1)} $ \\
%% \ \hspace*{2mm}  Check $\| L x^{(*, \ell)} \|$ against $\| L x^{(*,\ell-1)} \|$ and break if needed.    \\
%end
%
%Return $\sg{\xreg} = x^{(*,\ell-1)}$. 
%
%\JGN{Should the normalization of the gradient be explicitly written in the algorithm?}\\
%\SG{I think so (I fixed this, and also took the component-wise absolute value of the numerator). Also, to have coherent notations we may have to consider $g^{(\ell)}$ (instead of $g$), and so on...}\\
%
%\end{algorithm}

\begin{algorithm}[h]
\caption{Inner-Outer iterations for recovering $N_v \times N_h$ images} \label{alg:main}
Input an initial guess $x^{(\ast,0)}$ for the solution. Set $D^{(0)}=0$. Choose $p>0$ for (\ref{eq:GradientMap}).  %\\
%Set $L$ to multiple of discrete gradient operator \\

For $\ell=1,\ldots,$ until a stopping criterion is satisfied\\
\ \hspace*{2mm} Compute $g = |L x^{(*,\ell-1)}|/\|L x^{(*,\ell-1)}\|_{\infty}$, the rescaled gradient image of the current estimate. \\
\ \hspace*{2mm}  Set $d := 1 - g \dotexp p$ (produces map to [0,1])  \\
\ \hspace*{2mm} Set $D^{(\ell)} = \mbox{diag}(d) D^{(\ell-1)} $ \\
\ \hspace*{2mm} Set $M^{(\ell)} = D^{(\ell)} L$ \\
\ \hspace*{2mm} Run the hybrid algorithm for problem (\ref{eq:seqtrue}) to choose $\lambda_{\ast,\ell}$ and return $x^{(*,k_{\ell})}$. \\
%\ \hspace*{2mm} Set $\Lambda = \mbox{diag}(d)$ \\
% \ \hspace*{2mm}  Check $\| L x^{(*, \ell)} \|$ against $\| L x^{(*,\ell-1)} \|$ and break if needed.    \\
end

Return $\sg{\xreg} = x^{(*,\ell)}$. 

%\JGN{Should the normalization of the gradient be explicitly written in the algorithm?}\\
%\SG{I think so (I fixed this, and also took the component-wise absolute value of the numerator). Also, to have coherent notations we may have to consider $g^{(\ell)}$ (instead of $g$), and so on...}\\
%
\end{algorithm}

\sg{The outer iterations in Algorithm \ref{alg:main} need a starting guess $x^{(\ast,0)}$ to compute the first weighting matrix $D^{(1)}$. Although Algorithm \ref{alg:main} can accept any $x^{(\ast,0)}$ (even possibly determined by a different regularization method), a natural choice is to take $x^{(\ast,0)}$ as a constant vector: in this way, for $\ell=1$, we get $L x^{(*,\ell-1)}=0$ so that $D^{(1)}$ is the identity matrix of order $(N_h (N_v-1) + N_v(N_h-1))$, and the first problem in the sequence of quadratic problems (\ref{eq:seqtrue}) is actually a Tikhonov regularization problem with $R(x)=\|Lx\|_2^2$. A constant $x^{(\ast,0)}=0$ will be used for the numerical experiments in Section \ref{sec:results}. }
%This is the choice In this sense, Algorithm \ref{alg:main} is quite self-contained; this is 
%we propose a self-contained approach that starts  %, to keep it completely self contained, we 
%by simply computing the solution for $M = L$ (i.e., by taking the 2-norm of the gradient as a regularizer during the first outer iteration). 

Two appropriate stopping criteria should be set when implementing Algorithm \ref{alg:main}. The first one prescribes how to terminate the hybrid iterations, and should be devised in connection with the regularization parameter choice strategy, as discussed in Section \ref{ssec:choice} (basically, we monitor the stabilization of the approximate $\lambda_{k}$ along consecutive hybrid iterations). The second one prescribes how to terminate the outer iterations: since we would ideally like to iterate 
until there is no more real edge information to recover, we track 
%.   We need a stopping criteria to determine when we no longer have edge information to recover.   We track 
$ \| L x^{(*,\ell)} \|_2$. If we detect that this quantity is decreasing at step $\ell$, then we have likely oversmoothed our solution, so we break out of the (outer) loop. 
 
%   The resulting algorithm is \SGvar{sketched} in Algorithm \ref{alg:main}.     

%For $\ell = 1,\ldots$
%   \begin{itemize}
%     \item Let $v = L x_{\ell-1}$ represent the gradient image of current estimate.  
%      \item Set $v \leftarrow v / \| v \|_{\infty}$.   (normalize)
%      \item Let $d := 1 - v$ (or $(1 - v$.$^2)$)   (map to [0,1])
%      \item Define $\Lambda^{(\ell)} := \mbox{diag}(d)$
%      \item Set $D^{(\ell)} \leftarrow \Lambda^{(\ell)} D^{(\ell-1)}$.   
%    \end{itemize}

%Produces a diagonal weighting $D^{(\ell)}$ entries between 0 and 1.

% \subsection{Inner-Outer Iterations}
% We present the algorithm now assuming use of the 2nd hybrid scheme we discussed above; however, {\it any suitable hybrid solver can replace step one in the loop}.
% 
% If $x_{0}$ known, determine $D^{(0)}$ from $x_{0}$, else $D^{(0)} = I$.  
%For $\ell = 1,\ldots$ 
%    \begin{itemize}
%      \item Using $k$ steps of a hybrid iterative scheme, 
%     \[ x_{\ell,\lambda} = \arg \min_{x \in \mathcal{Z}_k} \| A x - b \|_2^2 + 
%                                                     \lambda_{\ell}^2 \| D^{(\ell-1)} L x \|_2^2  \]       
%       \item Generate $\Lambda^{(\ell)}$ from $x_{\ell,\lambda}$ as described in **
%       \item Set $D^{(\ell)} := \Lambda^{(\ell)} D^{(\ell - 1)}. $
%      \end{itemize}
%      
%   Further computational savings may be possible if any kind of `recycling' from each instance of step 1 can be used.   
 
 %%%%%%%%%%%%%%%%%%%%%%%%%%%%%%%%%%%%%%%
 \subsection{Analyzing the Main Algorithm} \label{ssec:analysis}
 
 %\SG{this was originally a different session; I think it makes sense to have it as a subsection of Sect. 4}
 
Directly from the definition of the diagonal weights given in Section \ref{ssec:weights}, it is immediate to state that:
%Because of the way the diagonal weights are defined, we first make the following observations:
\begin{enumerate}
\item If a diagonal entry of $D^{(\ell)}$ becomes 0, it stays 0 at future outer iterations.
\item The $i$th diagonal entry of two consecutive weight matrices is such that 
\begin{equation}\label{eq:consDer}
[D^{(\ell)}]_{ii} \le [D^{(\ell-1)}]_{ii}\,.
\end{equation}
\item Large values in the gradient images incur smaller weights, which are even smaller when $p\gg 1$ is selected.   
\end{enumerate}
%
%This tells us something
%
%
%A smaller number on the diagonal says {\it do not} penalize that pixel (too much) this step because it is a ``real'' edge that should not be smoothed.    
%For entries on the diagonal close to 1, we are enforcing smoothness (on those smaller jumps due to noise).    
To further analyze the behavior of Algorithm \ref{alg:main}, we ignore the subspace constraint $\Gamma^{(\ell)}_{k_\ell}$ in (\ref{eq:seqtrue}) and consider the sequence of minimization problems
 \begin{equation}\label{eq:seqnohybrid}
 x^{(\lambda,\ell)}=\arg\min_x \| A x - b \|_2^2 + \lambda_{\ell}^2 \| D^{(\ell)} L x \|_2^2 = \min_x J_{\lambda,\ell}(x),\quad\ell=1,2,3,\dots \,.
 \end{equation}
 %where $\lambda_{\ell}$ is given. 
% \SG{the previous version read $ \min_x \| A x^{(\lambda)} - b \|_2^2 + \lambda^2 \| D^{(\ell)} L x \|_2^2 = \min_x J_{\lambda,\ell}(x)$. I replaced $x^{(\lambda)}$ by $x$ and $\lambda$ by $\lambda_{\ell}$. I now consider $\arg\min$.} 
This assumption simplifies the following derivations, and it is also quite realistic because 
the hybrid solutions of (\ref{eq:seqtrue}) resemble the minimizers of (\ref{eq:seqnohybrid}) for $k$ sufficiently large (i.e., when enough inner hybrid iterations are performed). 
%(\SG{in practice, having iterative solutions can be a very rough approximation...in fact, our algorithm relies on two (nested) cycles of iterations: one for QR factorization, one for JBD solution of the Tikhonov (reweighted) problem -- we should observe this when there is some mismatch in the numerical results.})

\sg{Now let us give a preliminary qualitative overview of the role played by the regularization parameter $\lambda$ in (\ref{eq:seqnohybrid}). }
For smaller $\lambda$, the model fidelity through the operator $A^T A$ is the dominant term in the cost functional.  For larger $\lambda$ the regularization (derivative) term is dominant. Note that solving problem (\ref{eq:seqnohybrid}) is also mathematically equivalent to solving 
 \begin{equation}\label{eq:TikhNE}
 (A^T A + \lambda_{\ell}^2 L^T (D^{(\ell)})^2 L) x = A^T b.
 \end{equation}
In this reformulation it is easier to see that the second term plays the role of a diffusion-like operator (see, for instance, \cite{Bardsley, VogelOman}). 
%\\\MEK{if John observed this in his paper, we might need to cite}.\\  
 Diffusion operators are frequently used in denoising (see, for instance, \cite{ChenMacLachlanKilmer} and references therein), and so we may also state that %one interpretation is that 
a large $\lambda$ favors denoising. When $D^{(1)} = I$ (e.g., when Algorithm \ref{alg:main} runs with a constant initial guess $x^{(\ast,0)}$), \sg{the applied regularization just enforces smoothness on the derivative.}  
%constraint is just a smoothness constraint on the derivative.  
In our problems, we know that the image is not smooth, so overly smooth solution estimates result in a large data-fidelity mismatch.  Therefore, any reasonable parameter selection routine, would not select a $\lambda$ too large, since to do so would return a smooth solution at the price of a large residual; this is surely not the case when the discrepancy principle is adopted. However, if $x^{(*,0)}$ \sg{is already quite successful in revealing the edges of $\xtrue$}, then, for any vector $v$ of appropriate size,   
%\sg{Looking at the evolution of the parameter $\lambda$ (along the outer iterations) for $\ell=0,1,2,\dots$.}
%
%Suppose the image we wish to reconstruct is binary, and assume $x^{(*,0)}$ did reveal most of the transitions between the two values.    
%For any $v$, 
       \[ 
       \| D^{(1)} L v \|_{\sg{2}} < \| L v \|_{\sg{2}}, 
       \]   
%\JGN{did you mean to leave a subscript off the norm here, or should it be $\|\cdot\|_2^2$?}\\
%\SG{I replaced the subscript.}\\
so that, when running the second outer iteration of Algorithm \ref{alg:main}, a minimizer of (\ref{eq:seqnohybrid}) with larger $\lambda$ does not enforce smoothness uniformly across the entire solution, and could be tolerated.  Indeed,          
 the second outer iteration behaves more like denoising than deblurring, because %(constraint is only nonzero across) 
 \sg{the weights are large in} regions where the edges were not detected\sg{, and those are regions where a good amount of smoothing is still needed}. Thus, we expect the parameter selection \sg{strategy} to choose a larger value of $\lambda$ in the second step. \sg{We expect this to be a common trend also in the following outer iterations}, meaning that there is increased reliance on the denoising provided by the regularization term as the outer iterations progress.  \sg{We can prove that this is exactly the case when the discrepancy principle is used to select $\lambda_{\ell}$ for the $\ell$th problem in (\ref{eq:seqnohybrid}). This result (stated in Proposition \ref{prop:increaseLambda}) requires a fair amount of derivations, and it is preceded by a lemma (Lemma \ref{lem:proDiscrC}) that states that the discrepancy curves, i.e., the curves
\begin{equation}\label{eq:discrCurve}
(\lambda, \|b-Ax^{(\lambda,\ell)}\|_2)\,,\quad \lambda\geq 0
\end{equation}
decrease with respect to the outer iteration counter $\ell$.}

\begin{lemma}\label{lem:proDiscrC}
Assume that the null spaces of $A$ and $D^{(\ell)}L$ intersect trivially, for $\ell=1,2,\dots$. Let $r^{(\lambda,\ell)}=b-Ax^{(\lambda,\ell)}$ be the discrepancy associated to the $\ell$th problem in the sequence (\ref{eq:seqnohybrid}), with $\lambda_{\ell}=\lambda\geq 0$, $\ell=1,2,\dots$. 
%(possible notation clash: $\lambda\neq\lambda_{\ell}$ (i.e., the value satisfying the discrepancy principle, or the parameter selected)). 
Then 
\begin{equation}\label{eq:discrDecr}
\|r^{(\lambda,\ell)}\|_2\geq \|r^{(\lambda,\ell+1)}\|_2\,. %\,,\quad \mbox{for $\lambda\geq 0$ and $\ell=1,2,\dots\,$}.
\end{equation}
\end{lemma}

\begin{proof}
Consider the $\ell$th problem in the sequence (\ref{eq:seqnohybrid}). We first have to introduce some notations. Let
\[
\whA^{(\ell)}=A^TA+\lambda^2 L^T (D^{(\ell)})^2 L\,,
\]
i.e., the matrix appearing on the left-hand side of (\ref{eq:TikhNE}), with $\lambda_{\ell} = \lambda\geq 0$. Then
\[
\|r^{(\lambda,\ell)}\|_2^2=b^T(I - A(\whA^{(\ell)})^{-1}A^T)^2b\,.
\]
Assume that $[D^{(\ell)}]_{i_\ell,i_\ell}\neq [D^{(\ell+1)}]_{i_\ell,i_\ell}$, $i_{\ell}\in\{1_\ell,\dots,h_\ell\}\subset \{1,\dots,n\}$, i.e., the diagonal weighting matrices computed at two consecutive iterations differ at most for $h\geq 0$ entries. Then, thanks to (\ref{eq:consDer}), we can write
\[
(D^{(\ell)})^2=(D^{(\ell+1)})^2+\sum_{i=1}^{h}(d_{i_\ell})^2e_{i_\ell}e_{i_\ell}^T,\quad d_{i_\ell}>0\,.
\]
Let
\[
(D^{(\ell),j})^2=(D^{(\ell+1)})^2+\sum_{i=1}^{j}(d_{i_\ell})^2e_{i_\ell}e_{i_\ell}^T=(D^{(\ell),j-1})^2+(d_{j_\ell})^2e_{j_\ell}e_{j_\ell}^T,\quad j=0,\dots,h
\]
(so that $(D^{(\ell),0})^2=(D^{(\ell+1)})^2$, and $(D^{(\ell),h})^2=(D^{(\ell)})^2$). Correspondingly, let $r^{(\lambda,\ell),j}=b-Ax^{(\lambda,\ell),j}$, where $x^{(\lambda,\ell),j}=(\whA^{(\ell),j})^{-1}A^Tb$, $\whA^{(\ell),j}=A^TA+\lambda^2 L^T (D^{(\ell),j})^2 L$. Note that
\begin{equation}\label{hatAlj}
\whA^{(\ell),j} = \whA^{(\ell),j-1}+\lambda^2(d_{j_\ell})^2L^Te_{j_\ell}(L^Te_{j_\ell})^T\,,\quad j=1,\dots,h\,,
\end{equation}
i.e., the difference between two consecutive matrices of the form $\whA^{(\ell),j}$ is a symmetric rank-1 matrix. Let us fix $j\in\{0,\dots,h-1\}$. If we show that
\begin{equation}\label{newGoal}
\|r^{(\lambda,\ell),j+1}\|_2^2\geq\|r^{(\lambda,\ell),j}\|_2^2\,,
\end{equation}
we can conclude that $\|r^{(\lambda,\ell)}\|_2^2:=\|r^{(\lambda,\ell),h}\|_2^2\geq\|r^{(\lambda,\ell),h-1}\|_2^2\geq\dots\geq\|r^{(\lambda,\ell),0}\|_2^2=:\|r^{(\lambda,\ell+1)}\|_2^2$. Thanks to the Sherman-Morrison formula (see, e.g., \cite{GvL}),
\[
(\whA^{(\ell),j+1})\inv = (\whA^{(\ell),j})\inv - (\whA^{(\ell),j})\inv ww^T(\whA^{(\ell),j})\inv,\quad\mbox{where}\quad w=\frac{\lambda d_{j_\ell}L^Te_{j_\ell}}{\sqrt{1+\lambda^2 d_{j_\ell}^2(L^Te_{j_\ell})^T(\whA^{(\ell),j})\inv (L^Te_{j_\ell})} }
\]
and where we have used twice that $(\whA^{(\ell),j})\inv$ is symmetric positive definite. Then
\begin{eqnarray*}
\|r^{(\lambda,\ell),j+1}\|_2^2&=&b^T\left(I - A((\whA^{(\ell),j})\inv - (\whA^{(\ell),j})\inv ww^T(\whA^{(\ell),j})\inv)A^T\right)^2b\\
&=& b^T(I - A(\whA^{(\ell),j})\inv A^T+\whw\whw^T)^2 b\\
&=& \underbrace{b^T(I - A(\whA^{(\ell),j})\inv A^T)^2b}_{=\|r^{(\lambda,\ell),j}\|_2^2}%\\
%& & + \underbrace{b^T(I - A(\whA^{(\ell),j})\inv A^T)\whw\whw^Tb + b^T\whw\whw^T(I - A(\whA^{(\ell),j})\inv A^T)}_{=:\delta^{(1)}} \\
%& & 
+ \underbrace{2b^T\whw\whw^T(I - A(\whA^{(\ell),j})\inv A^T)b}_{=:\delta^{(1)}} %\\
%& & 
+ \underbrace{b^T\whw((\whw)^T(\whw))\whw^Tb}_{=:\delta^{(2)}}
\end{eqnarray*}
where we have defined $\whw = A(\whA^{(\ell),j})\inv w$. Here $\delta^{(2)}\geq 0$, since $\whw\whw^T$ is symmetric positive semi-definite. If we show that $\delta^{(1)}\geq 0$, then we have (\ref{newGoal}) and therefore (\ref{eq:discrDecr}). Let 
%$Q\lj R\lj=D\lj L$ be the economy-size QR factorization of $D\lj L$ and let 
\[
A=U\lj \Sigma\lj (X\lj)\inv,\quad D\lj L=V\lj \Gamma\lj (X\lj)\inv
\]
be the GSVD of the matrix pair $(A, D\lj L)$, \cite{GvL}. Then (exploiting the properties of the matrices involved in the decompositions above), (\ref{hatAlj}), and the definitions of $w$ and $\whw$, 
\begin{eqnarray*}
\delta^{(1)}&=&\underbrace{\frac{2\lambda^2(d_{j_\ell})^2}{1+\lambda^2 d_{j_\ell}^2(L^Te_{j_\ell})^T(\whA^{(\ell),j})\inv (L^Te_{j_\ell})}}_{=:\mu>0}
b^TA(\whA\lj)\inv L^Te_{j_\ell}e_{j_\ell}^TL(\whA\lj)\inv A^T(I - A(\whA\lj)\inv A^T)b\\
&=&\mu 
\overbrace{b^TU\lj}^{=:\whb^T} 
\underbrace{\Sigma\lj \left((\Sigma\lj)^2+\lambda^2(\Gamma\lj)^2\right)\inv}_{=:D\lj_1} 
\overbrace{X^TL^Te_{j_\ell}}^{=:\whd}
e_{j_\ell}^TLXD\lj_2(U\lj)^Tb
=\mu\whb^T D\lj_1\whd\, \whd^T D\lj_2\whb \,,
\end{eqnarray*}
where
\[
D\lj_2=((\Sigma\lj)^2+\lambda^2(\Gamma\lj)^2)\inv\Sigma\lj\left(I-(\Sigma\lj)^2((\Sigma\lj)^2+\lambda^2(\Gamma\lj)^2)\inv\right)\,.
\]
Note that both $D\lj_1$ and $D\lj_2$ are diagonal matrices with positive diagonal entries. One can easily show that $D\lj_1\whd\, \whd^T D\lj_2$ has at most rank-1, its only nonzero eigenvalue is $(D\lj_2\whd)^T(D\lj_1\whd)>0$ (since we are summing positive quantities), and is diagonalizable. Therefore it is positive semi-definite and $\delta^{(1)}\geq 0$, which concludes our proof.
\end{proof}

\sg{\begin{proposition}\label{prop:increaseLambda}
Assume that the null spaces of $A$ and $D^{(j)}L$ intersect trivially, for $j=1,2,\dots$, and assume that the regularization parameter $\lambda_\ell$ in (\ref{eq:seqnohybrid}) is determined according to the discrepancy principle. Then 
\begin{equation}\label{lambdavar}
\lambda_{\ell} \le \lambda_{\ell-1},\quad \ell=1,2,3,\dots\,.
\end{equation}
\end{proposition}}
\begin{proof} %\SG{very tentative!}\\
Take $\ell\geq 1$. Thanks to (\ref{eq:discrDecr}), the discrepancy curve (\ref{eq:discrCurve}) computed for the $\ell$th problem lays below the discrepancy curve computed for the $(\ell-1)$th problem. Also, one can easily prove that, for a fixed $\ell$ and a variable $\lambda\geq 0$, the discrepancy curve increases with respect to $\lambda$. Since applying the discrepancy principle consists in solving
\[
\|b-Ax^{(\lambda,\ell)}\|_2=\tau\|\eta\|_2\,,
\]
this implies (\ref{lambdavar}).
%Basic idea: we first prove that the discrepancy curve for the $\ell$th problem (i.e., a plot of $\|b-Ax_{\ell}^{\lambda}\|_2$ versus $\lambda$), lays below the . Since the discrepancy curves are increasing with respect to $\lambda$, and applying the discrepancy principle amounts to solving
%\[
%\|b-Ax_{\ell}^{\lambda}\|_2=\tau\|\eta\|_2\,,
%\]
%this implies (\ref{lambdavar}).  \\
%The tricky part is the first part of the proof. Possible approach (?): fix $\lambda$, and reformulate the discrepancy principle as
%\[
%b^T(I - A(A^TA+\lambda_{\ell}^2L^T(D^{(\ell)})^2L)^{-1}A^T )^2b=\tau^2\|\eta\|_2^2\,.
%\]
%Now assume, without loss of generality (?), that the $j$th diagonal entry is such that $[D^{(\ell)}]_{jj}=0$ and $[D^{(\ell-1)}]_{jj}> 0$. Then, we may use the Sherman-Morrison formula to link the expressions of \\$(A^TA+\lambda_{\ell}^2L^T(D^{(j)})^2L)^{-1}$ for $j=\ell-1,\ell$. 
\end{proof}

More generally, \sg{we experimentally find that the parameters $\lambda_{*,\ell}$ (often) satisfy $\lambda_{*,\ell} \le \lambda_{*,\ell-1}$ also when inner hybrid iterations are provided (i.e., when problem (\ref{eq:seqtrue}) is solved instead of (\ref{eq:seqnohybrid})), and that this also happens when the $\mathcal{L}$-curve criterion is used instead of the discrepancy principle (see Section \ref{sec:results})}.
%\\\MEK{I want to prove the sequence of selected parameters is non-increasing.  I can prove that the ${\mathcal L}$-curves, do not cross, and that the curve for $\ell$ is below the one for $\ell + 1$, but I don't know that that implies what I want it to.} \\   
%\JGN{I agree -- at least I don't see how it would imply the parameters are non-increasing. But still it might be nice
%   to include the proof in the paper.}\\
%\SG{I don't see it either. However, managing to prove that the curve for $\ell$ is below the one for $\ell + 1$ would be enough to prove that the discrepancy curve for $\ell$ lays below the discrepancy curve for $\ell-1$. In the previous draft it was assumed that the image $\xtrue$ is binary: do we need this in the proof?}\\
%\SG{I removed the figure (looked empty anyway, and we have similar ones in Section \ref{sec:results}).}  
 
% \begin{figure}
% %\includegraphics[height=5in,width=6in]{lcurves_peppers1.jpg}
% \caption{${\mathcal L}$-curves for example 1 which illustrates the relationship of the curves and corners as a function of outer iteration $\ell$.}
% \end{figure}
 
We conclude this section by stressing again that Algorithm \ref{alg:main} works well when the initial solution estimates $x^{(\ast,0)}$ or $x^{(\ast,1)}$ are able to recognize at least one decent edge: if the initial estimates are too smooth, which might be the case for example in deblurring with a large noise to signal ratio, our method can fail, as the mapping from the normalized gradient image to [0,1] (described in equation (\ref{eq:GradientMap})) may be totally inappropriate.   
The qualitative performance of Algorithm \ref{alg:main} also depends on the success of the parameter choice strategy but, as we show in the next section, it works quite well for many applications, providing significant enhancement in a completely automated fashion over the initial solution estimate.     
 
 %%%%%%%%%%%%%%%%%%%%%%%%%%%%%%%%%%%%%%%
  \section{Numerical Results}  \label{sec:results}
 All experiments were done in Matlab version \sg{9.1}. The test problems come from the MATLAB software packages IR Tools \cite{Gazzola2019} and
 AIR Tools II \cite{AIRToolsII}. 
 
In this section, we would like to assess the potentialities of our new algorithm, and we mainly do so by comparing it to other similar inner-outer iteration methods for Total Variation regularization. Namely, we consider the IRN method for Total Variation (IRN-TV) proposed in \cite{IRNtv}: similarly to our new method, IRN-TV updates some weights and defines a sequence of least squares problems (outer iterations), and solves each iteratively reweighted least squares problem iteratively (inner iterations). IRN-TV and our new method mainly differ in the way the weights are defined. The IRN-TV weights $D^{(\ell)}$ for the $\ell$th least-squares problem (\ref{eq:seqtrue}) have the following expression
\begin{equation}\label{eq:IRNTVweights}
D^{(\ell)} = W(x^{(\ast,\ell-1)})=I_2 \otimes (W_R^{(\ell)})^{1/2},\quad W_R^{(\ell)}=\mbox{diag}\left( \vvec(L_vX^{(\ast,\ell-1)})^2+ \vvec(X^{(\ast,\ell-1)}L_h^T)^2 \right)^{(q-2)/2},
\end{equation}
where $x^{(\ast,\ell-1)}$ is the solution of the $(\ell-1)$th least squares problem, $I_2$ is the identity matrix of order 2, and $1\leq q<2$. Taking $q=1$, fixing $y=\vvec(Y)$, $Y\in\R^{N_h\times N_v}$, and evaluating the above weights in $y$ leads to $\|W(y)Ly\|_2^2=\mbox{TV}(y)$, where $L$ is the discrete gradient operator (\ref{eq:gradient}), and $\mbox{TV}(\cdot)$ is defined as in (\ref{eq:TV}). Because of this, we can expect that the reweighted regularization term in IRN-TV is a good approximation of $\mbox{TV}(\xtrue)$ when $x^{(\ast,\ell-1)}$ is a good approximation of $\xtrue$. Moreover, when $ [\vvec(L_vX)^2+ \vvec(XL_h^T)^2]_i$ is big (typically when there is an edge between the $i$th and the $(i+1)$th pixel), then $[(W_R^{(\ell)})^{1/2}]_{ii}$ is small (so that the corresponding gradient component is not much penalized in the minimization process); viceversa, when $ [\vvec(L_vX)^2+ \vvec(XL_h^T)^2]_i\simeq 0$ (typically corresponding to smooth portions of $X$) then $[(W_R^{(\ell)})^{1/2}]_{ii}$ is huge (so that the corresponding gradient component is penalized even more in the minimization process). IRN-TV and our new method also differ in the way each least squares problem is solved: IRN-TV assumes that a good value of the regularization parameter is fixed, and employs CGLS for the inner iterations. \sg{We also make comparisons with two inner-outer iterative schemes that can be collocated somewhere in between IRN-TV and our new algorithm, namely: (a) we define the weights at the beginning of each outer iteration as in (\ref{eq:IRNTVweights}), and we solve each quadratic problem by the hybrid method described in Section \ref{ssec:hybridgen}, with the added benefit of adaptively choosing the regularization parameter; (b) we define the weights at the beginning of each outer iteration as in Section \ref{ssec:weights}, and we solve each quadratic problem by CGLS, having fixed a regularization parameter.}
%%% BASICALLY: (1) fix the linear solver, and compare the weight; (2) fix the weights, compare the linear solvers.

%% TV weights: we can see how the hybrid method performs with weights that reproduce total variation.

% This method is validated by performing comparisons to other TV solvers, and therefore we are happy in having comparisons with this IRN-TV method only. 
 
% \sg{For each test problem we display a number of results with the goal of assessing the behavior of our new algorithm, even with respect to different : specifically, we monitor the quality of the reconstructions at each outer iteration (i.e., when the regularizer is updated), as well as the progression of the inner iterations,  with respect to the value of p, parameter choice strategy...} the potentialities 
 
% \JGN{Some questions on what to display in the results:
% \begin{itemize}
% \item
% I'm not sure how to display the weights -- it would be hard to explain them, so I didn't include them.
% \item
% We may not want to include all of these examples, but if so, some of the information
% is redundant (like the true image phantoms).  And most of the text is cut-and-paste,
% so we'll need to make some modifications.
% \item
% I put in page breaks between examples for this draft so that it might be easier to see what goes
% with what -- we can change this later.
% \end{itemize}
% }
% 
%\sg{Just keep the first example with a lot of tests and displays; make sure that its behavior is sort-of-typical for other tests, too.}

%\newpage

\subsection{Examples from X-ray CT}

\begin{figure}[htbp]
\begin{center}
\begin{tabular}{cc}
\includegraphics[width=5cm]{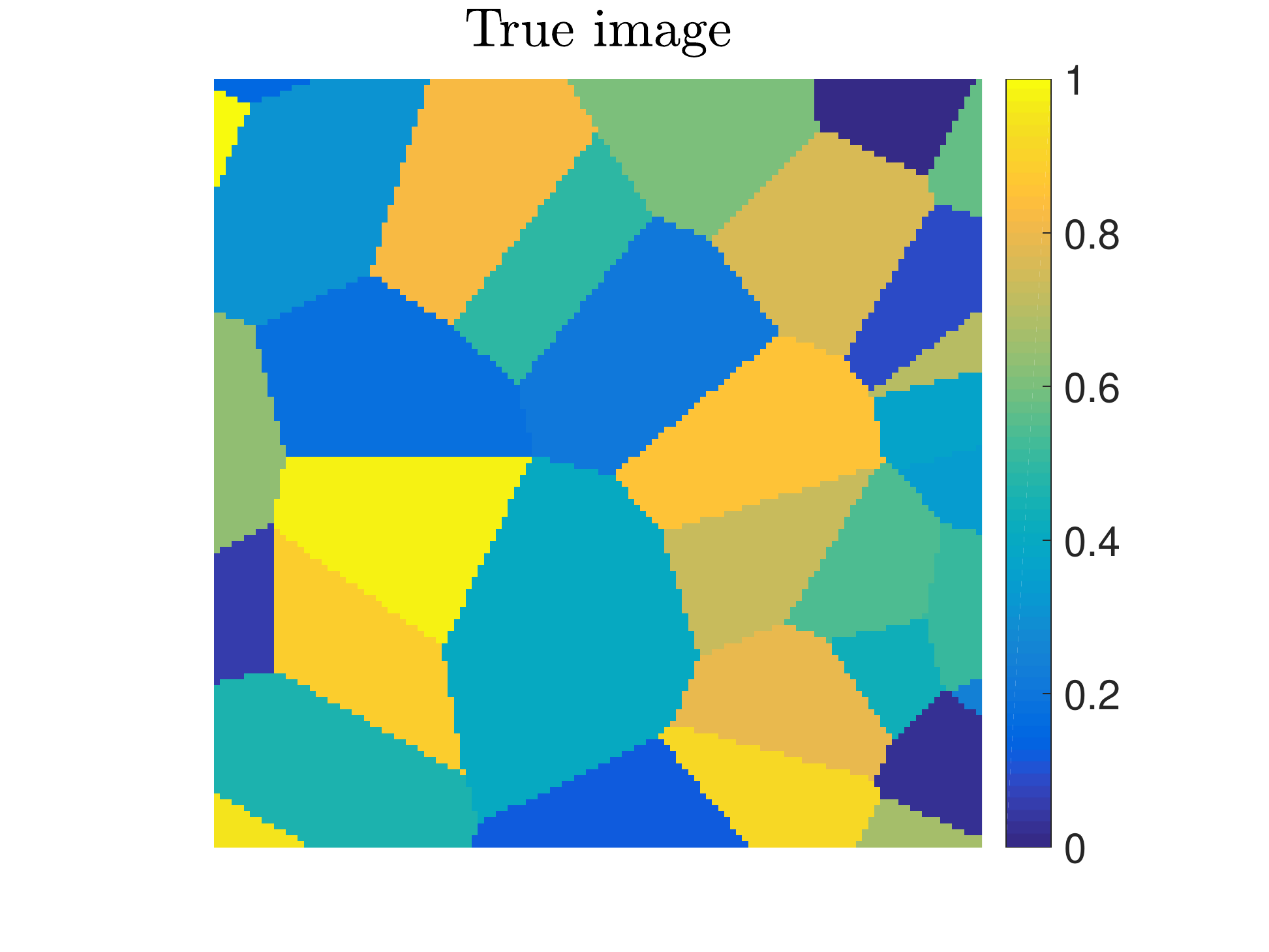} &
\includegraphics[width=5cm]{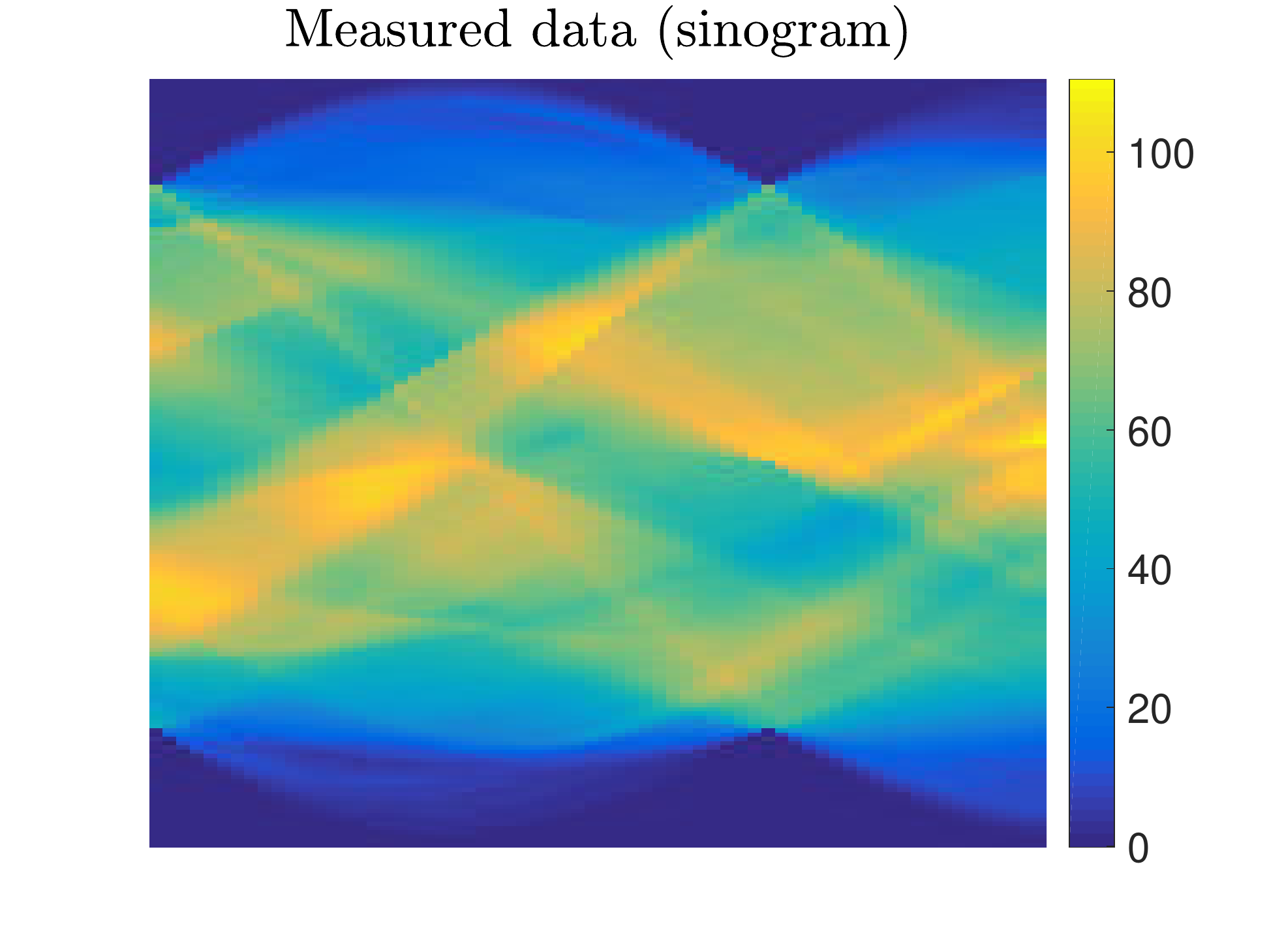}
\end{tabular}
\end{center}
\caption{\emph{Parallel-beam CT} test problem. Left frame: true image phantom $x$. Right frame: the measured data (sinogram) $b$.}
\label{fig:Radon1aData}
\end{figure}

\paragraph{Parallel-beam CT} %: `grains' phantom}
%An extensive analysis for this test problems; not so many plots will be displayed in the following. The plot displayed here are representative of the ones obtained for the other test problems. 
We use IR Tools to setup the following X-ray tomography simulation:
\begin{itemize}
\item
Generate a simulated true phantom image called ``grains" from AIR Tools II, of size $128 \times 128$ pixels.
\item
Construct noise-free projections, along with the matrix $A$ that simulates the ray trace forward operator,
assuming a parallel beam X-ray transmission, with data collected at angles $0, 2, \ldots, 130$ degrees.
\item
Add 0.1\% normally distributed (white Gaussian) noise.
\end{itemize}
More specifically, the test problem is generated with the following MATLAB statements:
\begin{verbatim}
     ProblemOptions = PRset('angles', 0:2:130, 'phantomImage', 'grains');
     [A, b_true, x_true, ProblemInfo] = PRtomo(128, ProblemOptions);
     b = PRnoise(b_true, 1e-3);
\end{verbatim}
Figure~\ref{fig:Radon1aData} shows the true phantom image, along with the
measured data $b$.

\begin{figure}[htbp]
\begin{center}
\begin{tabular}{cc}
\includegraphics[width=6cm]{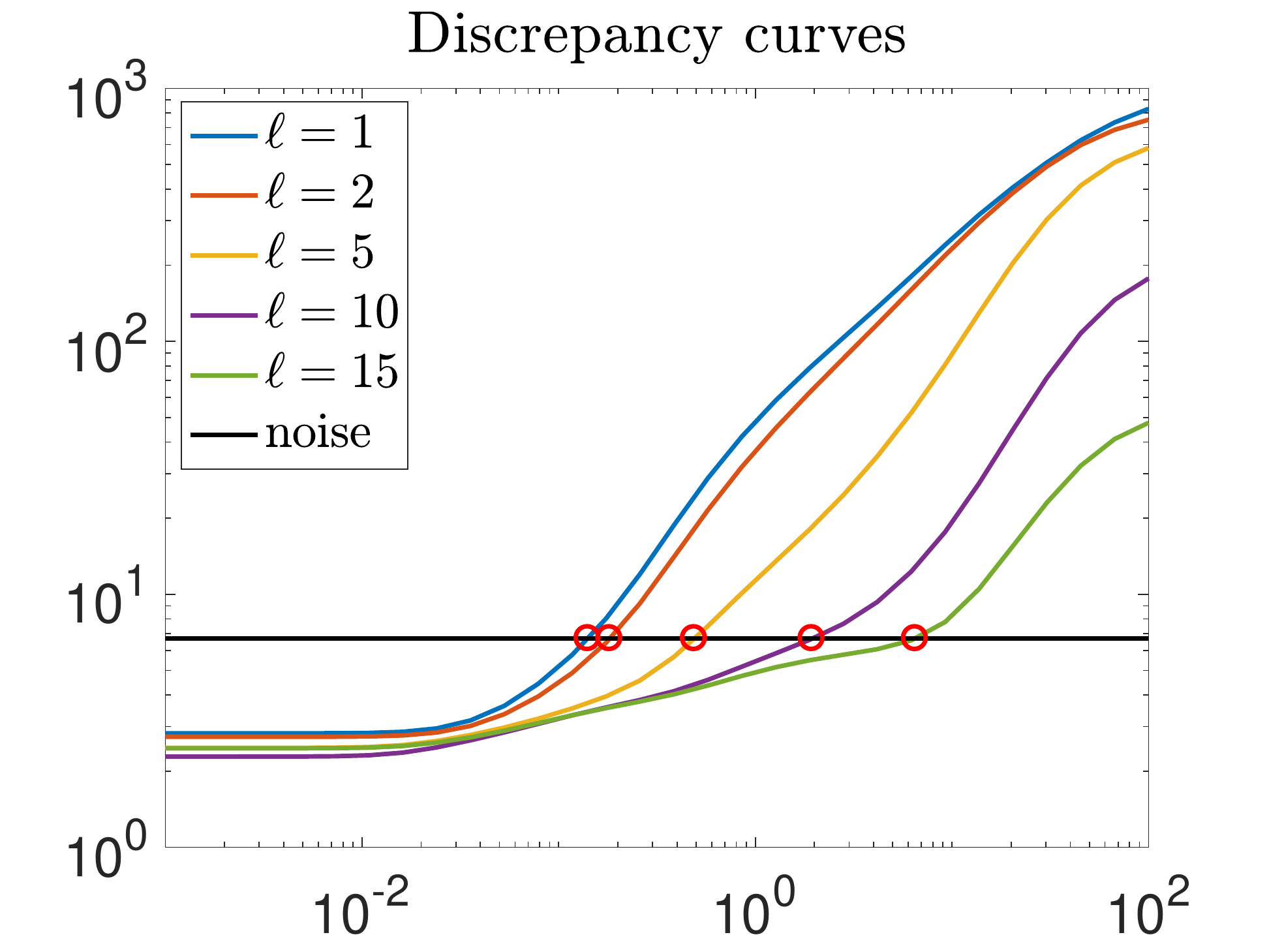} &
\includegraphics[width=6cm]{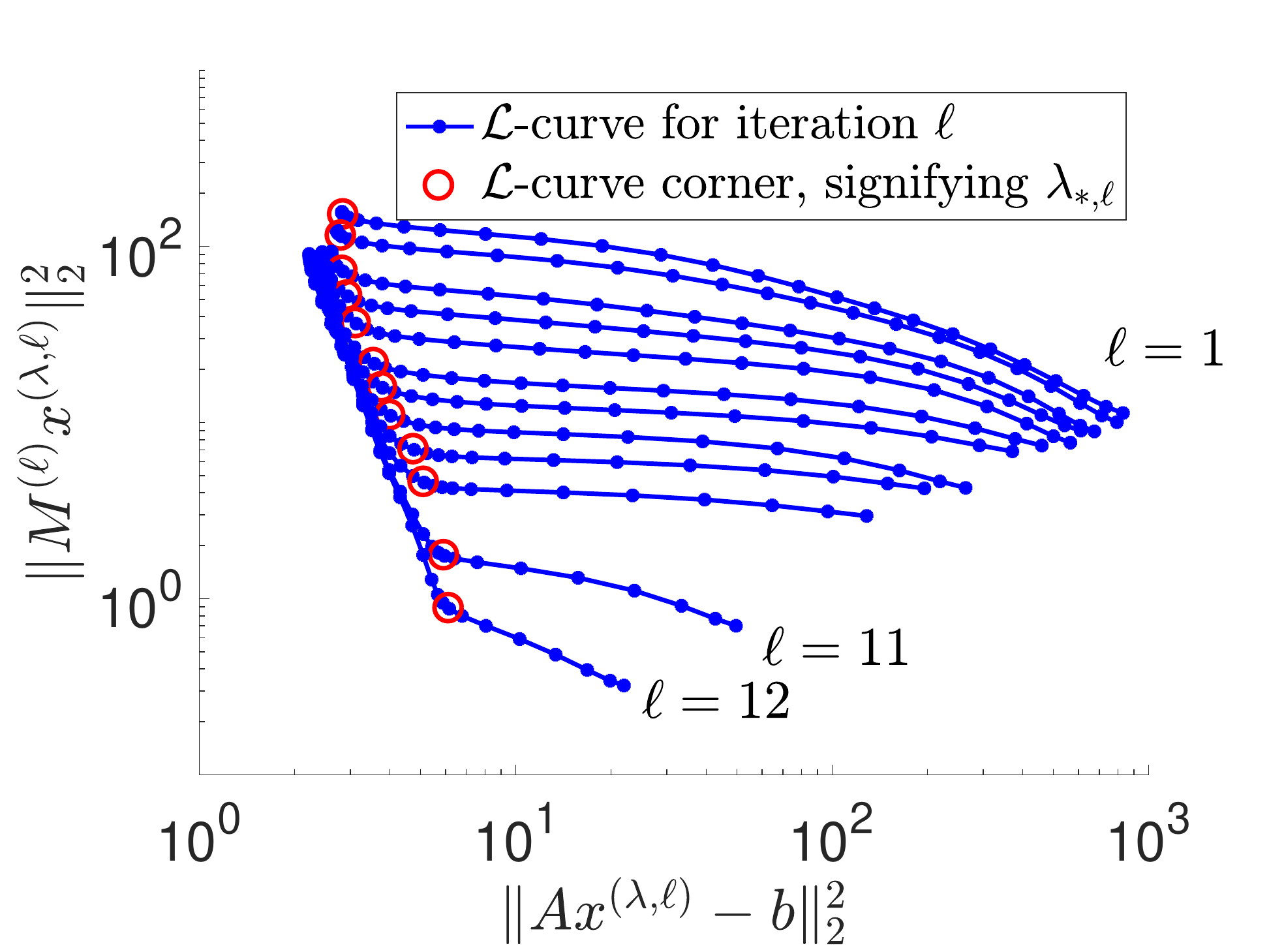} 
\end{tabular}
\end{center}
\caption{\emph{Parallel-beam CT} test problem. Left frame: logarithmic plot of the discrepancy curves, i.e., $\|Ax^{(\lambda,k_{\ell})}-b\|_2$ versus $\lambda$, at the end of some selected cycles of inner iterations (as implied from the text in the legend, these correspond to the outer iterations $\ell=1,2,5,10,15$); the red circles denote the intersection between each discrepancy curve and the noise level line, which corresponds to the chosen
regularization parameter, $\lambda_{*,\ell}$, for the particular outer iteration. Right frame: logarithmic plot of the ${\mathcal L}$-curves for each outer iteration (as implied from the text in the plot, the top curve corresponds to the first outer iteration, $\ell = 1$, and the curves
below this correspond sequentially to iterations $\ell = 2, 3, \ldots, 11, 12$); the red circles
denote corners of each ${\mathcal L}$-curve, which corresponds to the chosen
regularization parameter, $\lambda_{*,\ell}$, for the particular outer iteration.}
\label{fig:Radon1aLCurves}
\end{figure}

First of all, we test our new algorithm for different parameter choice strategies employed within the inner hybrid scheme for generalized Tikhonov regularization (Section \ref{ssec:hybridgen}). More precisely, we consider the discrepancy principle and the ${\mathcal L}$-curve criterion. The left frame of Figure~\ref{fig:Radon1aLCurves} shows a plot of different discrepancy curves, defined analogously to (\ref{eq:discrCurve}) as the graphs of the function $\|b-Ax^{(\lambda,k_{\ell})}\|_2$ versus $\lambda$, and corresponding to different outer iterations $\ell$. The right frame of Figure~\ref{fig:Radon1aLCurves} shows a plot of the different ${\mathcal L}$-curves obtained at the end of each inner iteration cycle. Looking at this graphs, we can see that both the discrepancy curves and the ${\mathcal L}$-curves are nested, in agreement with the analysis performed in Section \ref{ssec:analysis}. 

%When employing the 
%
%: the former behavior is analysed in for each iteration. Observe the nesting
%property of the different ${\mathcal L}$-curves, and how the chosen regularization
%parameter (corresponding to the corner of each ${\mathcal L}$-curve) increases.

When our algorithm is implemented with the stopping criterion described in Section \ref{ssec:mainalgo} (i.e., stopping as soon as $\|Lx^{(\ast,\ell)}\|_2$ decreases during two consecutive outer iterations) we have termination after $\ell=15$ iterations if the discrepancy principle is adopted in the inner iterations, and after $\ell=12$ iterations if the ${\mathcal L}$-curve criterion is adopted in the outer iterations. 

\begin{figure}[htbp]
\begin{center}
\begin{tabular}{cc}
\includegraphics[width=6cm]{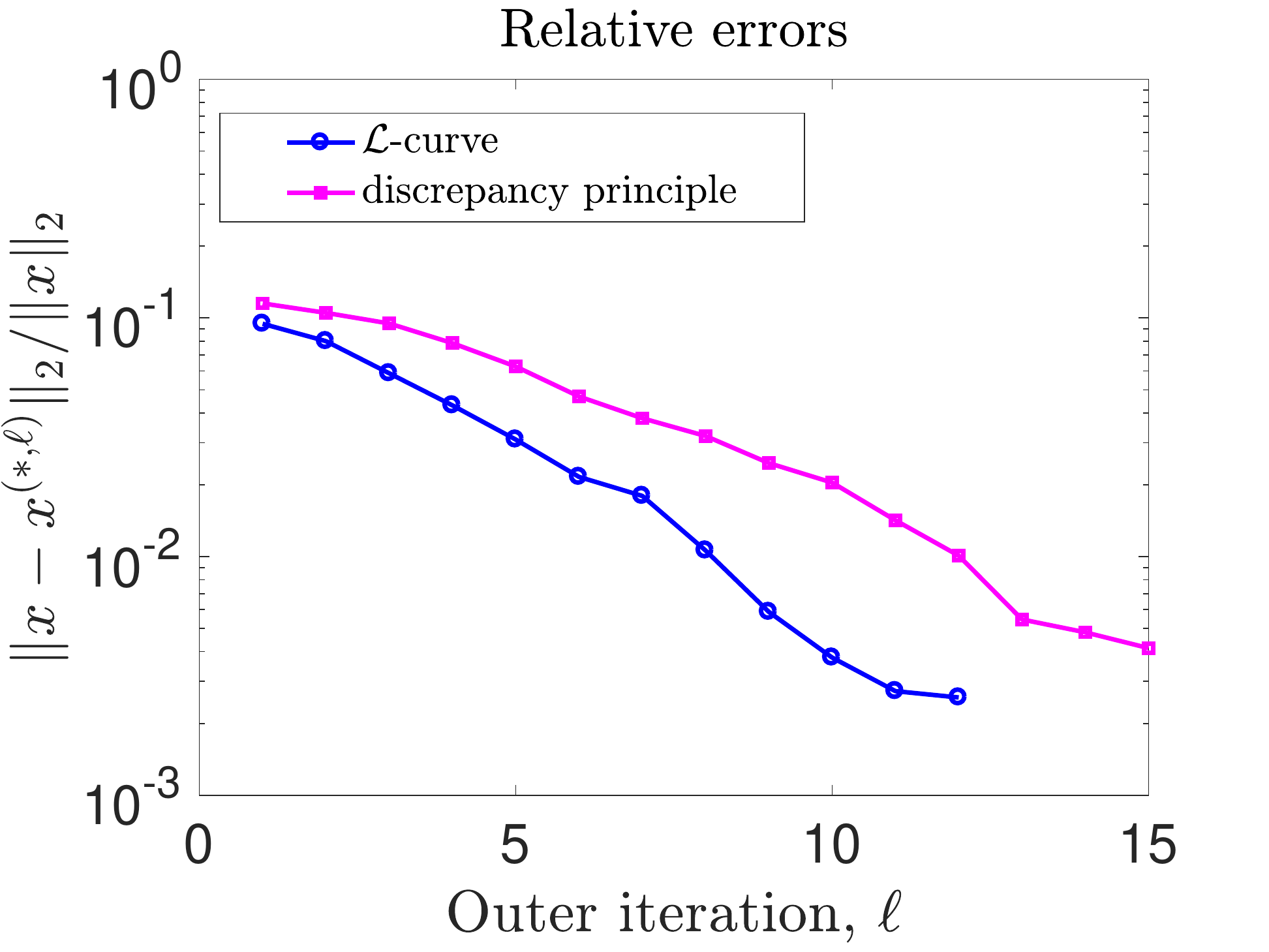} &
\includegraphics[width=6cm]{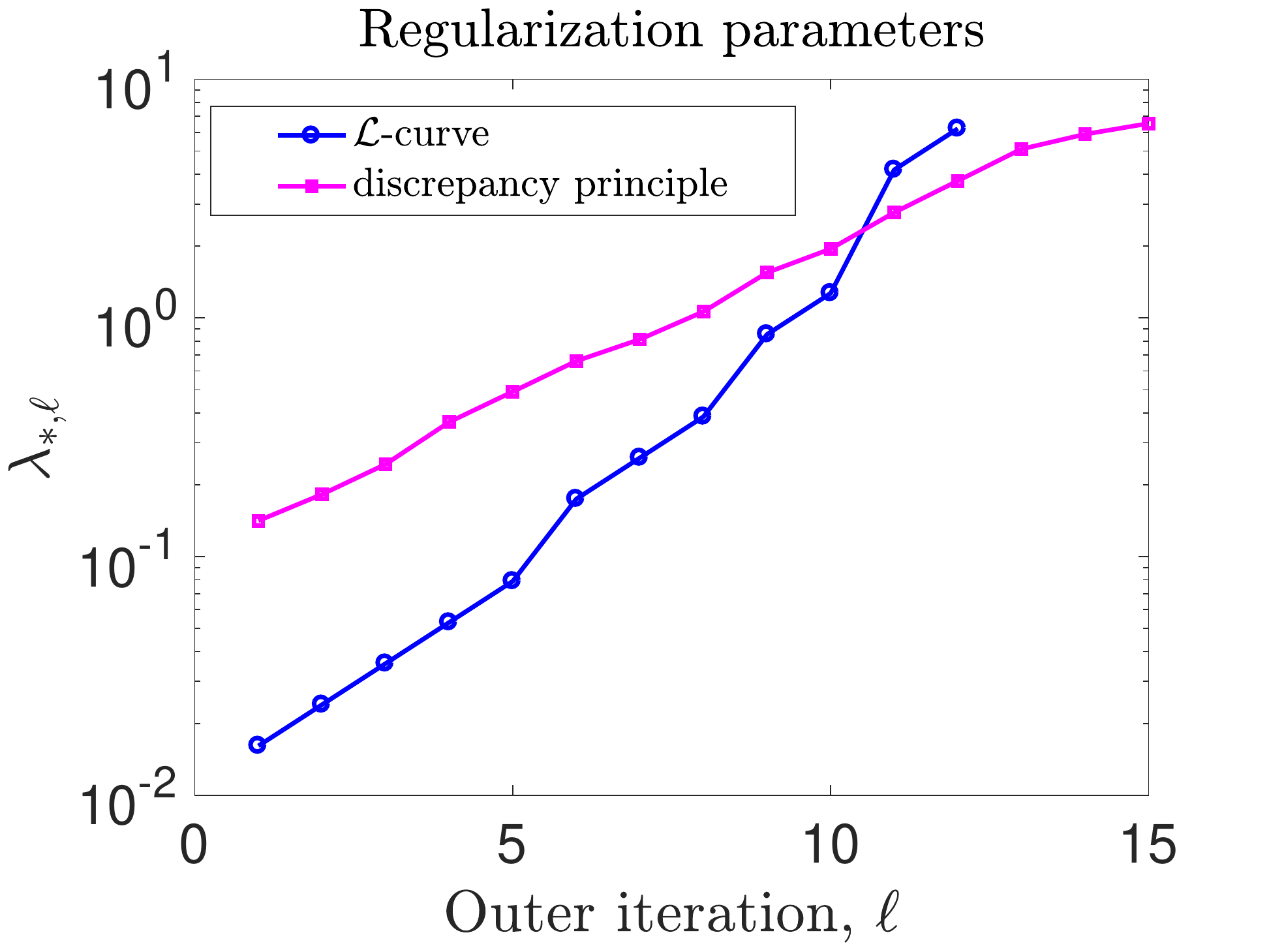}
\end{tabular}
\end{center}
\caption{\emph{Parallel-beam CT} test problem. Relative errors and regularization parameters values at each outer iterations $\ell$,until the stopping criterion is satisfied. Both the discrepancy principle and the ${\mathcal L}$-curve criterion are considered.}
%Basically: effect of different parameter choice methods.
\label{fig:Radon1aIterations}
\end{figure}

Figure~\ref{fig:Radon1aIterations} shows a plot of the relative errors
and chosen regularization parameters at each outer iteration, when the discrepancy principle and the ${\mathcal L}$-curve criterion are used to select the regularization parameter during the inner hybrid iterations. Note that, as expected, the regularization parameters increase as
the outer iteration proceeds: this illustrates the property derived in Section \ref{ssec:analysis} for the discrepancy principle, which experimentally holds also for the ${\mathcal L}$-curve criterion. This behavior of the regularization parameter is meaningful and desirable: indeed, as the outer iterations increase, we are recovering approximate solutions of enhanced quality that result in improved weights for the regularization term, which in turn should be weighted more to achieve reconstructions of even higher quality. 
%This can also explain why the reconstructions obtained by the ${\mathcal L}$-curve criterion are of better quality than the discrepancy principle one: namely, they are obtained selecting a bigger regularization parameter during the final outer iterations. 

%
%terminated after $\ell=12$ outer iterations, where at iteration $\ell$ the hybrid 
%method determined an estimate of the regularization parameter $\lambda_{*,\ell}$,
%and a corresponding reconstructed image, $x^{(*,\ell)}$.  

%\begin{figure}[htbp]
%\begin{center}
%\includegraphics[width=7cm]{FigsExampleRadon1a/LCurves1a} 
%\end{center}
%\caption{${\mathcal L}$-curves for each iteration, for the test problem. As implied from
%the text in the plot, the top curve corresponds to the first outer iteration, $\ell = 1$, and the curves
%below this correspond sequentially to iterations $\ell = 2, 3, \ldots, 11, 12$.  The red circles
%denote corners of each ${\mathcal L}$-curve, which correspond to the chosen
%regularization parameter, $\lambda_{*,\ell}$ for the particular iteration.}
%\label{fig:Radon1aLCurves}
%\end{figure}

\begin{figure}[htbp]
\begin{center}
\begin{tabular}{cc}
%\footnotesize{$x^{(\ast,1)}$, discrepancy p., $\lambda_{\ast,1}=0.1413$} &
%\footnotesize{$x^{(\ast,13)}$, discrepancy p., $\lambda_{\ast,15}=6.5309$} \\
\includegraphics[width=5cm]{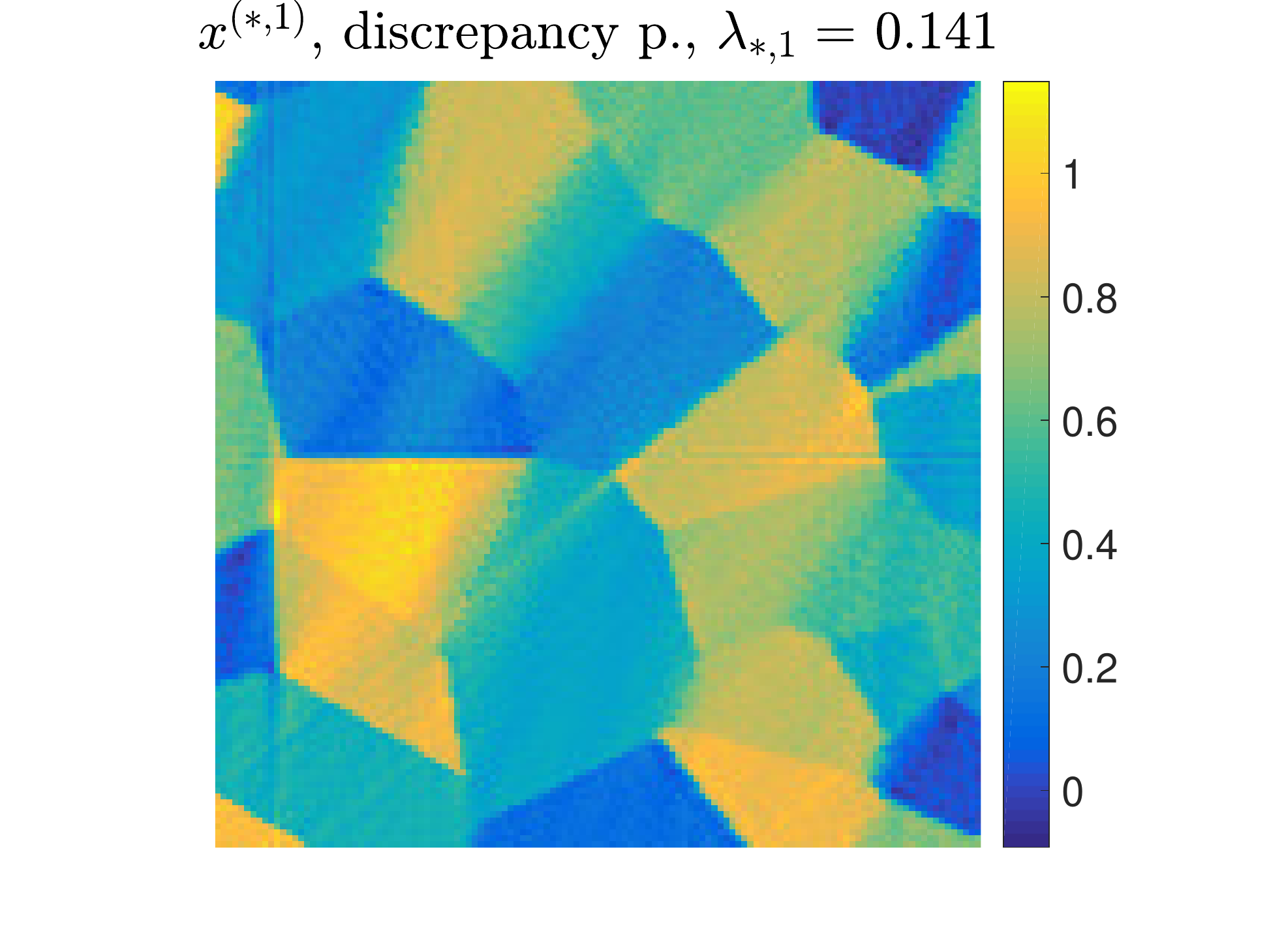} &
\includegraphics[width=5cm]{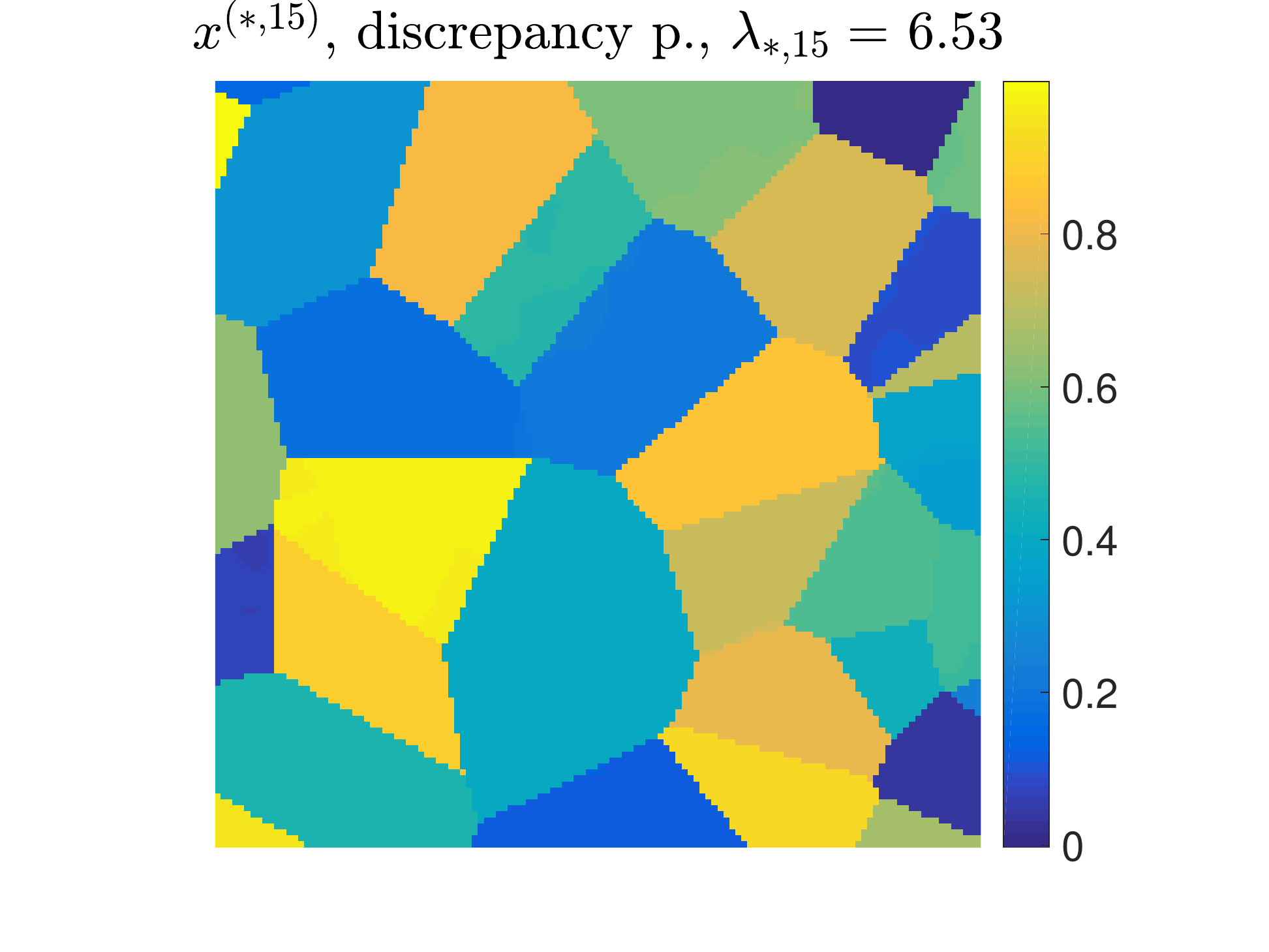}\\
%\footnotesize{$x^{(\ast,1)}$, ${\mathcal L}$-curve, $\lambda_{\ast,1}=0.0161$} &
%\footnotesize{$x^{(\ast,13)}$, ${\mathcal L}$-curve, $\lambda_{\ast,12}=6.2102$} \\
\includegraphics[width=5cm]{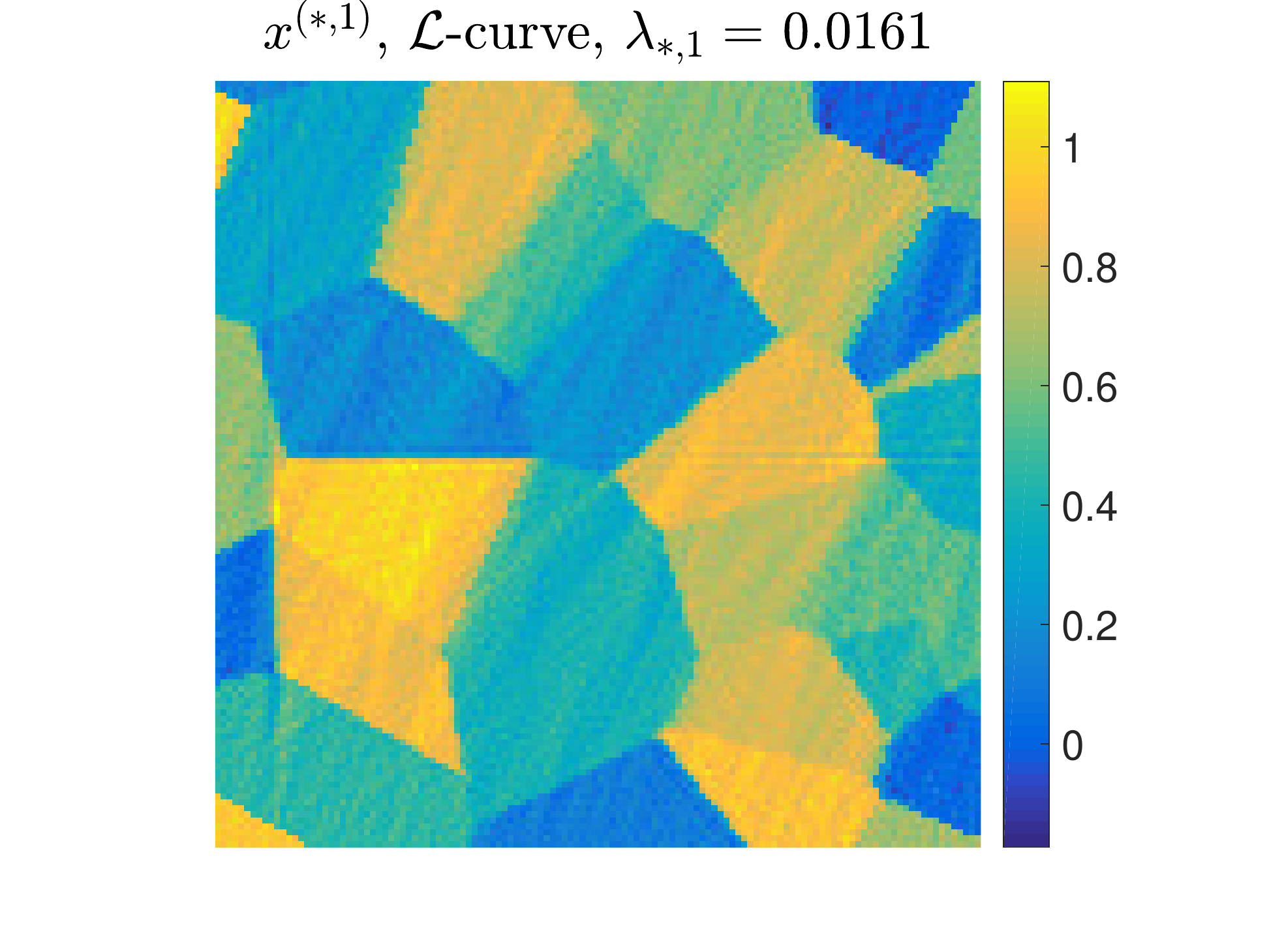} &
\includegraphics[width=5cm]{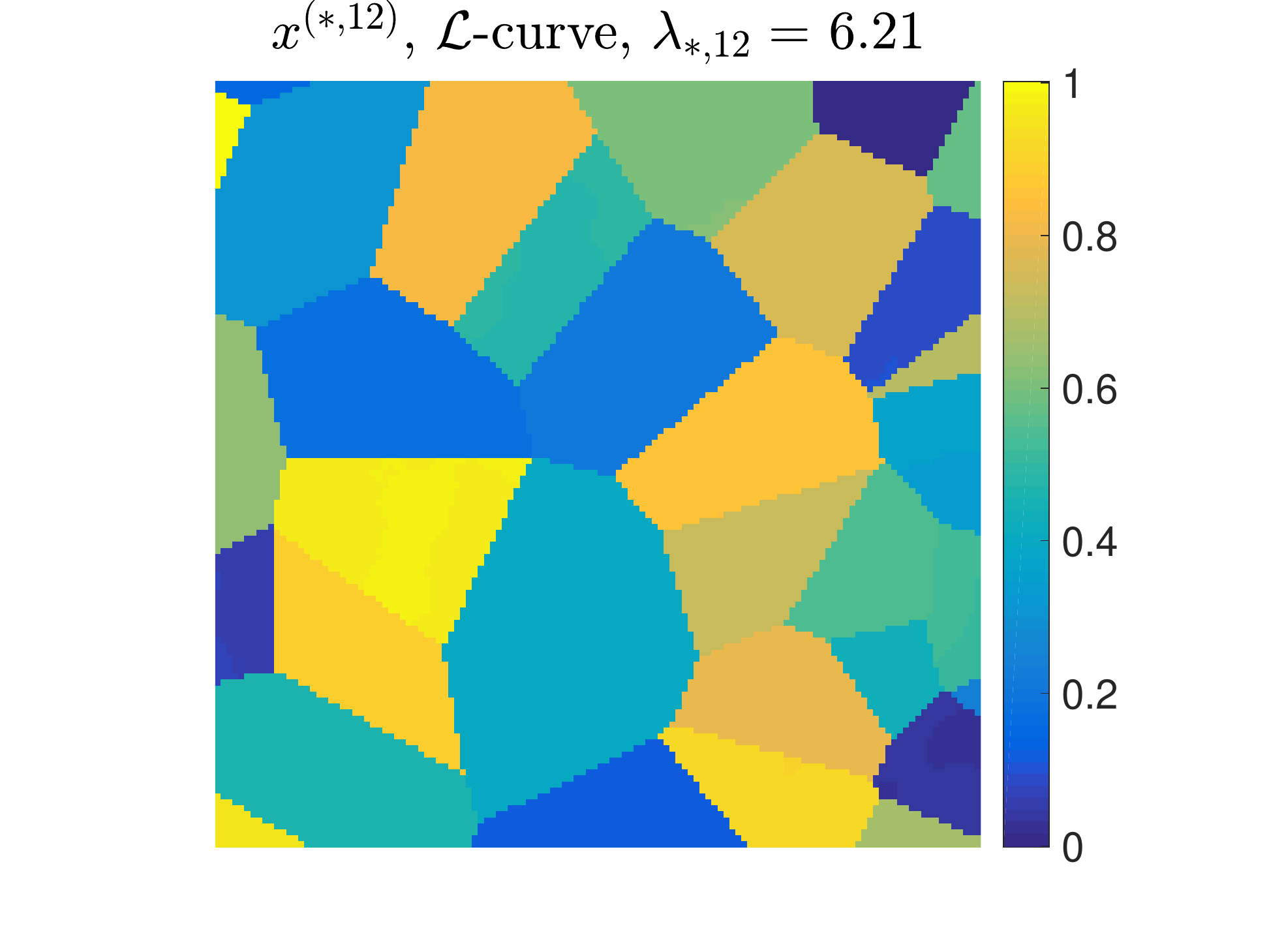}
\end{tabular}
\end{center}
\caption{\emph{Parallel-beam CT} test problem. On the left are the reconstructions 
at the initial outer iteration, and on the right are the reconstructions at the
final outer iterations. The corresponding regularization parameters chosen
by the hybrid method are displayed above each image. Solutions computed by the discrepancy principle and the ${\mathcal L}$-curve are displayed above and below, respectively.}
\label{fig:Radon1aSolutions}
\end{figure}

Computed reconstructions for the first outer iteration (that is, $x^{(*,1)}$), and for the
final outer iteration (that is, $x^{(*,15)}$ or $x^{(*,12)}$, depending on the parameter choice strategy) are shown 
in Figure~\ref{fig:Radon1aSolutions}. As we can see from these plots, there is a 
significant improvement in the reconstructions, and in particular 
the edges at the final outer iteration are much sharper than in the initial outer iteration. Also, pixel intensities are very close to those of the true image.

\begin{figure}[htbp]
\begin{center}
\begin{tabular}{cc}
\hspace{-0.2cm}\includegraphics[width=8.8cm]{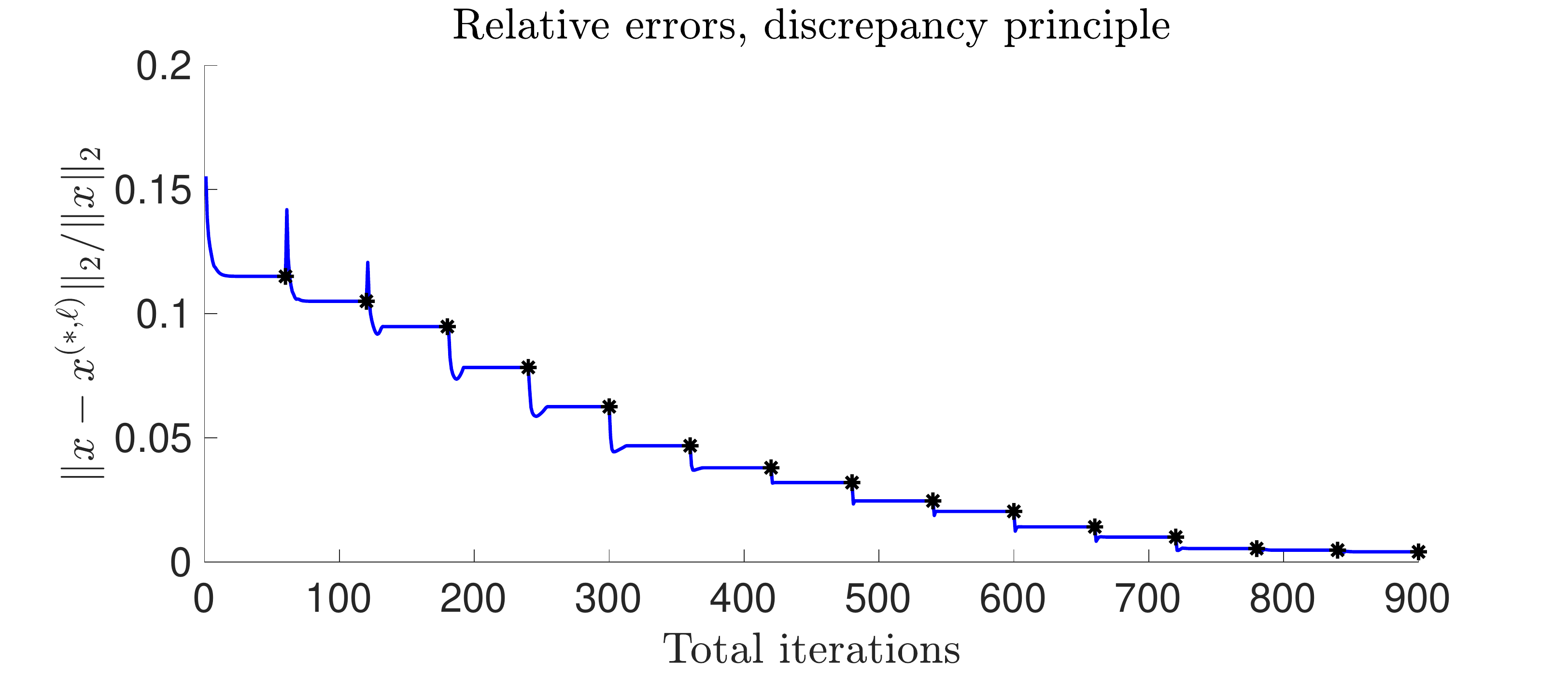} & 
\hspace{-1cm}\includegraphics[width=8.8cm]{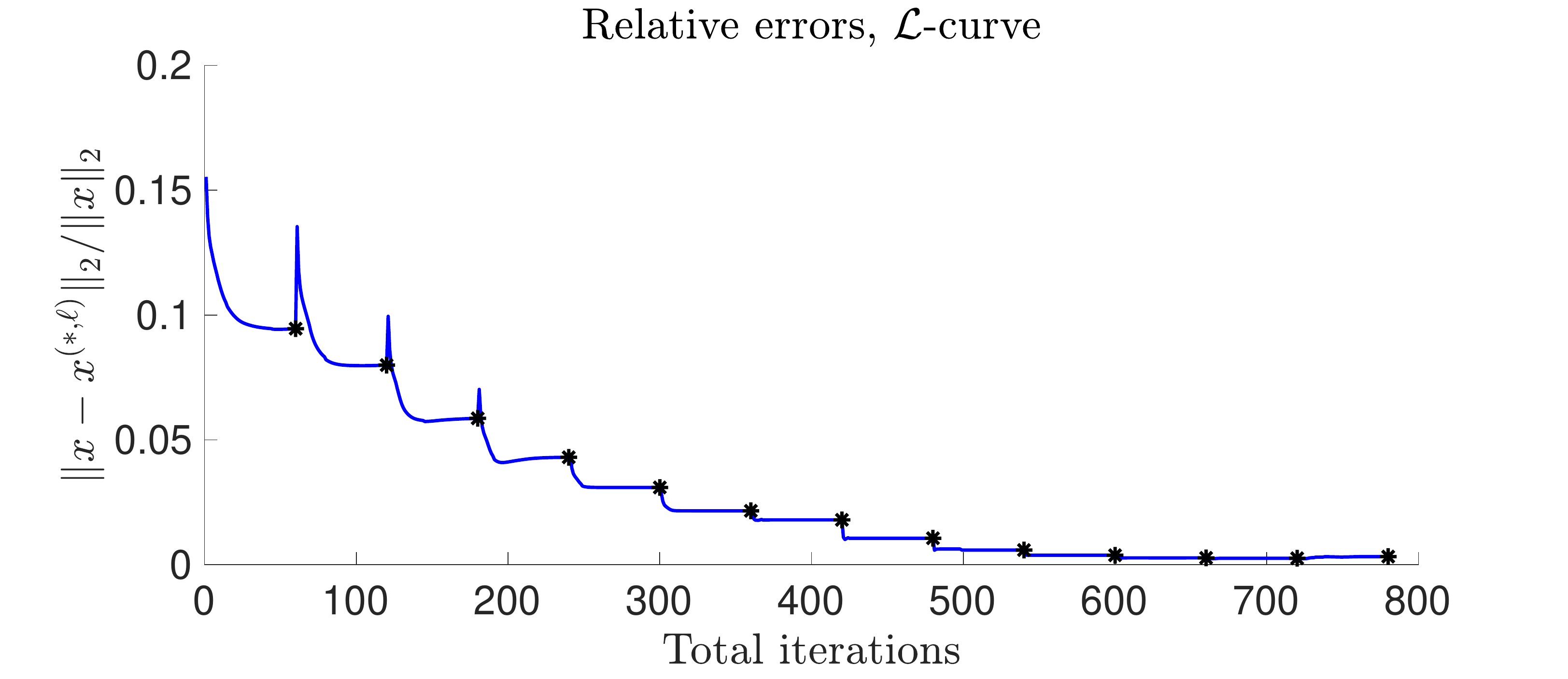}
\end{tabular}
\end{center}
\caption{\emph{Parallel-beam CT} test problem. Relative errors versus number of total iterations, when the discrepancy principle (left frame) and the ${\mathcal L}$-curve (right frame) are used to select the regularization parameter at each inner iteration. Black asterisks mark the last iteration of each inner cycle.}
\label{fig:Radon1total}
\end{figure}

Since our new algorithm is an inner-outer iterative strategy, and since so far only the behavior of the solution and regularization parameter across the outer iterations has been displayed, Figure \ref{fig:Radon1total} displays the relative errors against the number of total iterations. The end of each inner iteration cycle is marked by an asterisk. We can see that each inner iteration cycle consists of 60 iterations, i.e., the maximum allowed number of inner iterations: despite the inner stopping criterion for the inner iterations not being very effective for this example, we can observe that the quality of the reconstruction stabilizes and is almost optimal after some inner iterations (although it seems to slightly deteriorate for the discrepancy principle): this means that both adaptive parameter choice strategies are quite effective. 

\begin{figure}[htbp]
\begin{center}
\begin{tabular}{cc}
\includegraphics[width=6cm]{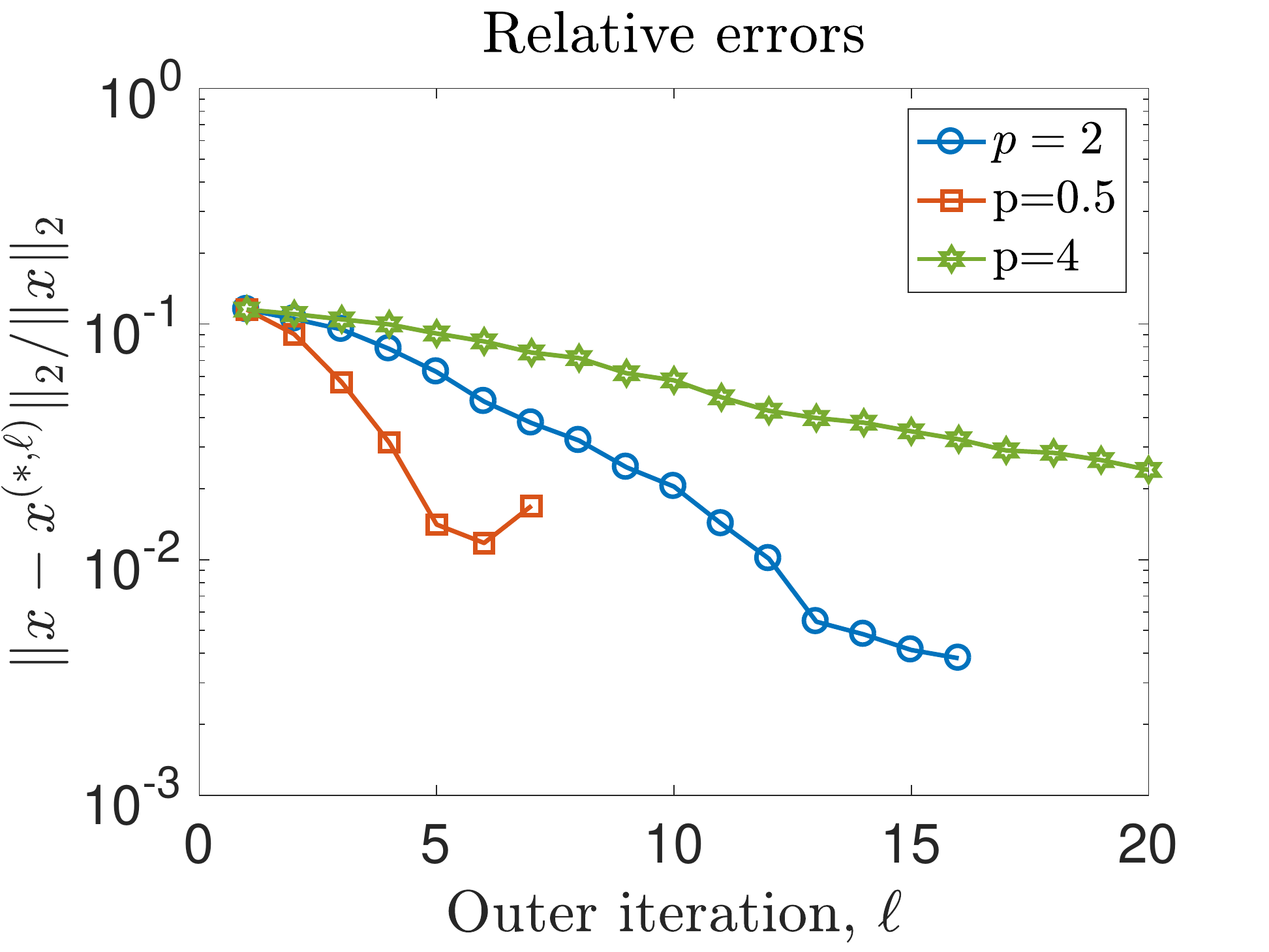} &
\includegraphics[width=6cm]{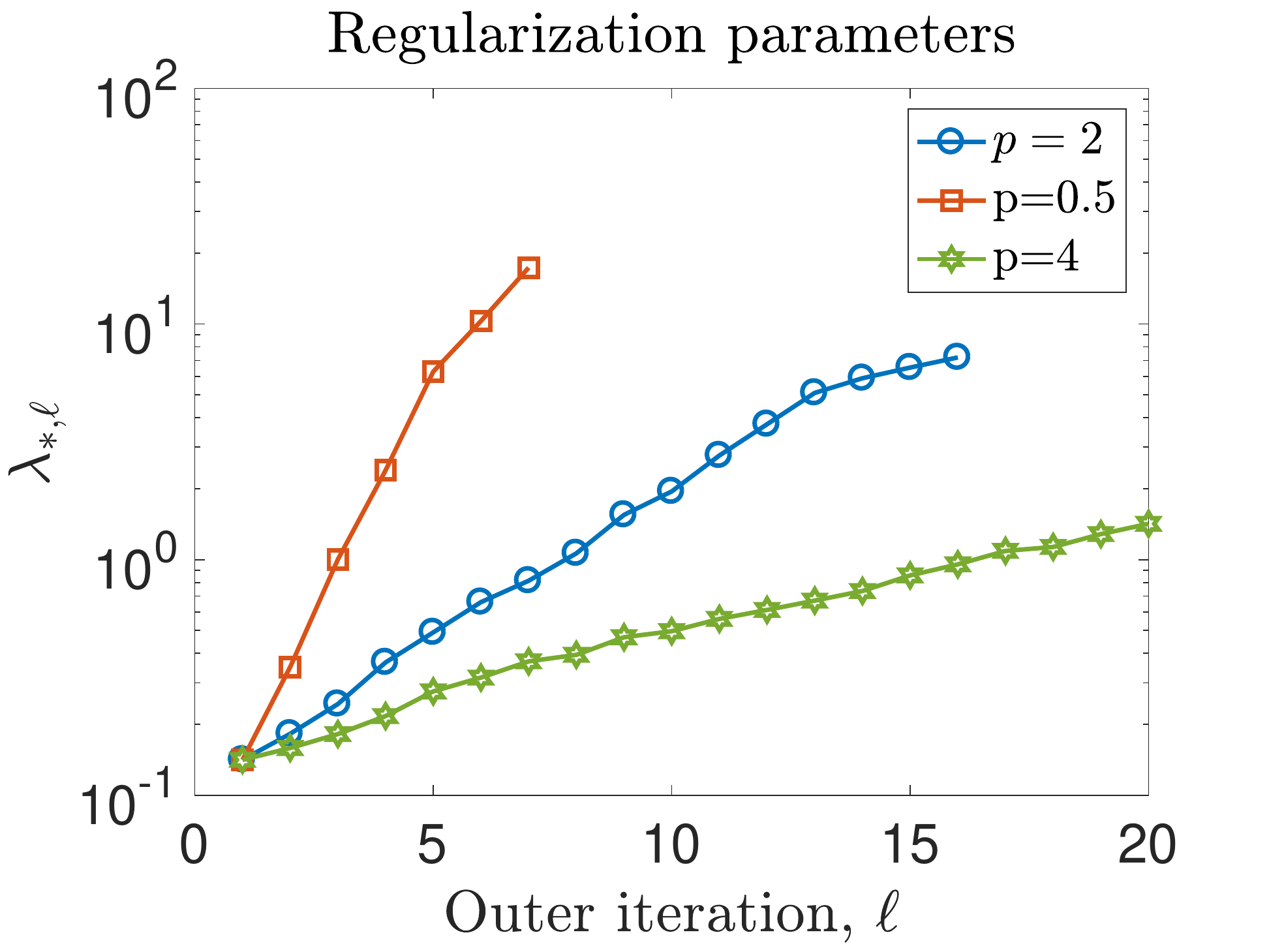}
\end{tabular}
\end{center}
\caption{\emph{Parallel-beam CT} test problem. Relative errors and regularization parameter versus number of outer iterations for three different choices of the power $p$ in (\ref{eq:GradientMap}).}
\label{fig:Radon1p}
\end{figure}

\begin{figure}[htbp]
\begin{center}
\begin{tabular}{cccc}
\small{exact} & \small{$p=2$} & \small{$p=0.5$} & \small{$p=4$}\\ 
\includegraphics[width=4cm]{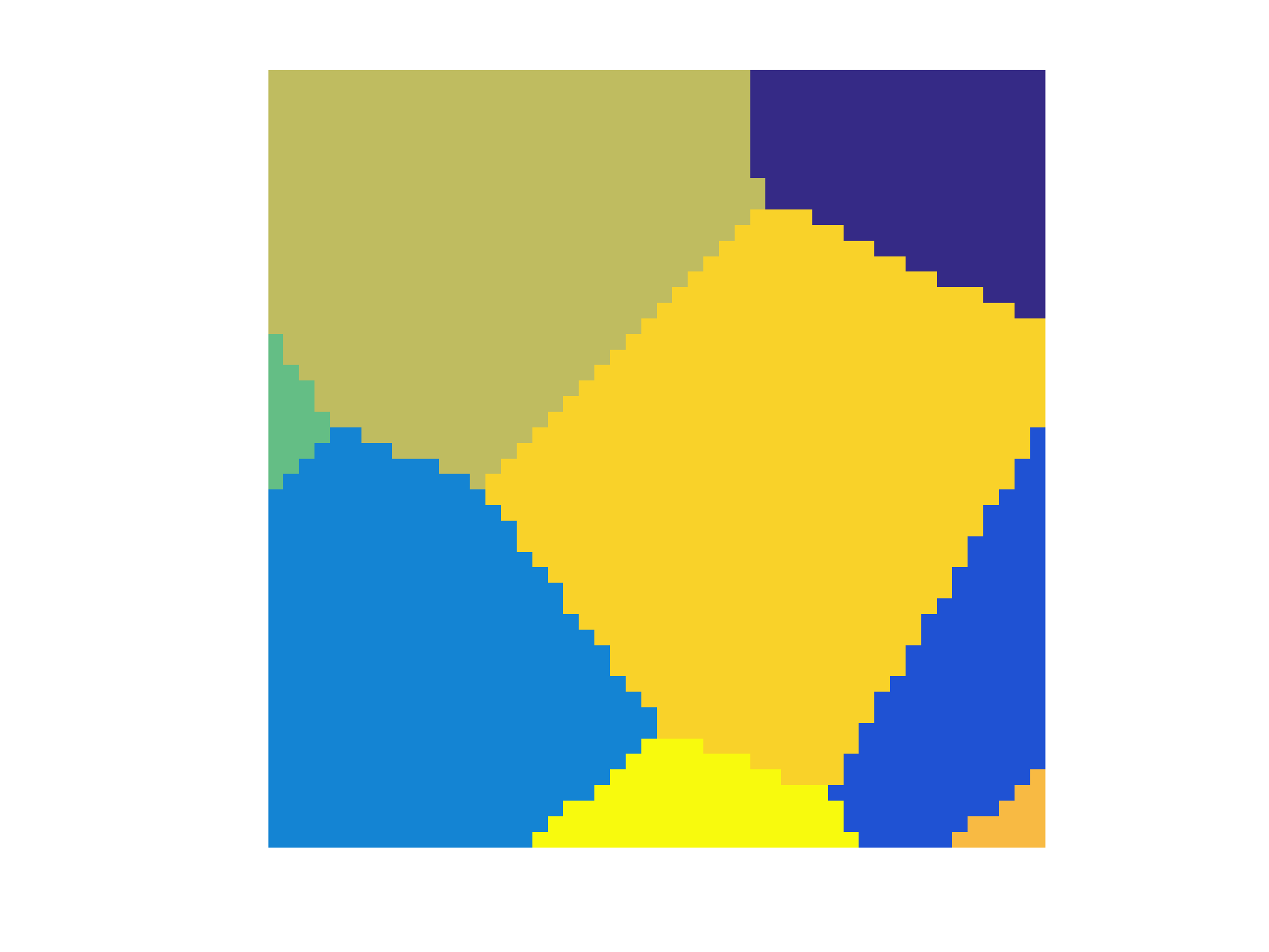} &
\includegraphics[width=4cm]{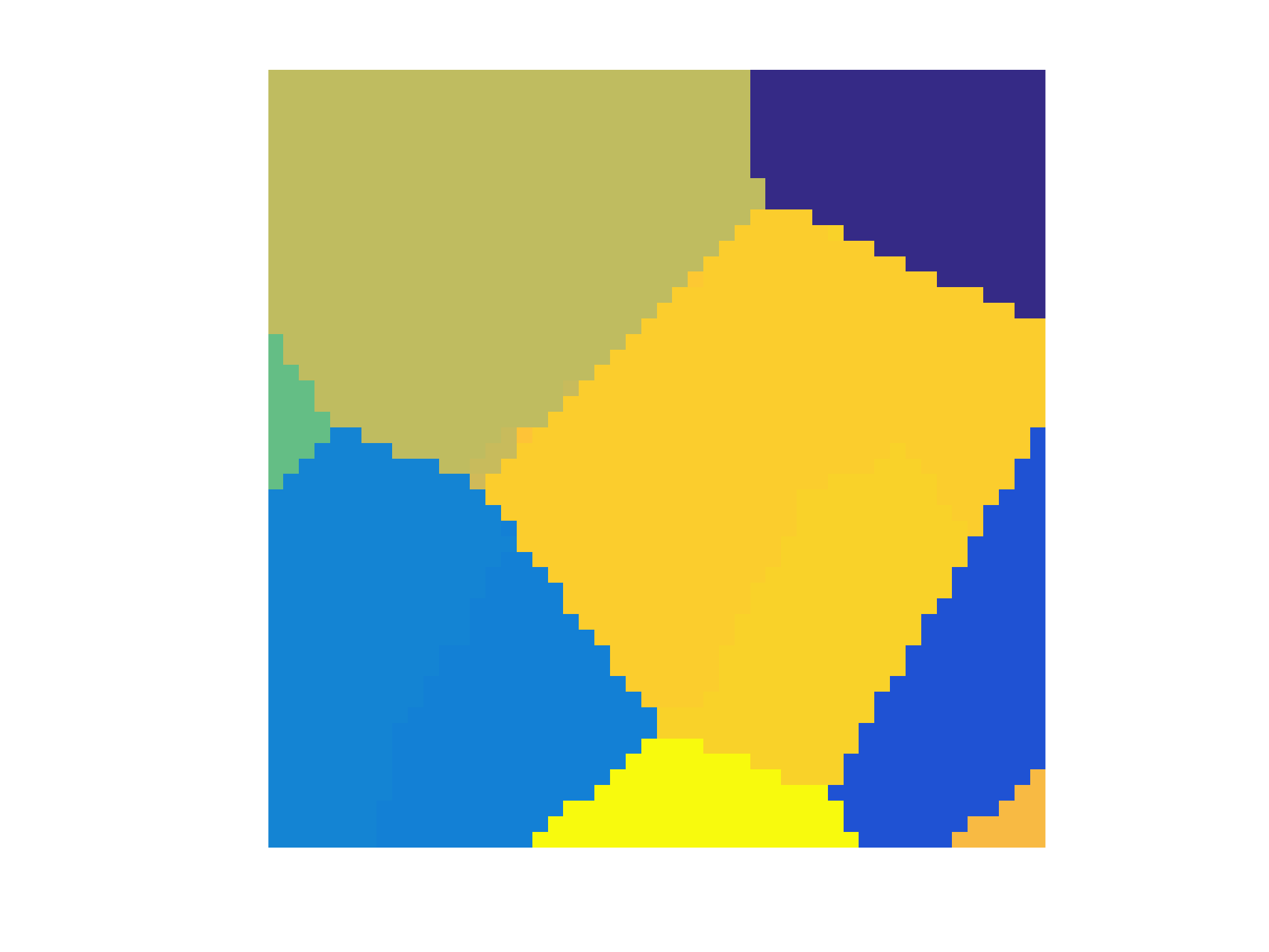} & 
\includegraphics[width=4cm]{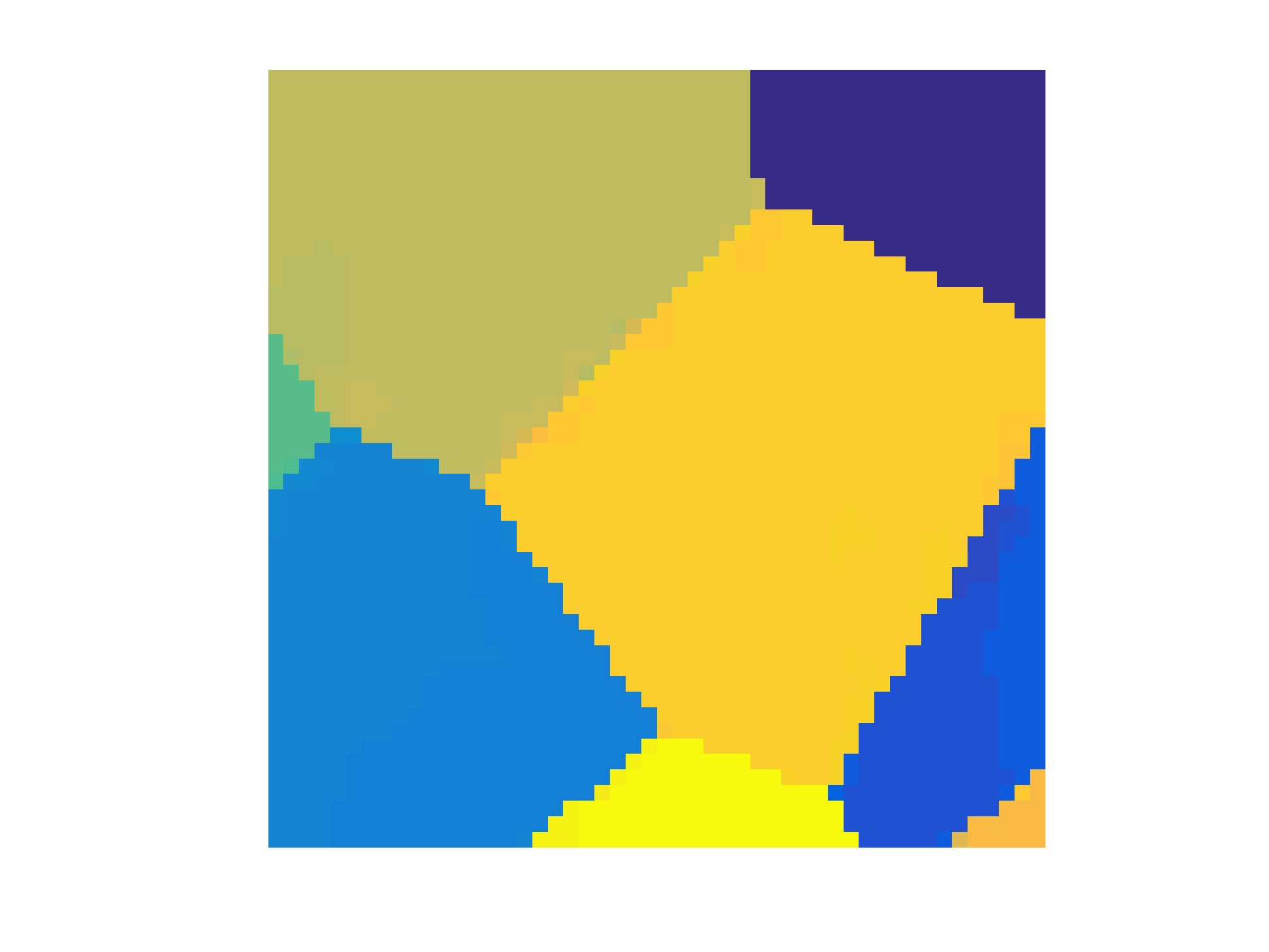} & 
\includegraphics[width=4cm]{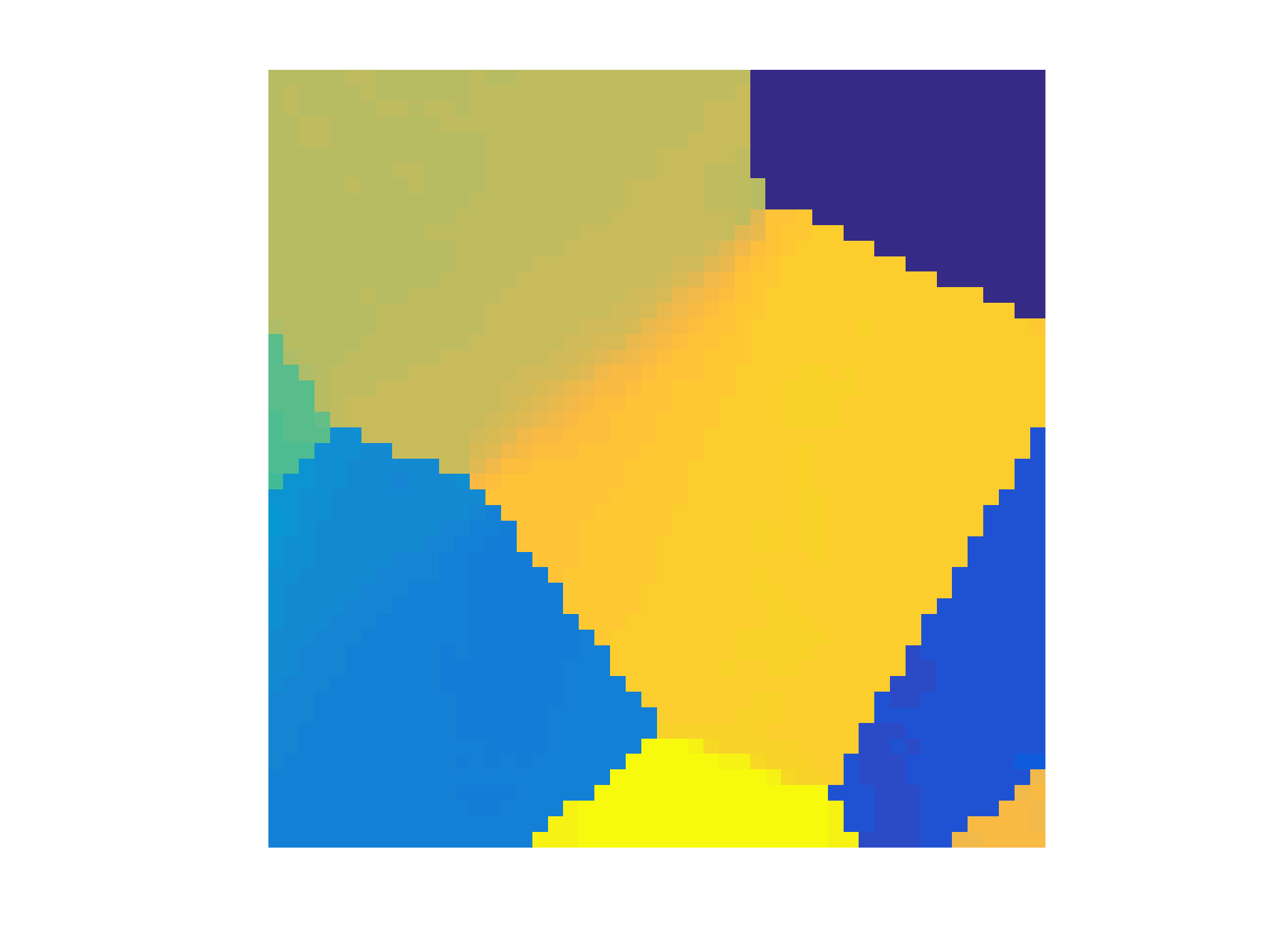} \\
\includegraphics[width=4cm]{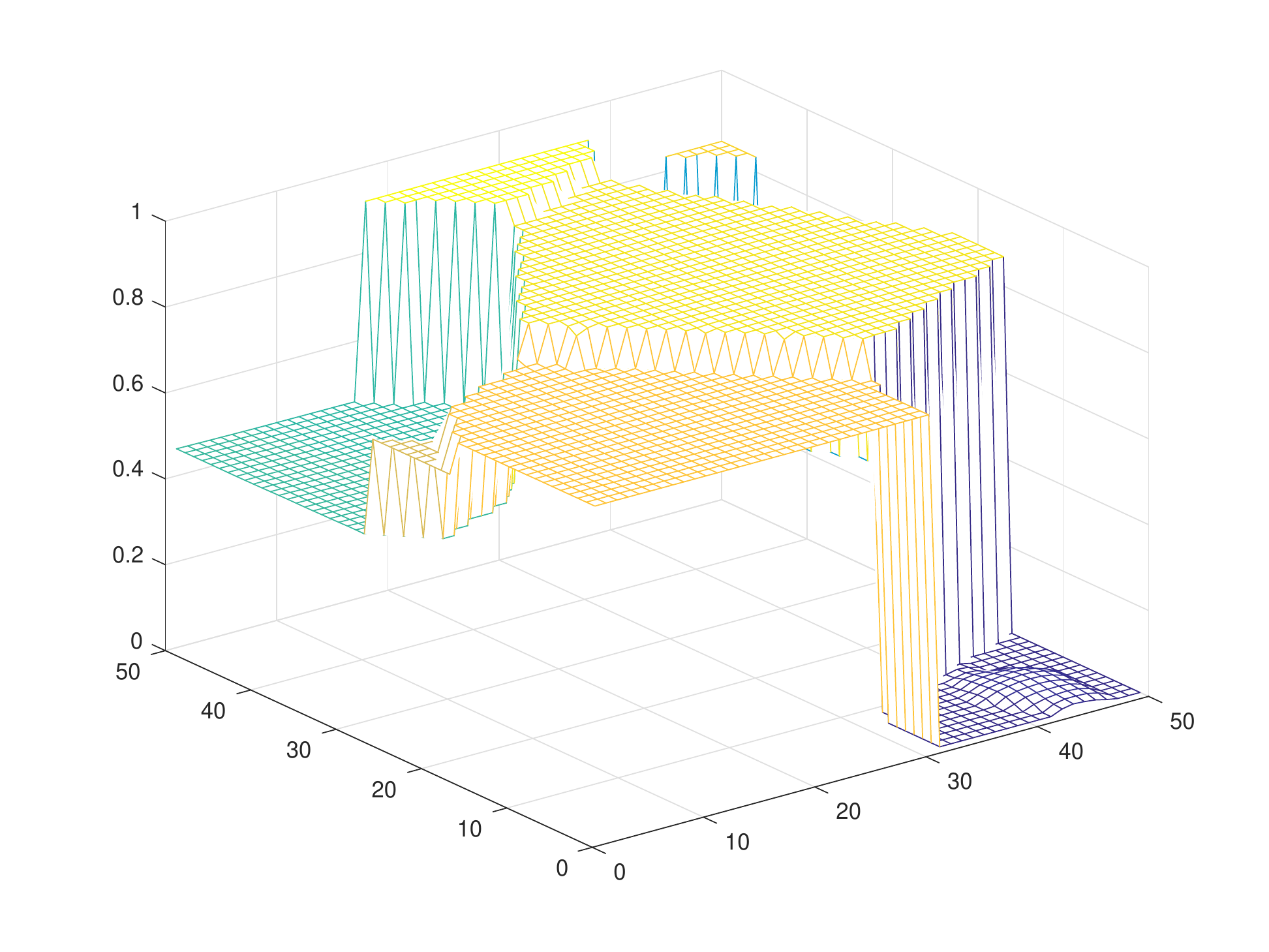} &
\includegraphics[width=4cm]{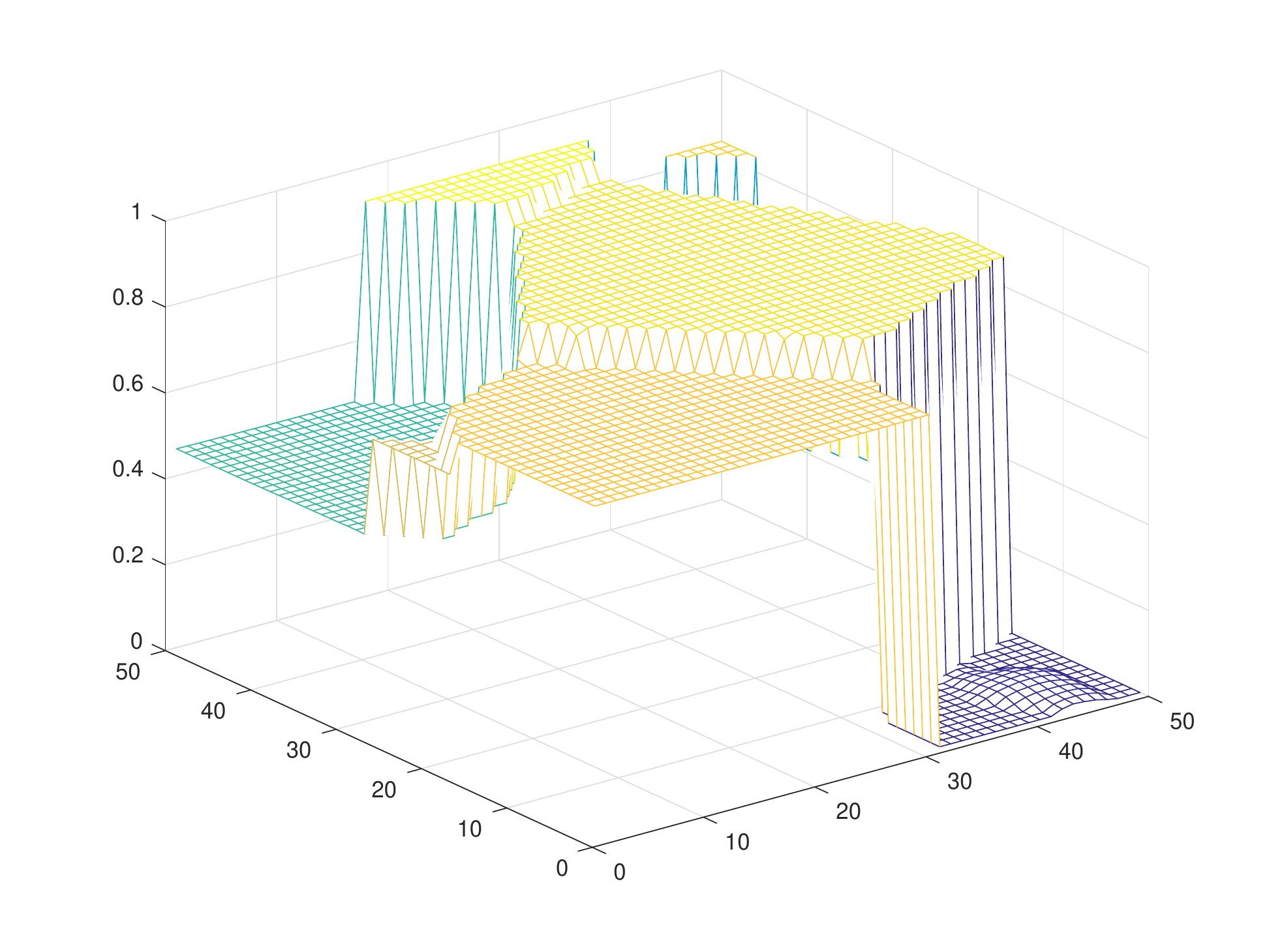} & 
\includegraphics[width=4cm]{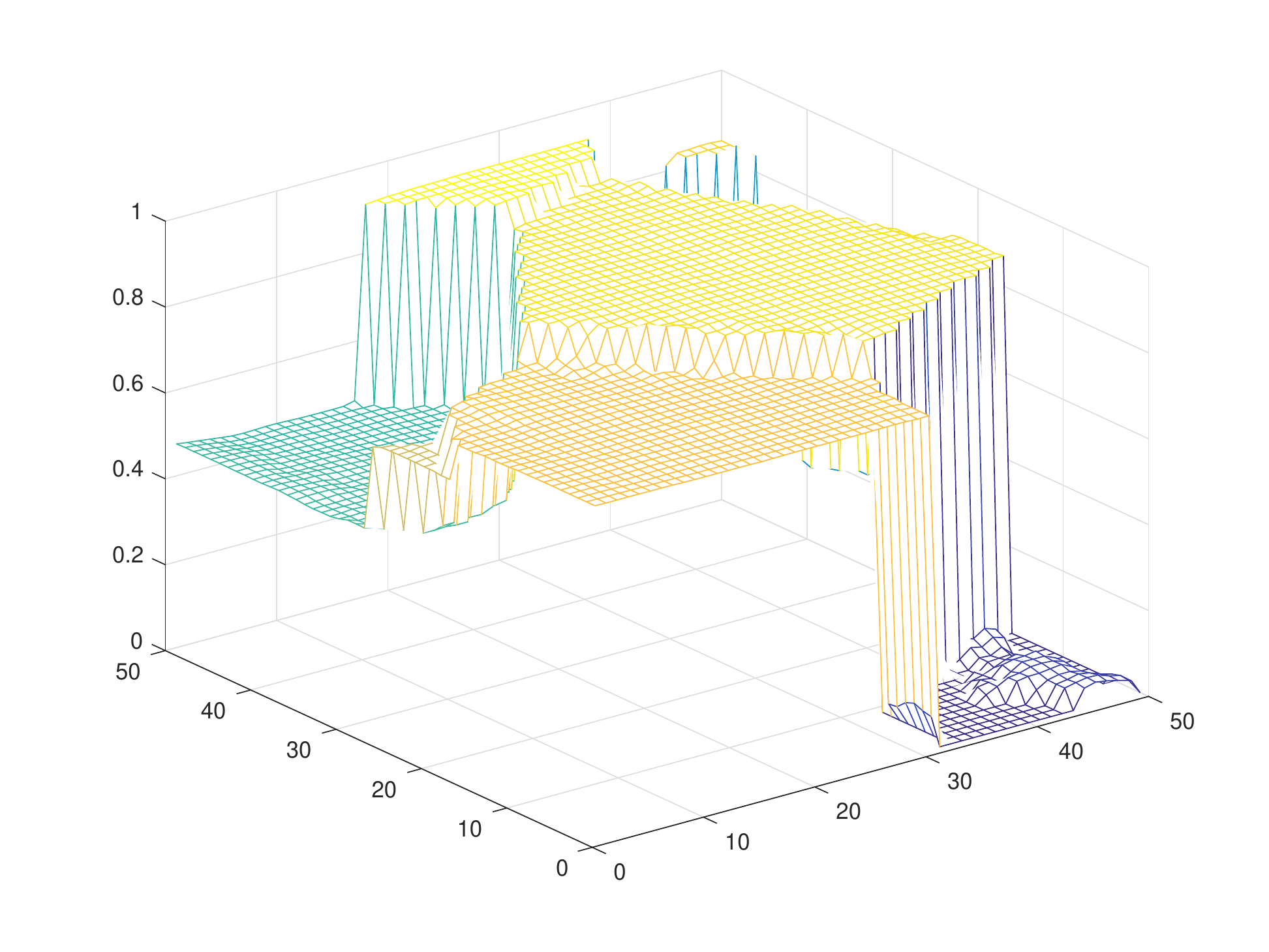} & 
\includegraphics[width=4cm]{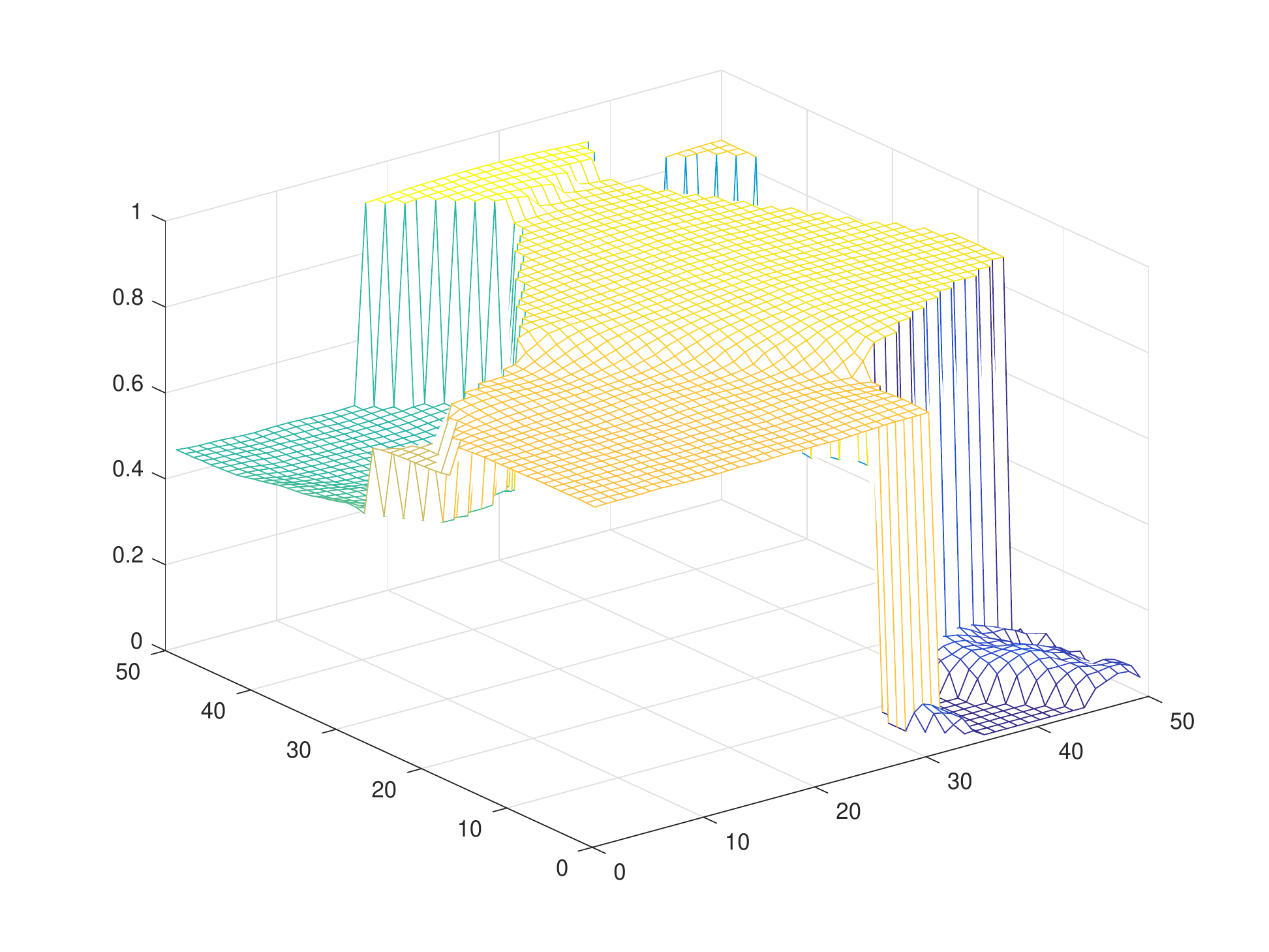} \\
\end{tabular}
\end{center}
\caption{\emph{Parallel-beam CT} test problem. Zoom-ins of the reconstructions obtained when the outer iterations terminate, also displayed as surface plots.}
\label{fig:Radon1zoom}
\end{figure}

The definition of the weights at each outer iteration depends on a power $p>0$ (see, for instance, equation (\ref{eq:GradientMap})). In Section \ref{ssec:weights} we made the argument that the choices $p\ll 1$ and $p\gg 1$ correspond to less and more smooth reconstructions, respectively. We now experimentally assess how the value of $p$ affects the quality of the reconstructions when the discrepancy principle is employed to adaptively set the regularization parameter at each inner iteration. 
%, and when the latter is adaptively selected according to the discrepancy principle. 
Figure \ref{fig:Radon1p} shows the values of the relative errors and the regularization parameter versus the number of outer iterations for three values of $p$; $p=2$ is the value selected to display the previous graphs.
Zoom-ins of the reconstructed images are displayed in Figure \ref{fig:Radon1zoom}, also as surfaces. Although all the considered values of $p$ deliver reconstructions of excellent quality, we can see some slight spurious oscillations in the supposedly constant patches reconstructed taking $p=0.5$, and we can clearly see that some of the edges reconstructed using $p=4$ are washed out.
Interestingly enough, still comparing the choices $p=0.5$ and $p=4$ in Figure \ref{fig:Radon1p}, we can see that, since $p=0.5$ intrinsically enforces less smoothness through the weights, a larger regularization parameter is automatically set to compensate for this, and less outer iterations are performed; viceversa, since $p=4$ intrinsically enforces more smoothness through the weights, a smaller regularization parameter is automatically set to compensate for it, and more outer iterations are performed. To generate the following graphs and when considering all the following test problems, we will always take $p=2$.
%In the following tests and examples we will always choose  for it and the weights already enforce a large amount of smoothness, see``intrinsically'' results in an overall smoother reconstruction, so that a smaller regularization parameter will do; the choice $p\ll 1$ ``intrinsically'' results in an overall less smooth reconstruction, so that a bigger regularization parameter should be selected. In the following we will stick to $p=2$. 

\begin{figure}[htbp]
\begin{center}
\begin{tabular}{cc}
\includegraphics[width=6cm]{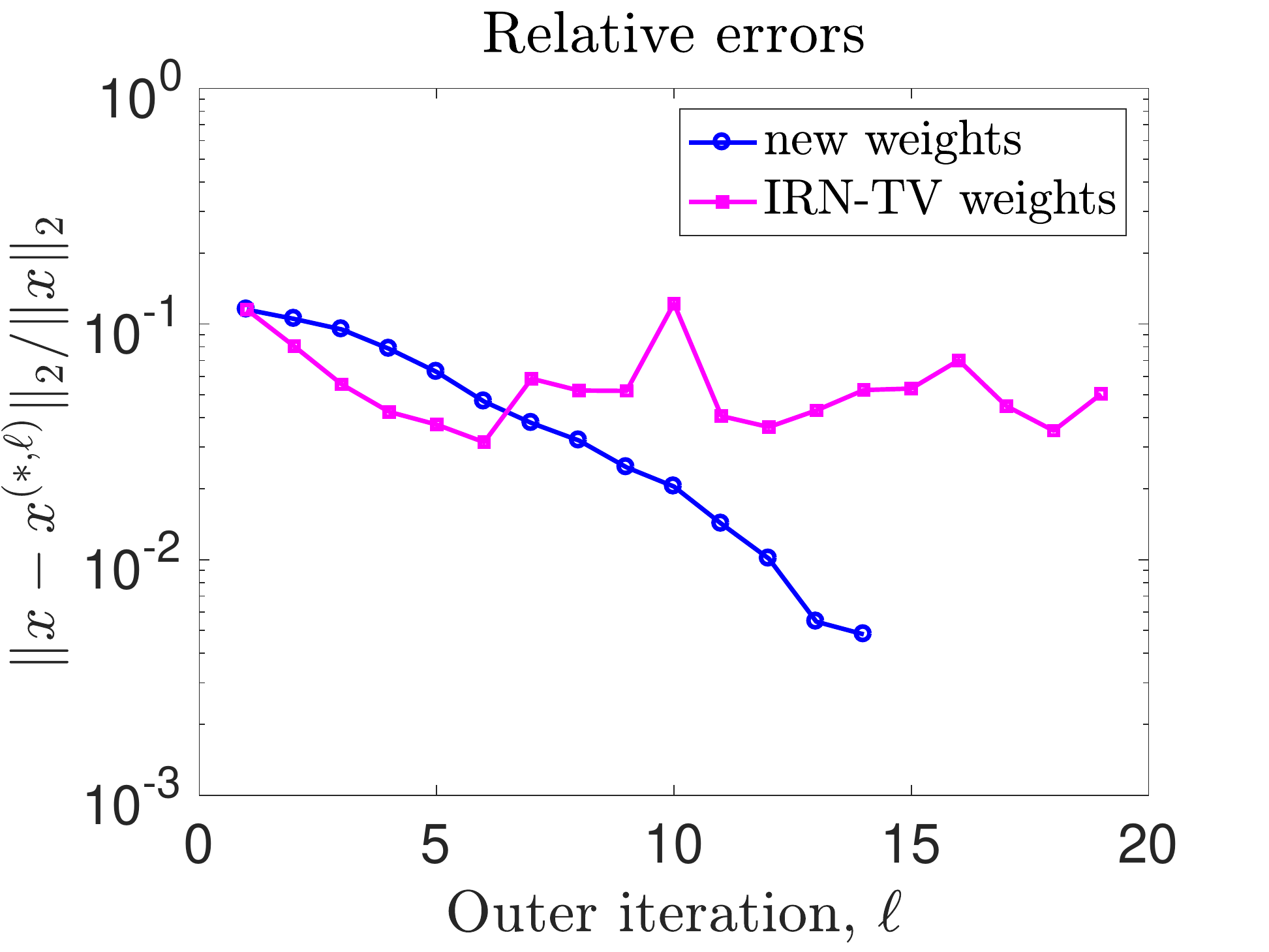} &
\includegraphics[width=6cm]{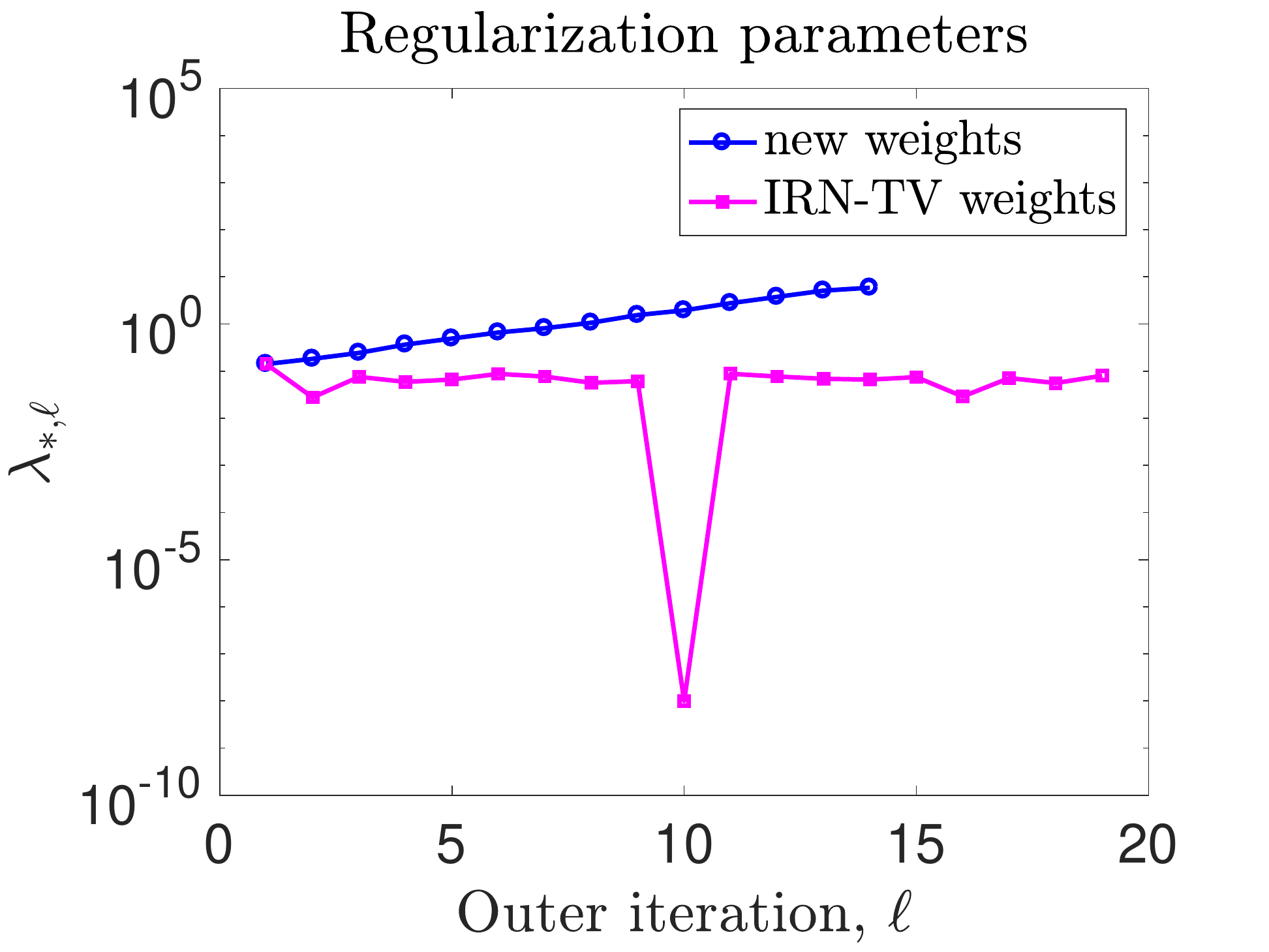}
\end{tabular}
\end{center}
\caption{\emph{Parallel-beam CT} test problem. Relative errors and regularization parameter versus number of (outer) iterations; both methods use the hybrid method for general form Tikhonov during the inner iterations, and adaptively select the regularization parameter according to the discrepancy principle.}
\label{fig:Radon1TV}
\end{figure}

\begin{figure}[htbp]
\begin{center}
\begin{tabular}{ccc}
\hspace{-0.5cm}\includegraphics[width=5.5cm]{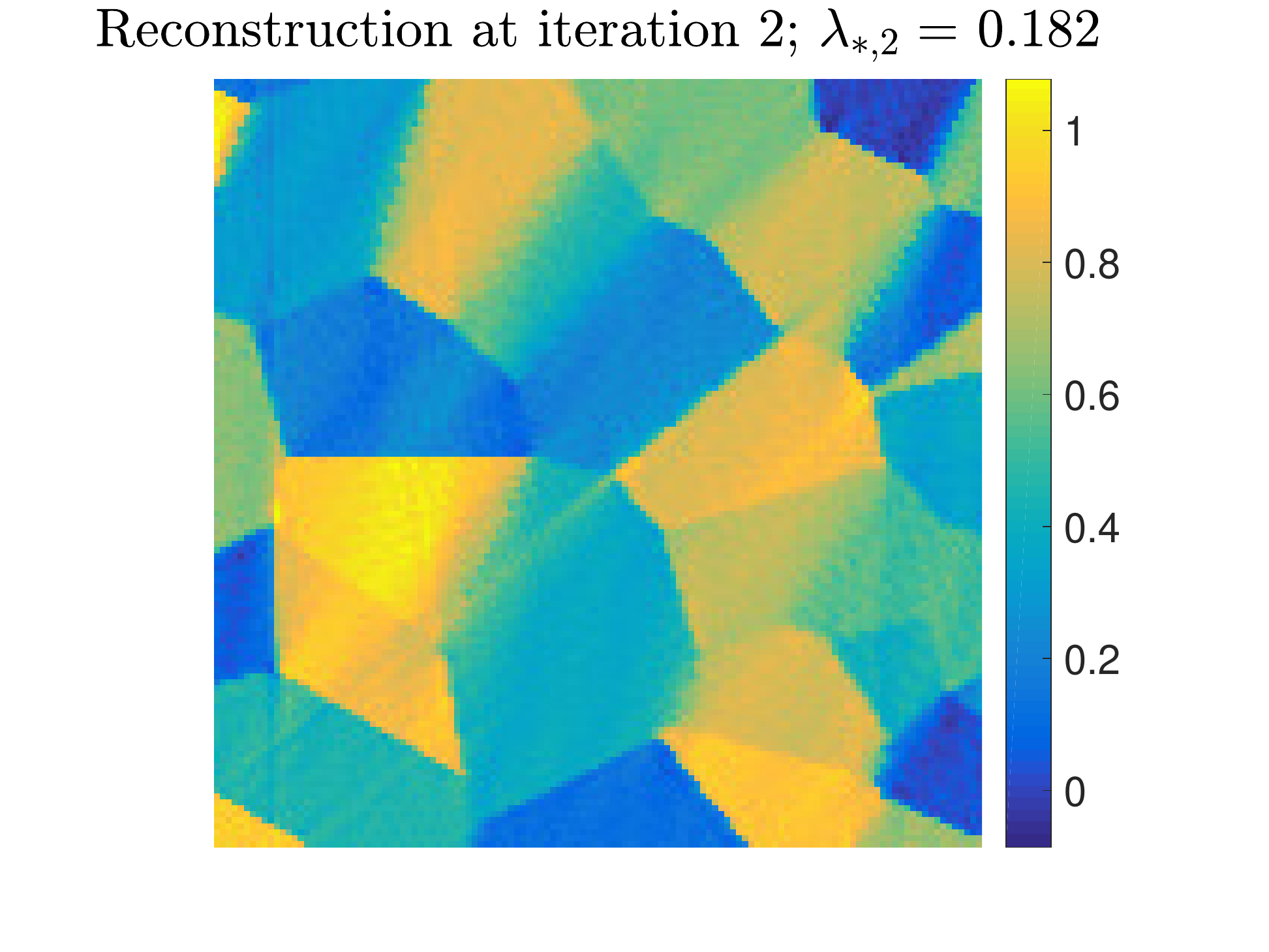} &
\hspace{-0.5cm}\includegraphics[width=5.5cm]{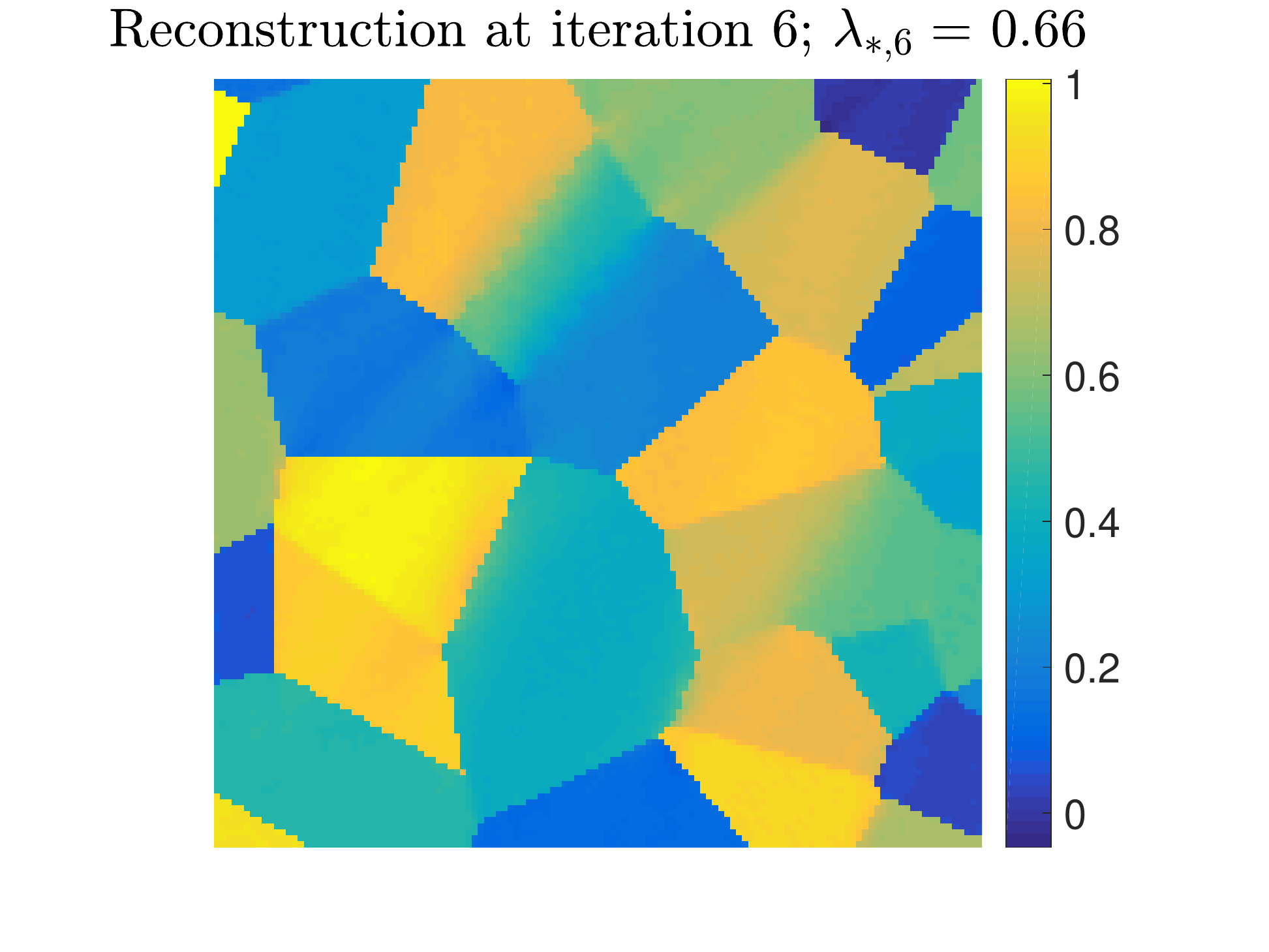} & 
\hspace{-0.5cm}\includegraphics[width=5.5cm]{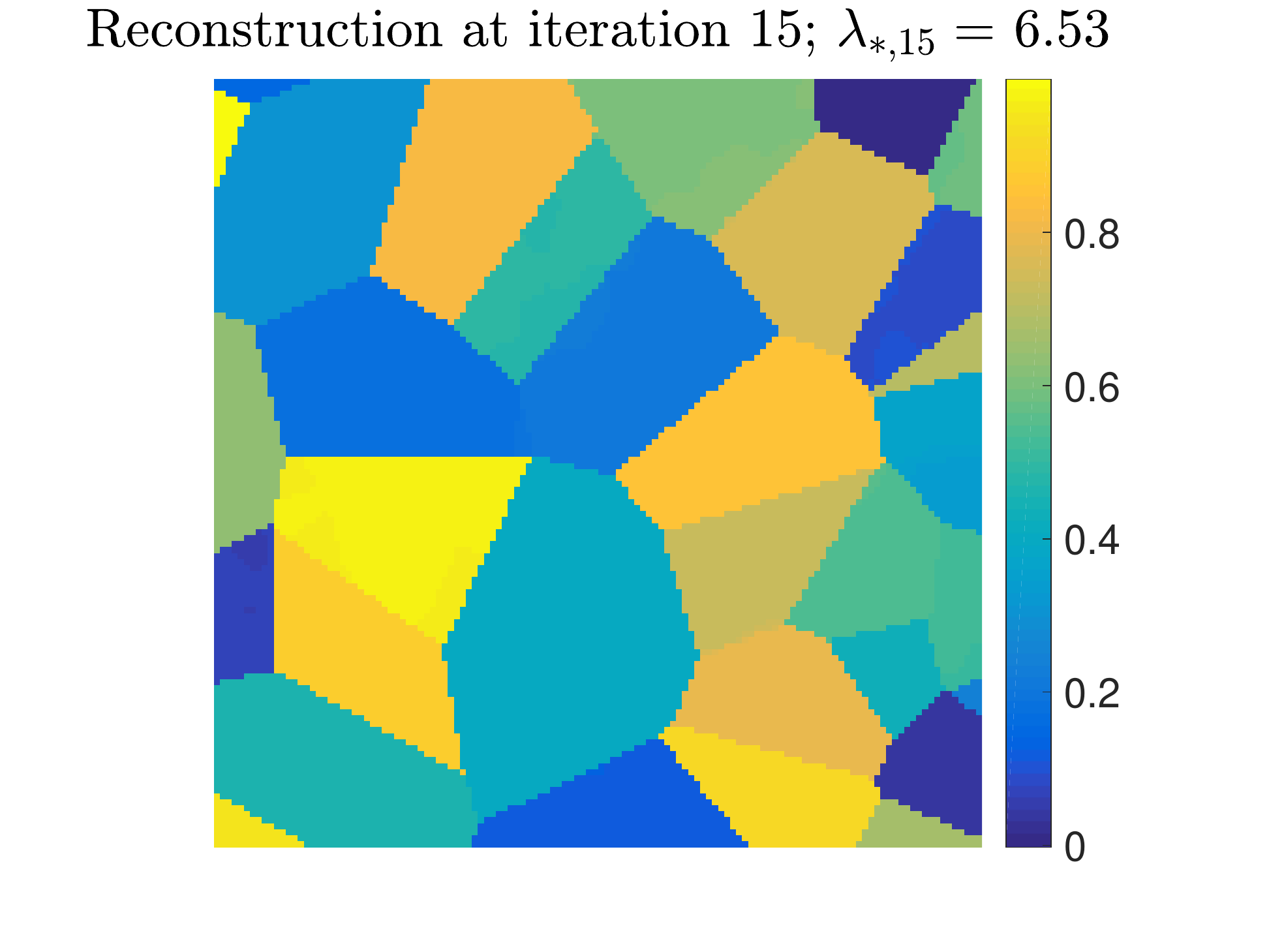}\\ 
\hspace{-0.5cm}\includegraphics[width=5.5cm]{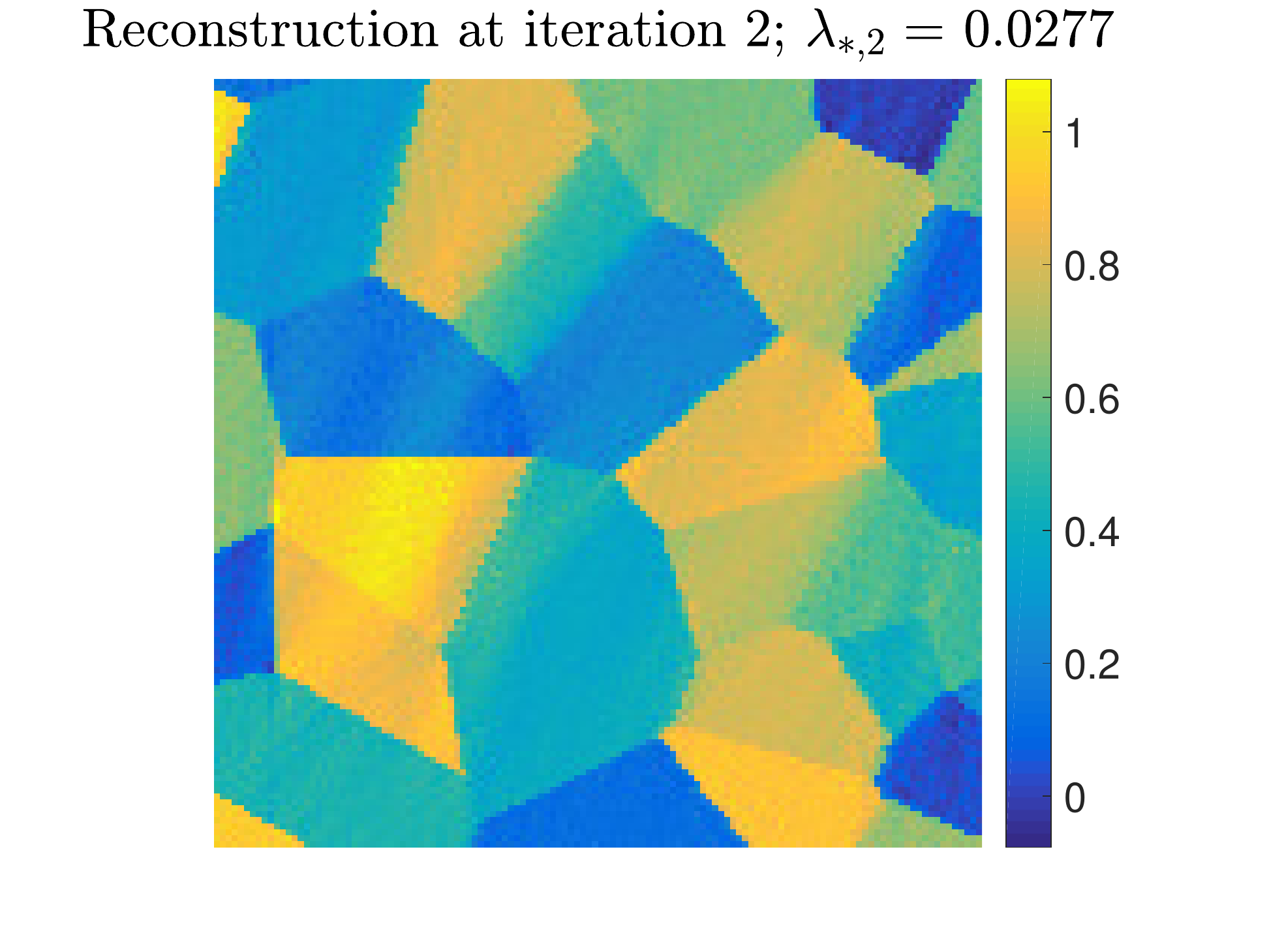} &
\hspace{-0.5cm}\includegraphics[width=5.5cm]{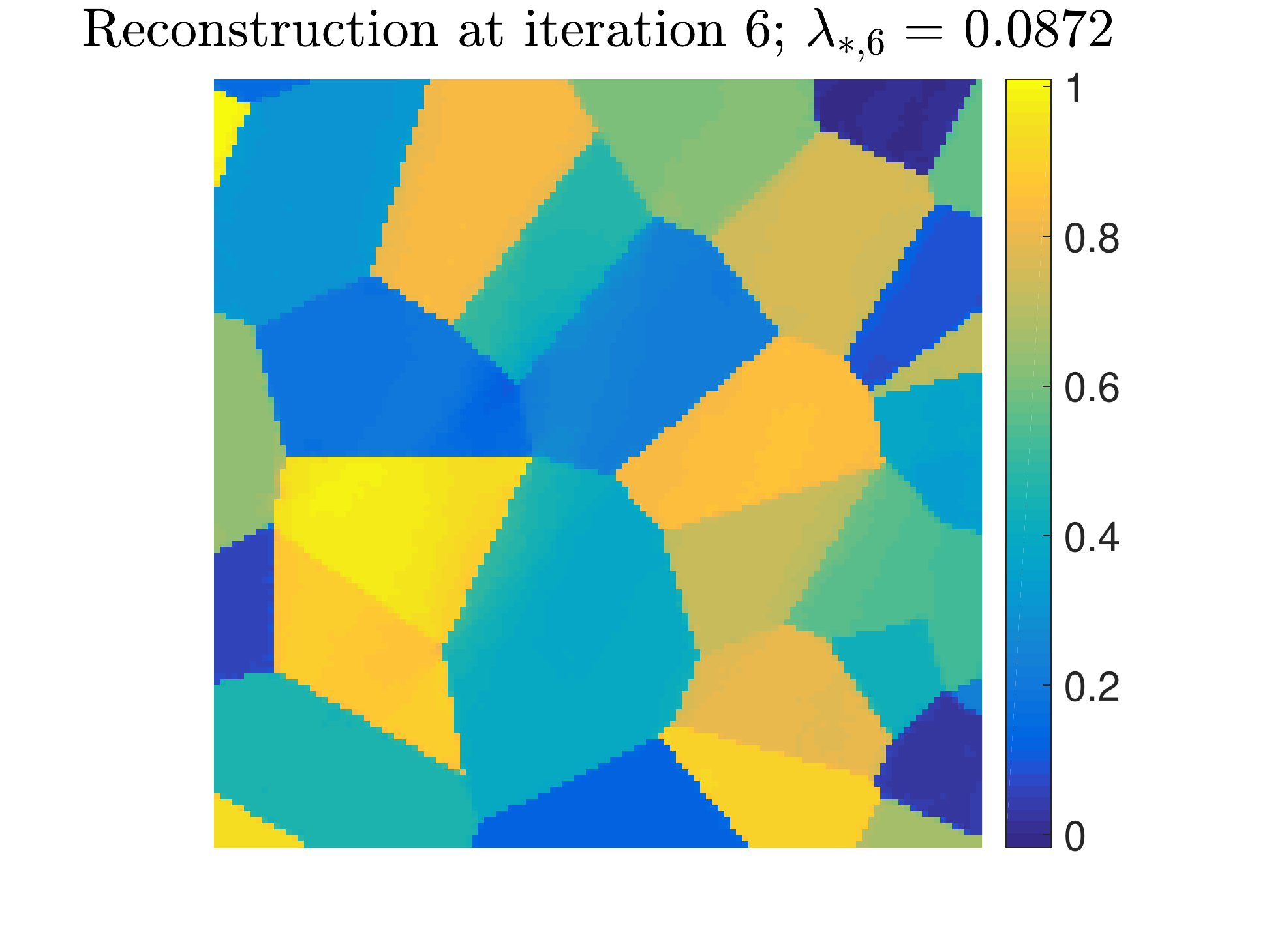} & 
\hspace{-0.5cm}\includegraphics[width=5.5cm]{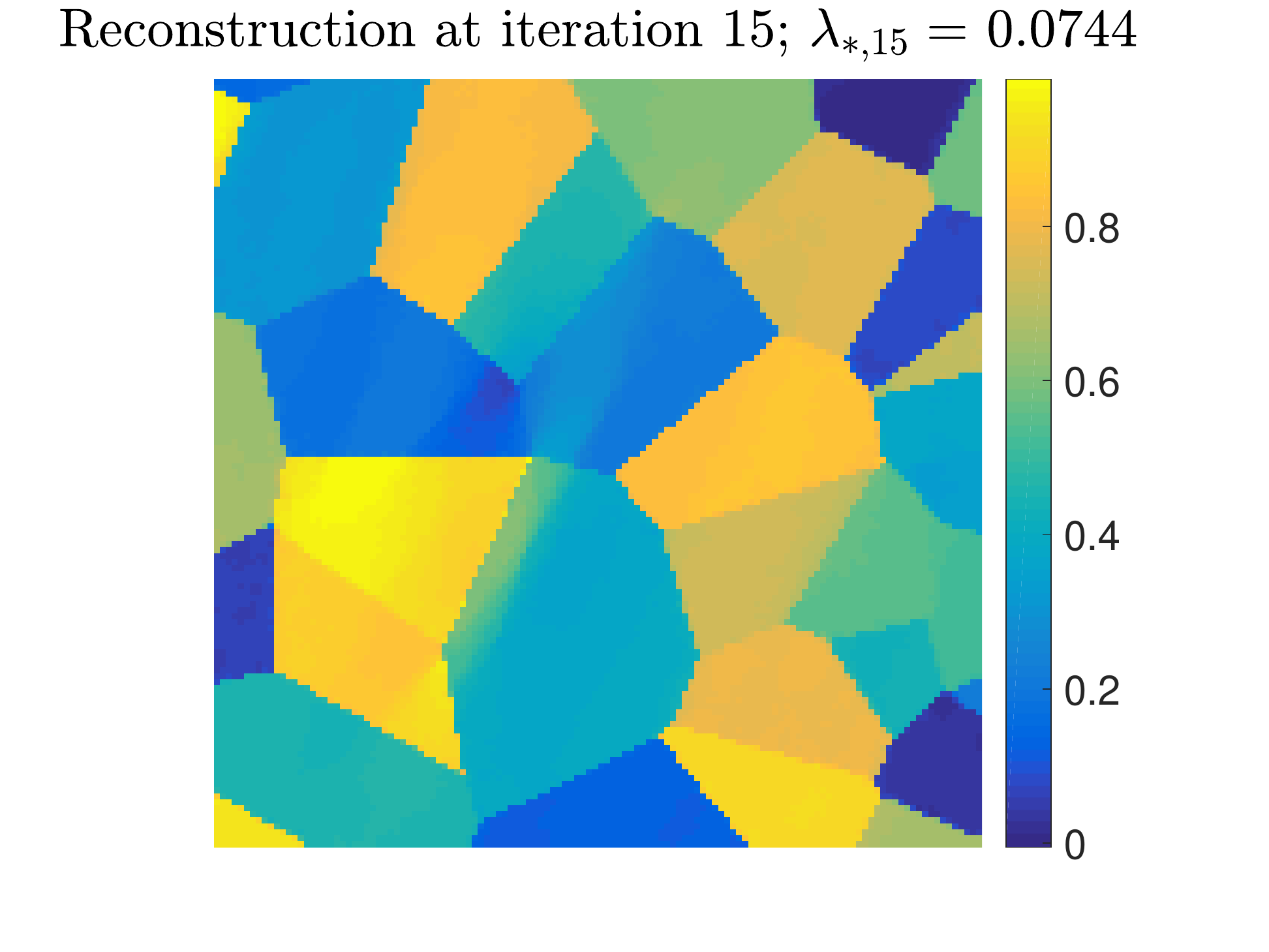}
\end{tabular}
\end{center}
\caption{\emph{Parallel-beam CT} test problem. Upper row: reconstructions obtained employing the new weights. Lower row: reconstructions obtained employing the IRN-TV weights. The reconstructions obtained at the outer iterations $\ell=2,6,15$ are displayed, reporting also the value of the regularization parameter $\lambda_{\ast,\ell}$.}
\label{fig:Radon1TVrec}
\end{figure}

\begin{figure}[htbp]
\begin{center}
\begin{tabular}{ccc}
\includegraphics[width=5cm]{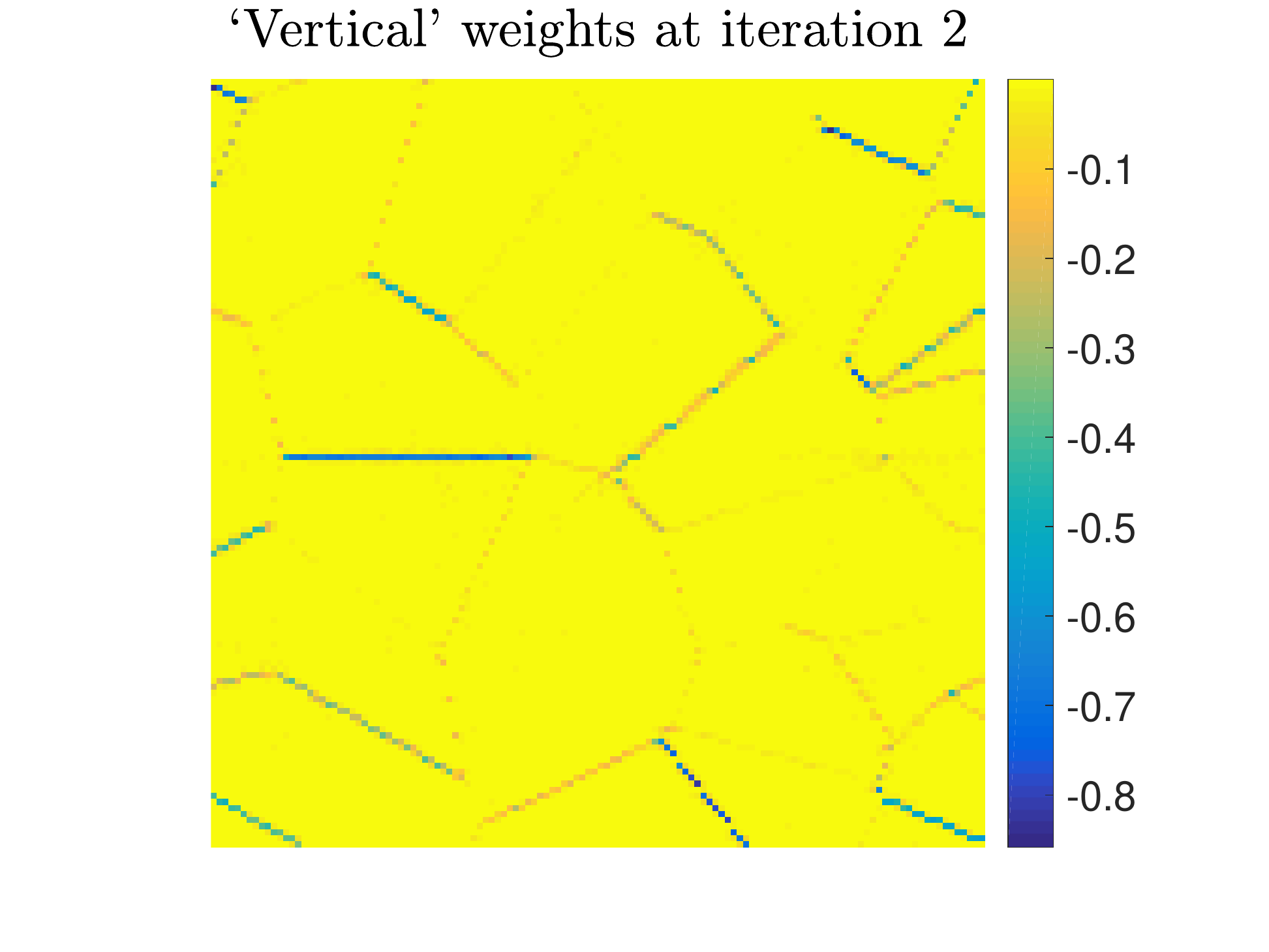} &
\includegraphics[width=5cm]{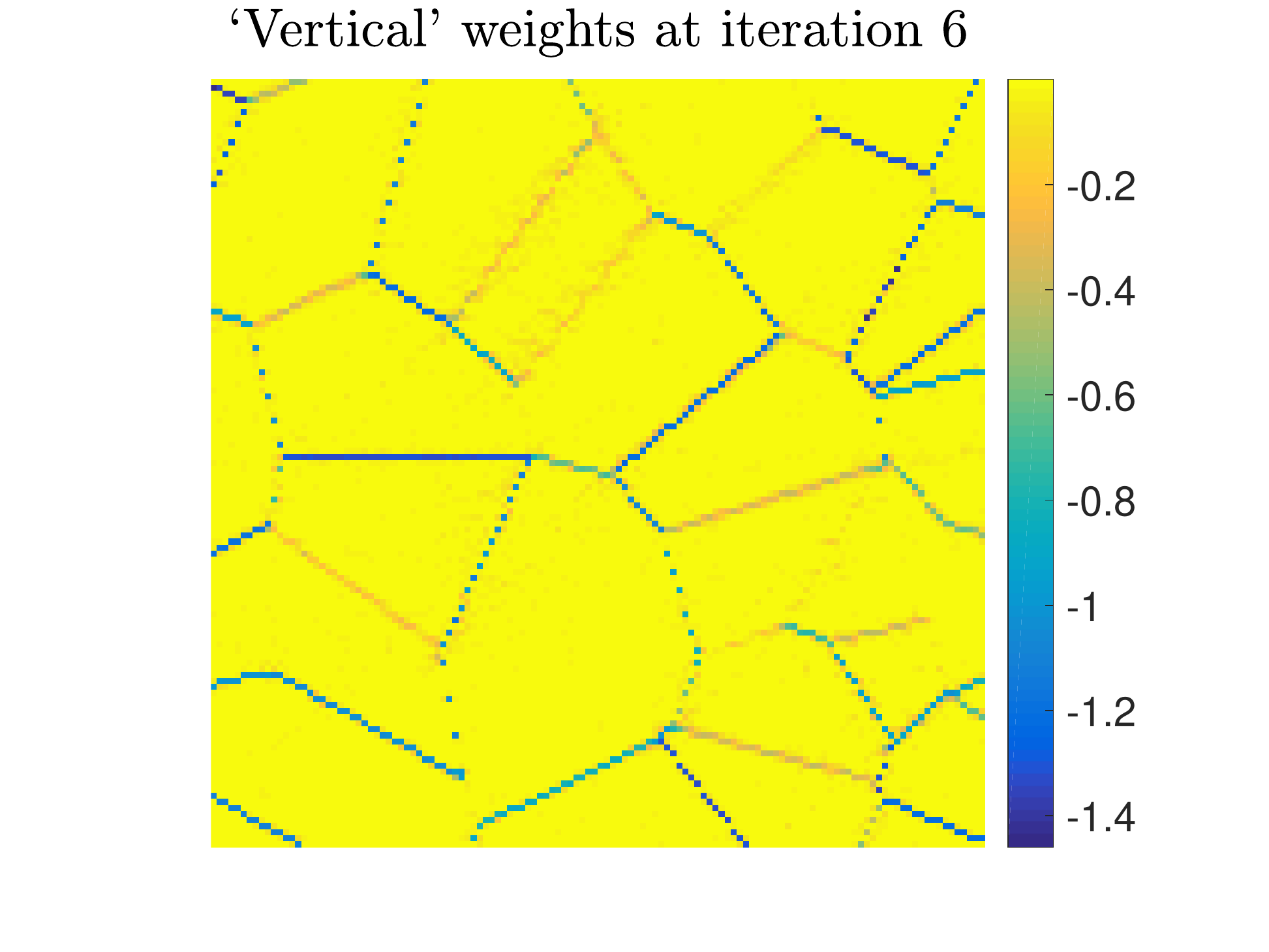} & 
\includegraphics[width=5cm]{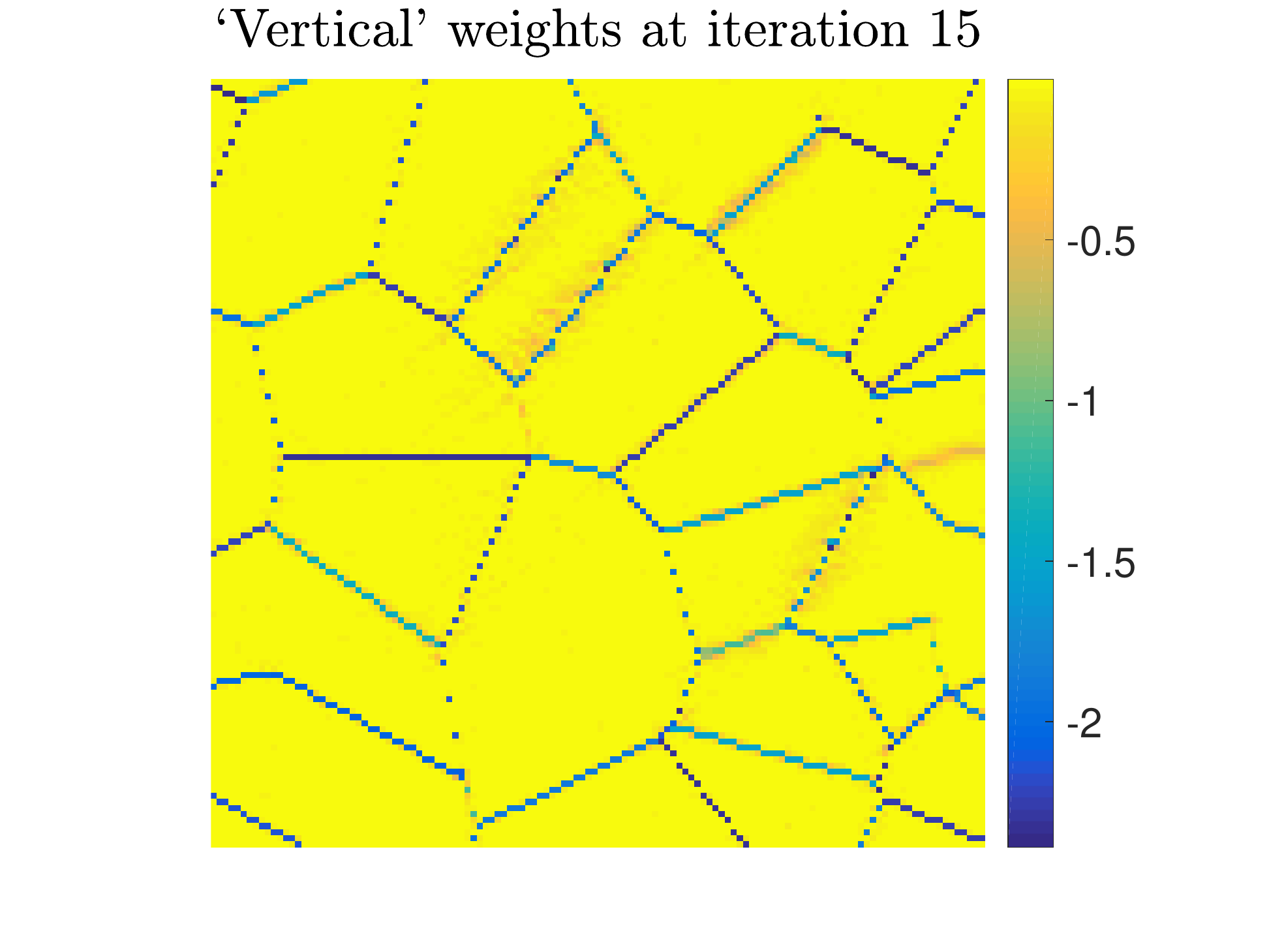}\\ 
\includegraphics[width=5cm]{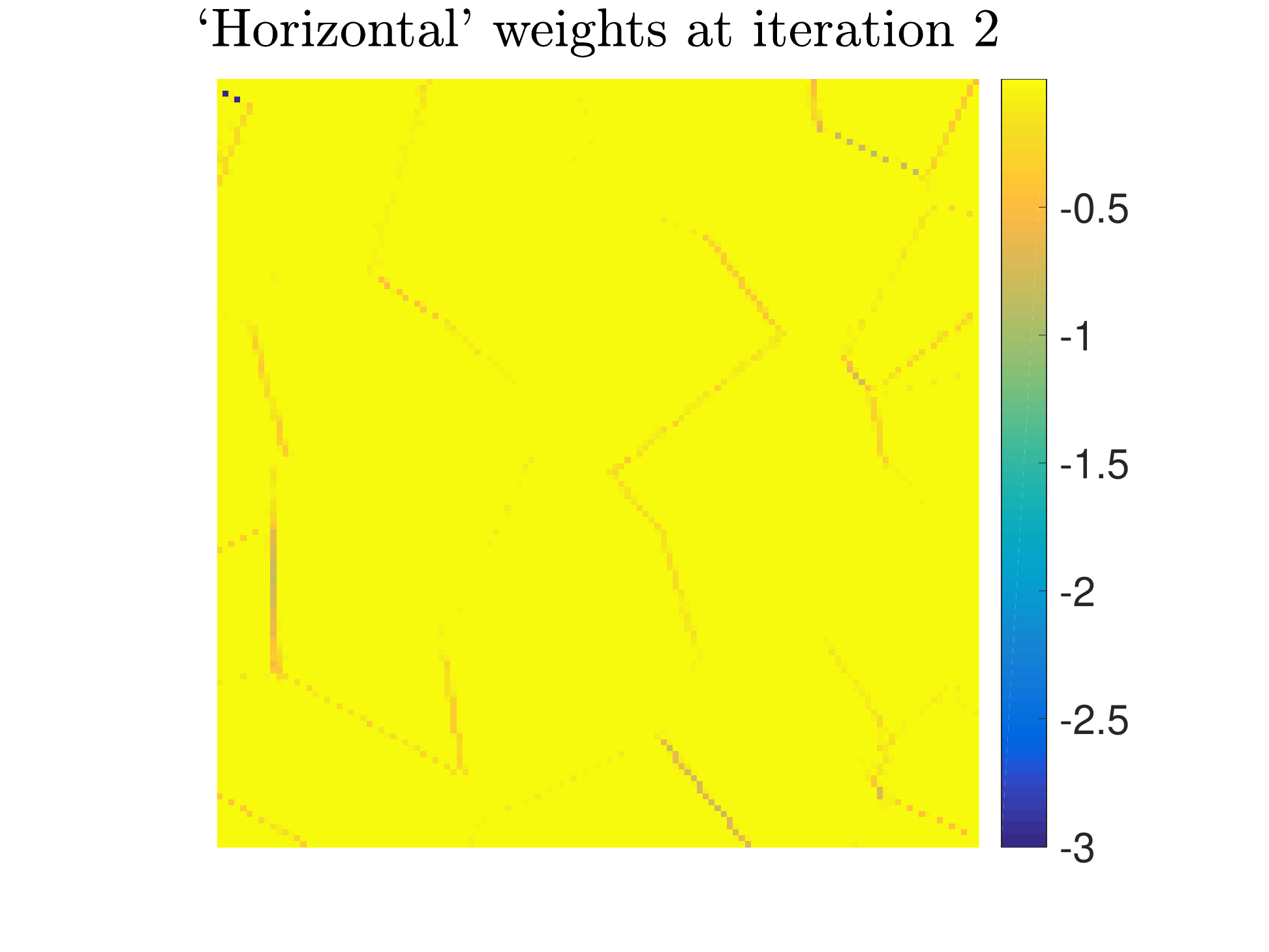} &
\includegraphics[width=5cm]{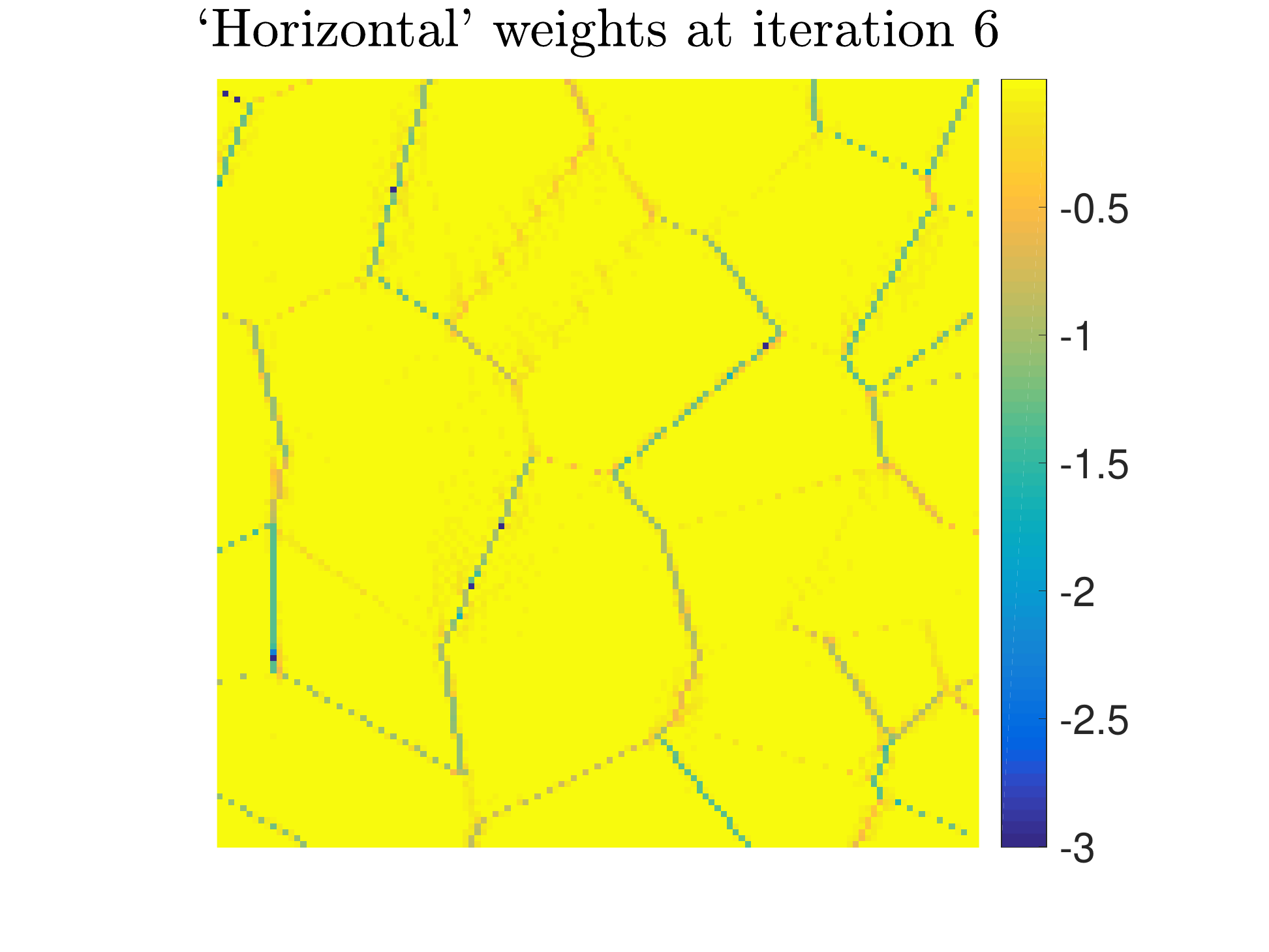} & 
\includegraphics[width=5cm]{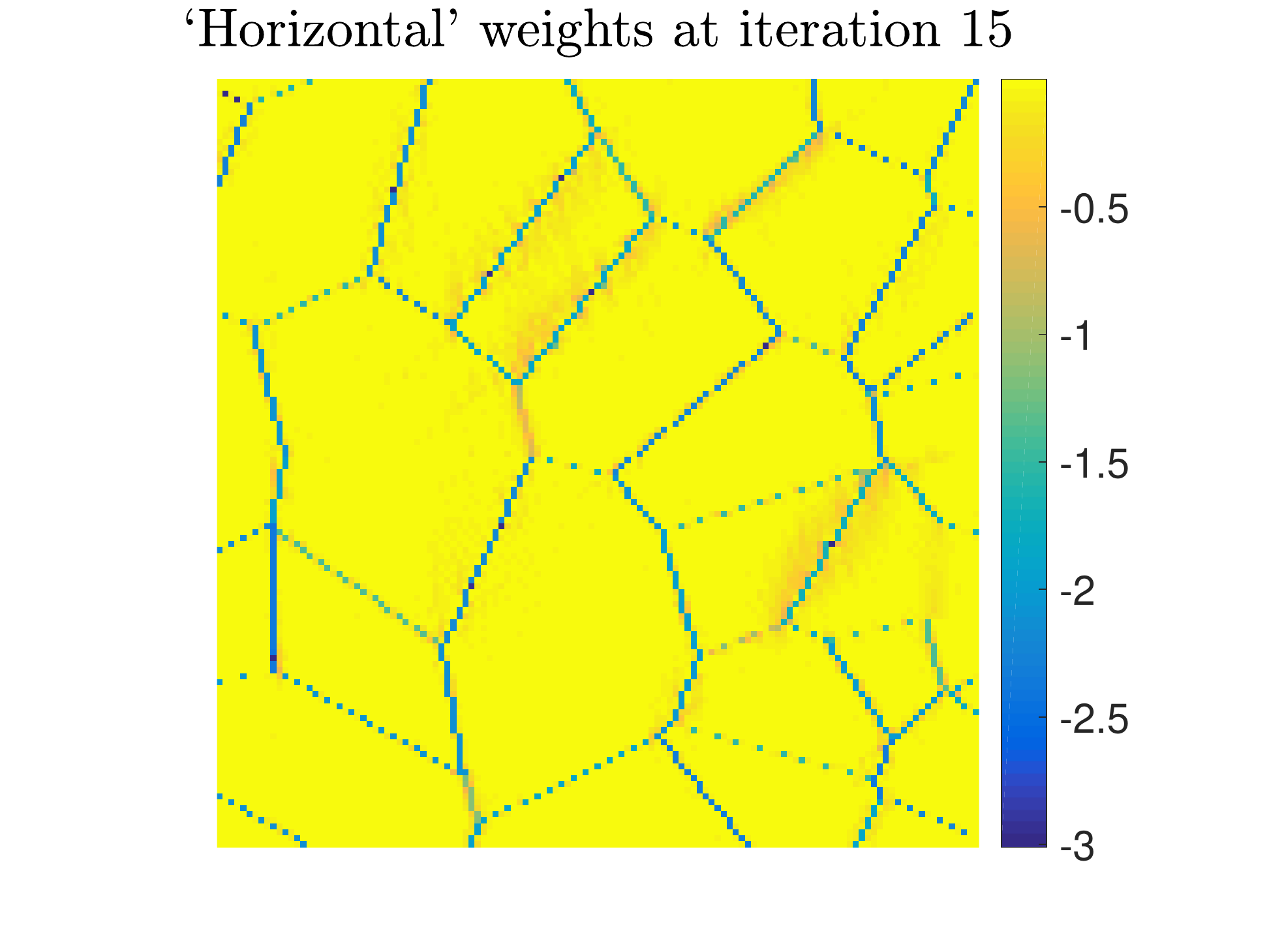}\\
\includegraphics[width=5cm]{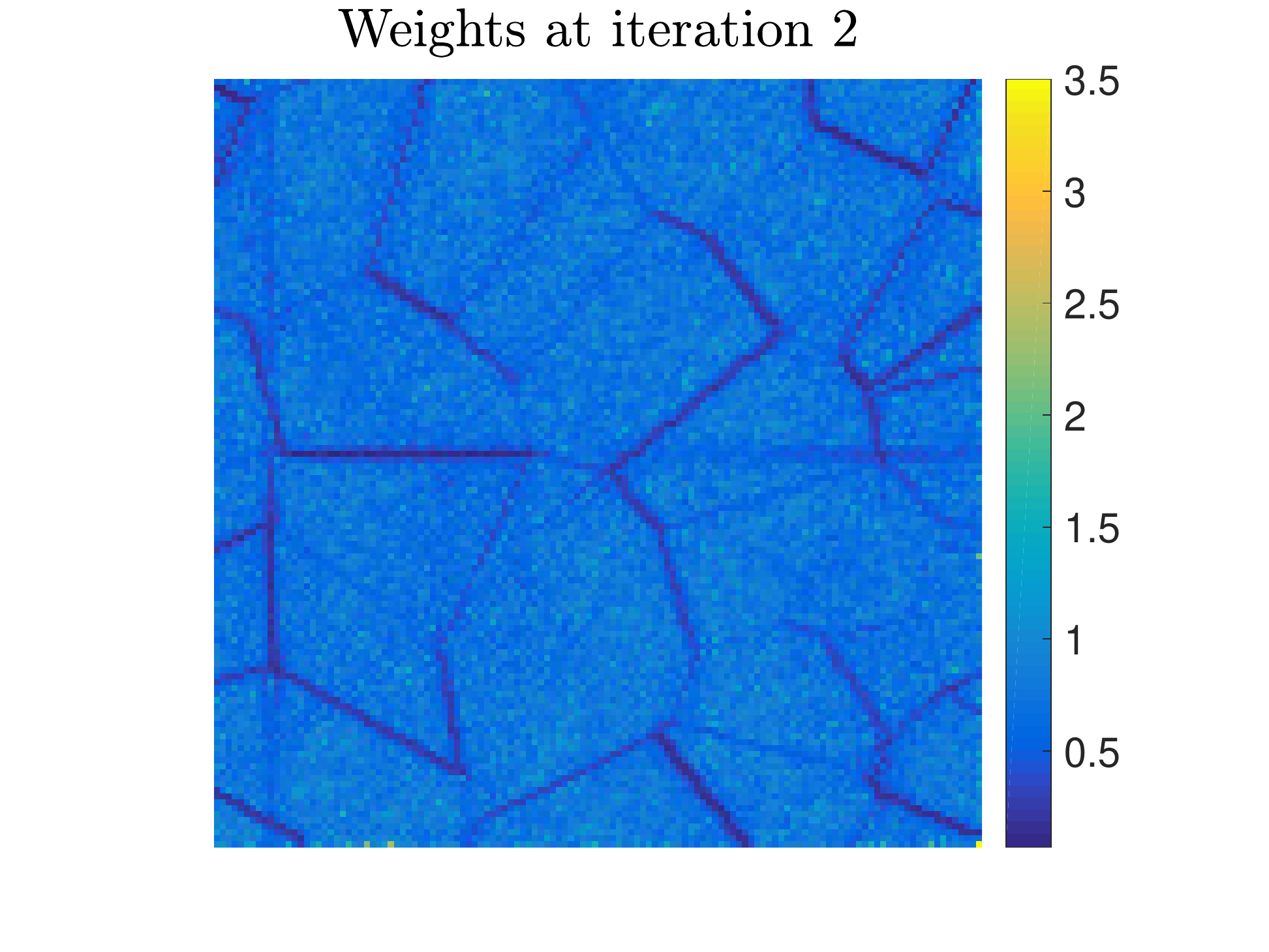} &
\includegraphics[width=5cm]{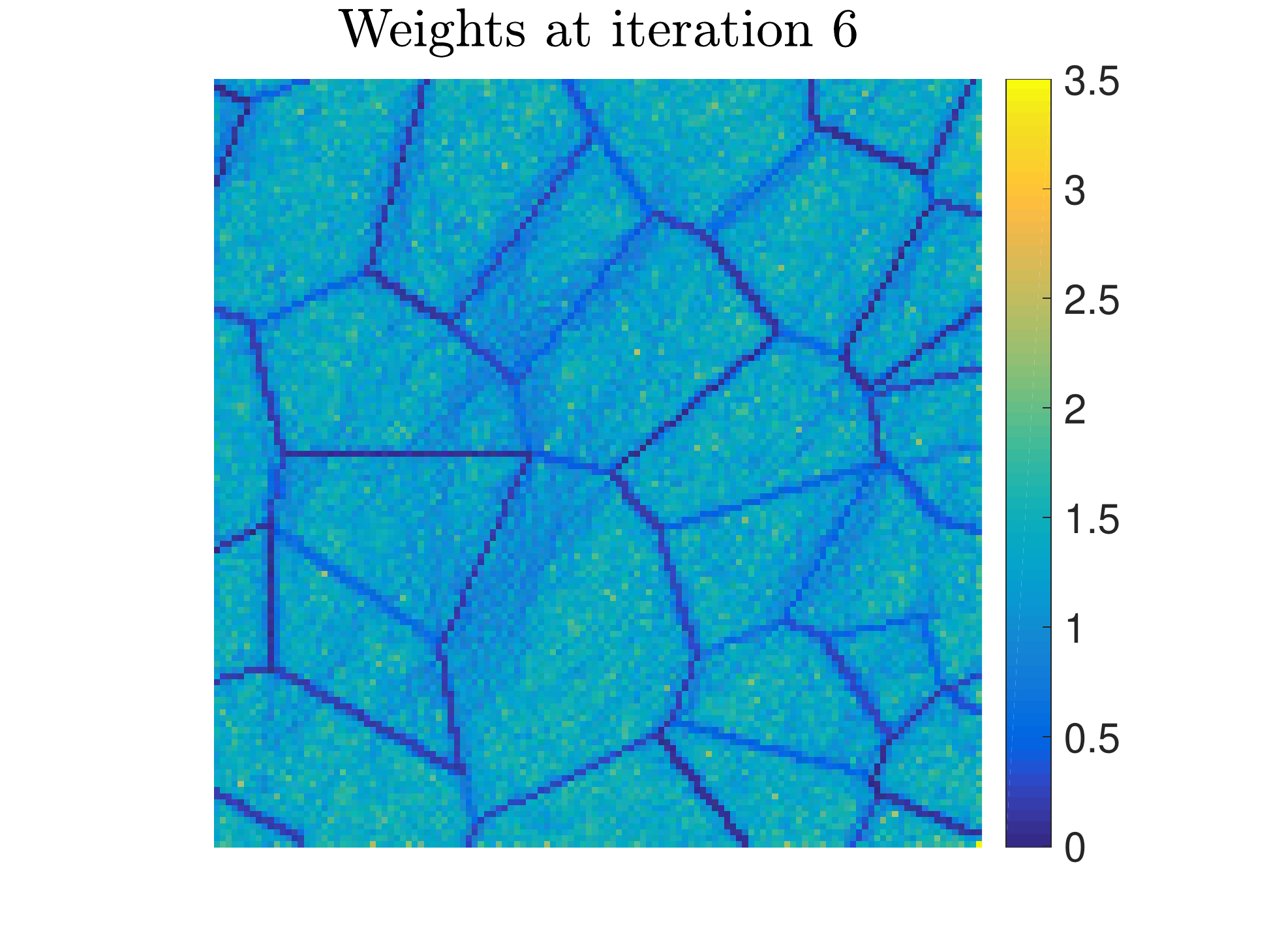} & 
\includegraphics[width=5cm]{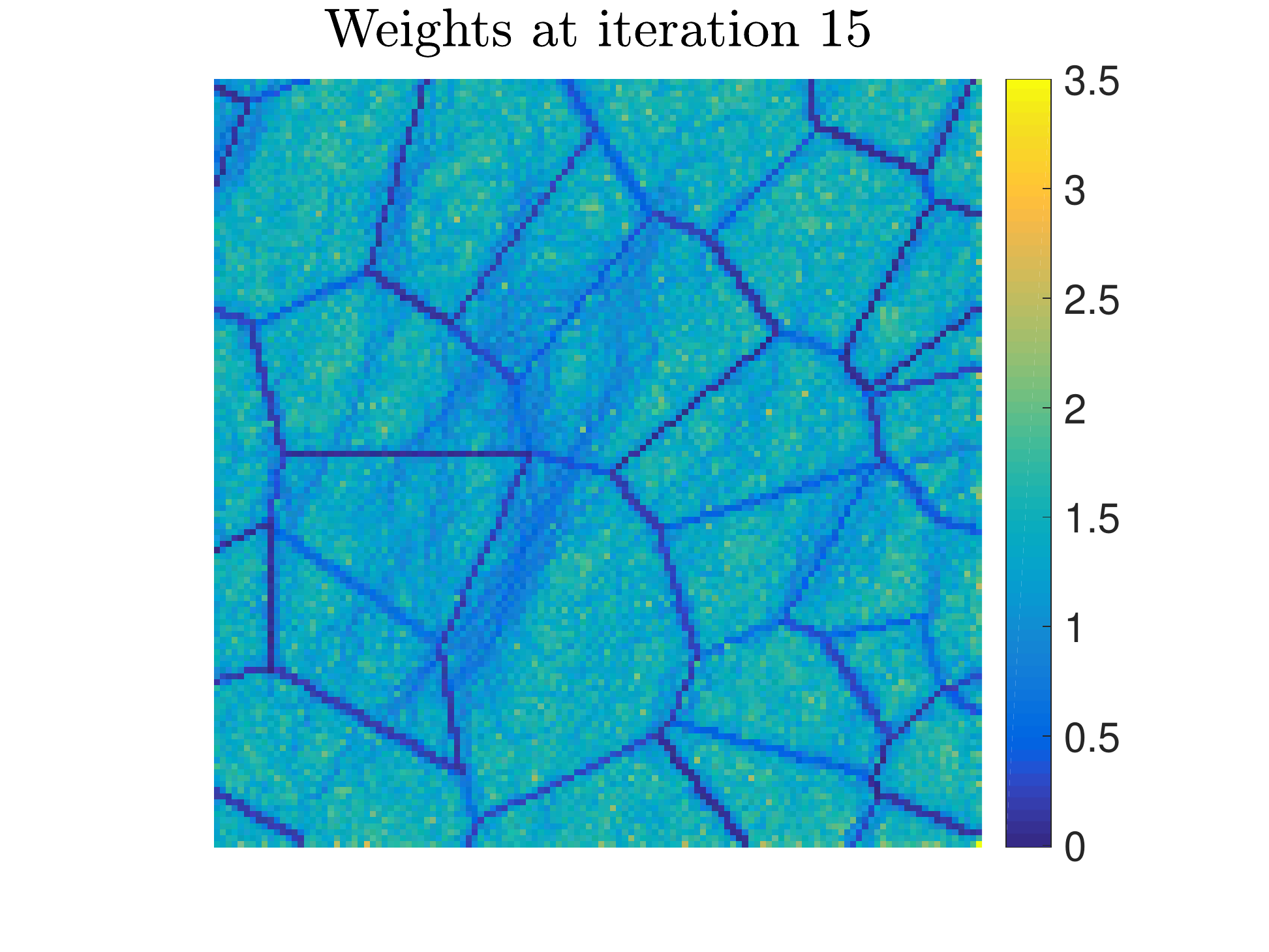}
\end{tabular}
\end{center}
\caption{\emph{Parallel-beam CT} test problem. Upper row: new weights, to be applied to the vertical derivatives. Middle row: new weights, to be applied to the horizontal derivatives. Lower row: IRN-TV weights. The weights used at the outer iterations $\ell=2,6,15$ are displayed in logarithmic scale.}
\label{fig:Radon1TVweight}
\end{figure}

We now consider some comparisons between the new solver and different inner-outer iterative methods for edge enhancement in imaging. First of all, we assess the effect of a different choice of the weight matrix: namely, we keep an inner iteration scheme based on the hybrid solver for general form Tikhonov described in Section \ref{ssec:hybridgen}, together with the discrepancy principle to adaptively set the regularization parameter, and we choose the weights as in (\ref{eq:IRNTVweights}), i.e., the weights associated to the IRN-TV methods. In Figure \ref{fig:Radon1TV} we display the behavior of the relative errors and regularization parameters versus the number of outer iterations, considering the new weights and the IRN-TV ones. 
%In both cases, joint bidiagonalization, with adaptive regularization parameter selection (according to the discrepancy principle) is used. 
We can clearly see that the IRN-TV method seems to perform well during the early (outer) iterations, but then the quality of the reconstructed solution rapidly stagnates and deteriorates, while the new method keeps improving. Also, the regularization parameter selected by the discrepancy principle seems to be almost constant right from the early outer iterations when the IRN-TV weights are used, while the regularization parameter selected by the same rule keeps increasing when the new weights are used.
% computes  the new method is amazing at the end.
Correspondingly, Figure \ref{fig:Radon1TVrec} shows the reconstructions computed by the two methods at selected outer iterations: namely, at iterations $\ell=2$ (i.e., as soon as the reweighting becomes effective; see Section \ref{ssec:mainalgo}), $\ell=6$ (i.e., when the relative error associated to the IRN-TV weights is lower than the one obtained employing the new weights), and $\ell=15$ (i.e., when the stopping criterion for the outer iterations is satisfied employing the new weights). We can clearly see that, contrarily to the new method, the reconstructions obtained using IRN-TV weights do not improve much when the outer iterations proceed, and piecewise constant features are never completely recovered.
In order to better understand the reasons behind such a different performances of these two inner-outer iterative schemes, in Figure \ref{fig:Radon1TVweight} we display the entries of the weighting diagonal matrices at selected outer iterations, rearranged as images. Note that, directly from the definition of the weights (see, for instance, (\ref{eq:GradientMap}) and (\ref{eq:IRNTVweights})), when considering the new weights it is meaningful to display two images (one for the weights applied to the vertical derivatives, and one for the weights applied to the horizontal derivatives), while only one image suffices when considering the IRN-TV weights (because the same weights containing information about both the vertical and the horizontal derivatives are applied to both the vertical and the horizontal derivatives). Looking at Figure \ref{fig:Radon1TVweight} it is evident that the new weights associated to both the vertical and horizontal derivatives correctly recover most of the edge locations (which are mapped to the smallest values $\geq 0$), and the smooth regions (which are mapped to the largest values $\leq 1$): in this way, appropriate penalization happens and the reconstructions as well as the weights improve along the outer iterations. The same is not true for the IRN-TV weights: while the locations of the edges are somewhat recovered (and the smallest $\geq 0$ weights are assigned to them), the smooth regions do not properly show up (and very oscillating small weights are assigned to them too): these weights are clearly not effective in enforcing piecewise smoothness, and almost no improvement can be seen in both reconstructions and weights as the outer iterations proceed. 

%%% GOOD SENTENCE!
% It is clear that the weights associated to the IRN-RW method are quite poor in revealing the structure of the image; the new ones are much better. 

\begin{figure}[htbp]
\begin{center}
\begin{tabular}{cc}
\includegraphics[width=5.9cm]{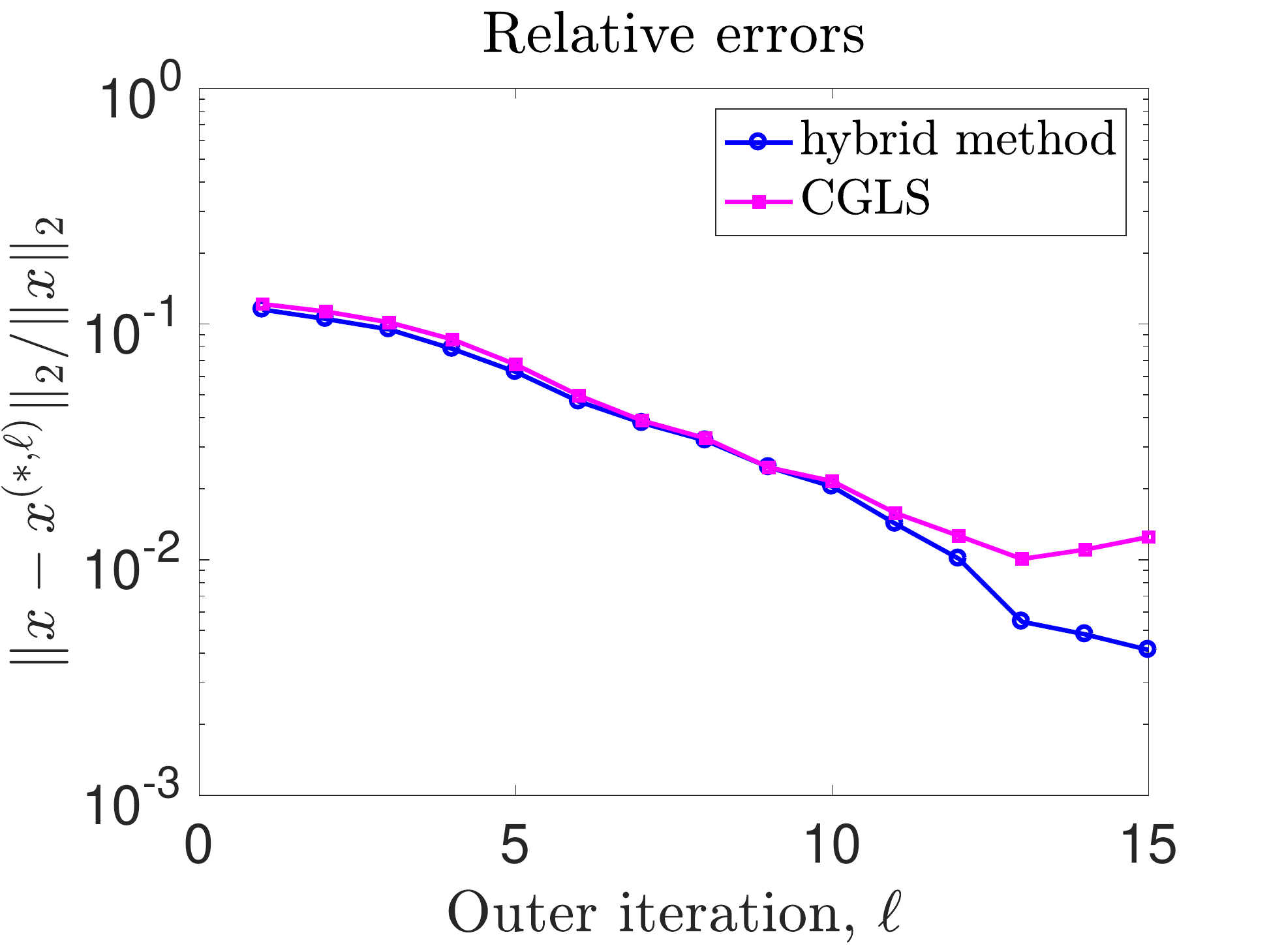} &
\includegraphics[width=5.9cm]{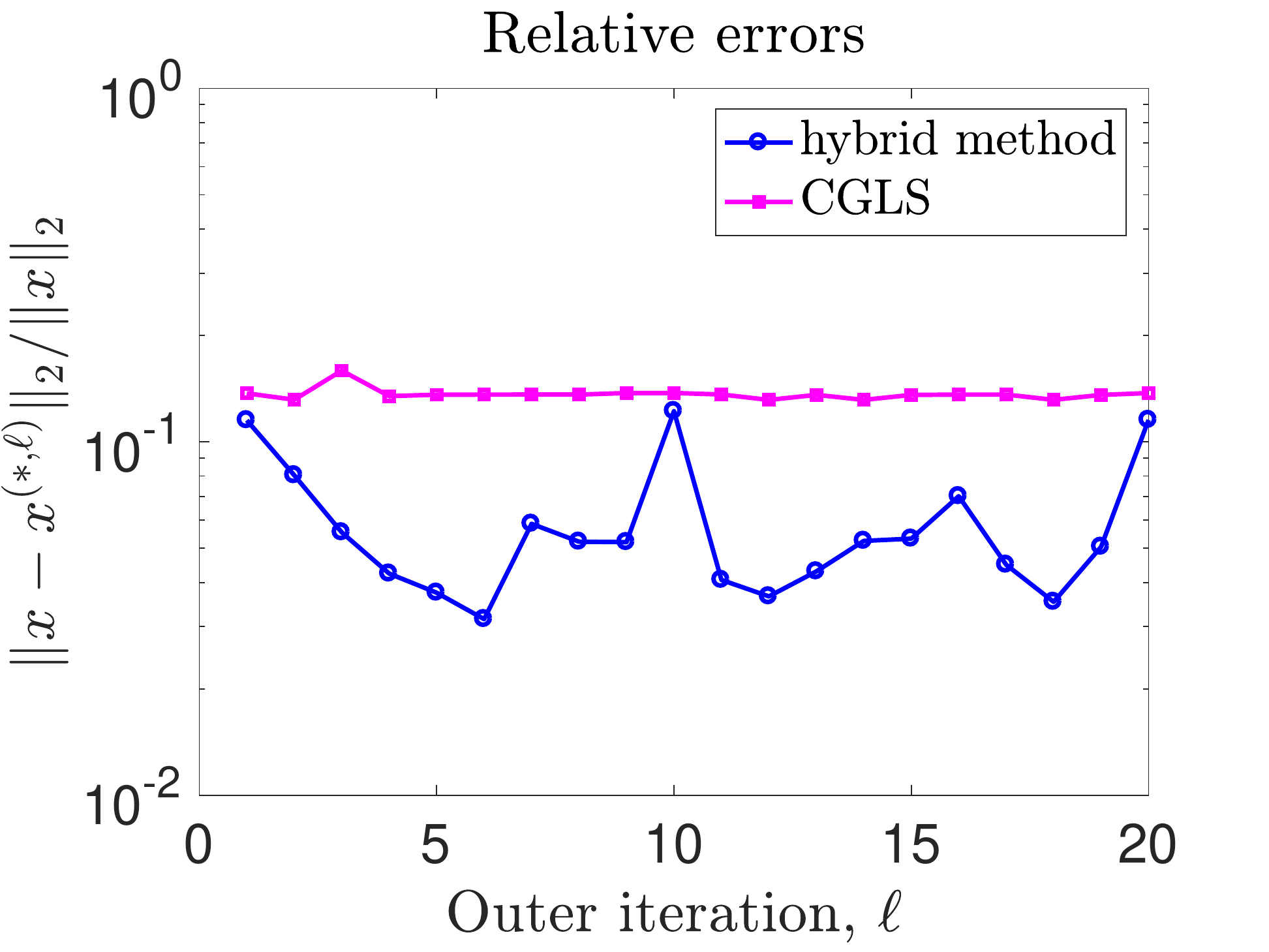}
\end{tabular}
\end{center}
\caption{\emph{Parallel-beam CT} test problem. Relative errors versus number of (outer) iterations for different inner solvers. Left frame: the new weights are used. Right frame: the IRN-TV weights are used.}
\label{fig:Radon1TV_comparisons}
\end{figure}

Finally, we show two more comparisons, that assess the influence of the inner iterative solver on the overall behavior of the method. Namely, we consider the inner-outer iterative schemes implemented with both the new and the IRN-TV weights, and with both the hybrid and the CGLS methods as inner solvers. Since CGLS requires a fixed value of the regularization parameter to be available in advance of each cycle of inner iterations, we first run the methods based on the hybrid solver for general form Tikhonov that adaptively chooses the regularization parameter at each inner iteration according to the discrepancy principle: in this way, a parameter $\lambda_{\ast,\ell}$ is eventually set at the $\ell$th outer iteration. When running the methods based on CGLS, we take $\lambda_{\ast,\ell}$ as regularization parameter for the $\ell$th outer iteration. Figure \ref{fig:Radon1TV_comparisons} displays the relative errors versus the number of outer iterations for the four instances of inner-outer iterative methods just described.
In both the frames displayed in Figure \ref{fig:Radon1TV_comparisons} we can see that the results obtained using CGLS as inner iterative solver are not as good as the ones obtained using the hybrid strategy; this is mostly evident when the IRN-TV weights are used. This could be because $\lambda_{\ast,\ell}$ is supposed to be a nearly-optimal regularization parameter when the joint bidiagonalization algorithm is used to project the $\ell$th quadratic problem (\ref{eq:seqtrue}); however CGLS projects the $\ell$th quadratic problem into a different subspace, and the same regularization parameter $\lambda_{\ast,\ell}$ might not be such a nearly-optimal choice for CGLS. Moreover, adaptively choosing the regularization parameter for the CGLS inner iterations might help improving the quality of the reconstructions. 
%When comparing with the classical IRN approach, both with the new weights and with the ``traditional'' TV weights, the results may be not so good also because of a different choice of the regularization parameter: we keep the regularization parameter that would have been chosen at the end of each JBD cycle, but this mat not be optimal for the following reasons. (1) The parameter is fixed for each linear problem, and not adaptively chosen during the iterations. (2) The parameter satisfies the discrepancy principle for a problem that is projected in a different approximation subspace (different from the CGLS one). 

%Perhaps change phantom immediately. 

%\newpage

%\subsection{Example from X-ray CT, Limited Angle}
\paragraph{Parallel-beam CT with limited angle}

\begin{figure}[htbp]
\begin{center}
\begin{tabular}{cc}
%\footnotesize{True image} & \footnotesize{Measured data (sinogram)}\\
\includegraphics[height=4cm]{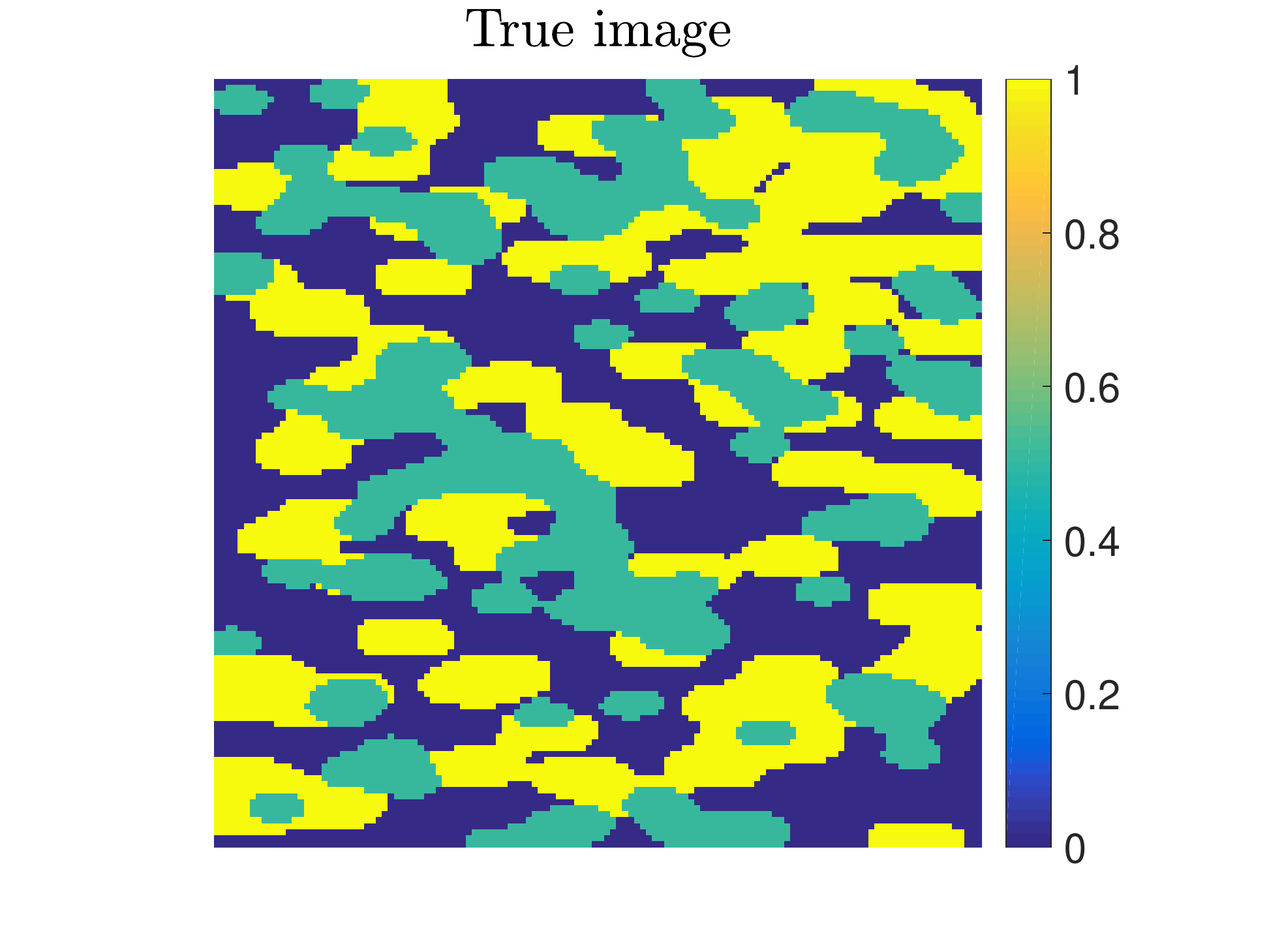} &
\includegraphics[height=4cm]{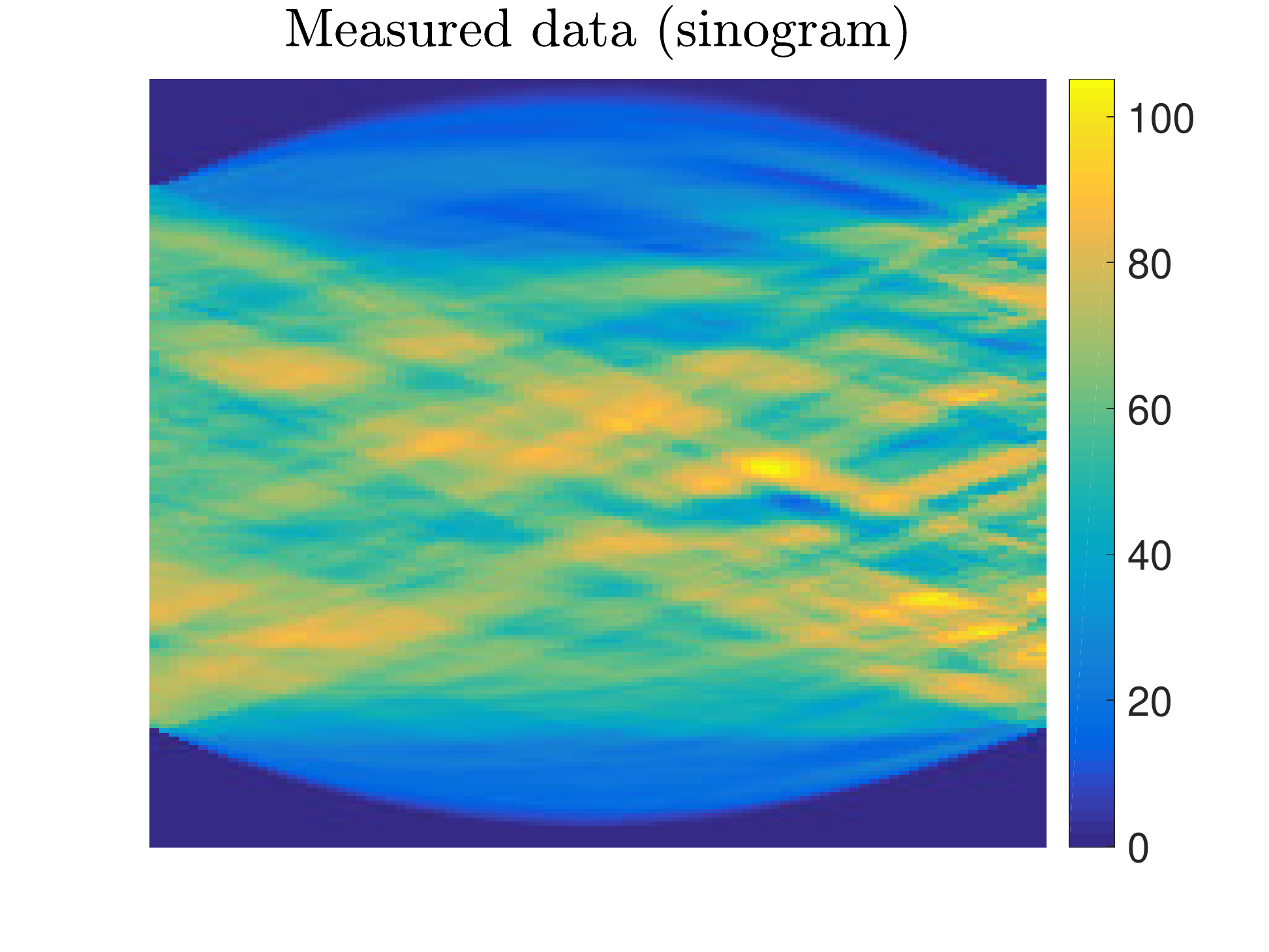}
\end{tabular}
\end{center}
\caption{\emph{Parallel-beam CT with limited angle} test problem. Left frame: true image phantom $x$. Right frame: the measured data (sinogram) $b$, where the projections are taken only at angles $0, 1, \ldots, 90$.}
\label{fig:Radon1bData}
\end{figure}

We use IR Tools to vary some of the options defining the previous X-ray tomography simulation. In particular, we still wish to consider a parallel beam X-ray transmission problem, but we would like to change the phantom image to the so-called ``three-phases'' one from AIR Tools II of size $128\times 128$ pixels, which has many piecewise constant patterns, and we would like to make the reconstruction problem more challenging by using (limited) projection angles $0, 1, \ldots, 90$. Within IR Tools, this problem can be generated by just replacing the first line of the MATLAB statements specific of the previous example, i.e., we should define a new \texttt{ProblemOptions} structure in the following way:
\begin{verbatim}
     ProblemOptions = PRset('angles', 0:90, 'phantomImage', 'threephases');
\end{verbatim}
%This is a slightly
%more limited angle case than the previous example, making the problem more
%ill-posed, and more difficult to obtain good reconstructed images.
%
%Using IR Tools, the problem can be generated with the following MATLAB statements:
%\begin{verbatim}
%     ProblemOptions = PRset('angles', 0:90, 'phantomImage', 'grains');
%     [A, b_true, x_true, ProblemInfo] = PRtomo(128, ProblemOptions);
%     b = PRnoise(b_true, 1e-3);
%\end{verbatim}
%
Figure~\ref{fig:Radon1bData} shows the true phantom image, along with the
measured data $b$ (which are again corrupted by Gaussian white noise of level $10^{-3}$). 

\begin{figure}[htbp]
\begin{center}
\begin{tabular}{cc}
\includegraphics[width=6cm]{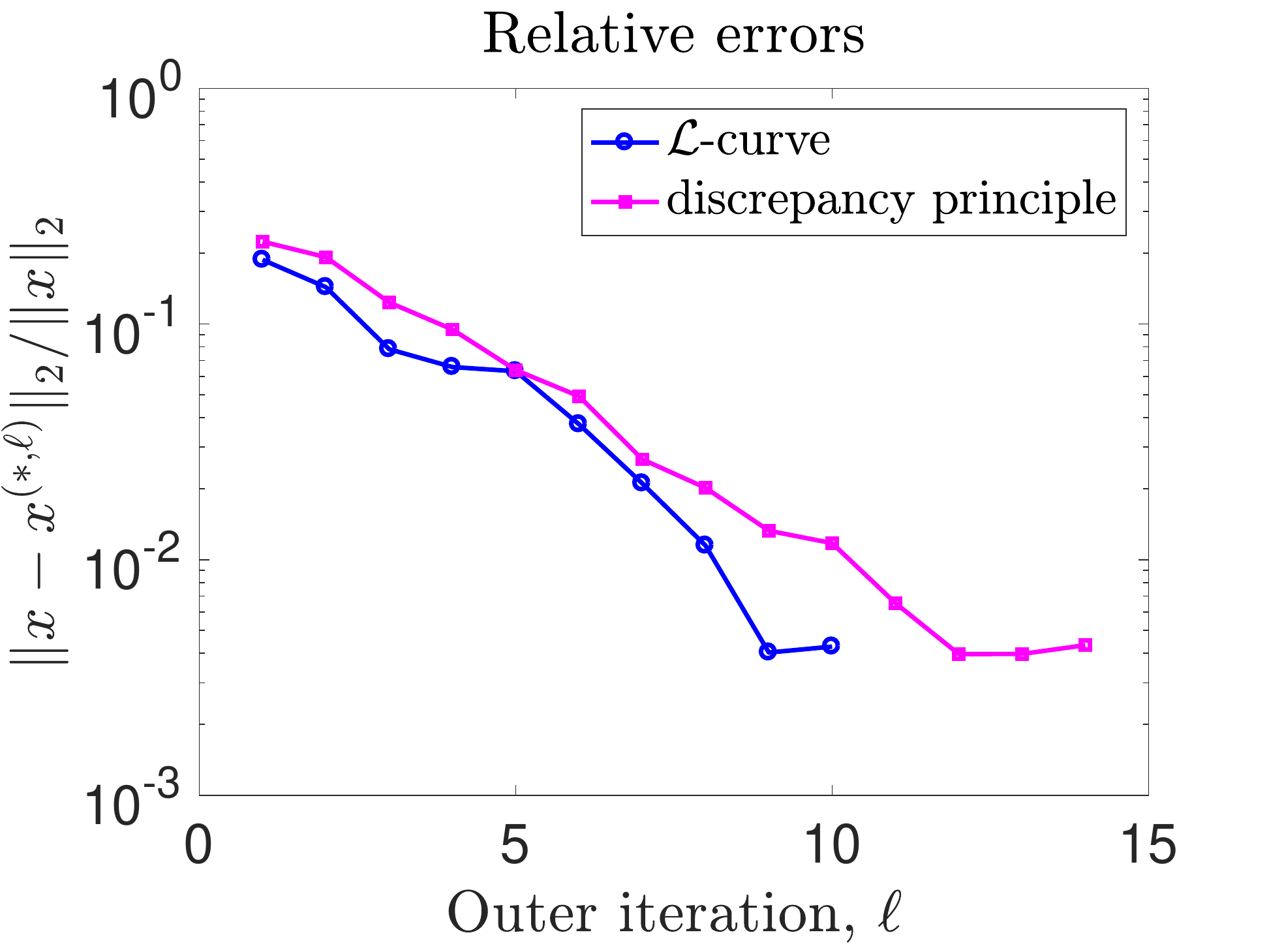} &
\includegraphics[width=6cm]{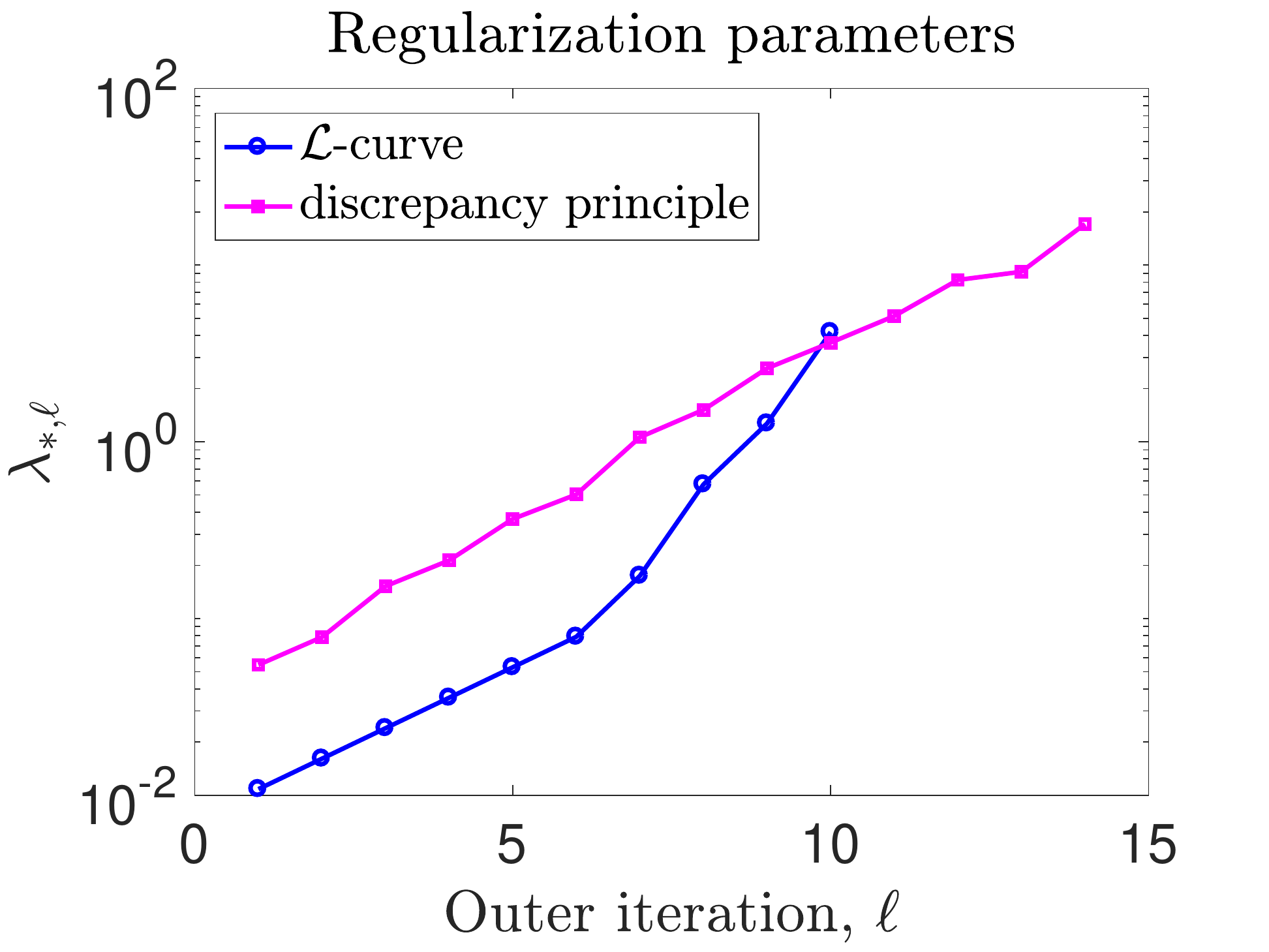}
\end{tabular}
\end{center}
\caption{\emph{Parallel-beam CT with limited angle} test problem. Relative errors and regularization parameters values at each outer iteration $\ell$,until the stopping criterion is satisfied. Both the discrepancy principle and the ${\mathcal L}$-curve criterion are considered.}
\label{fig:Radon1bIterations}
\end{figure}

As in the previous example, we first run our algorithm with different parameter choice strategies within the inner hybrid scheme for generalized Tikhonov regularization: when the discrepancy principle and the $\mathcal{L}$-curve criterion are used, the stopping rule of Section \ref{ssec:mainalgo} (decrease of $\|Lx^{(\ast,\ell)}\|_2$) is satisfied after $\ell=14$ and $\ell=10$ outer iterations, respectively. Figure~\ref{fig:Radon1bIterations} shows a plot of the relative errors
and chosen regularization parameters at each outer iteration. We can clearly see that, when considering both the discrepancy principle and the $\mathcal{L}$-curve criterion, the reconstructions greatly improve as the outer iterations proceed, with the latter strategy being more efficient (it employs less outer iterations to achieve a relative error comparable to the one obtained by the discrepancy principle at the end of the iterations). Note that, as expected (see again the theory presented in Section \ref{ssec:analysis}), the regularization parameters increase as
the outer iterations proceed and, because of this and the
fact that identification of the edges is also improved, the relative errors generally decrease. The behavior of the relative error versus the number of total iterations, and the dependency of the results on the power $p>0$ appearing within the weights, are analogous to the ones displayed within the previous test problem. 
%We can also see that the regularization parameter selected at the end of each inner iteration cycle consistently increases when both parameter rules are employed: this behavior is expected and, in . 
%he terminated after 4 outer iterations, where at iteration $\ell$ the hybrid 
%method determined an estimate of the regularization parameter $\lambda_{*,\ell}$,
%and a corresponding reconstructed image, $x^{(*,\ell)}$.  
%Figure~\ref{fig:Radon1bLCurves}
%shows a plot of all ${\mathcal L}$-curves for each iteration. Observe the nesting
%property of the ${\mathcal L}$-curves, and 
%how the chosen regularization
%parameter (corresponding to the corner of each ${\mathcal L}$-curve) increases.
%
%
%

\begin{figure}[htbp]
\begin{center}
\begin{tabular}{cc}
\includegraphics[width=6cm]{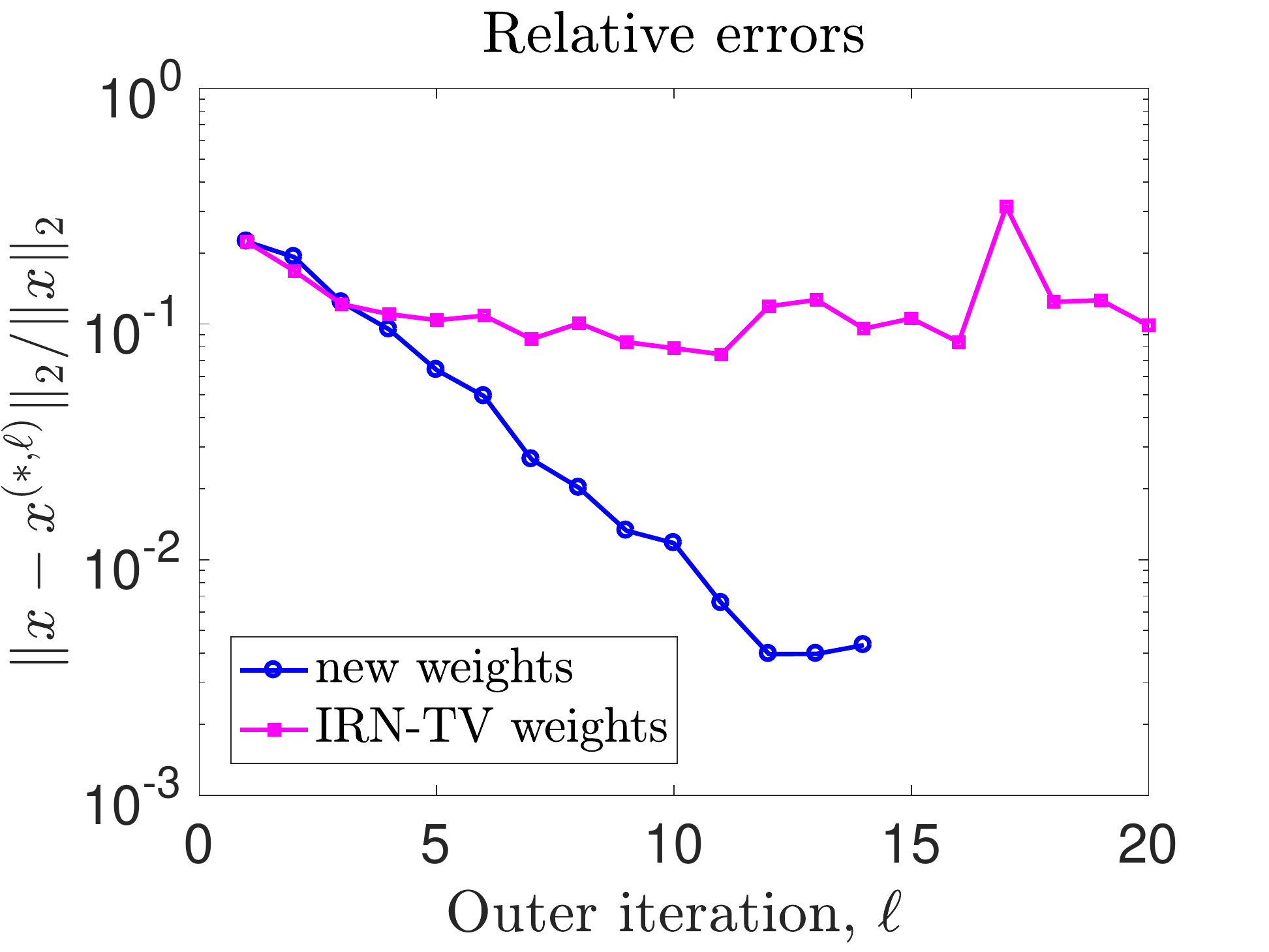} &
\includegraphics[width=6cm]{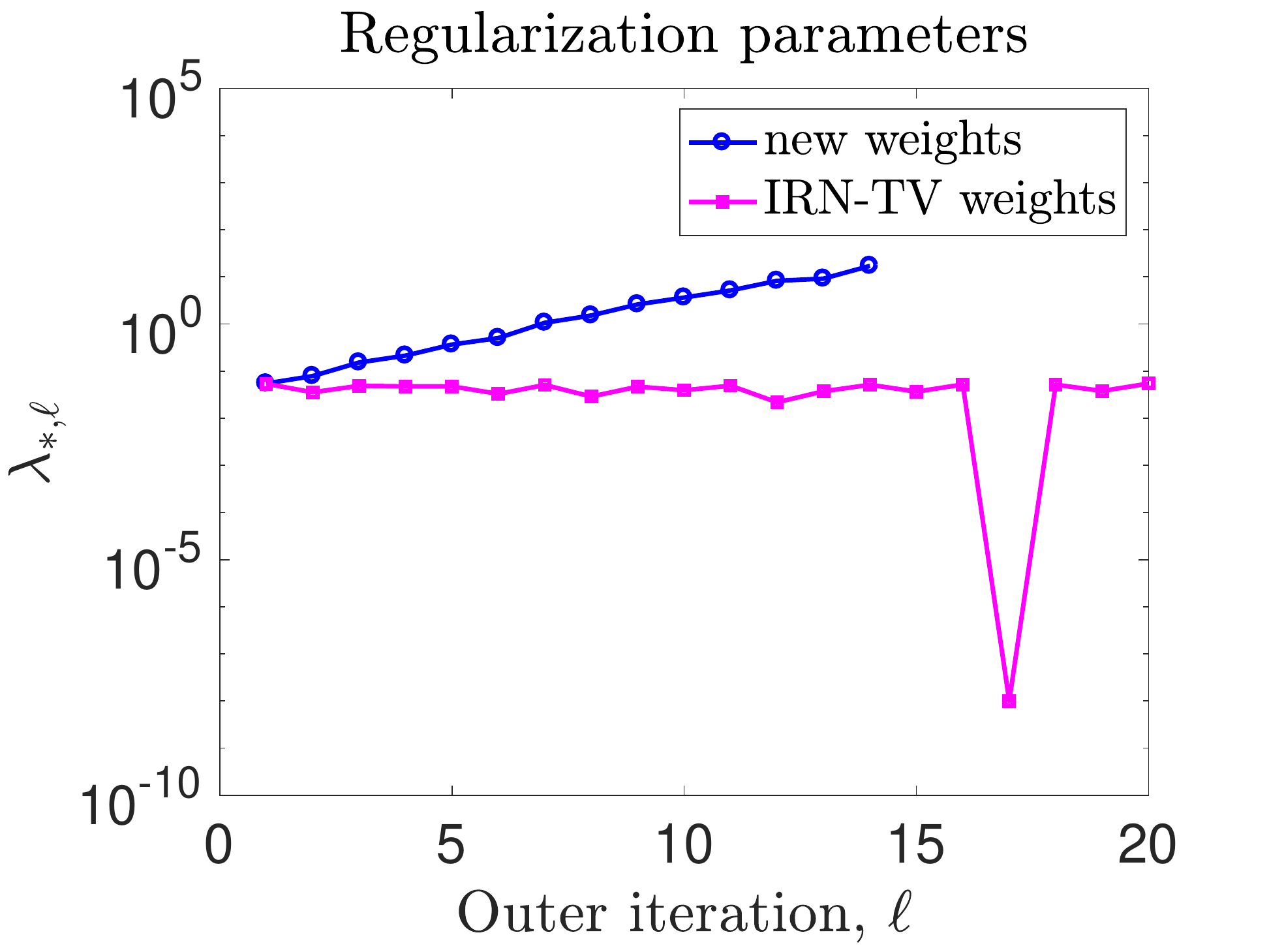}
\end{tabular}
\end{center}
\caption{\emph{Parallel-beam CT with limited angle} test problem. Relative errors and regularization parameter versus number of (outer) iterations; both methods use the hybrid method for general form Tikhonov during the inner iterations, and adaptively select the regularization parameter according to the discrepancy principle.}
\label{fig:Radon2TV}
\end{figure}

\begin{figure}[htbp]
\begin{center}
\begin{tabular}{cc}
\includegraphics[width=6cm]{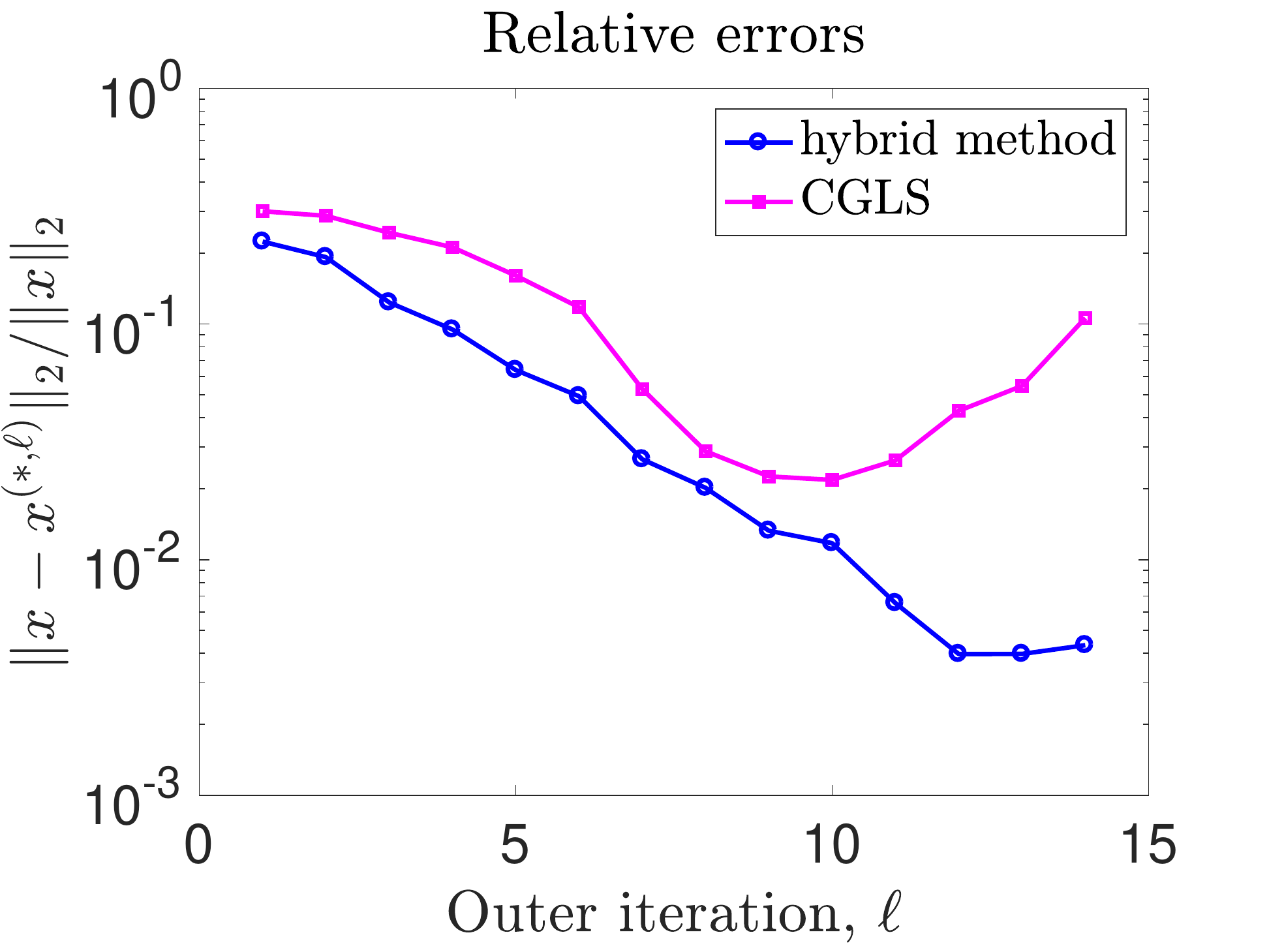} &
\includegraphics[width=6cm]{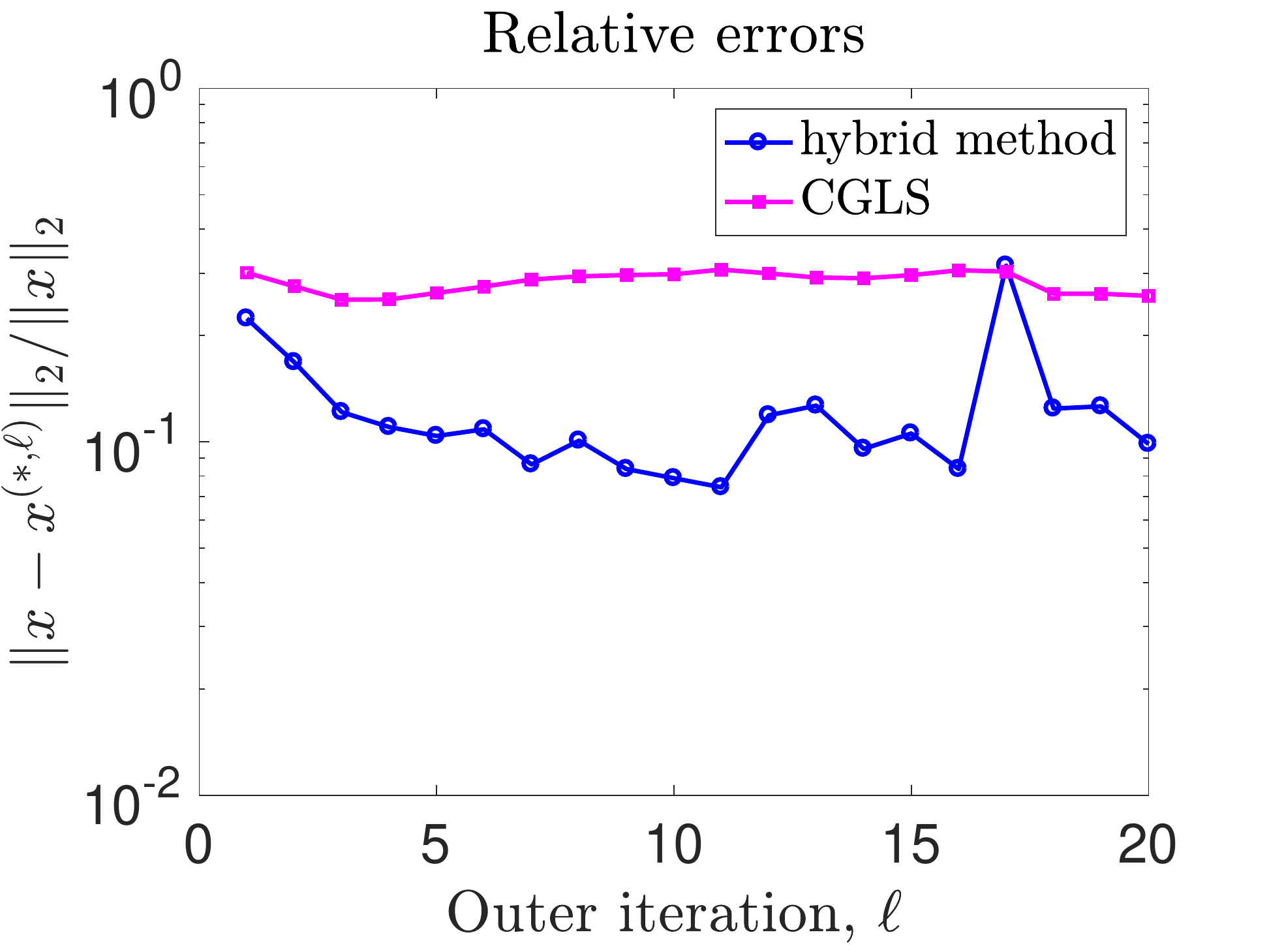}
\end{tabular}
\end{center}
\caption{\emph{Parallel-beam CT with limited angle} test problem. Relative errors versus number of (outer) iterations. Left frame: the new weights are used. Right frame: the IRN-TV weights are used.}
\label{fig:Radon2TV_comparisons}
\end{figure}

Next, we compare the new method to other inner-outer iterative methods for edge enhancement in imaging. Figure \ref{fig:Radon2TV} displays the relative error and regularization parameter values versus the number of outer iterations for the new method, and for a method that still employs the hybrid solver for general form Tikhonov regularization to handle the inner iterations while updating the IRN-TV weights (\ref{eq:IRNTVweights}) at each outer iteration. For both methods, the discrepancy principle is employed to adaptively choose the regularization parameter at each inner iteration. As for the previous test problem, both methods perform similarly during the first (outer) iterations, but then the quality of the IRN-TV solutions rapidly stagnates and eventually worsens, while the new method keeps improving. With respect to the previous test problem, in this setting IRN-TV noticeably performs better. Also, the value of the automatically selected regularization parameter is always quite small for IRN-TV, with a steep decrease during one of the final outer iterations: this could explain why the method based on IRN-TV stagnates and worsens at some point: namely, the method fails to give more weight to the regularization term when the latter should become more accurate. 
% computes  the new method is amazing at the end.
Figure \ref{fig:Radon2TV_comparisons} assesses the influence of the inner iterative solver on the overall behavior of the method. Namely, we consider the inner-outer iterative schemes implemented with both the new and the IRN-TV weights, and with both the hybrid and the CGLS methods as inner solvers. The hybrid method chooses the regularization parameter adaptively according to the discrepancy principle as the iterations proceed, and the value $\lambda_{\ast,\ell}$ selected when the $\ell$th inner iteration cycle terminates is taken to be the fixed regularization parameter to be set in advance of the $\ell$th CGLS inner iteration cycle. 
% requires a fixed value of the regularization parameter to be available in advance of each cycle of inner iterations, we first run the methods based on the hybrid solver for general form Tikhonov that adaptively chooses the regularization parameter at each inner iteration according to the discrepancy principle: in this way, a parameter 
% $\lambda_{\ast,\ell}$ is eventually set at the $\ell$th outer iteration. When running the methods based on CGLS, we take $\lambda_{\ast,\ell}$ as regularization parameter for the $\ell$th outer iteration.
In both the frames displayed in Figure \ref{fig:Radon2TV_comparisons} we can see that the results obtained using CGLS as inner iterative solver are not as good as the ones obtained using the hybrid strategy, especially when the IRN-TV weights are used. Moreover, the CGLS solution obtained using the new weights seems to deteriorate as the outer iterations proceed (as a sort of ``semiconvergence'' appears): this can be probably avoided if the fixed regularization parameter is chosen differently. 

\begin{figure}[htbp]
\begin{center}
\begin{tabular}{ccc}
\hspace{-0.3cm}\includegraphics[width=5.5cm]{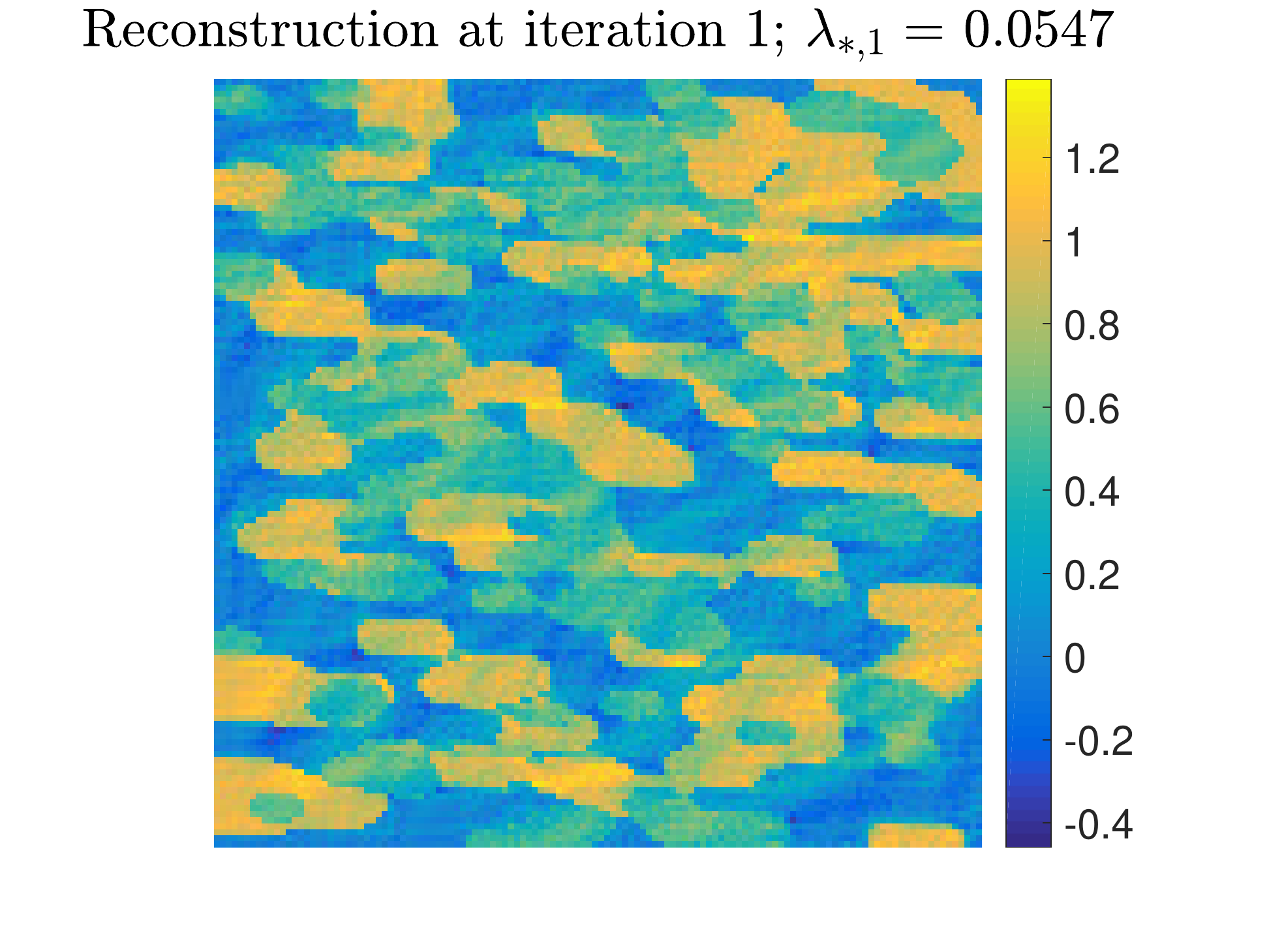} &
\hspace{-0.3cm}\includegraphics[width=5.5cm]{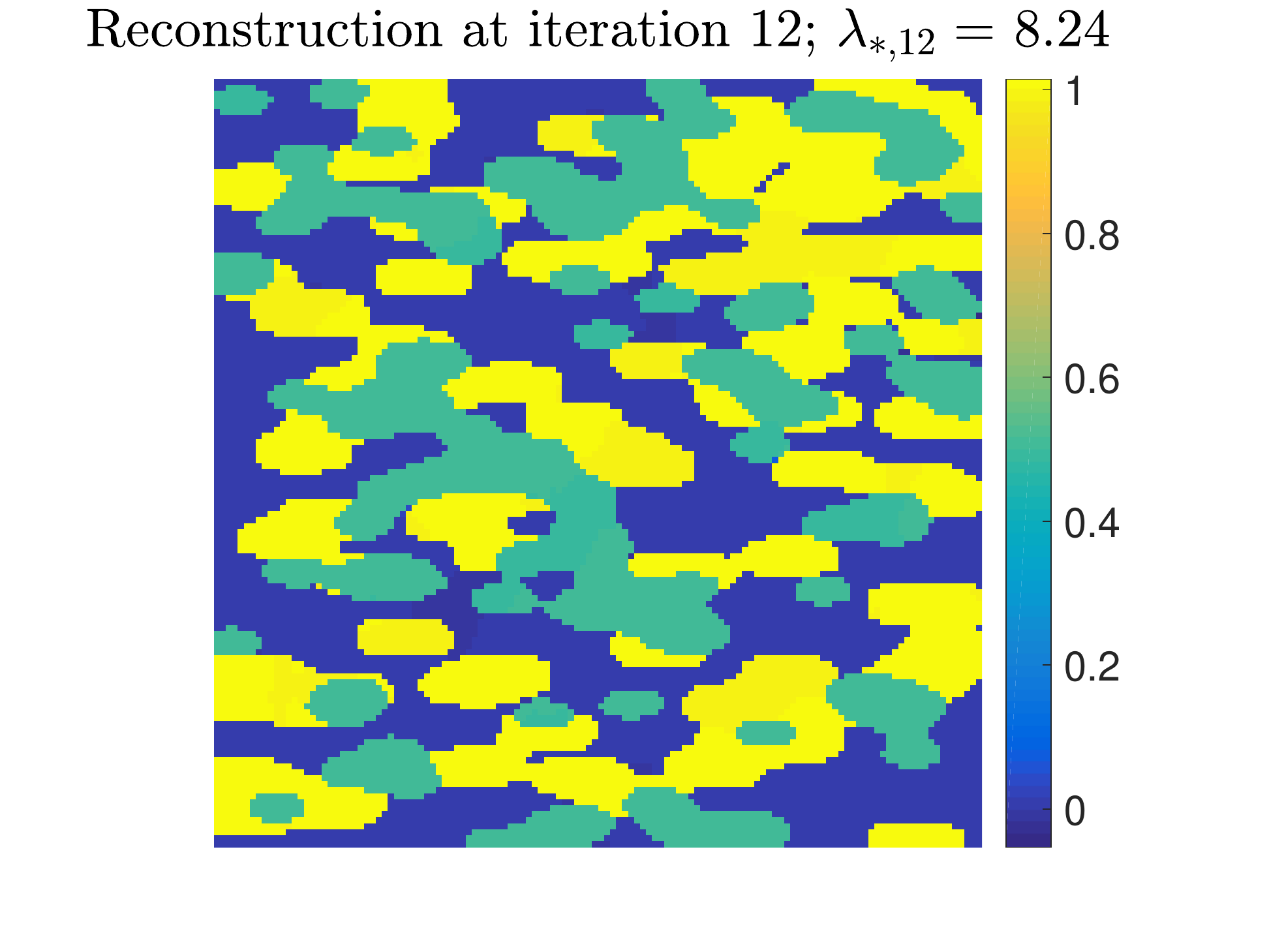} &
\hspace{-0.3cm}\includegraphics[width=5.5cm]{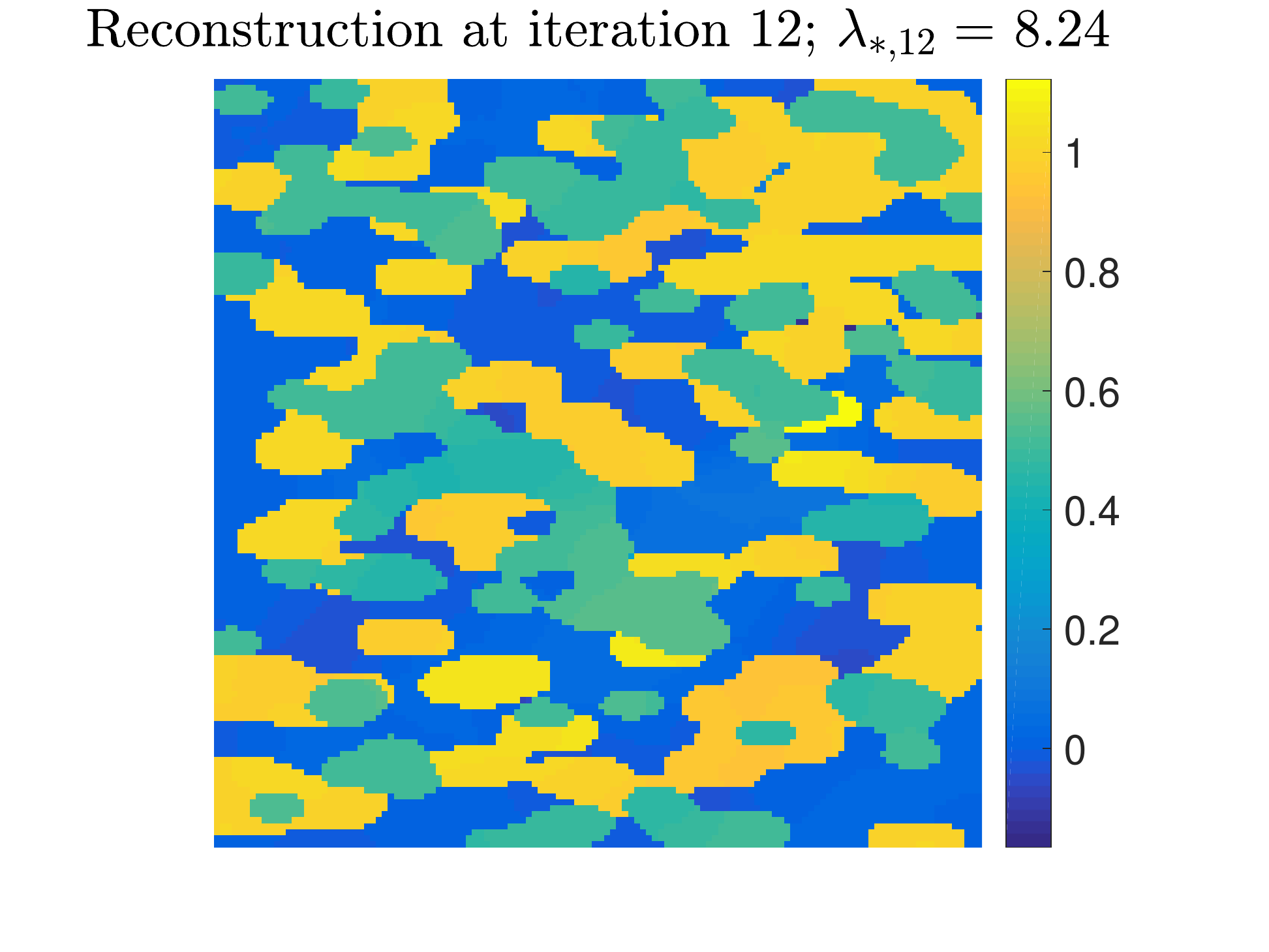}\\
\hspace{-0.3cm}\includegraphics[width=5.5cm]{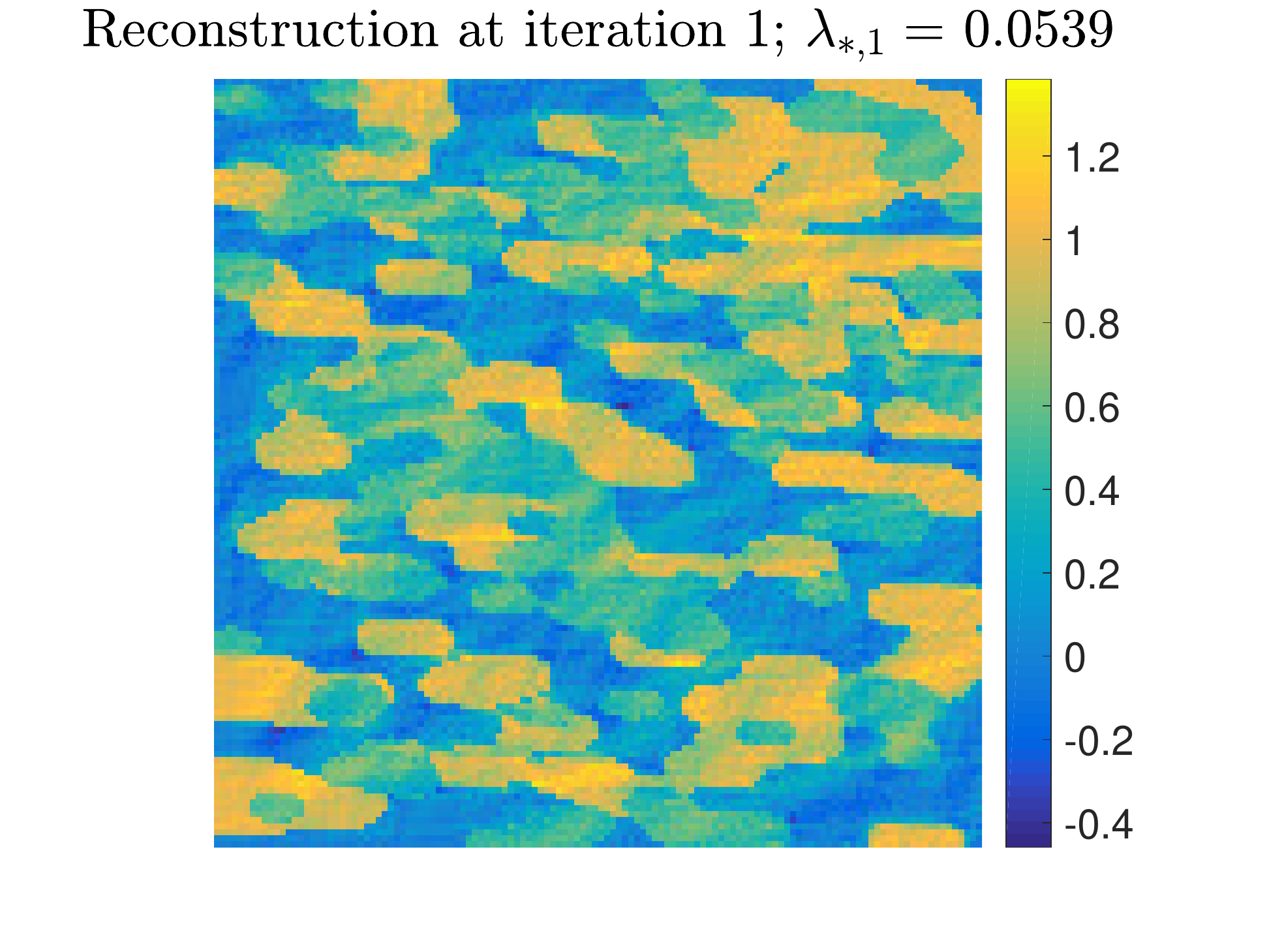} &
\hspace{-0.3cm}\includegraphics[width=5.5cm]{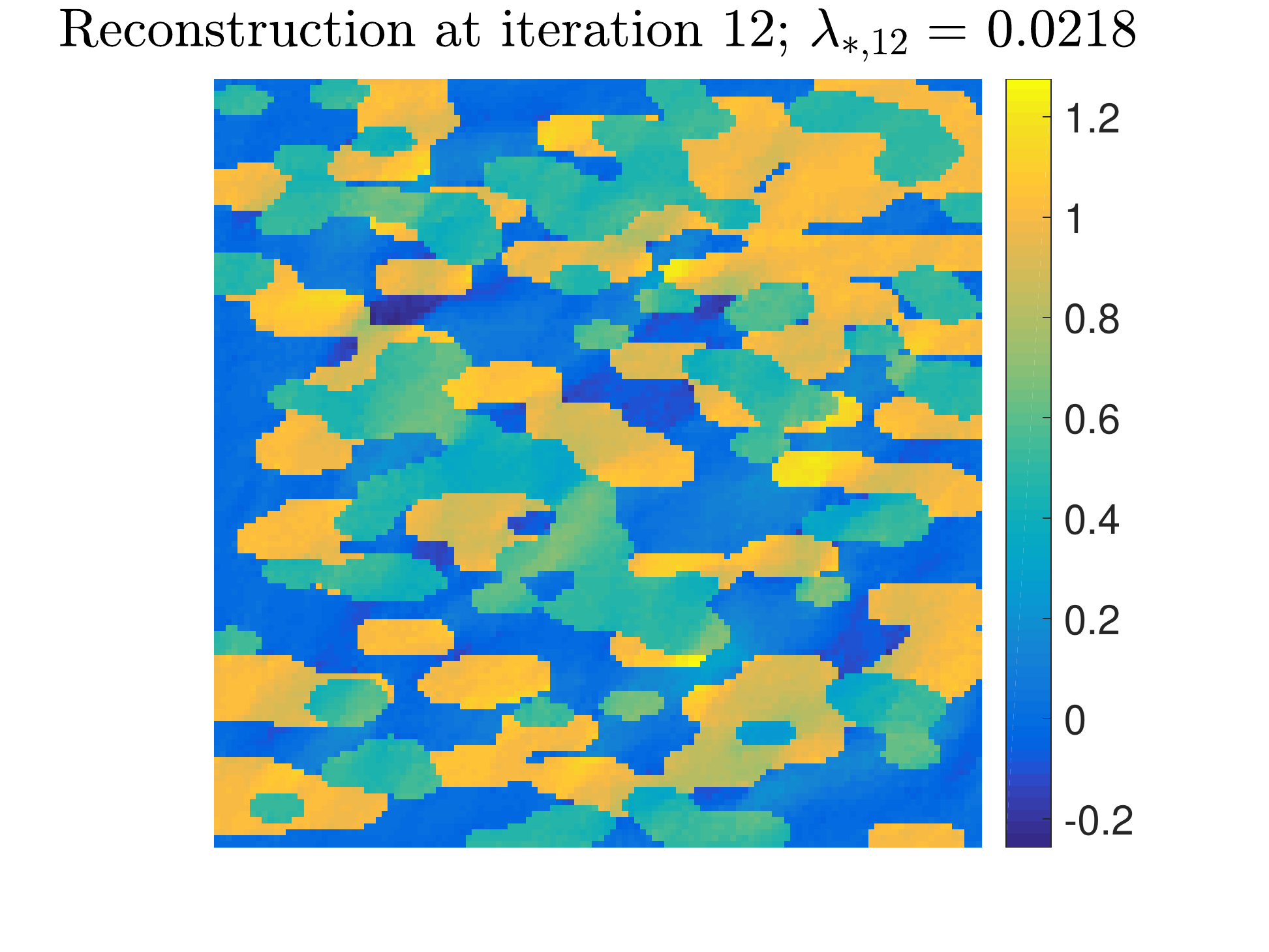} &
\hspace{-0.3cm}\includegraphics[width=5.5cm]{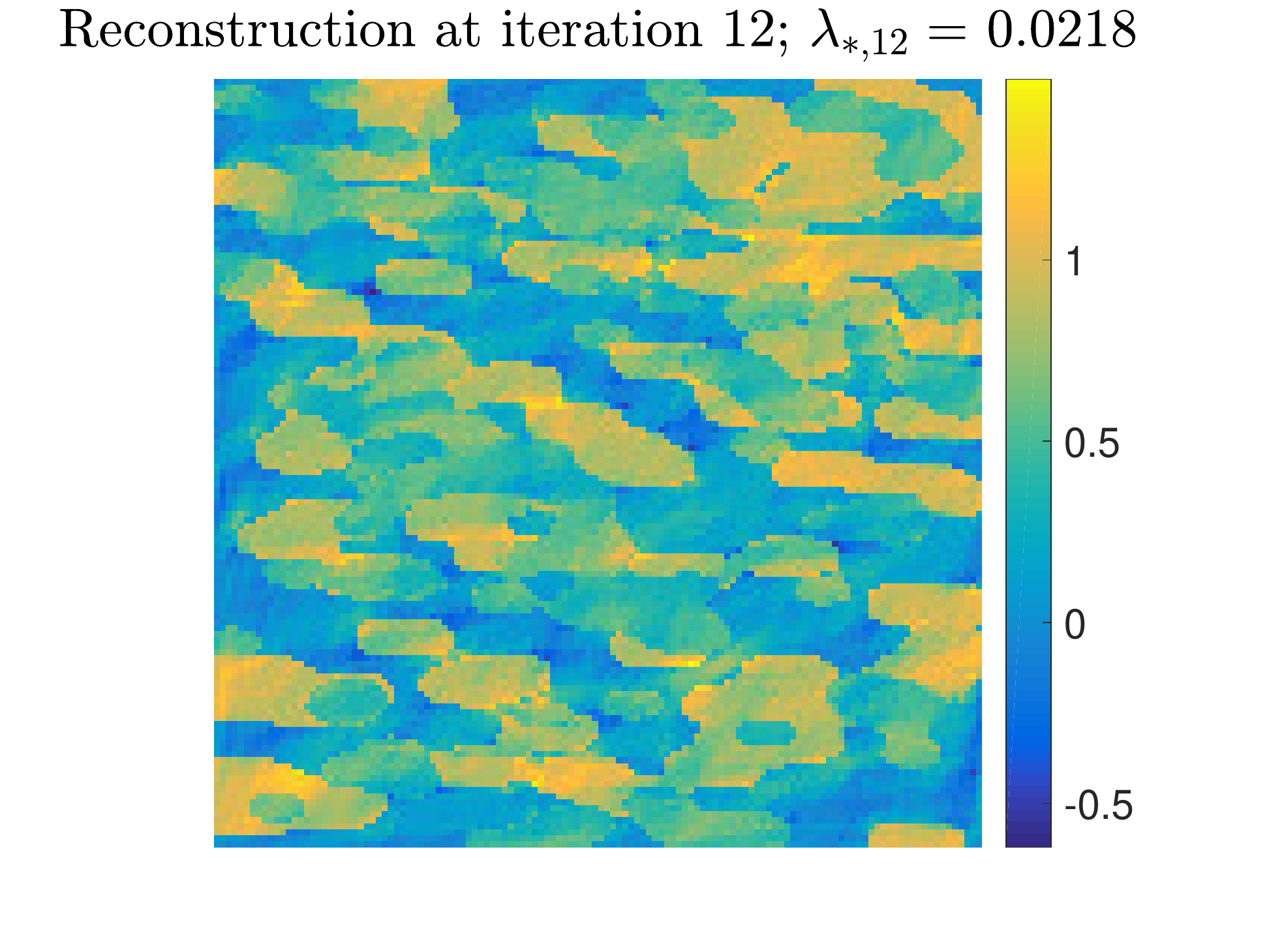}
\end{tabular}
\end{center}
\caption{\emph{Parallel-beam CT with limited angle} test problem. Upper row: reconstructions obtained using the new weights, at different outer iterations, and by different inner linear solvers. Lower row: reconstructions obtained using the IRN-TV weights, at different outer iterations, and by different inner linear solvers. The corresponding regularization parameters chosen
by the hybrid method, and used within CGLS, are displayed above each image.}
\label{fig:Radon1bSolutions}
\end{figure}

\begin{figure}[htbp]
\begin{center}
\begin{tabular}{ccc}
\includegraphics[width=5cm]{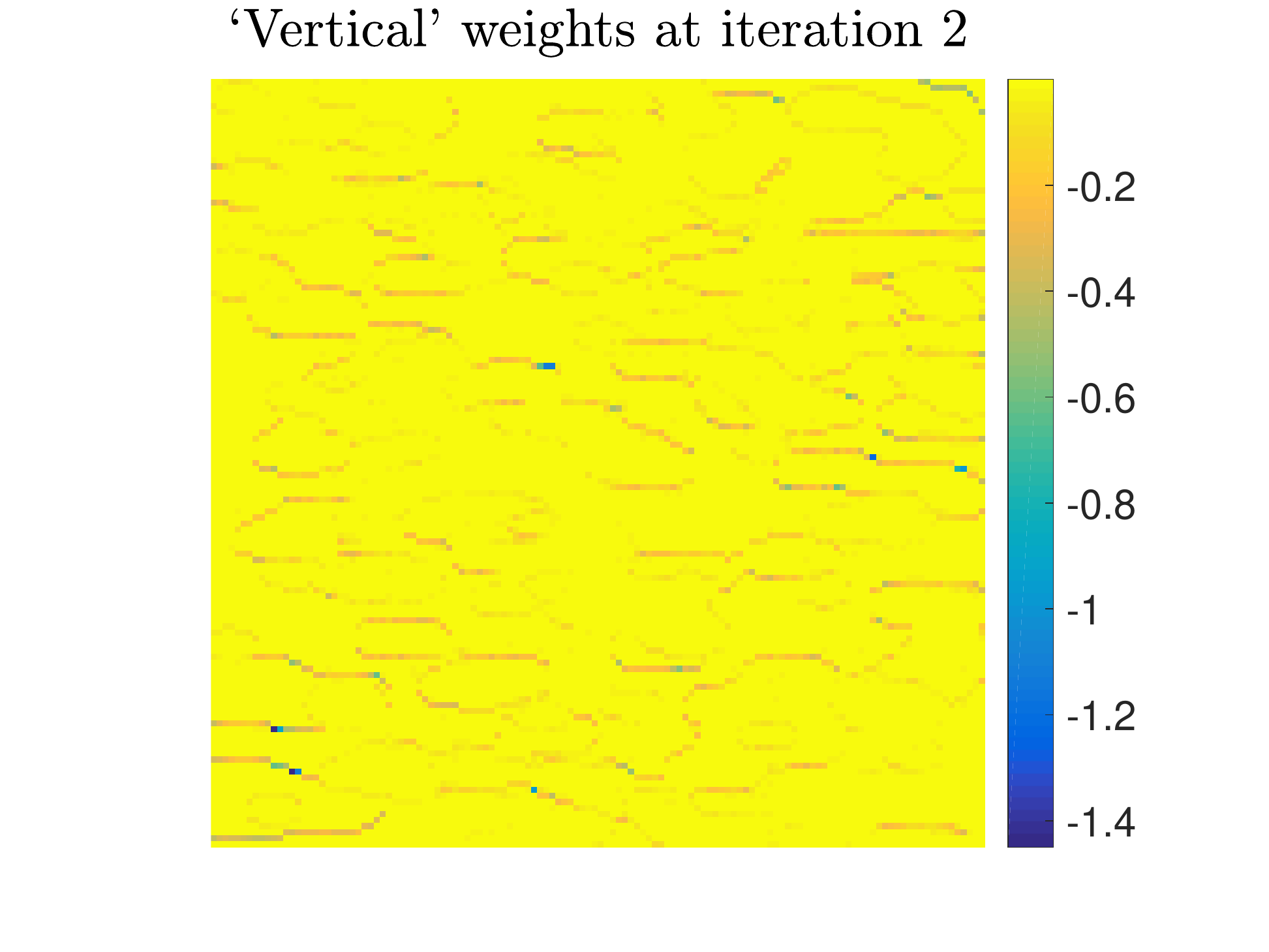} &
\includegraphics[width=5cm]{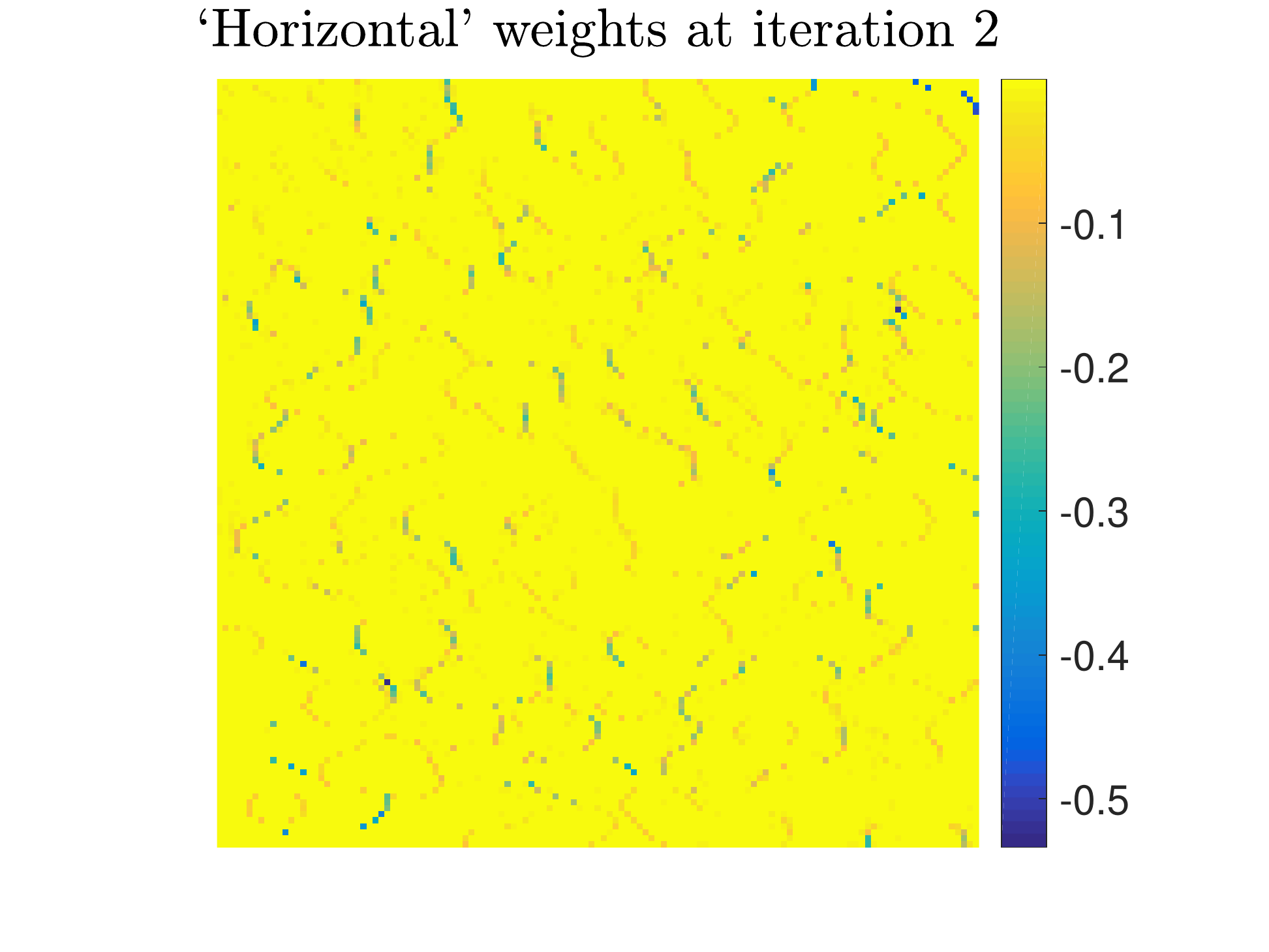} & 
\includegraphics[width=5cm]{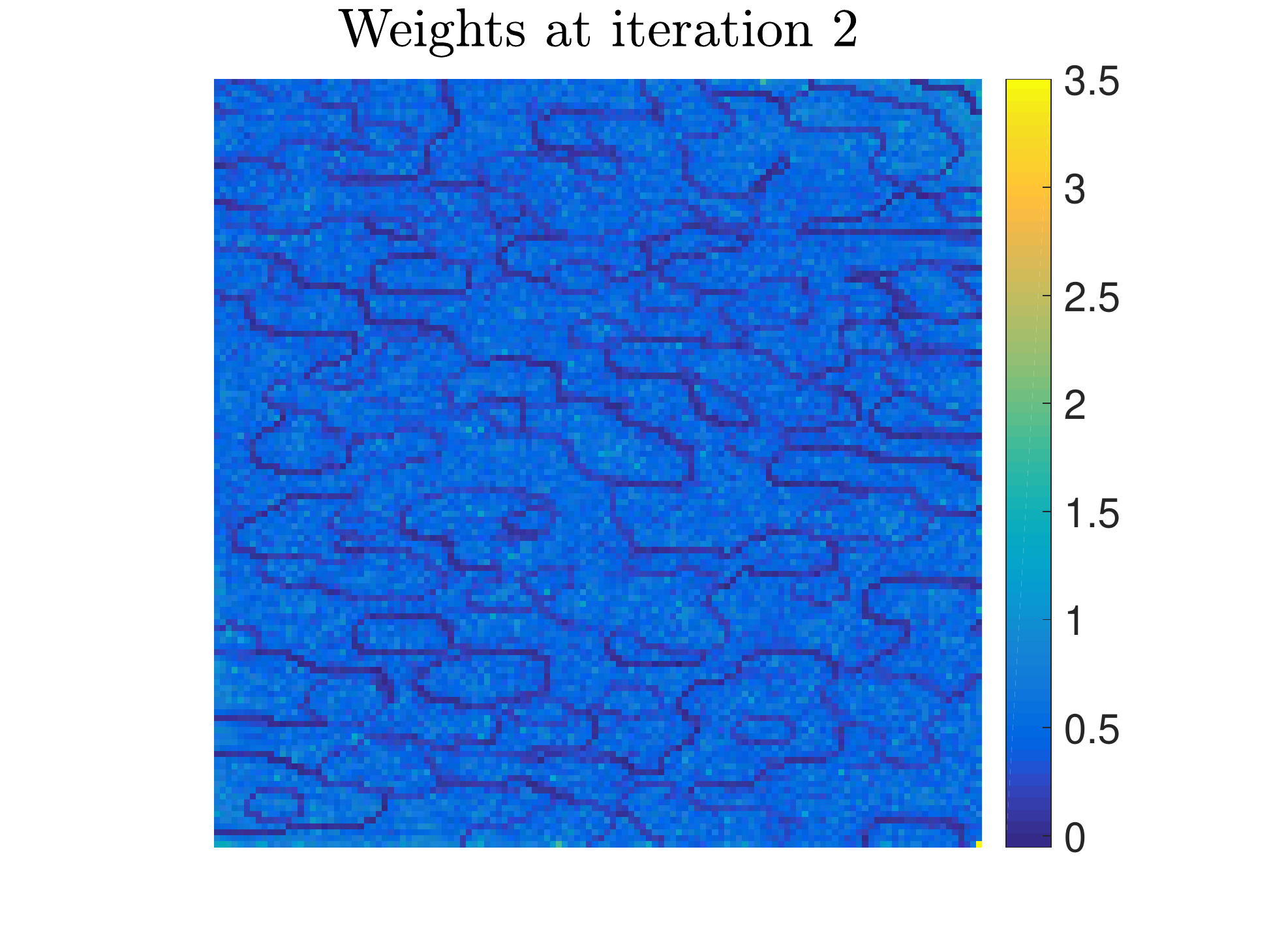}\\ 
\includegraphics[width=5cm]{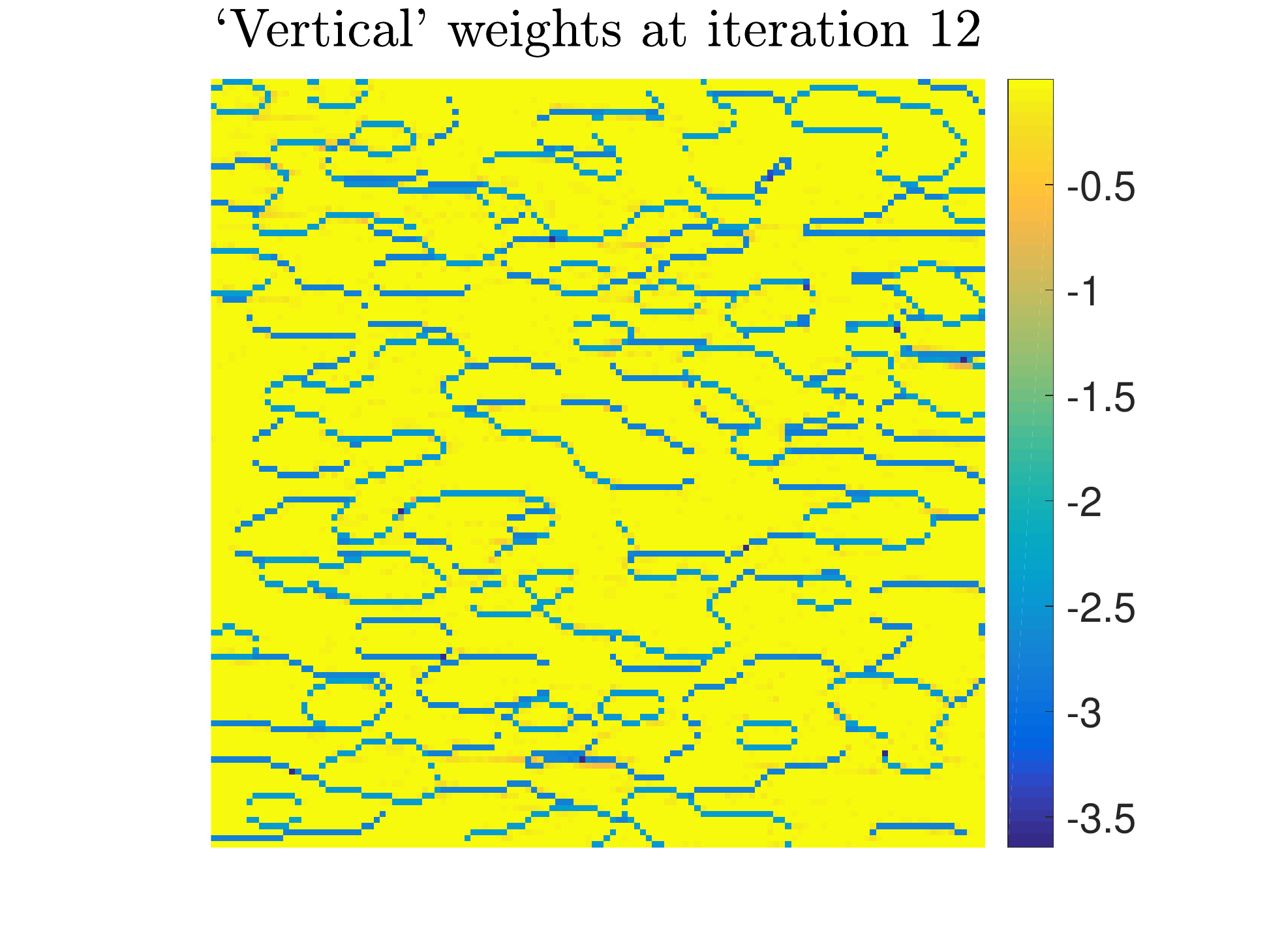} &
\includegraphics[width=5cm]{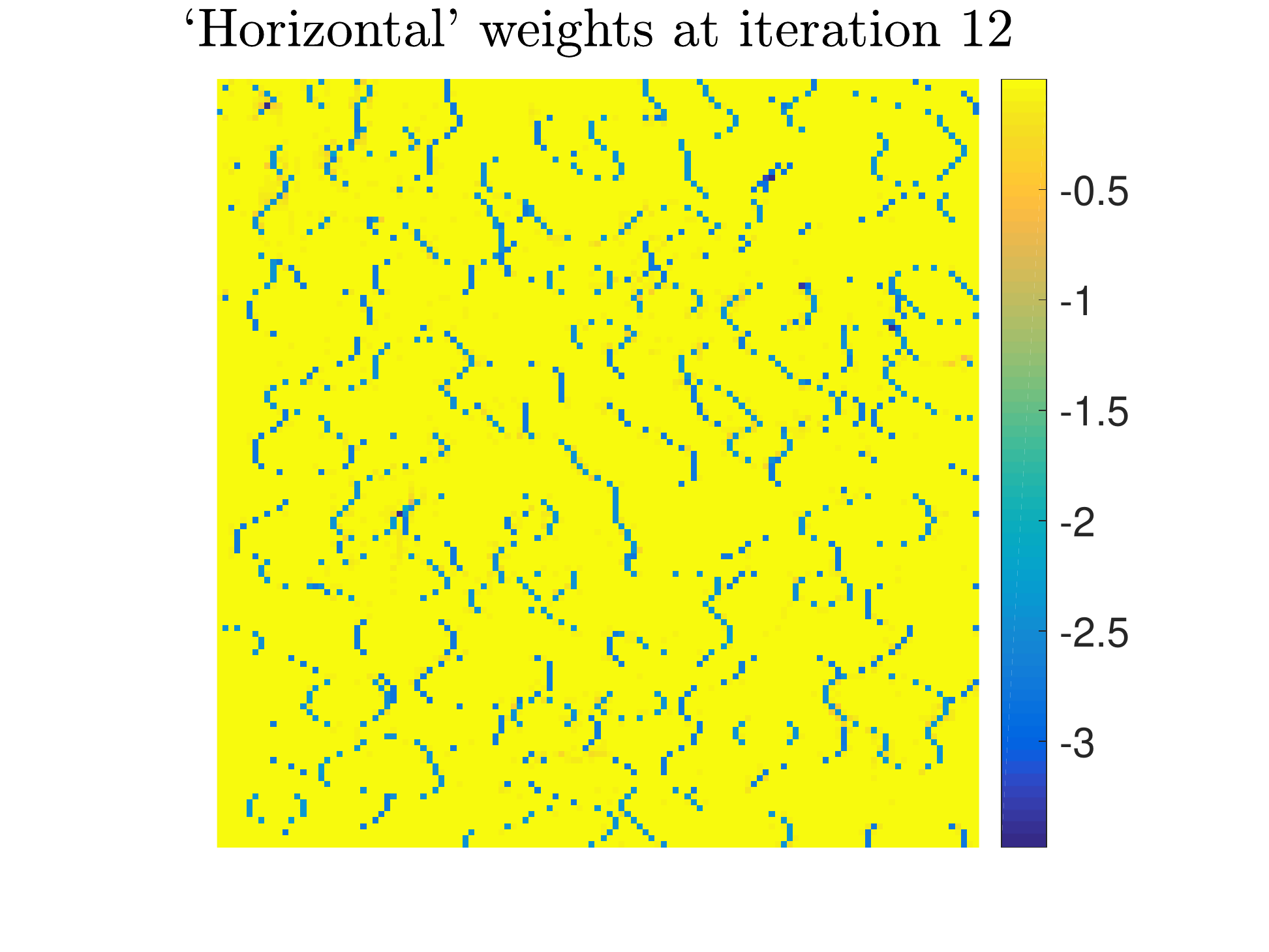} & 
\includegraphics[width=5cm]{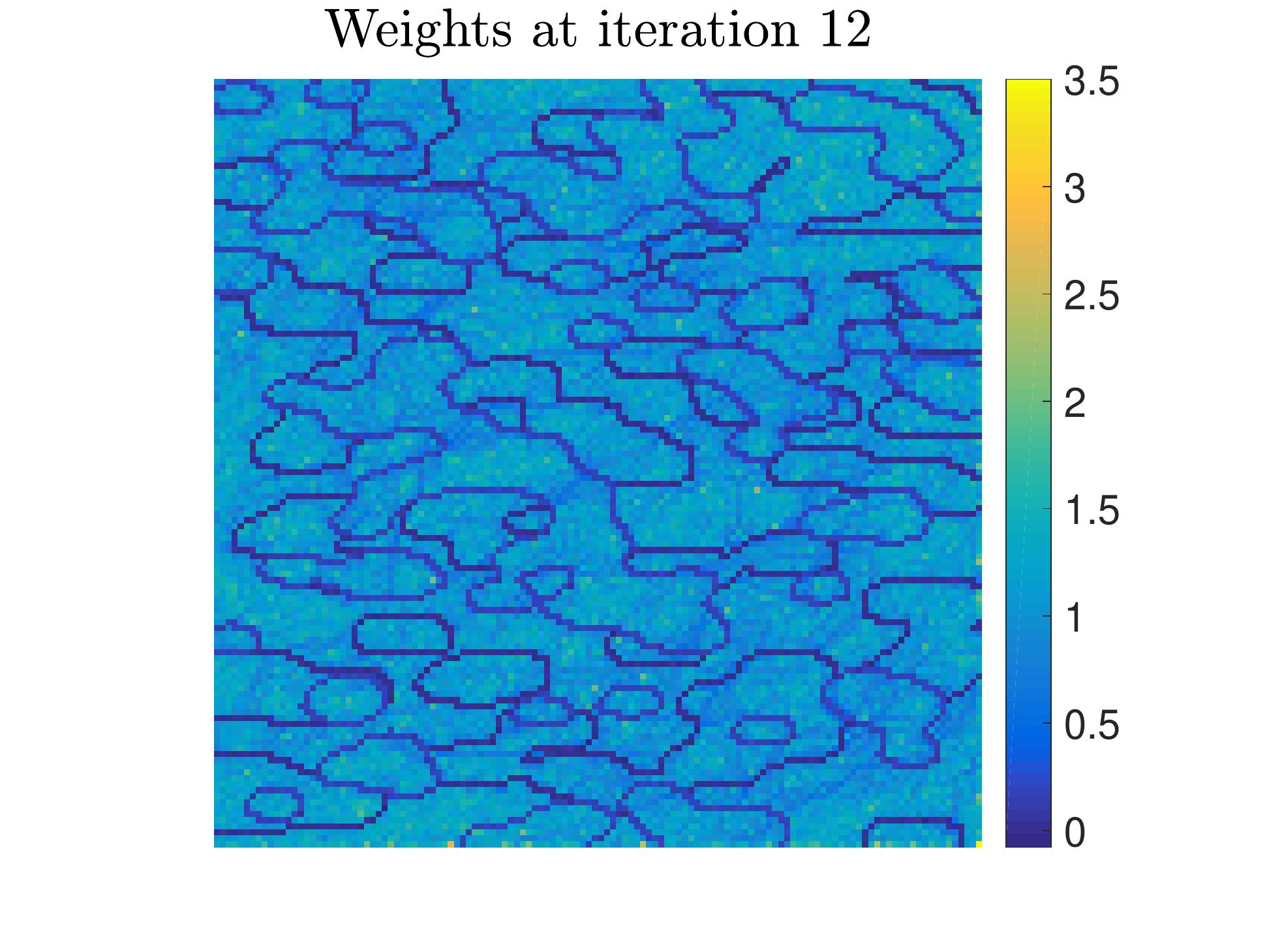}\\ \end{tabular}
\end{center}
\caption{\emph{Parallel-beam CT with limited angle} test problem. Left column: new weights, to be applied to the vertical derivatives. Middle column: new weights, to be applied to the horizontal derivatives. Right column: IRN-TV weights. All the values are displayed in logarithmic scale.}
\label{fig:Radon2TVweight}
\end{figure}

Finally, Figure \ref{fig:Radon1bSolutions} displays some relevant reconstructions. We show the initial reconstructions $x^{(\ast,1)}$ obtained at the end of the first inner iteration cycle, when a Tikhonov-regularized problem with regularization term $R(x) = \|Lx\|_2^2$ is used, and the reconstructions $x^{(\ast,12)}$. We consider the inner-outer iterative solvers that employ both the new and the IRN-TV weights, and both the hybrid method (with discrepancy principle) and CGLS. We can clearly see that, except for when the IRN-TV weights are used together with CGLS, there is a 
significant improvement in the reconstructions as the outer iterations proceed, and in particular 
the edges at the final outer iteration are much sharper than at the initial outer iteration. Interestingly enough, when the new weights are used together with CGLS, the edges and piecewise features of the solution seem fully recovered in the final reconstruction, but the pixel values are different from the original ones. 
%
%\begin{figure}[htbp]
%\begin{center}
%\includegraphics[width=7cm]{FigsExampleRadon1b/LCurves1b} 
%\end{center}
%\caption{${\mathcal L}$-curves for each iteration, for the second test problem. As implied from
%the text in the plot, the top curve corresponds to the first outer iteration, $\ell = 1$, and the curves
%below this correspond sequentially to iterations $\ell = 2, 3, 4$.  The red circles
%denote corners of each ${\mathcal L}$-curve, which correspond to the chosen
%regularization parameter, $\lambda_{*,\ell}$ for the particular iteration.}
%\label{fig:Radon1bLCurves}
%\end{figure}
%
%Computed reconstructions for the first outer iteration (that is, $x^{(*,1)}$), and for the
%final outer iteration (that is, $x^{(*,4)}$) are shown 
%in Figure~\ref{fig:Radon1bSolutions}. As we can see from these plots, there is a 
%significant improvement in the reconstructions, and in particular 
%the edges at the final outer iteration are much sharper than in the initial outer iteration.
%
Figure \ref{fig:Radon2TVweight} displays the entries of the weight diagonal matrices at outer iteration $\ell=2$ (i.e., as soon as reweighting becomes effective) and at iteration $\ell=12$, for both the new reweighting strategy (\ref{eq:GradientMap}) and the IRN-TV reweighting strategy; when considering the new reweighting strategy, both the weights applied to the vertical and horizontal derivatives are displayed. Edges and piecewise constant portions of the image to be recovered remarkably show up in the new weights: although these features are not so evident when $\ell=2$ (indeed, they are defined with respect to the reconstruction $x^{(\ast,1)}$ that is still very inaccurate, see Figure \ref{fig:Radon1bSolutions}), they can be recovered as the outer iterations progress. The IRN-TV weights are not so effective in revealing the structure of the image: although the edges seem to be better recovered when comparing the weights used at outer iteration $\ell=2$ and $\ell=12$, the piecewise constant parts are mainly weighted by small weights, too, and therefore are not so effectively smoothed in the reconstruction process.

%\newpage

%\subsection{Example from X-ray CT, More Limited Angle}
\paragraph{Parallel-beam CT with more limited angle}

\begin{figure}[htbp]
\begin{center}
\begin{tabular}{cc}
%\footnotesize{True image} & \footnotesize{Measured data (sinogram)}\\
\includegraphics[height=4cm]{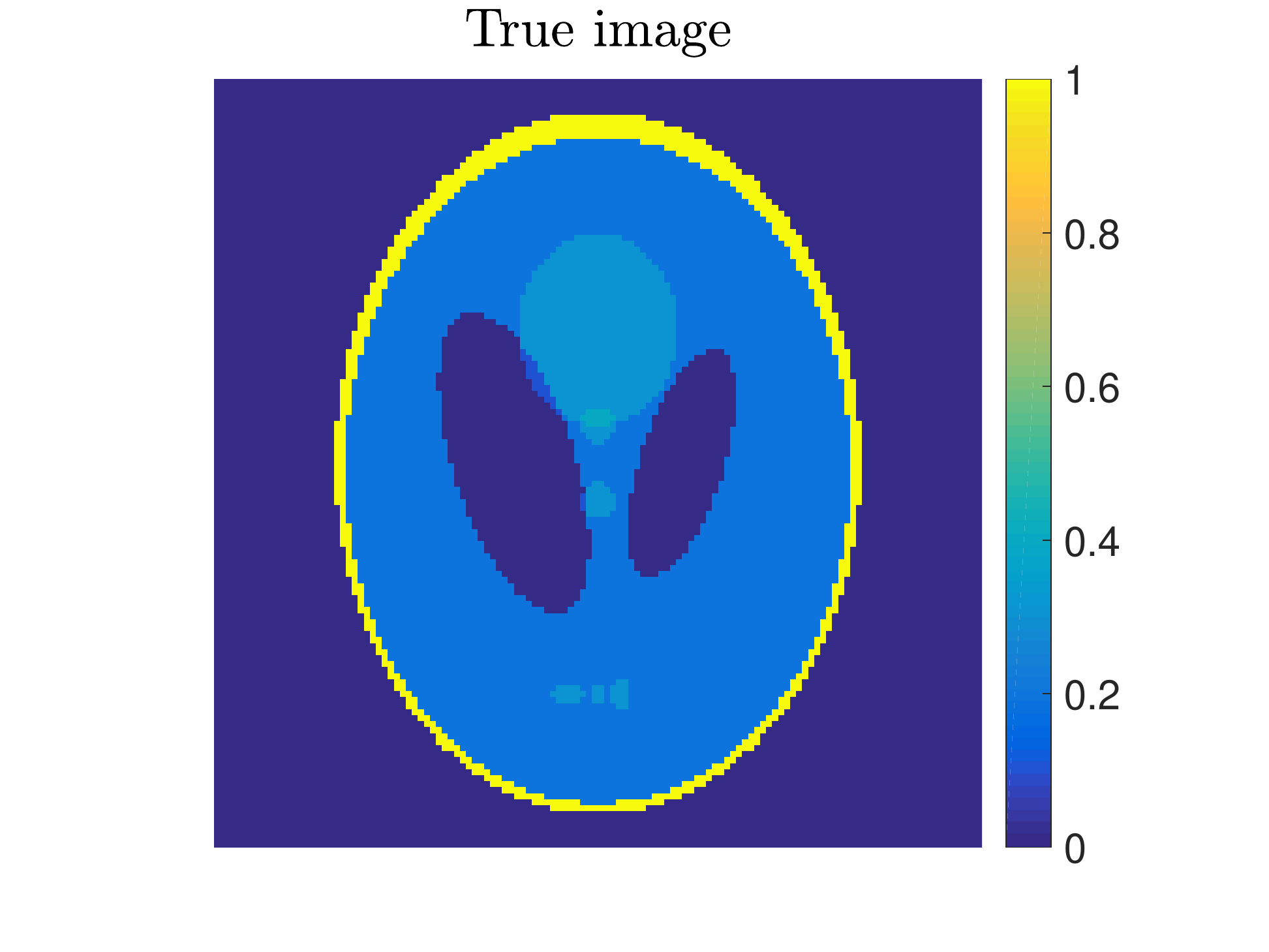} &
\includegraphics[height=4cm]{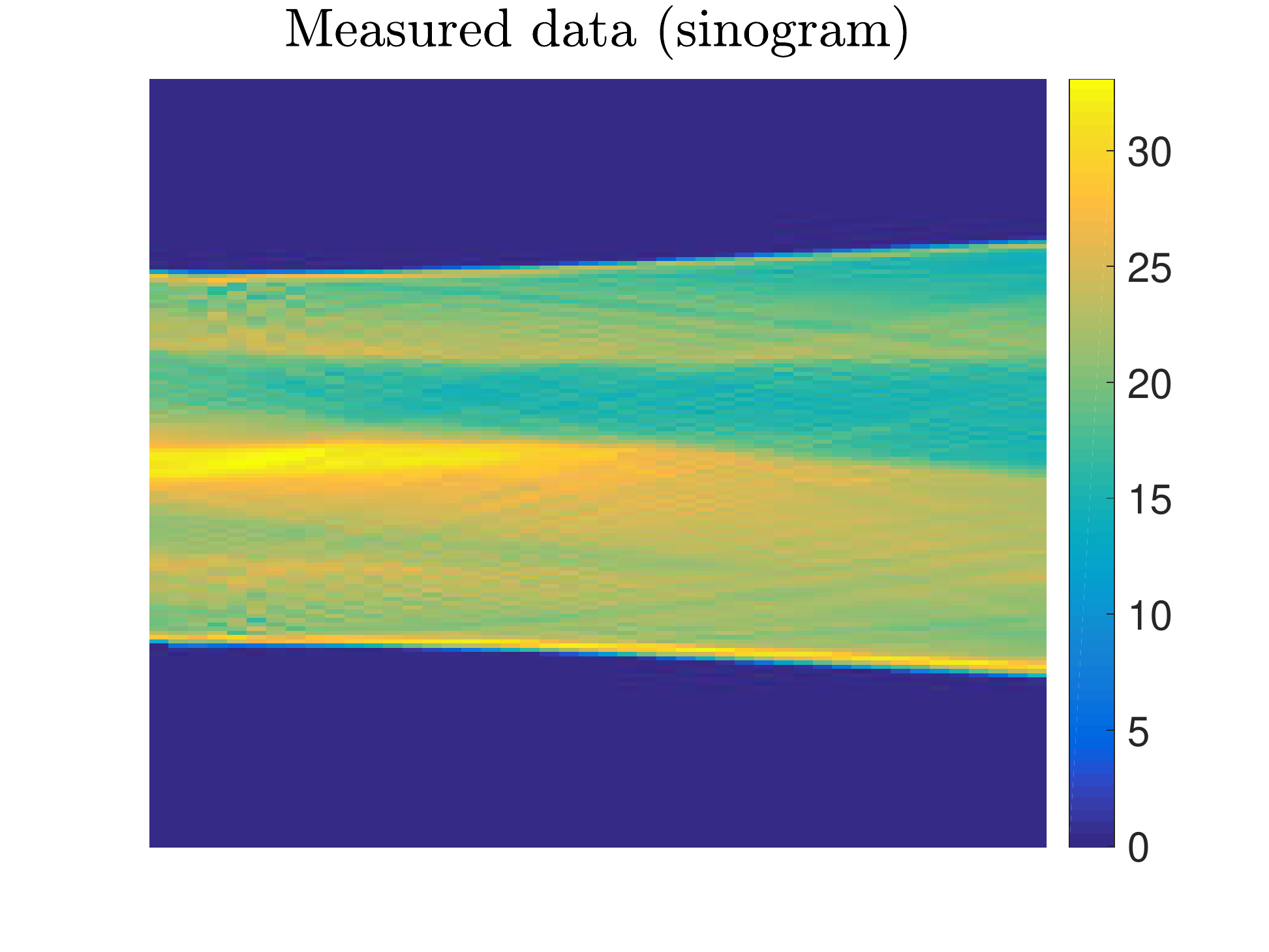}
\end{tabular}
\end{center}
\caption{\emph{Parallel-beam CT with more limited angle} test problem. Left frame: true image phantom $x$. Right frame: the measured data (sinogram) $b$, where the projections are taken only at angles $0, 1, \ldots, 45$.}
\label{fig:Radon1cData}
\end{figure}

We use IR Tools to vary some of the options defining the previous X-ray tomography simulations. We still wish to consider a parallel beam X-ray transmission problem, but we would like to change the phantom image to the well-known Shepp-Logan one, of size $128\times 128$ pixels. Moreover, we would like to make the reconstruction problem even more challenging by using even more limited projection angles $0, 1, \ldots, 45$. Again, we generate this test problem within IR Tools, using instructions similar to the previous ones, where a new \texttt{ProblemOptions} structure is defined in the following way:
\begin{verbatim}
     ProblemOptions = PRset('angles', 0:45, 'phantomImage', 'shepplogan');
\end{verbatim}
%This is a slightly
%more limited angle case than the previous example, making the problem more
%ill-posed, and more difficult to obtain good reconstructed images.
%
%Using IR Tools, the problem can be generated with the following MATLAB statements:
%\begin{verbatim}
%     ProblemOptions = PRset('angles', 0:90, 'phantomImage', 'grains');
%     [A, b_true, x_true, ProblemInfo] = PRtomo(128, ProblemOptions);
%     b = PRnoise(b_true, 1e-3);
%\end{verbatim}
%
Figure~\ref{fig:Radon1cData} shows the true phantom image, along with the
measured data $b$ (which are again corrupted by Gaussian white noise of level $10^{-3}$). 

As in the previous examples, we first run our algorithm with different parameter choice strategies within the inner hybrid scheme for generalized Tikhonov regularization: when the discrepancy principle is used, the method performs the maximum allowed number of outer iterations $\ell=20$; when  the $\mathcal{L}$-curve criterion is used, the stopping rule of Section \ref{ssec:mainalgo} (decrease of $\|Lx^{(\ast,\ell)}\|_2$) is satisfied after only $\ell=7$ outer iterations. Figure~\ref{fig:Radon1bIterations} shows a plot of the relative errors
and chosen regularization parameters at each outer iteration. We can clearly see that, when considering both the discrepancy principle and the $\mathcal{L}$-curve criterion, the reconstructions greatly improve as the outer iterations proceed, with the latter strategy being very efficient; however, the discrepancy principle eventually achieves a lower relative error (although the $\mathcal{L}$-curve criterion may still achieve the same relative error, if not that the stopping criterion is satisfied quite early). Note that, as expected, the regularization parameters increase with the outer iterations.
%We can also see that the regularization parameter selected at the end of each inner iteration cycle consistently increases when both parameter rules are employed: this behavior is expected and, in . 
%he terminated after 4 outer iterations, where at iteration $\ell$ the hybrid 
%method determined an estimate of the regularization parameter $\lambda_{*,\ell}$,
%and a corresponding reconstructed image, $x^{(*,\ell)}$.  
%Figure~\ref{fig:Radon1bLCurves}
%shows a plot of all ${\mathcal L}$-curves for each iteration. Observe the nesting
%property of the ${\mathcal L}$-curves, and 
%how the chosen regularization
%parameter (corresponding to the corner of each ${\mathcal L}$-curve) increases.
%
%
%
\begin{figure}[htbp]
\begin{center}
\begin{tabular}{cc}
\includegraphics[width=6cm]{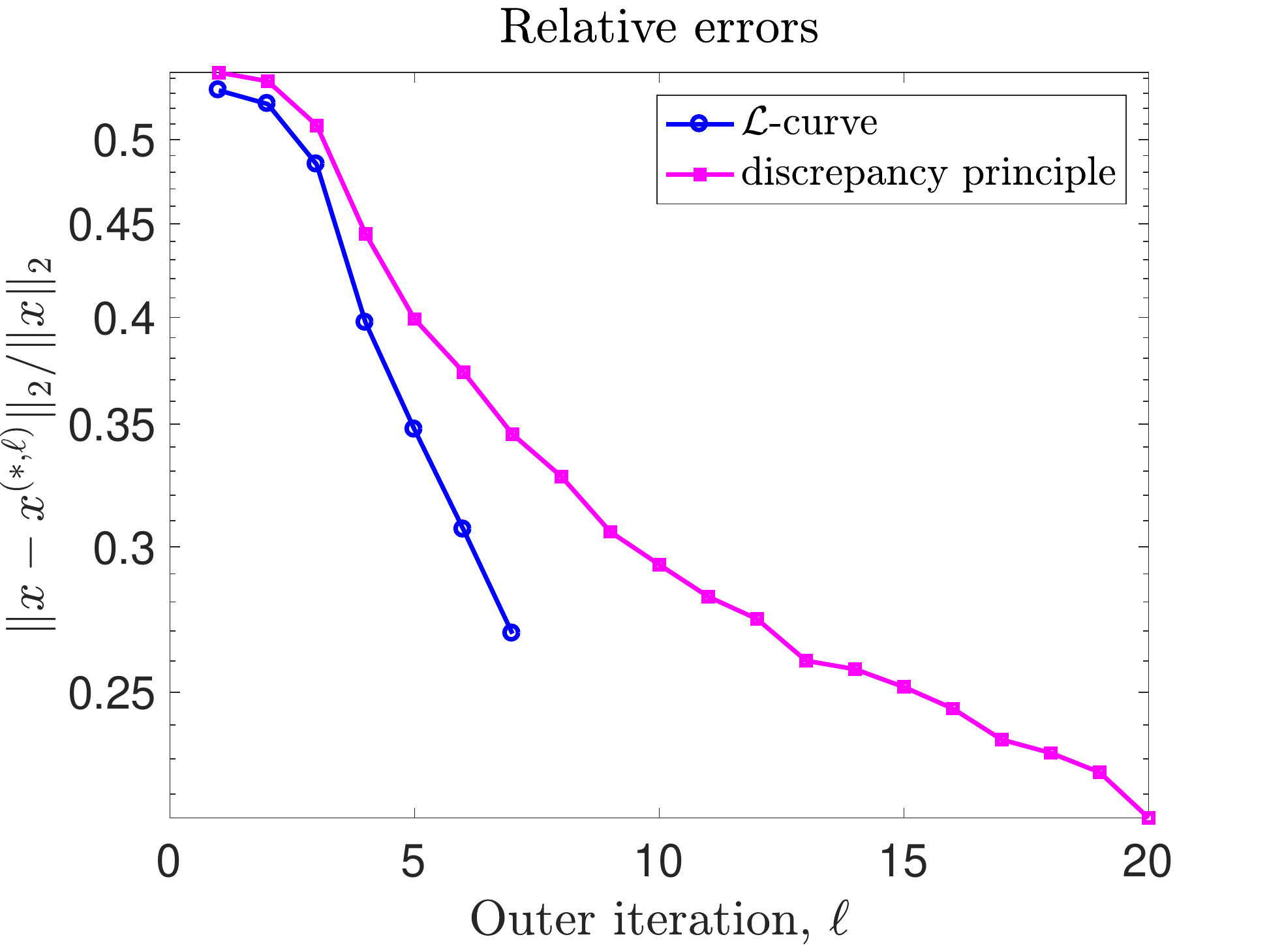} &
\includegraphics[width=6cm]{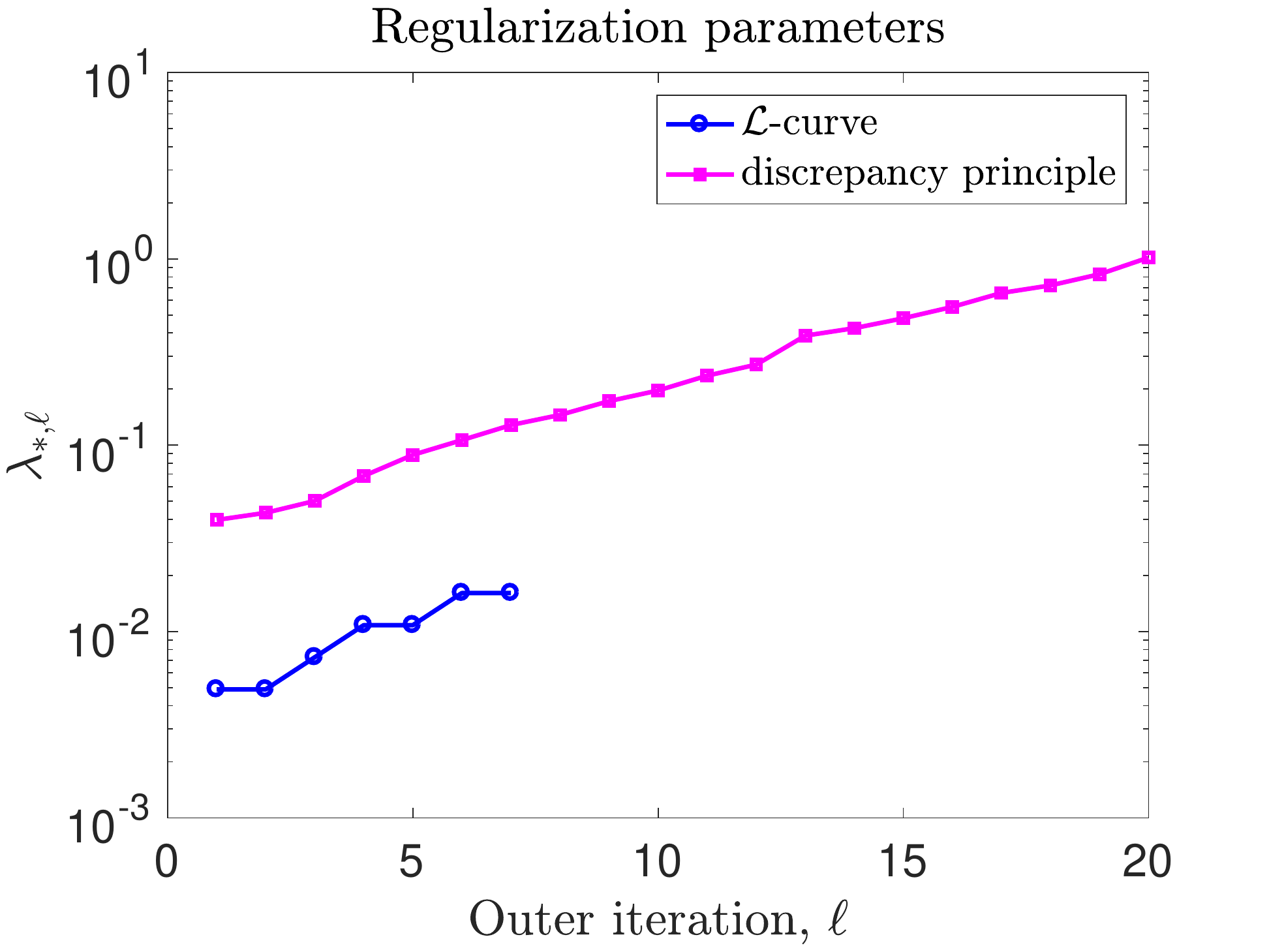}
\end{tabular}
\end{center}
\caption{\emph{Parallel-beam CT with more limited angle} test problem. Relative errors and regularization parameters values at each outer iteration $\ell$,until the stopping criterion is satisfied. Both the discrepancy principle and the ${\mathcal L}$-curve criterion are considered.}
\label{fig:Radon1cIterations}
\end{figure}

\begin{figure}[htbp]
\begin{center}
\begin{tabular}{cc}
\includegraphics[width=6cm]{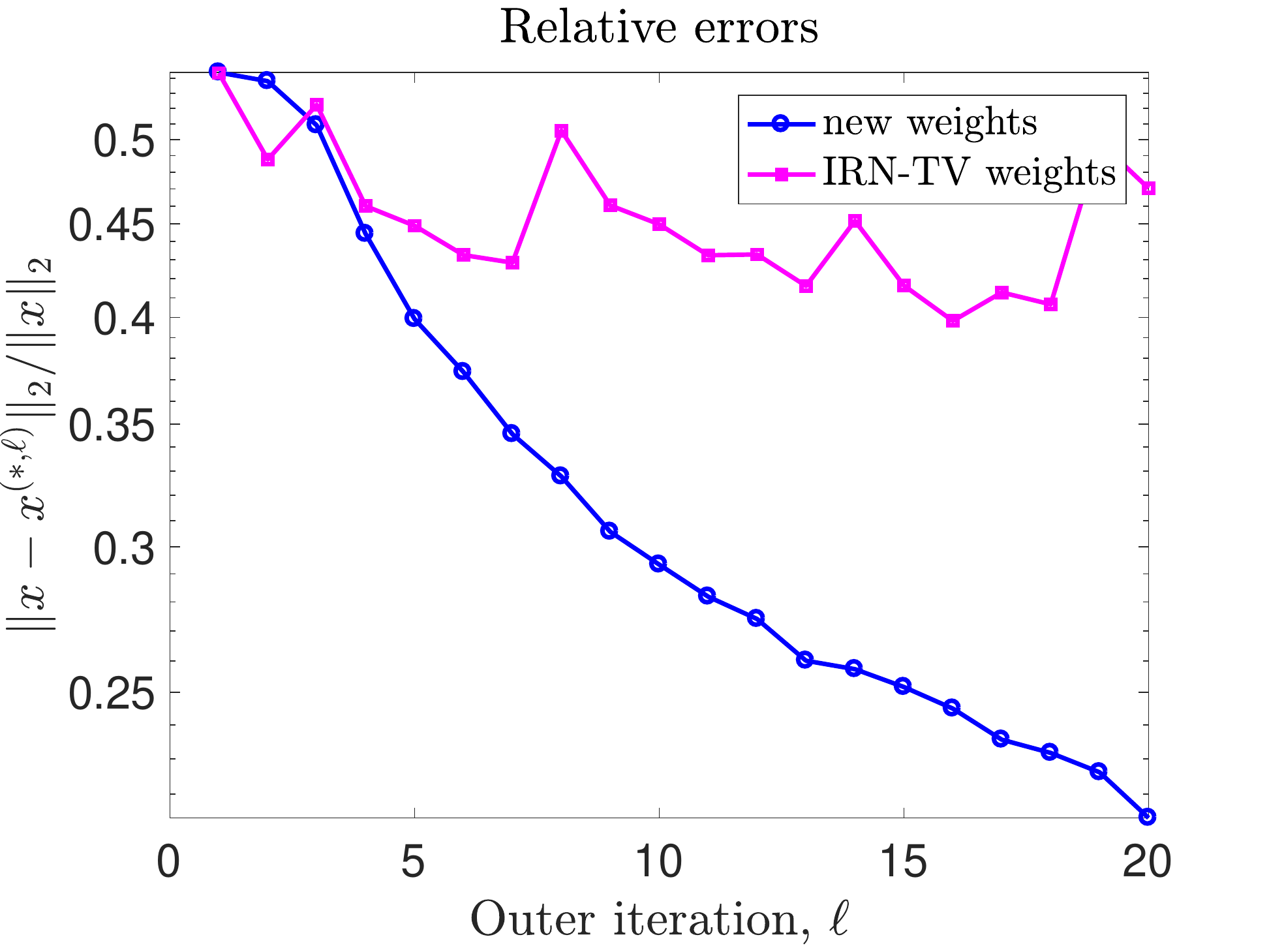} &
\includegraphics[width=6cm]{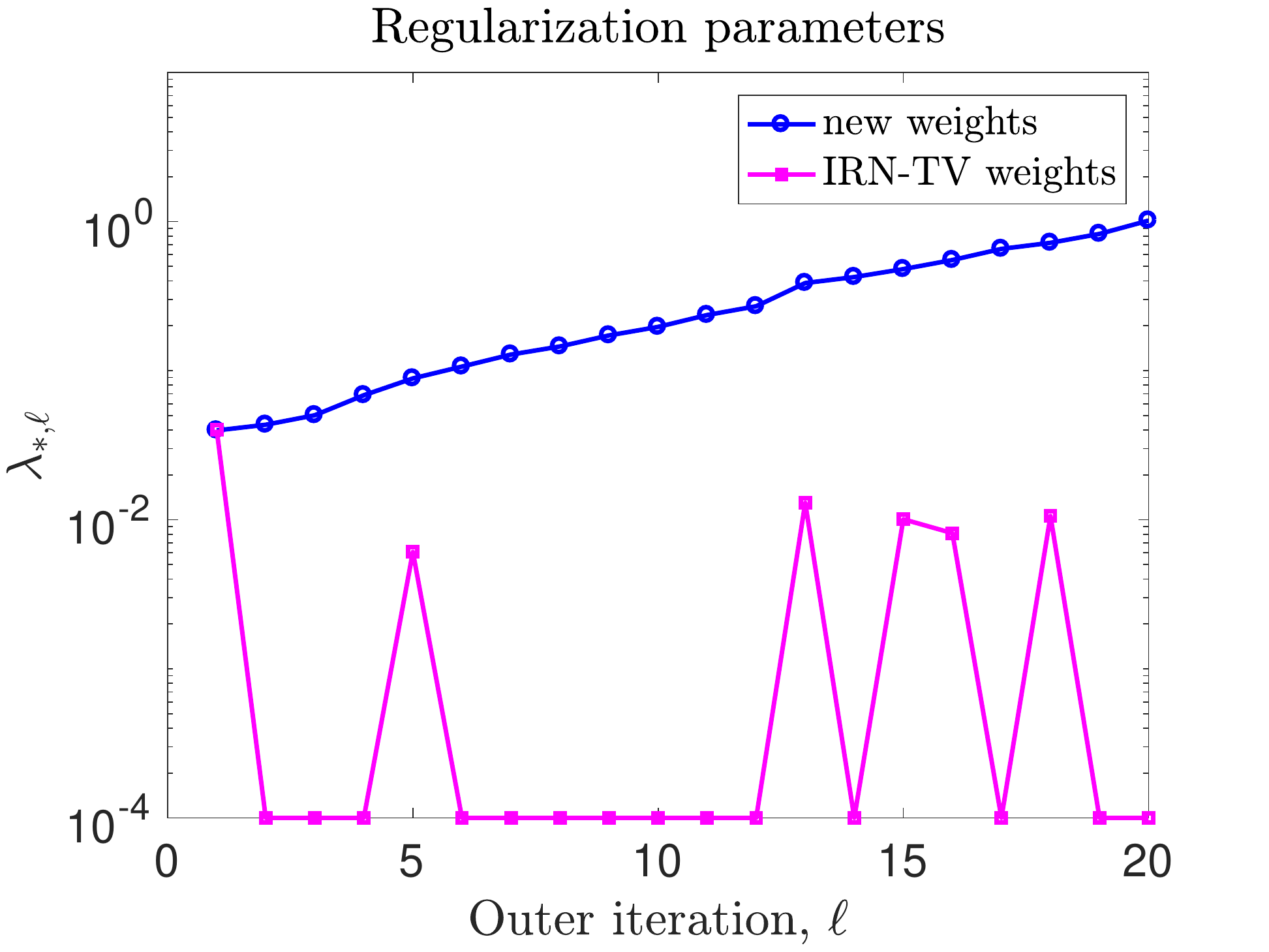}
\end{tabular}
\end{center}
\caption{\emph{Parallel-beam CT with more limited angle} test problem. Relative errors and regularization parameter versus number of (outer) iterations; both methods use the hybrid method for general form Tikhonov during the inner iterations, and adaptively select the regularization parameter according to the discrepancy principle.}
\label{fig:Radon3TV}
\end{figure}

\begin{figure}[htbp]
\begin{center}
\begin{tabular}{cc}
\includegraphics[width=6cm]{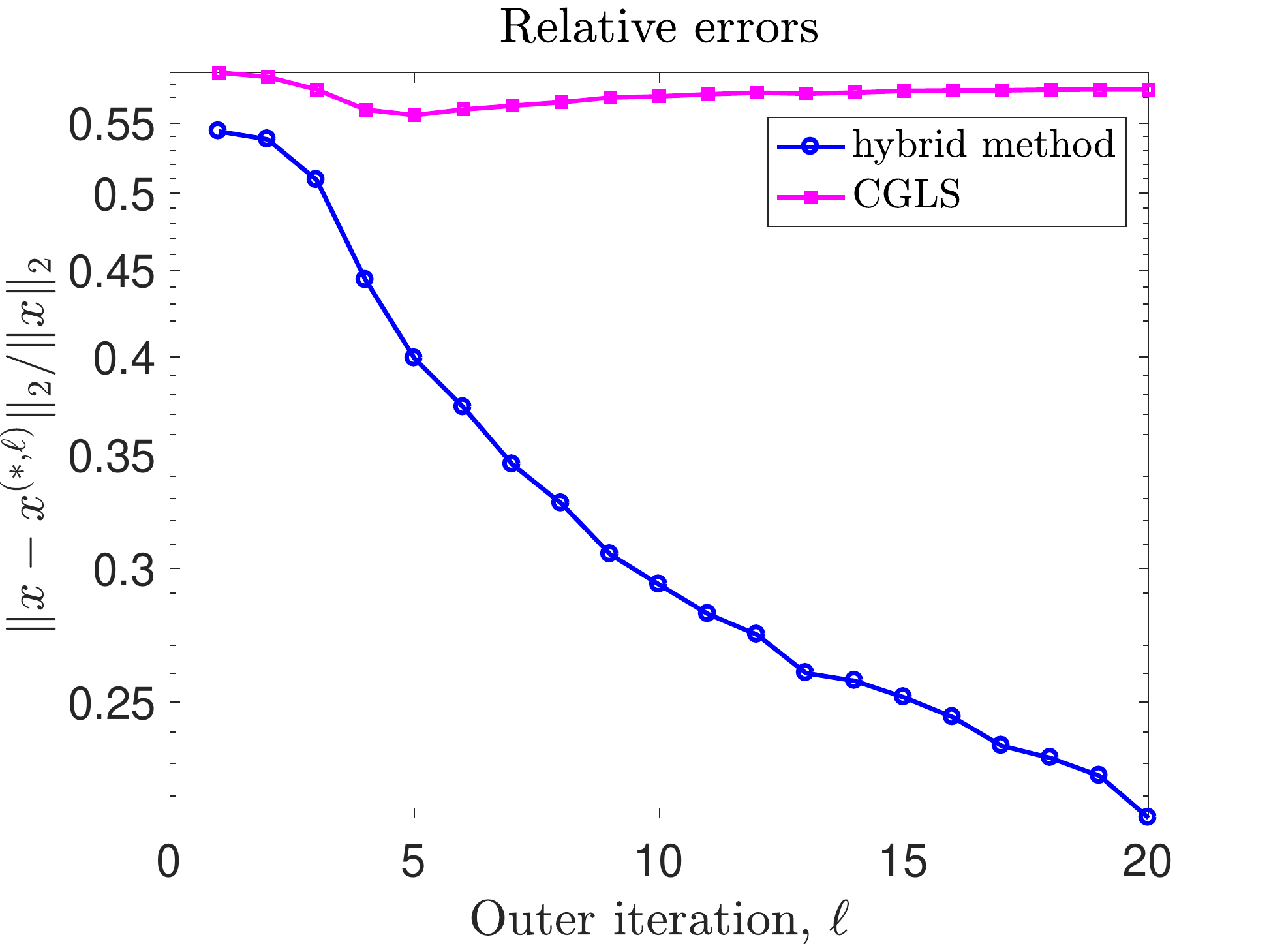} &
\includegraphics[width=6cm]{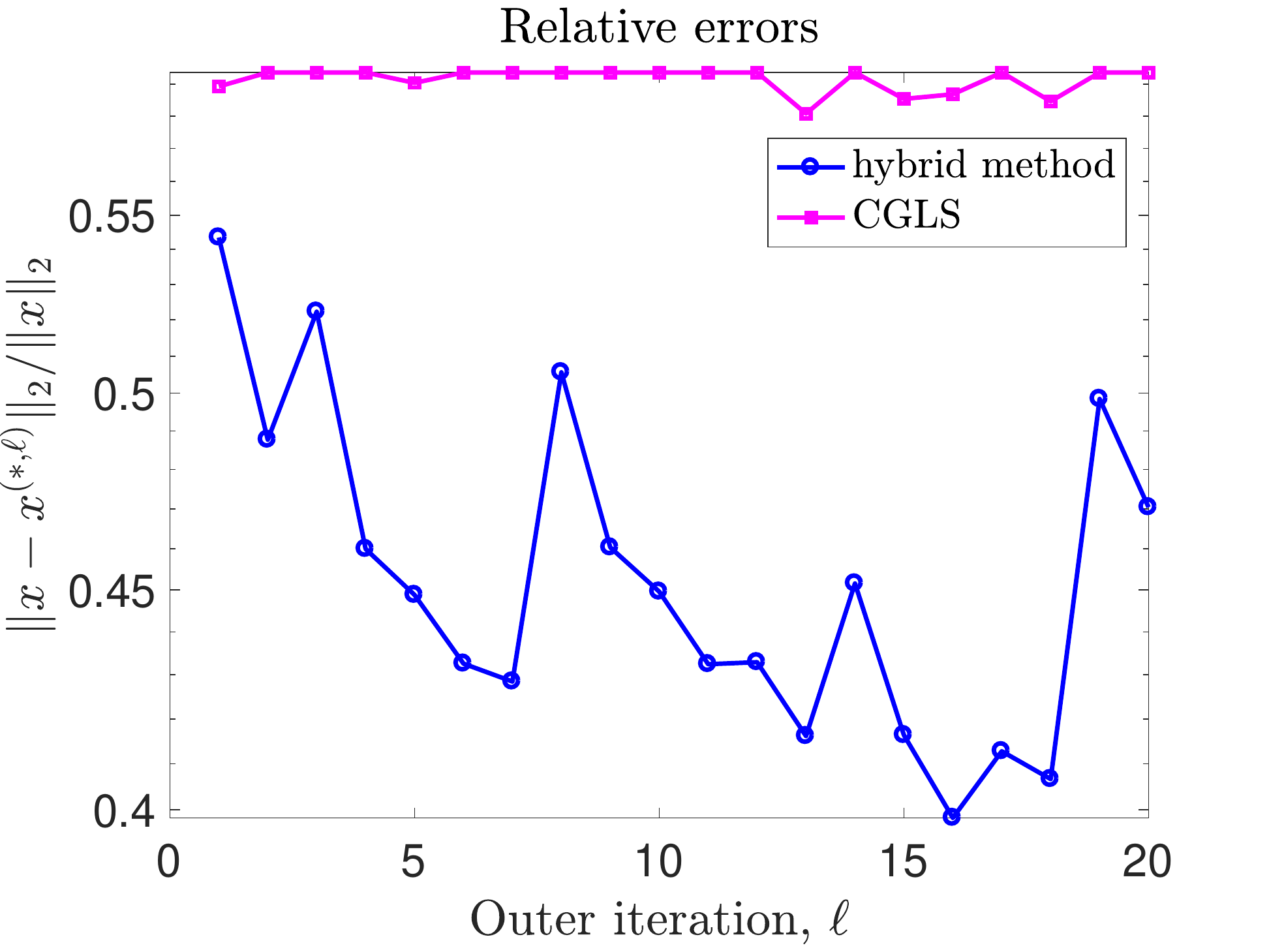}
\end{tabular}
\end{center}
\caption{\emph{Parallel-beam CT with more limited angle} test problem. Relative errors versus number of (outer) iterations. Left frame: the new weights are used. Right frame: the IRN-TV weights are used.}
\label{fig:Radon3TV_comparisons}
\end{figure}

Next, we compare the new method to other inner-outer iterative methods for edge enhancement in imaging. Figure \ref{fig:Radon3TV} displays the relative error and regularization parameter values versus the number of outer iterations for the new method, and for a method that still employs the hybrid solver for general form Tikhonov regularization to handle the inner iterations while udating the IRN-TV weights (\ref{eq:IRNTVweights}) at each outer iteration. For both methods, the discrepancy principle is employed to adaptively choose the regularization parameter at each inner iteration. 
We can clearly see that, when the IRN-TV weights are used, the behavior of the relative errors is quite hectic and, although the quality of the solution computed at the second outer iteration is good, this trend is not maintained during the following outer iterations. This may be because the automatically selected regularization parameter for IRN-TV keeps oscillating between values of the order of $10^{-2}$ and values of the order of $10^{-4}$. 
%%%; indeed, for this test problem, setting the regularization parameter through the $\mathcal{L}$-curve results in a more accurate reconstructions
%% Similarly to the previous test problems, both methods perform well during the first (outer) iterations, with the IRN-TV weights being quite effective in reducing the error; however, the quality of the IRN-TV solutions eventually stagnates, while the new method keeps improving. 
%%%Similarly to the previous test problems, both methods perform well during the first (outer) iterations, with the IRN-TV weights being quite effective in reducing the error; however, the quality of the IRN-TV solutions eventually stagnates, while the new method keeps improving. 
%%The automatically selected regularization parameter for IRN-TV keeps oscillating between values of the order of $10^{-2}$ and values of the order of $10^{-4}$. 
% computes  the new method is amazing at the end.
Figure \ref{fig:Radon3TV_comparisons} assesses the influence of the inner iterative solver on the overall behavior of the method. Namely, we consider the inner-outer iterative schemes implemented with both the new and the IRN-TV weights, and with both the hybrid and the CGLS methods as inner solvers. The hybrid method chooses the regularization parameter adaptively according to the discrepancy principle as the iterations proceed, and the value $\lambda_{\ast,\ell}$ selected when the $\ell$th inner iteration cycle terminates is taken to be the fixed regularization parameter to be set in advance of the $\ell$th CGLS iteration cycle. 
% requires a fixed value of the regularization parameter to be available in advance of each cycle of inner iterations, we first run the methods based on the hybrid solver for general form Tikhonov that adaptively chooses the regularization parameter at each inner iteration according to the discrepancy principle: in this way, a parameter 
% $\lambda_{\ast,\ell}$ is eventually set at the $\ell$th outer iteration. When running the methods based on CGLS, we take $\lambda_{\ast,\ell}$ as regularization parameter for the $\ell$th outer iteration.
In both the frames displayed in Figure \ref{fig:Radon3TV_comparisons} we can see that the results obtained using CGLS hardly improve when the outer iterations proceed, especially when the IRN-TV weights are used. Again, this can probably be avoided if the fixed regularization parameter is chosen differently. 

\begin{figure}[htbp]
\begin{center}
\begin{tabular}{ccc}
\hspace{-0.3cm}\includegraphics[width=5.5cm]{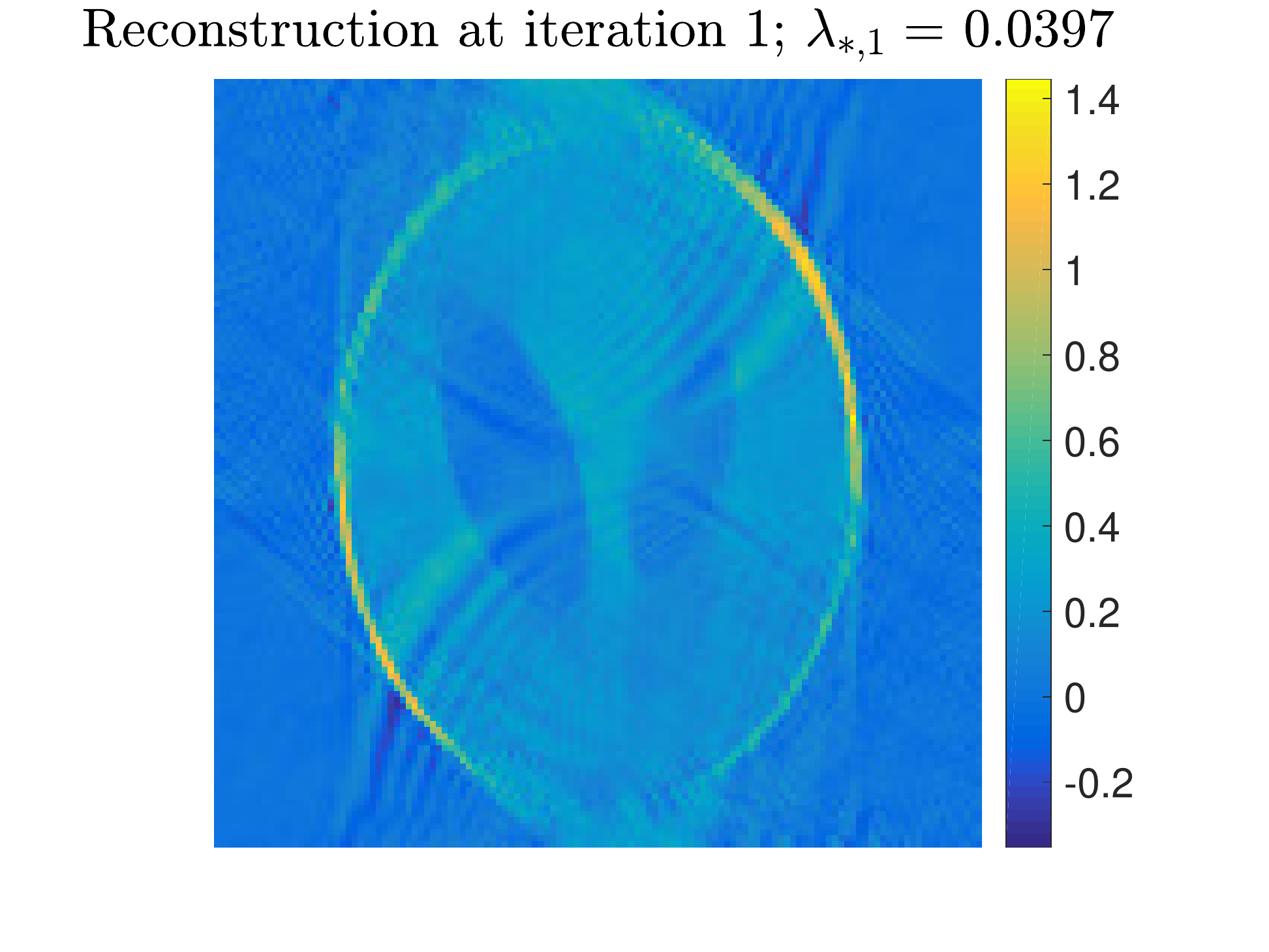} &
\hspace{-0.3cm}\includegraphics[width=5.5cm]{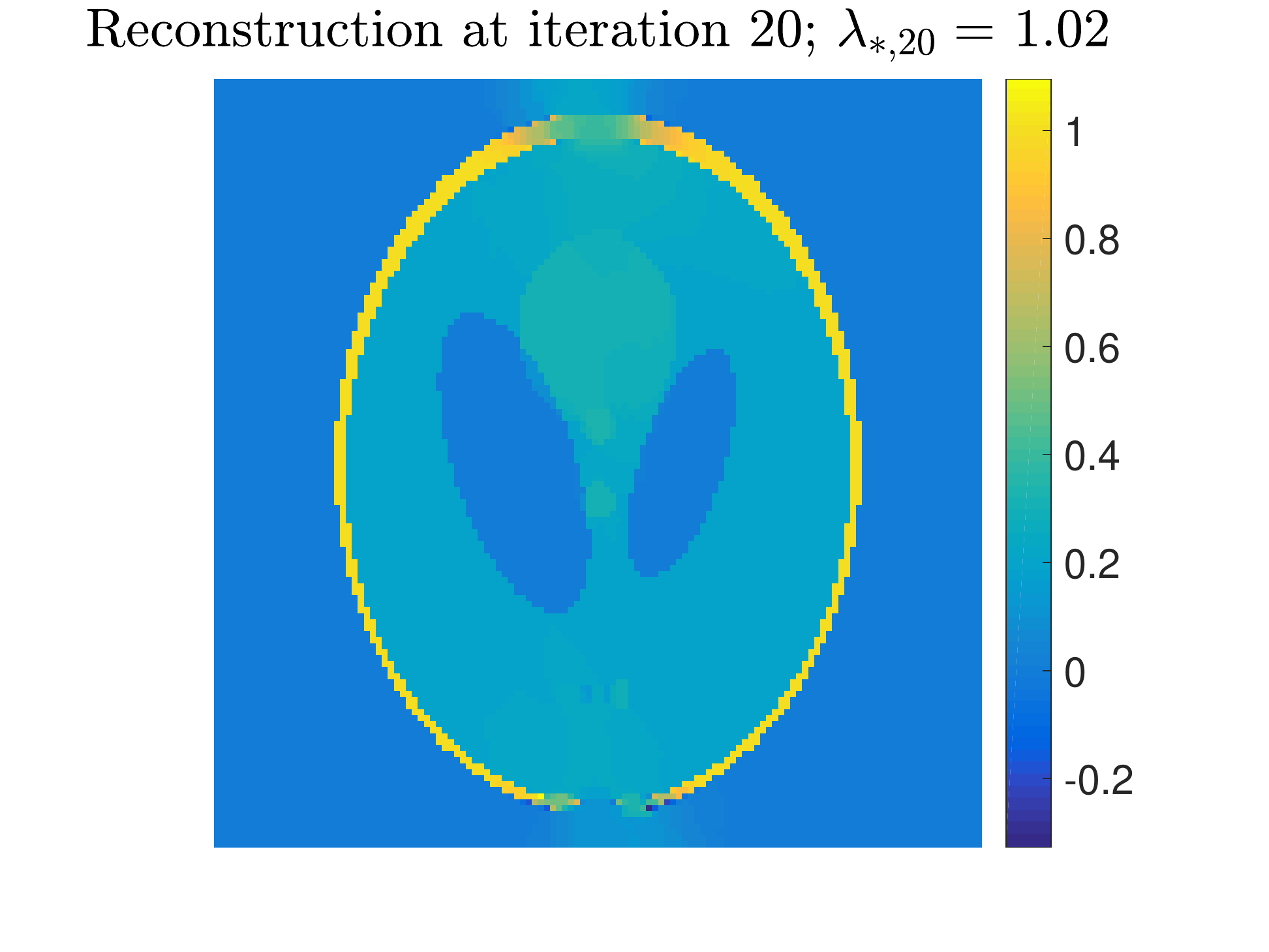} &
\hspace{-0.3cm}\includegraphics[width=5.5cm]{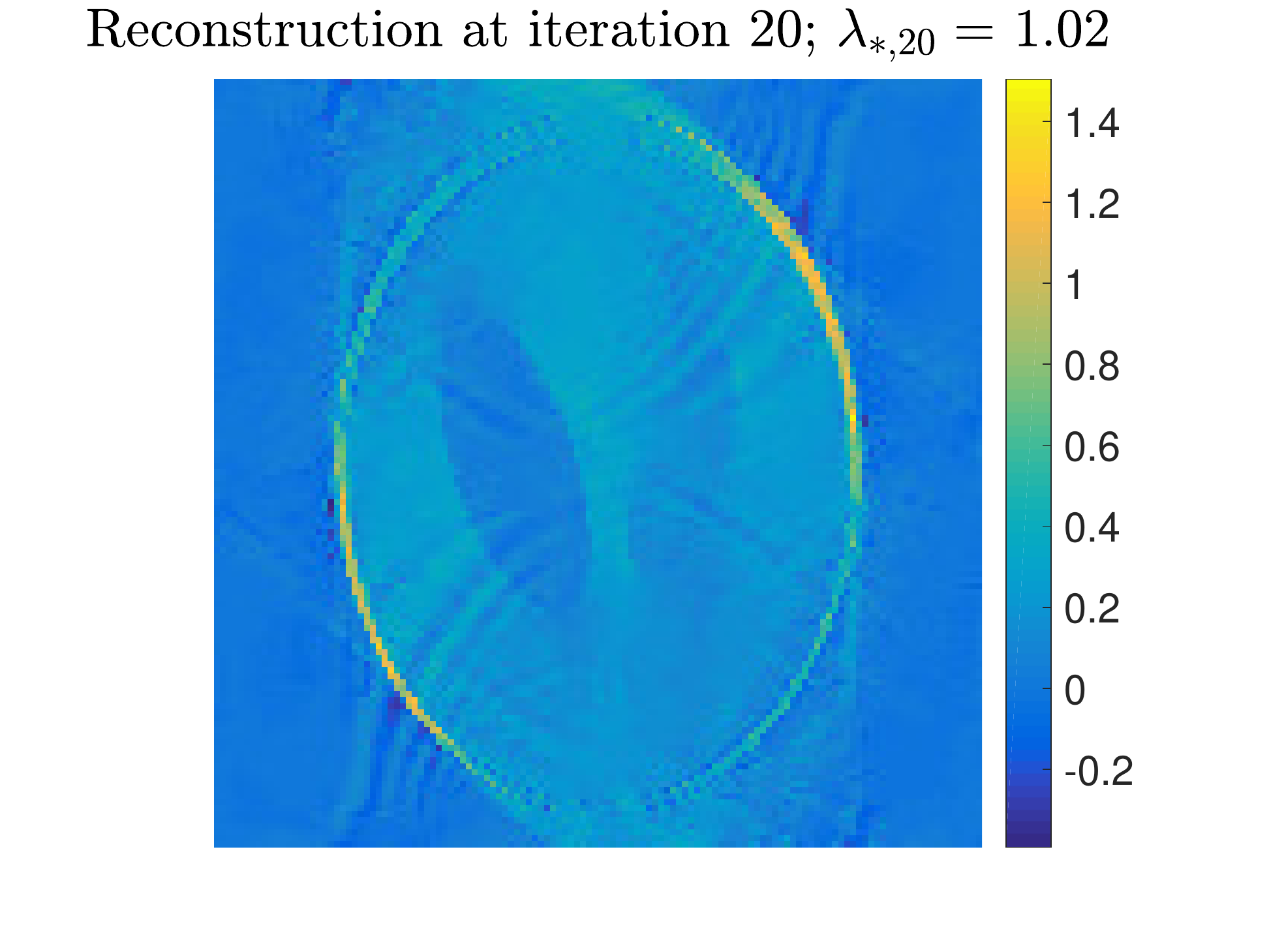}\\
\hspace{-0.3cm}\includegraphics[width=5.5cm]{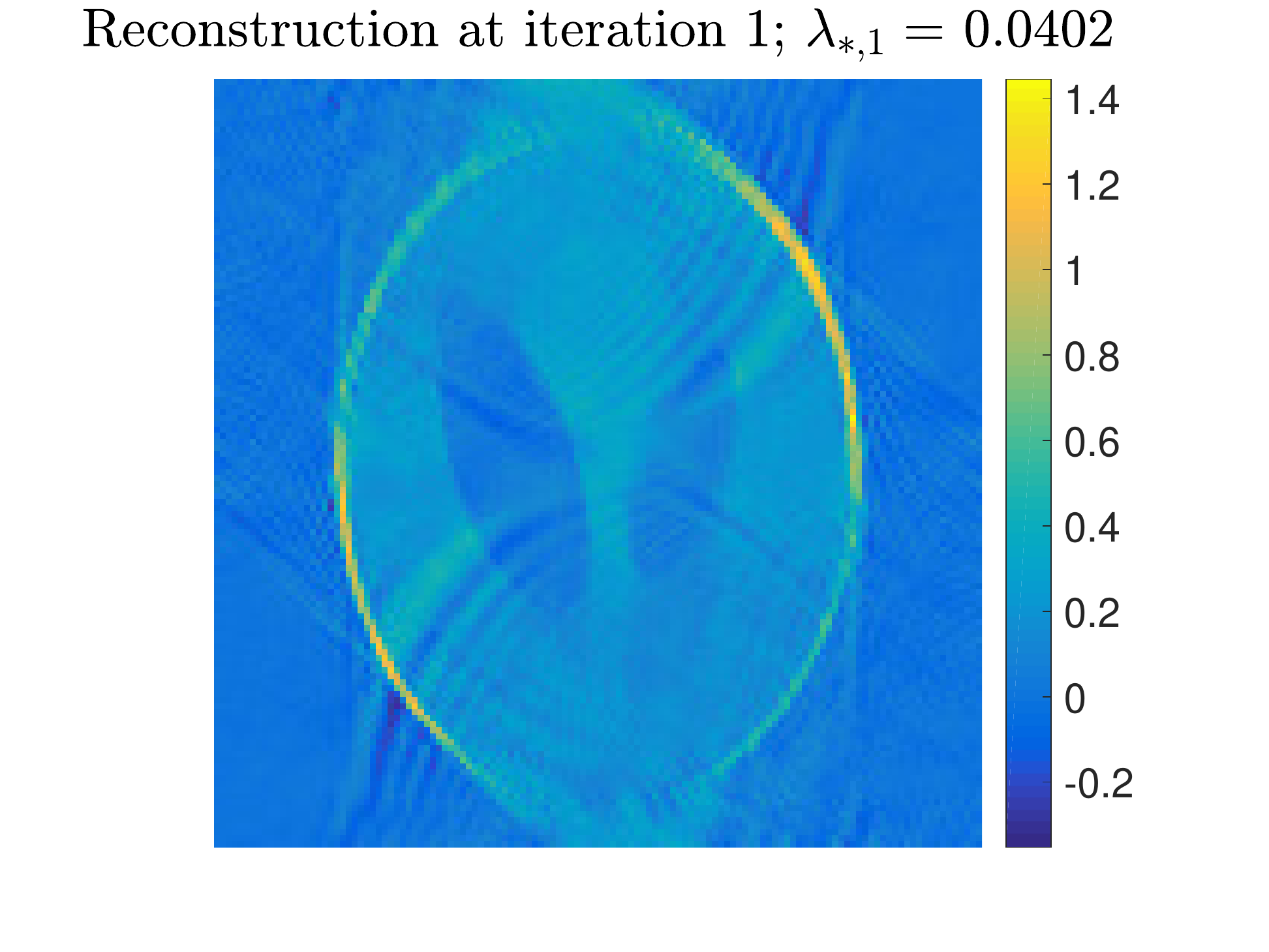} &
\hspace{-0.3cm}\includegraphics[width=5.5cm]{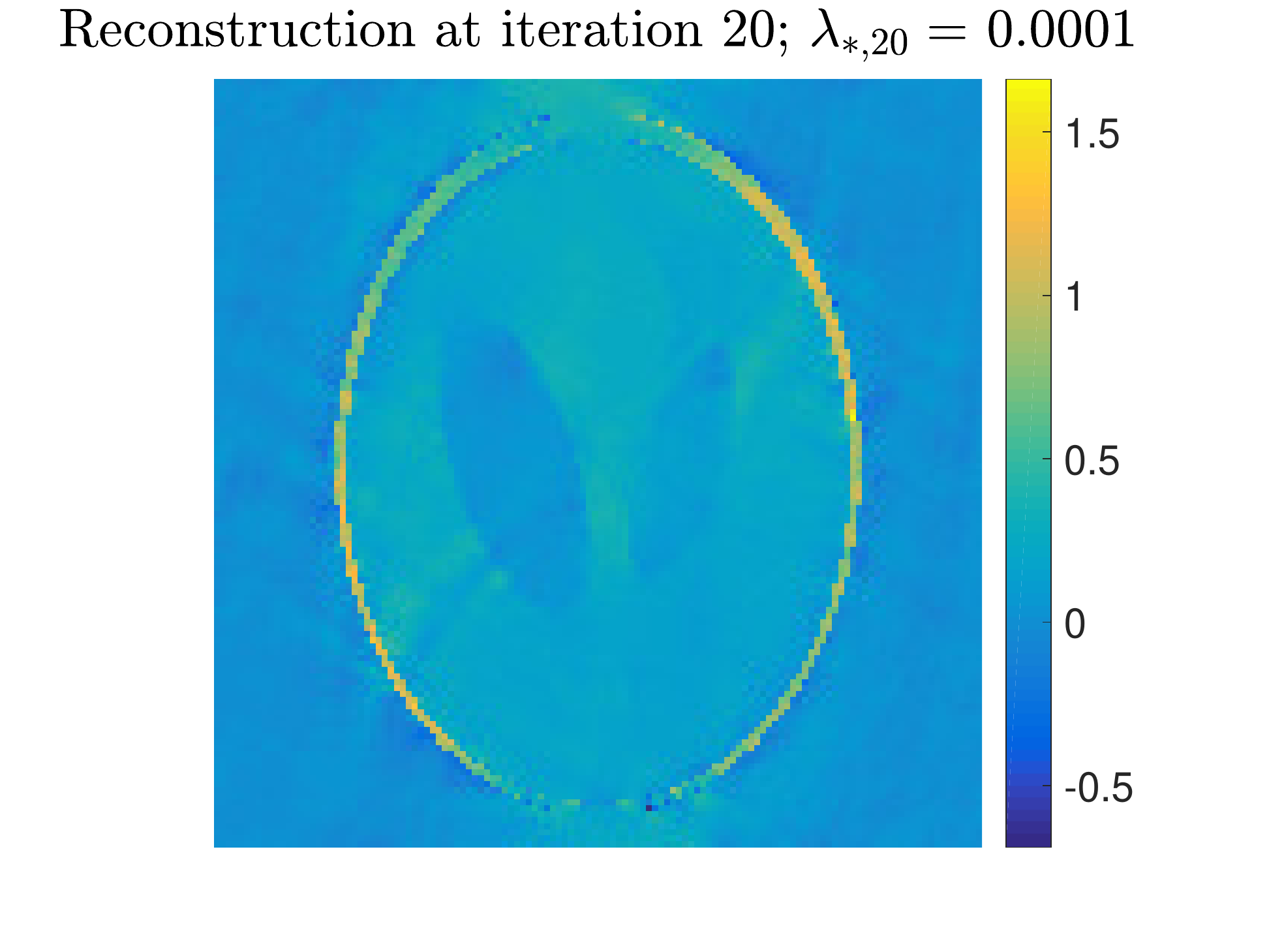} &
\hspace{-0.3cm}\includegraphics[width=5.5cm]{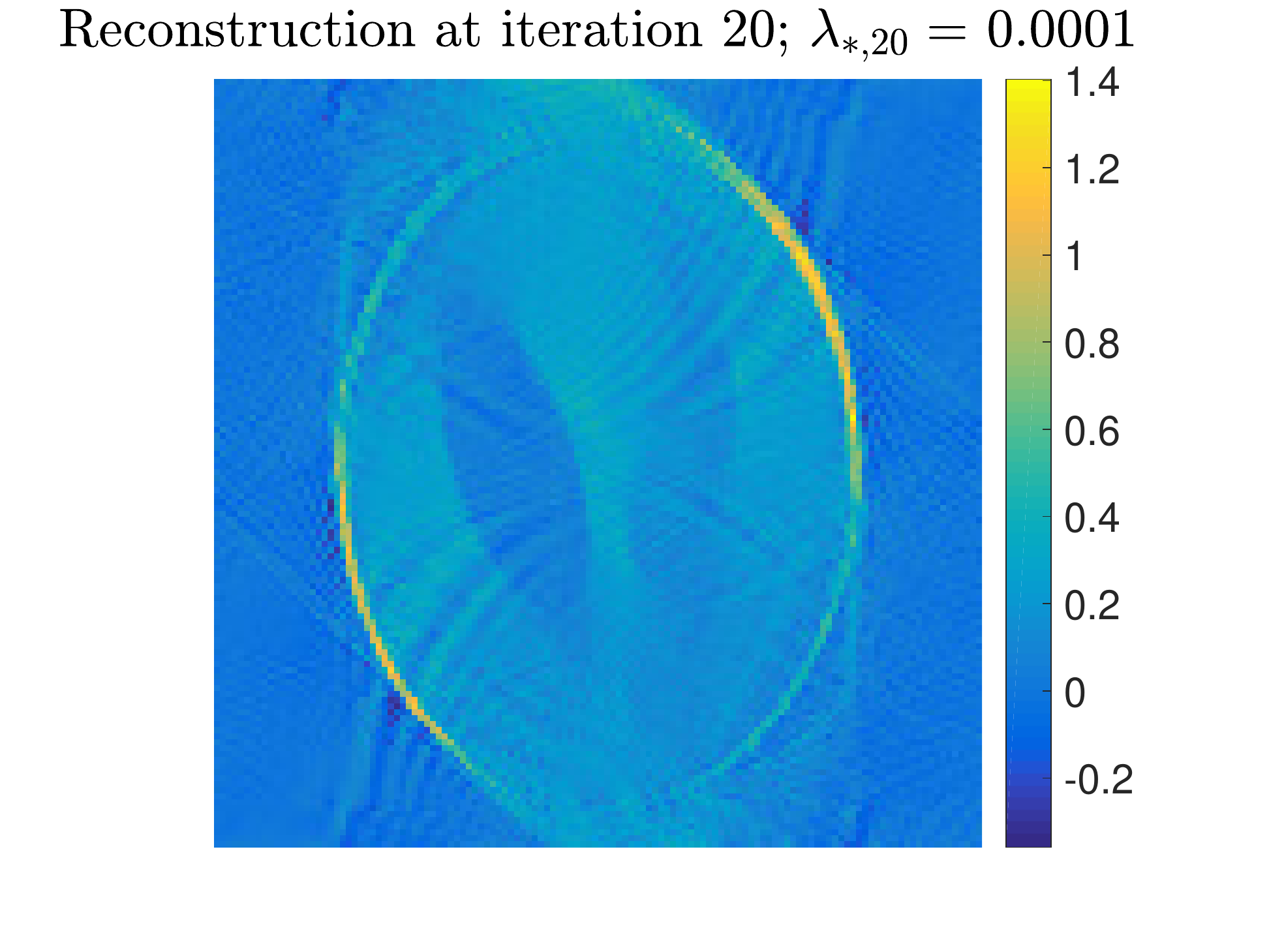}
\end{tabular}
\end{center}
\caption{\emph{Parallel-beam CT with more limited angle} test problem. Upper row: reconstructions obtained using the new weights, at different outer iterations, and by different inner linear solvers. Lower row: reconstructions obtained using the IRN-TV weights, at different outer iterations, and by different inner linear solvers. The corresponding regularization parameters chosen
by the hybrid method, and used within CGLS, are displayed above each image.}
\label{fig:Radon1cSolutions}
\end{figure}

Finally, Figure \ref{fig:Radon1cSolutions} displays some relevant reconstructions. We show the initial reconstructions $x^{(\ast,1)}$ obtained at the end of the first inner iteration cycle, when a Tikhonov-regularized problem with regularization term $R(x) = \|Lx\|_2^2$ is used, and the reconstructions $x^{(\ast,20)}$ obtained when the maximum number of outer iterations is performed. We consider the inner-outer iterative solvers that employ both the new and the IRN-TV weights, and both the hybrid method (with the discrepancy principle as a parameter choice rule) and CGLS. We can clearly see that, when the hybrid method is used as inner solver, there is a 
significant improvement in the reconstructions as the outer iterations proceed and, in particular, 
the edges at the final outer iteration are much sharper than at the initial outer iteration: this is especially true for the new weights, while some artifacts are evident in the IRN-TV reconstructions. As expected, the CGLS reconstructions are still very corrupted, even at the end of the full cycle of outer iterations. 
%
%\begin{figure}[htbp]
%\begin{center}
%\includegraphics[width=7cm]{FigsExampleRadon1b/LCurves1b} 
%\end{center}
%\caption{${\mathcal L}$-curves for each iteration, for the second test problem. As implied from
%the text in the plot, the top curve corresponds to the first outer iteration, $\ell = 1$, and the curves
%below this correspond sequentially to iterations $\ell = 2, 3, 4$.  The red circles
%denote corners of each ${\mathcal L}$-curve, which correspond to the chosen
%regularization parameter, $\lambda_{*,\ell}$ for the particular iteration.}
%\label{fig:Radon1bLCurves}
%\end{figure}
%
%Computed reconstructions for the first outer iteration (that is, $x^{(*,1)}$), and for the
%final outer iteration (that is, $x^{(*,4)}$) are shown 
%in Figure~\ref{fig:Radon1bSolutions}. As we can see from these plots, there is a 
%significant improvement in the reconstructions, and in particular 
%the edges at the final outer iteration are much sharper than in the initial outer iteration.
%
Figure \ref{fig:Radon2TVweight} displays the entries of the weight diagonal matrices at outer iterations $\ell=2$ and $\ell=20$, for both the new reweighting strategy (\ref{eq:GradientMap}) and the IRN-TV reweighting strategy. When considering the new reweighting strategy, we can clearly see that both the weights to be applied to the vertical and horizontal derivatives are appropriate and improve with  increasing outer iterations. 
%; the same is not true for the weights to be applied to the horizontal  derivatives, which seem to detect some spurious edges in the first reconstruction $x^{(\ast,1)}$, which are enhanced during the subsequent outer iterations. Despite this, spurious edges are hardly visible in the final reconstruction. 
Similarly to the previous test problems, the IRN-TV weights are not so effective in revealing the structure of the image; some true and spurious edges seem to be assigned very small weights at the end of the iterations.
\begin{figure}[htbp]
\begin{center}
\begin{tabular}{ccc}
\includegraphics[width=5cm]{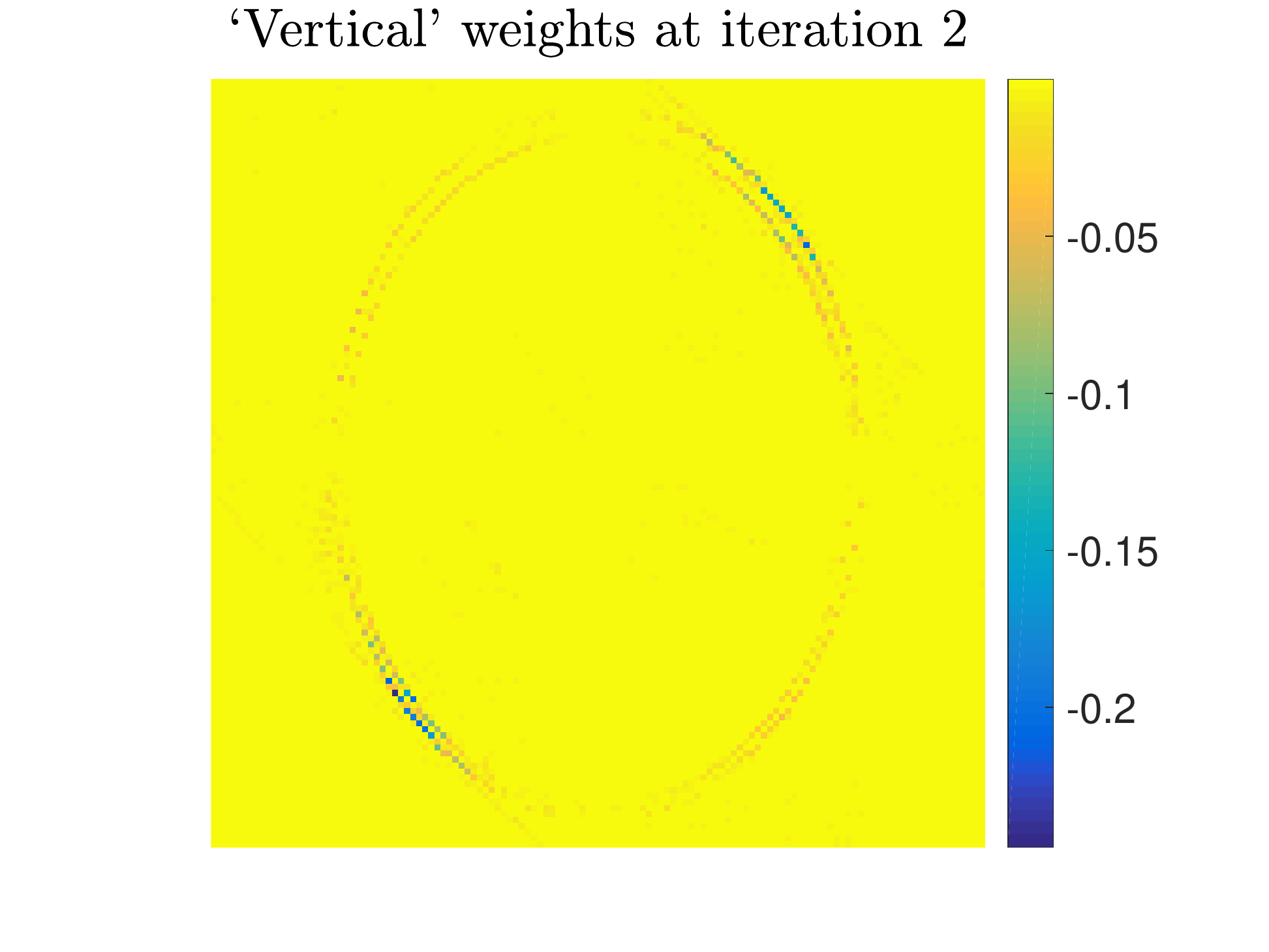} &
\includegraphics[width=5cm]{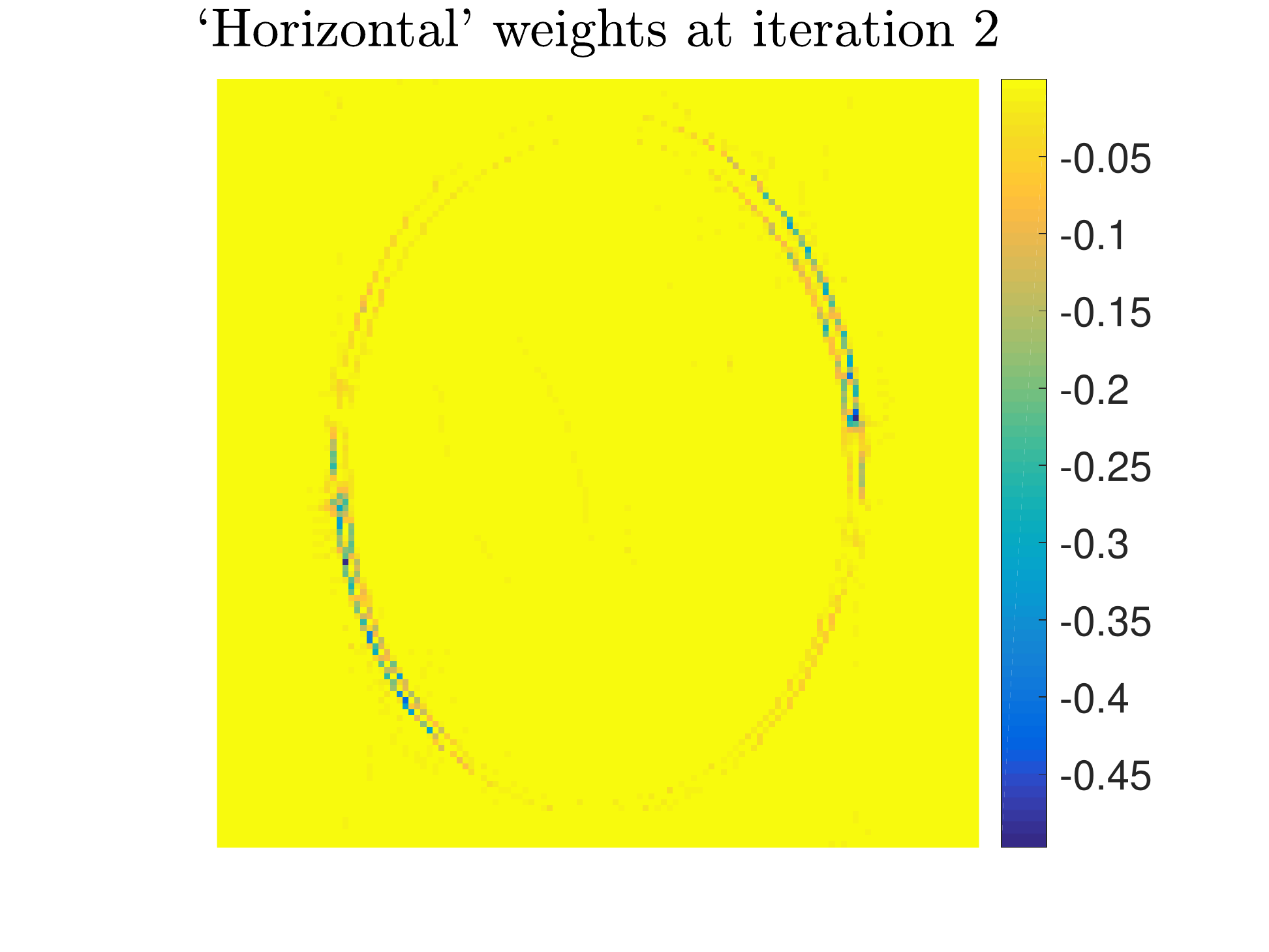} & 
\includegraphics[width=5cm]{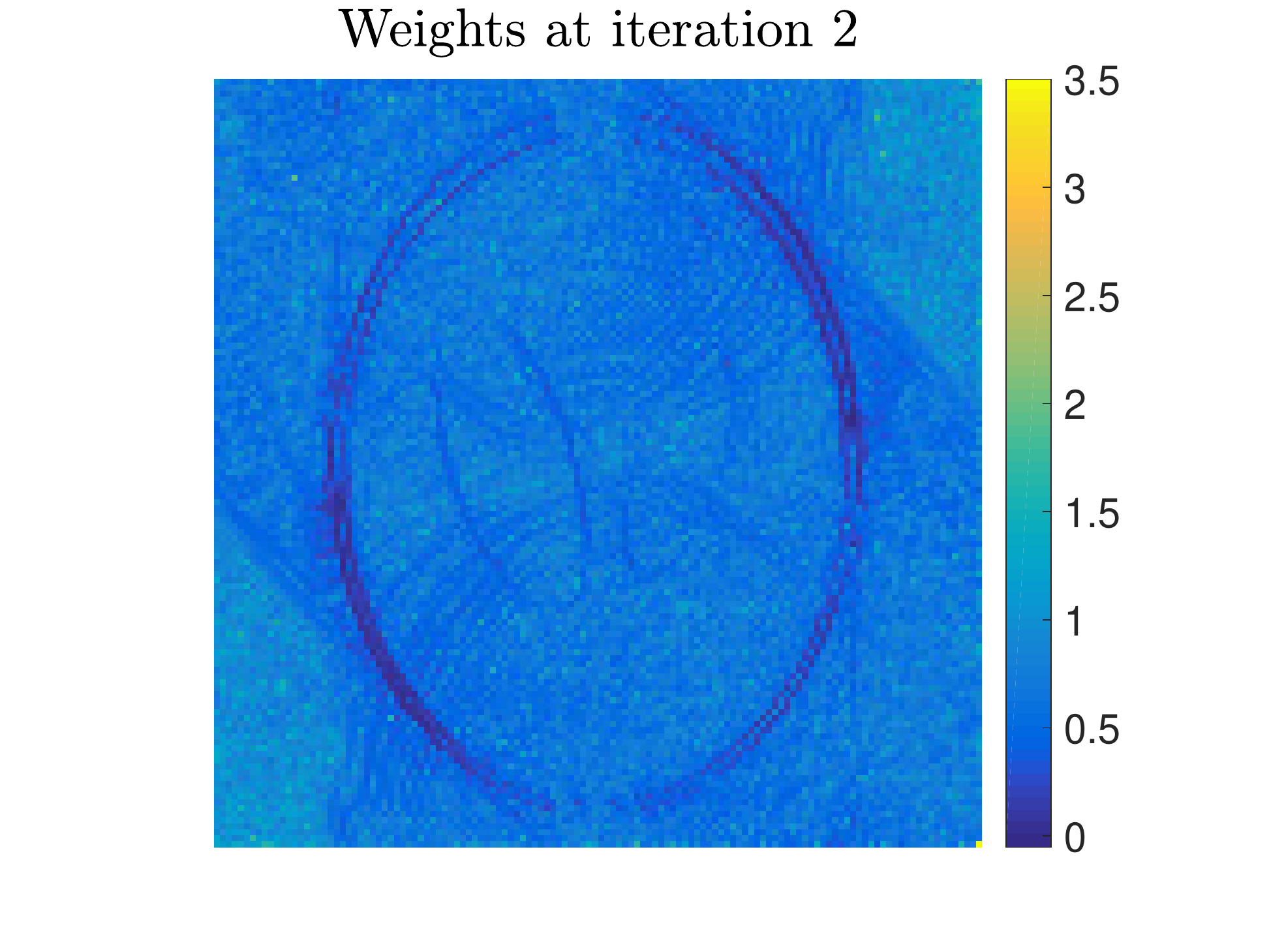}\\ 
\includegraphics[width=5cm]{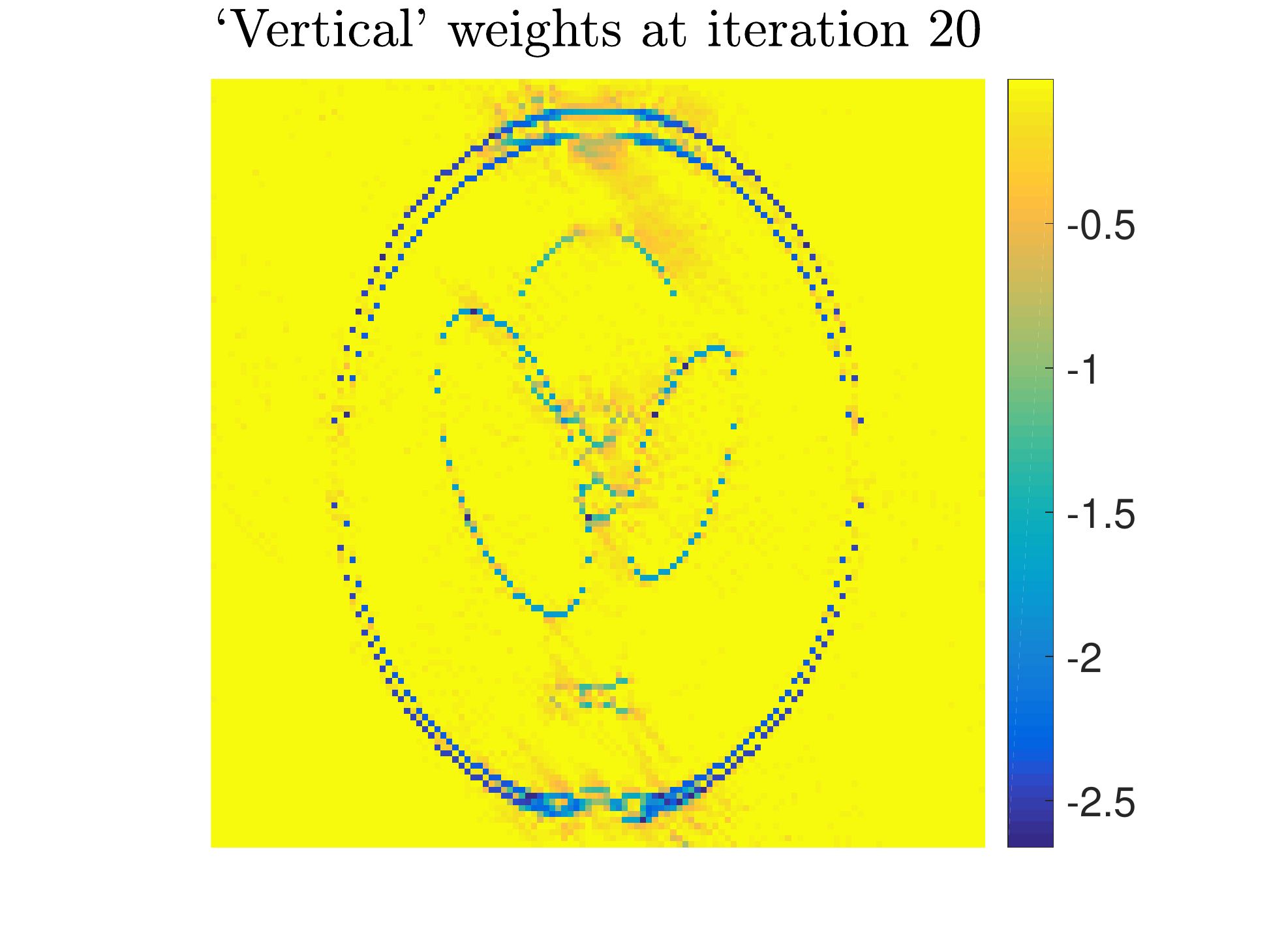} &
\includegraphics[width=5cm]{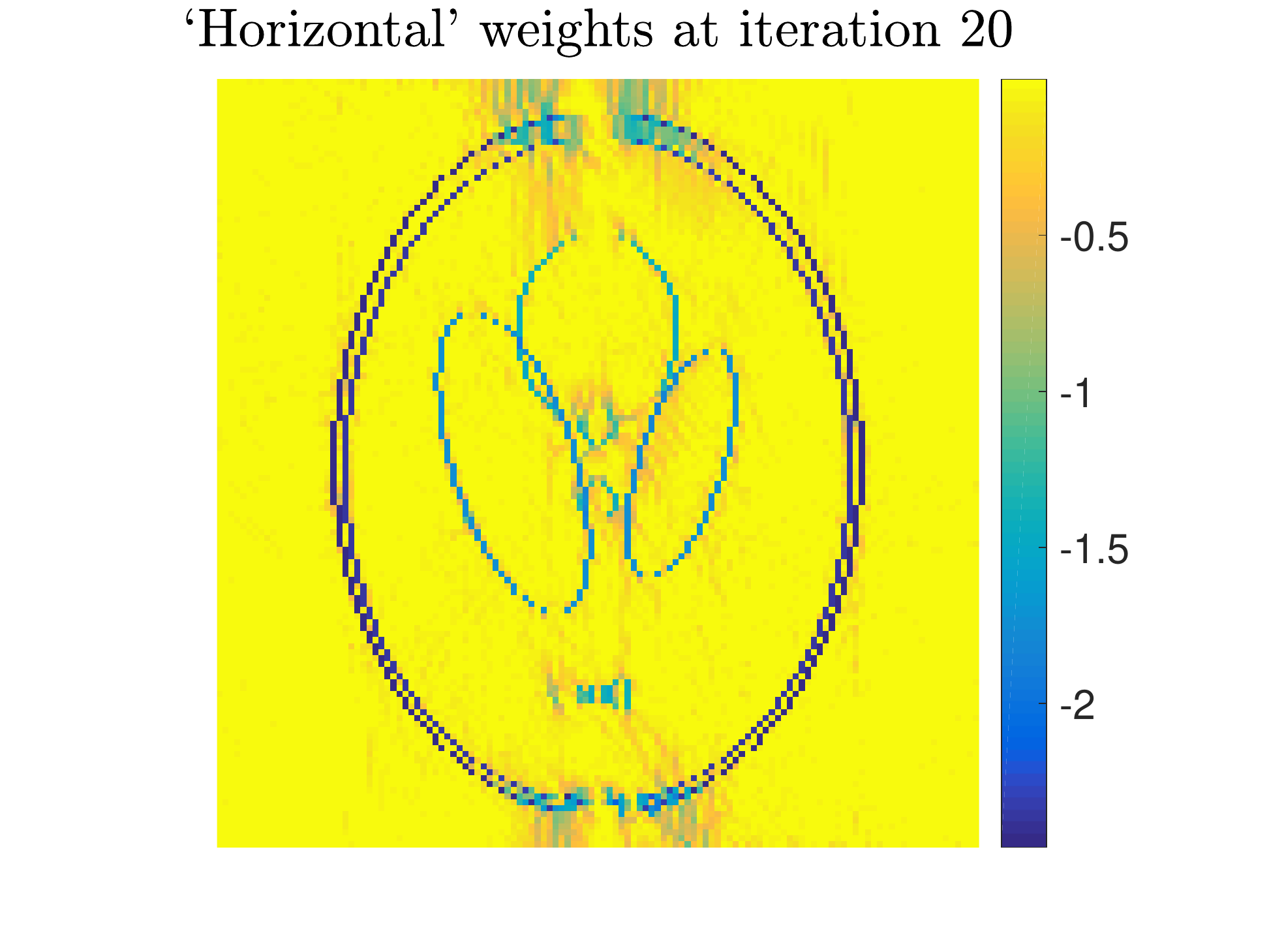} & 
\includegraphics[width=5cm]{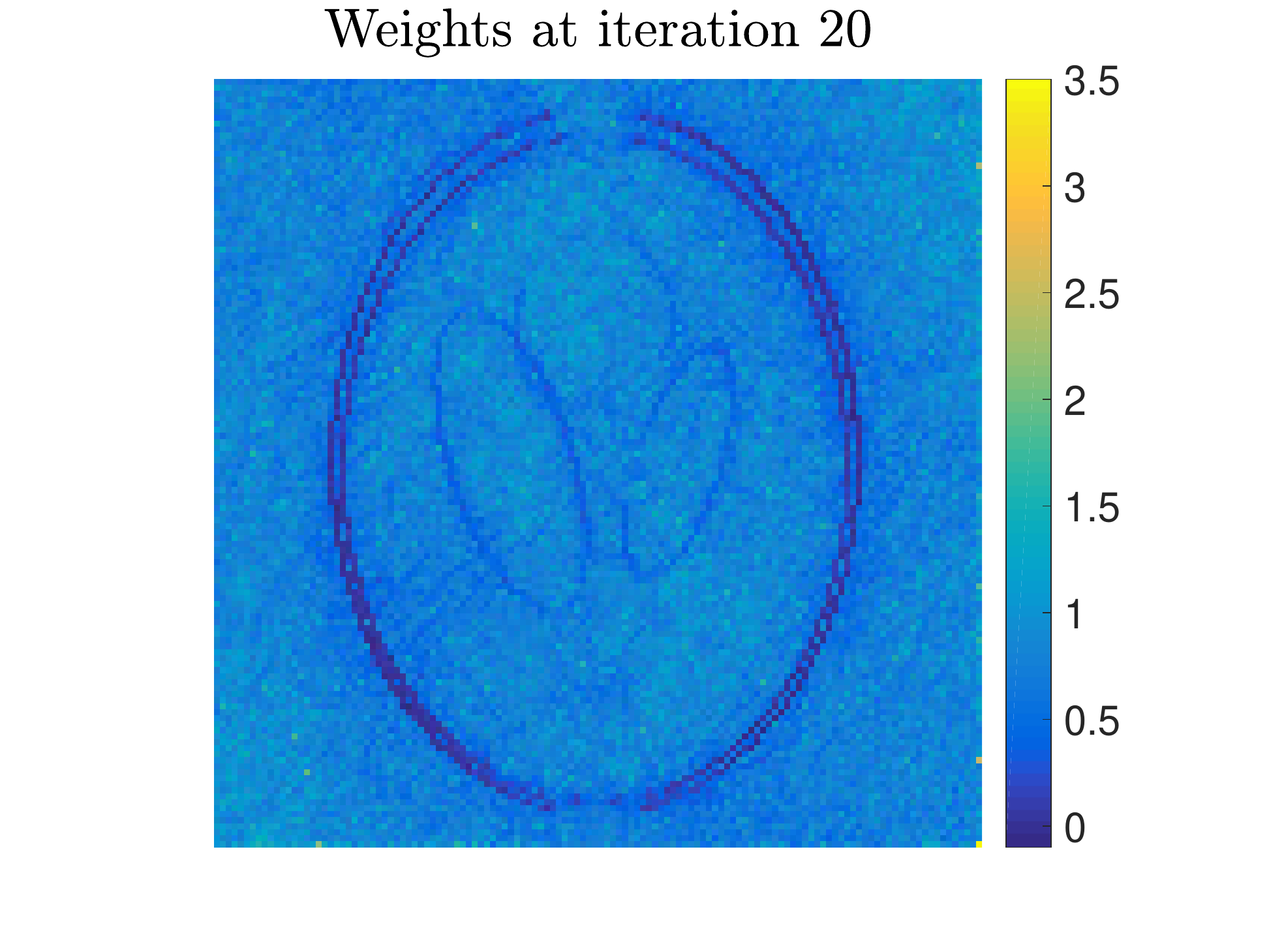}\\ \end{tabular}
\end{center}
\caption{\emph{Parallel-beam CT with more limited angle} test problem. Left column: new weights, to be applied to the vertical derivatives. Middle column: new weights, to be applied to the horizontal derivatives. Right column: IRN-TV weights. The pixel values are displayed in logarithmic scale.}
\label{fig:Radon3TVweight}
\end{figure}

\subsection{Examples from image deblurring}

\paragraph{Shaking blur} 

\begin{figure}[htbp]
\begin{center}
\begin{tabular}{cc}
\footnotesize{True image} & \footnotesize{Corrupted image}\\
\includegraphics[height=4cm]{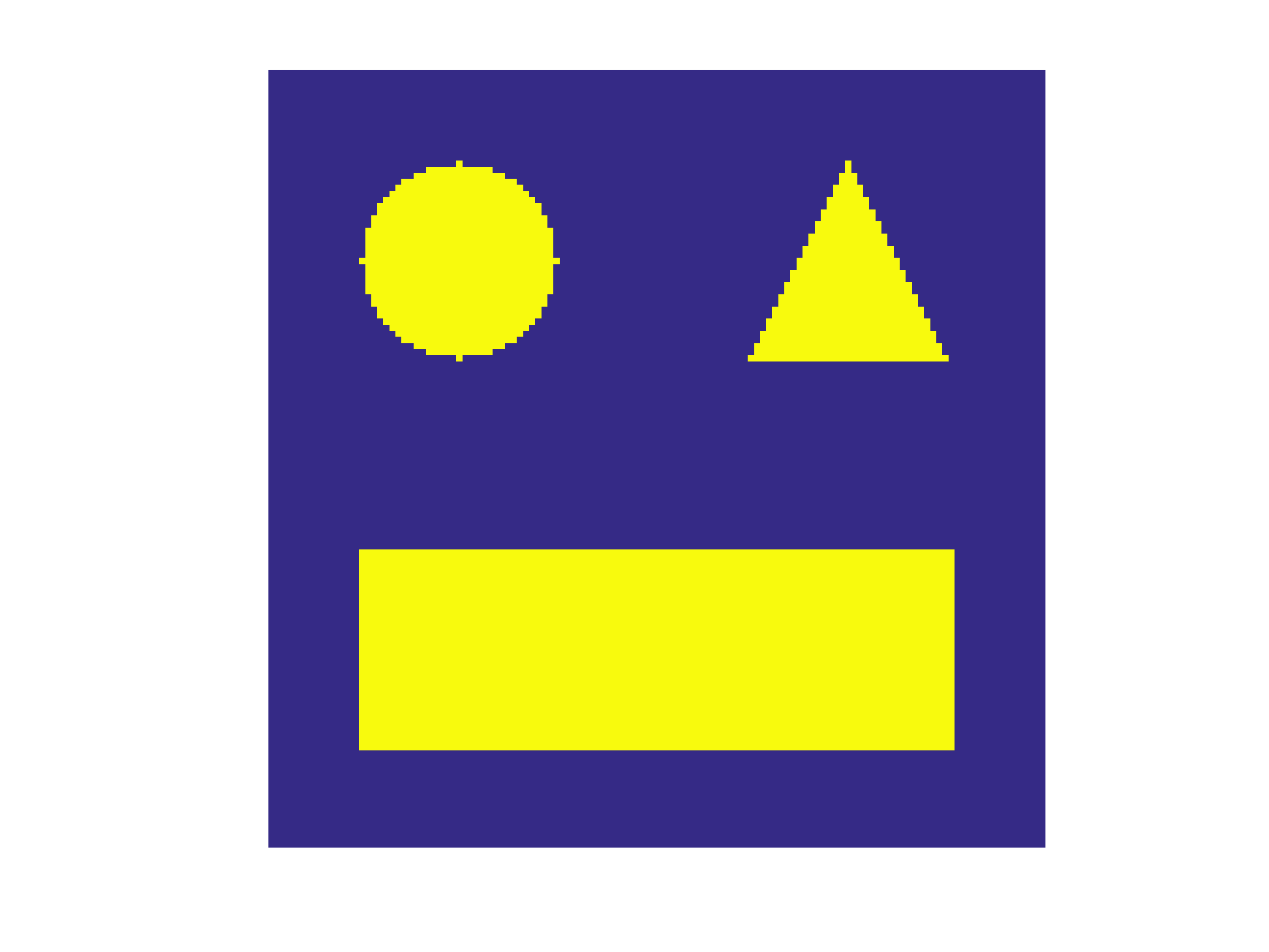} &
\includegraphics[height=4cm]{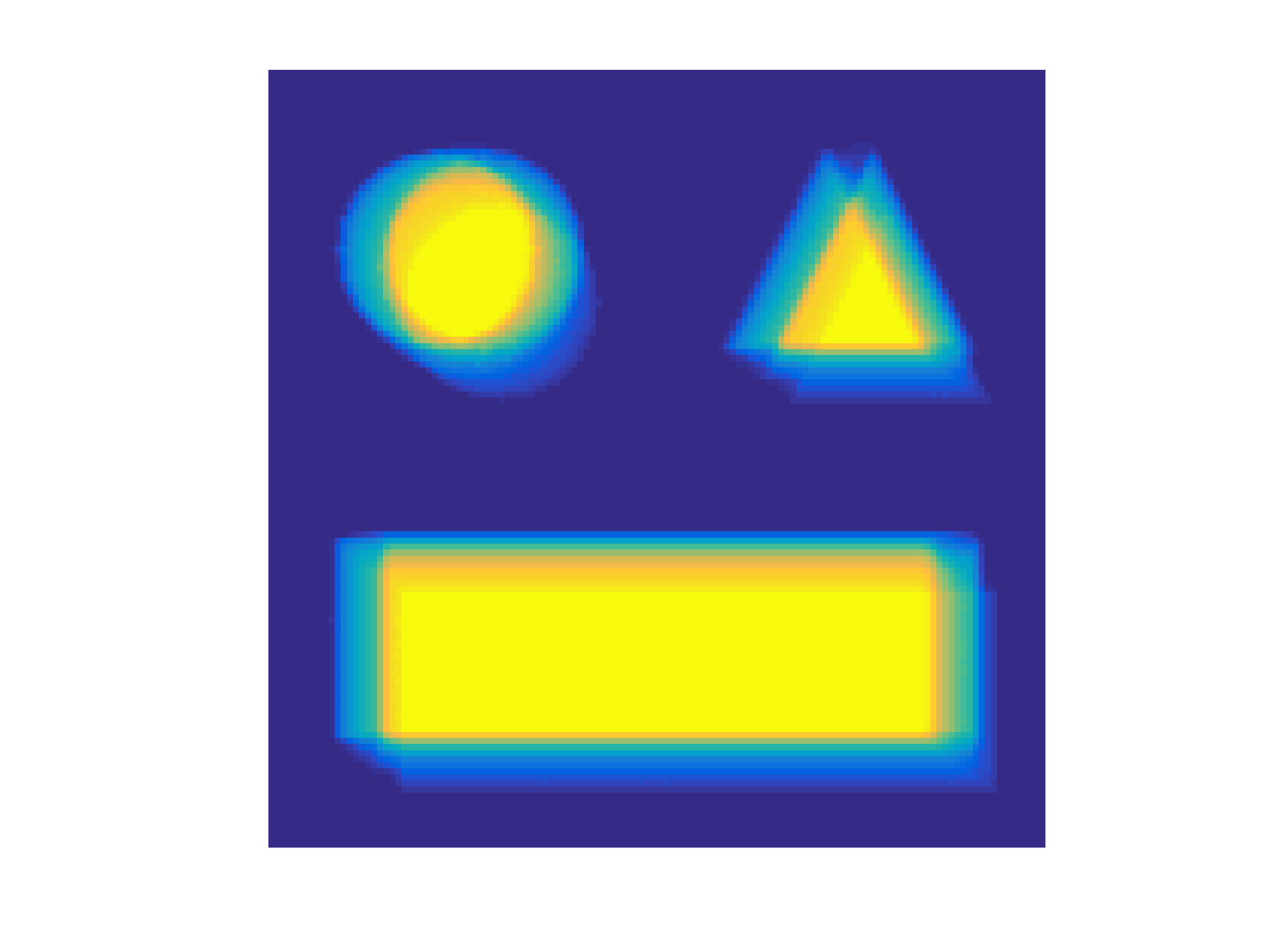}
\end{tabular}
\end{center}
\caption{\emph{Shaking blur} test problem. Left frame: true sharp image $x$. Right frame: the measured data (corrupted image) $b$.}
\label{fig:Blur1Data}
\end{figure}

We use IR Tools to setup the following shaking blur simulation 
\begin{itemize}
\item 
Generate a simulated true sharp geometric image called ``pattern1'', of size $128\times 128$ pixels.
\item
Construct noise-free blurred image, along with the matrix $A$ that simulates the shaking blur forward operator of medium intensity
\item
Add 0.1\% normally distributed (white Gaussian) noise.
\end{itemize}
More specifically, the test problem is generated with the following MATLAB statements:
\begin{verbatim}
     ProblemOptions = PRset('trueImage', 'pattern1', 'BlurLevel', medium); 
     [A, b_true, x_true, ProblemInfo] = PRblurshake(128, ProblemOptions);
     b = PRnoise(b_true, NoiseLevel);
\end{verbatim}

Figure~\ref{fig:Blur1Data} shows the true sharp image, along with the
measured data $b$.

\begin{figure}[htbp]
\begin{center}
\begin{tabular}{cc}
\includegraphics[width=6cm]{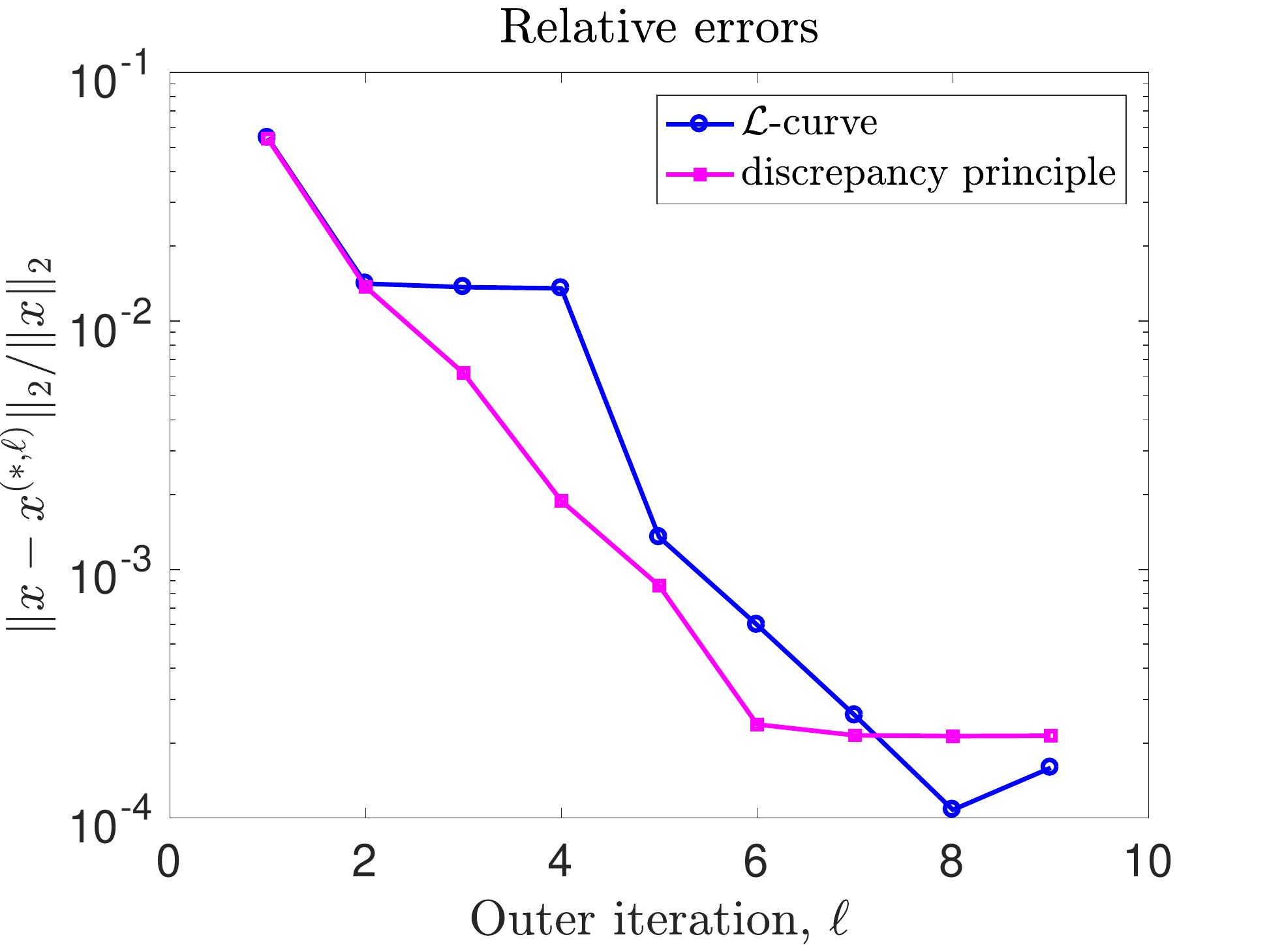} &
\includegraphics[width=6cm]{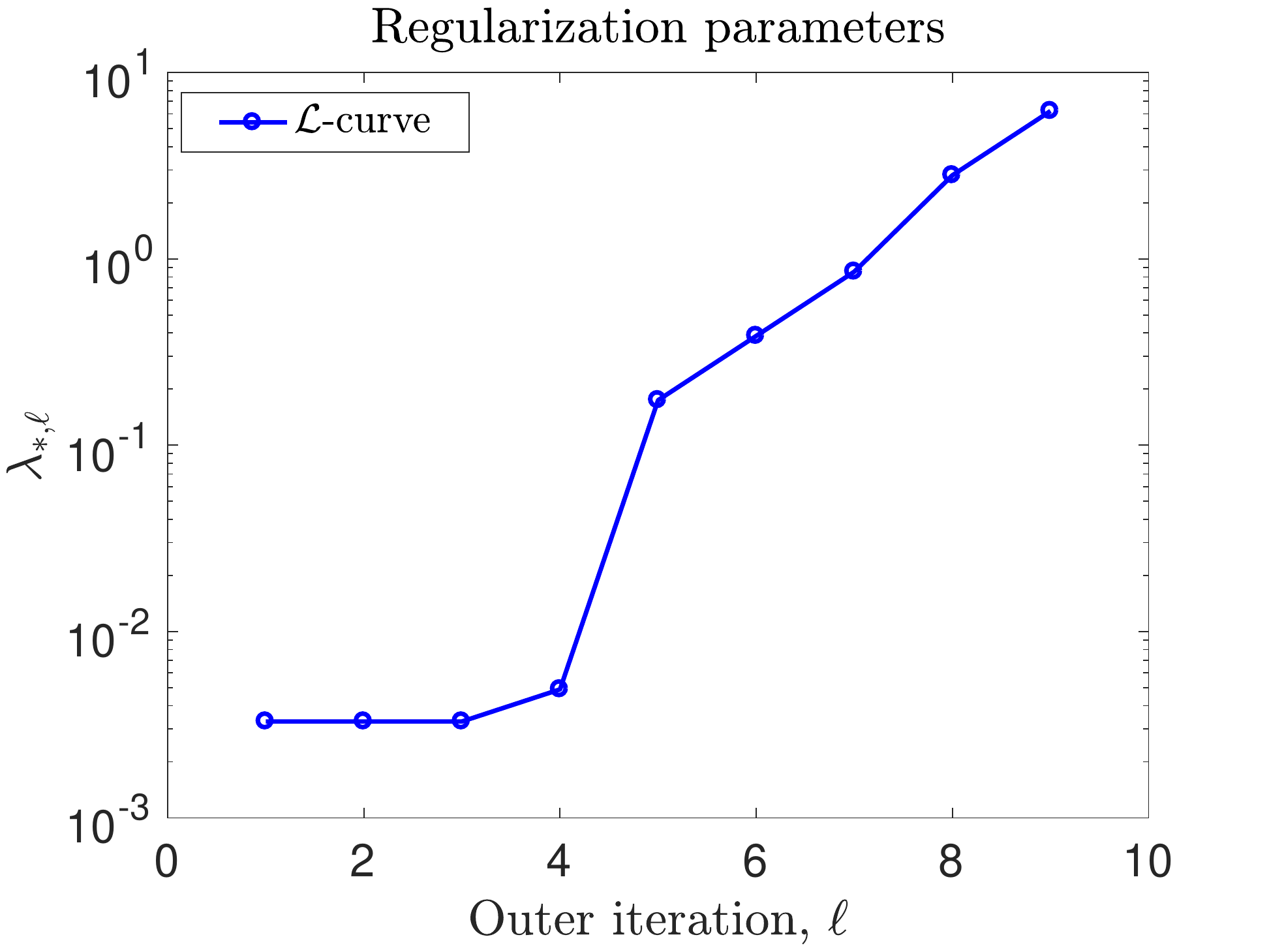}
\end{tabular}
\end{center}
\caption{\emph{Shaking blur} test problem. Relative errors and regularization parameters values at each outer iteration $\ell$,until the stopping criterion is satisfied. Both the discrepancy principle and the ${\mathcal L}$-curve criterion are considered; the regularization parameter computed by the former is always numerically zero, and therefore it is not displayed in the rightmost frame.}
\label{fig:Blur1Iterations}
\end{figure}

As in the previous examples, we first run our algorithm for 9 outer iterations with different parameter choice strategies within the inner hybrid scheme for generalized Tikhonov regularization. Figure~\ref{fig:Blur1Iterations} shows a plot of the relative errors
and chosen regularization parameters at each outer iteration. We can clearly see that, when considering both the discrepancy principle and the $\mathcal{L}$-curve criterion, the reconstructions greatly improve as the outer iterations proceed, with the former strategy being very efficient; however, the $\mathcal{L}$-curve criterion eventually achieves a slightly lower relative error. The regularization parameter selected by the $\mathcal{L}$-curve criterion is quite small and almost stagnates during the first outer iterations (where the corners of the $\mathcal{L}$-curves may be hard to distinguish), but considerably increase during the final outer iterations: as remarked in the previous subsection, this behavior is meaningful and desirable. On the contrary, the regularization parameters computed by the discrepancy principle are always either zero or numerically zero (and therefore they are not displayed in the rightmost frame of Figure ~\ref{fig:Blur1Iterations}). This means that the discrepancy curves almost never cross the noise level line: this could be because the noise level $\|\eta\|_2/\|\btrue\|_2=10^{-3}$ is not very high, and/or this test problem is not very ill-conditioned. It is interesting to note that, despite of this, the regularization term favorably affects the approximation subspace for the solution, and the stopping criterion for the inner iterations prevents the quality of the solutions to degenerate. 

\begin{figure}[htbp]
\begin{center}
\begin{tabular}{cc}
\includegraphics[width=6cm]{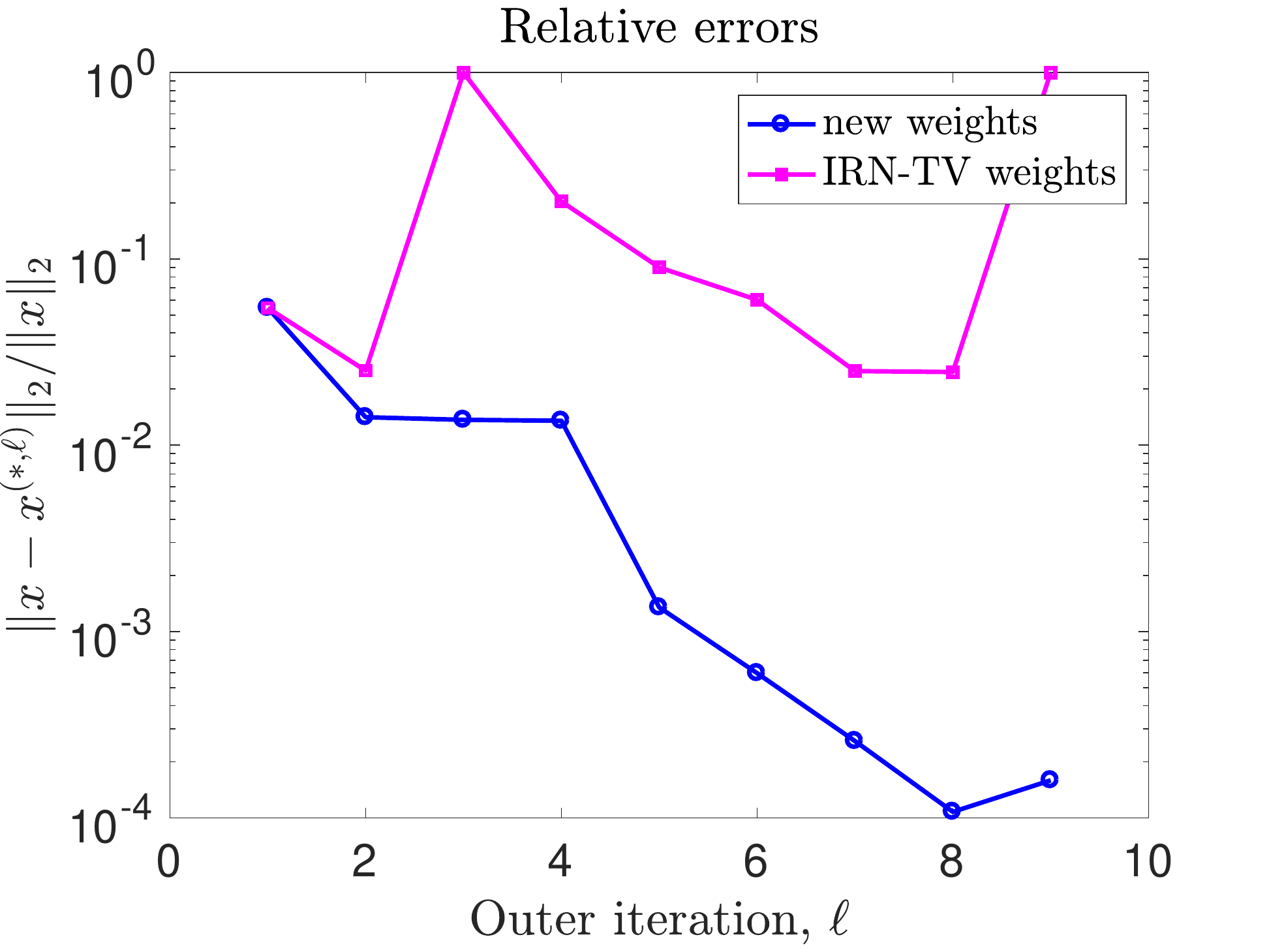} &
\includegraphics[width=6cm]{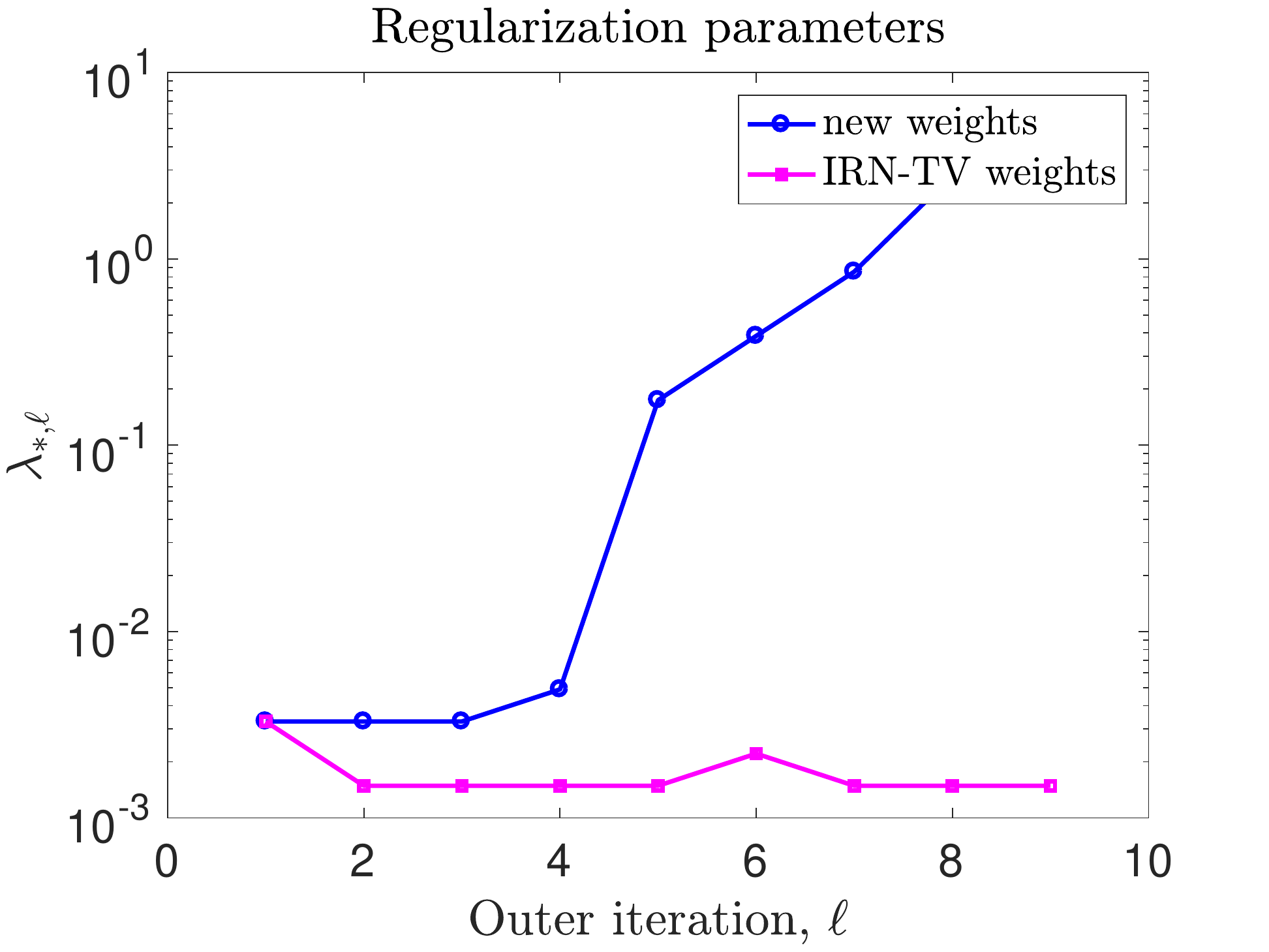}
\end{tabular}
\end{center}
\caption{\emph{Shaking blur} test problem. Relative errors and regularization parameter versus number of (outer) iterations; both methods use the hybrid method for general form Tikhonov during the inner iterations, and adaptively select the regularization parameter according to the ${\mathcal L}$-curve.}
\label{fig:Blur1TV}
\end{figure}

\begin{figure}[htbp]
\begin{center}
\begin{tabular}{cc}
\includegraphics[width=6cm]{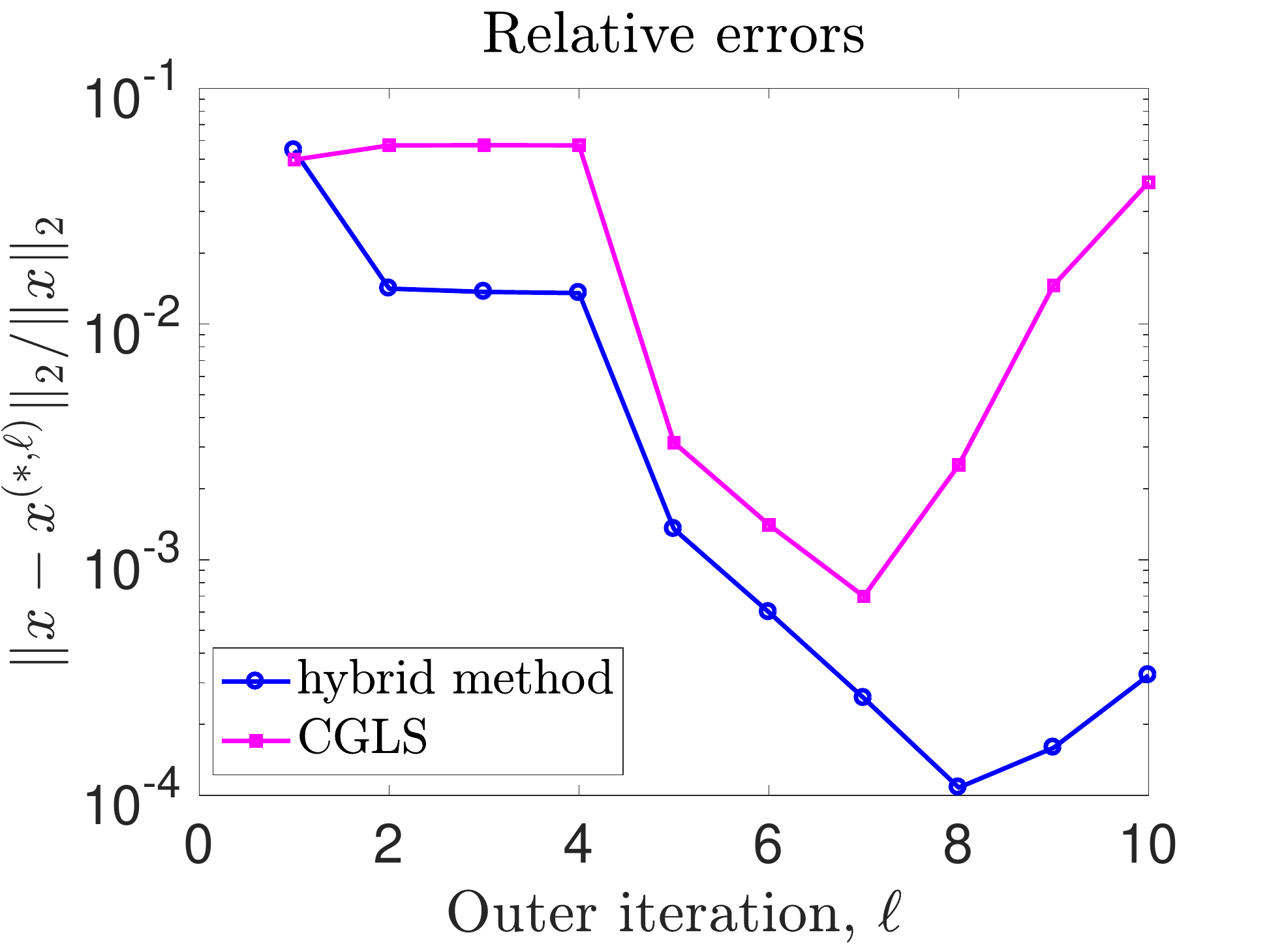} &
\includegraphics[width=6cm]{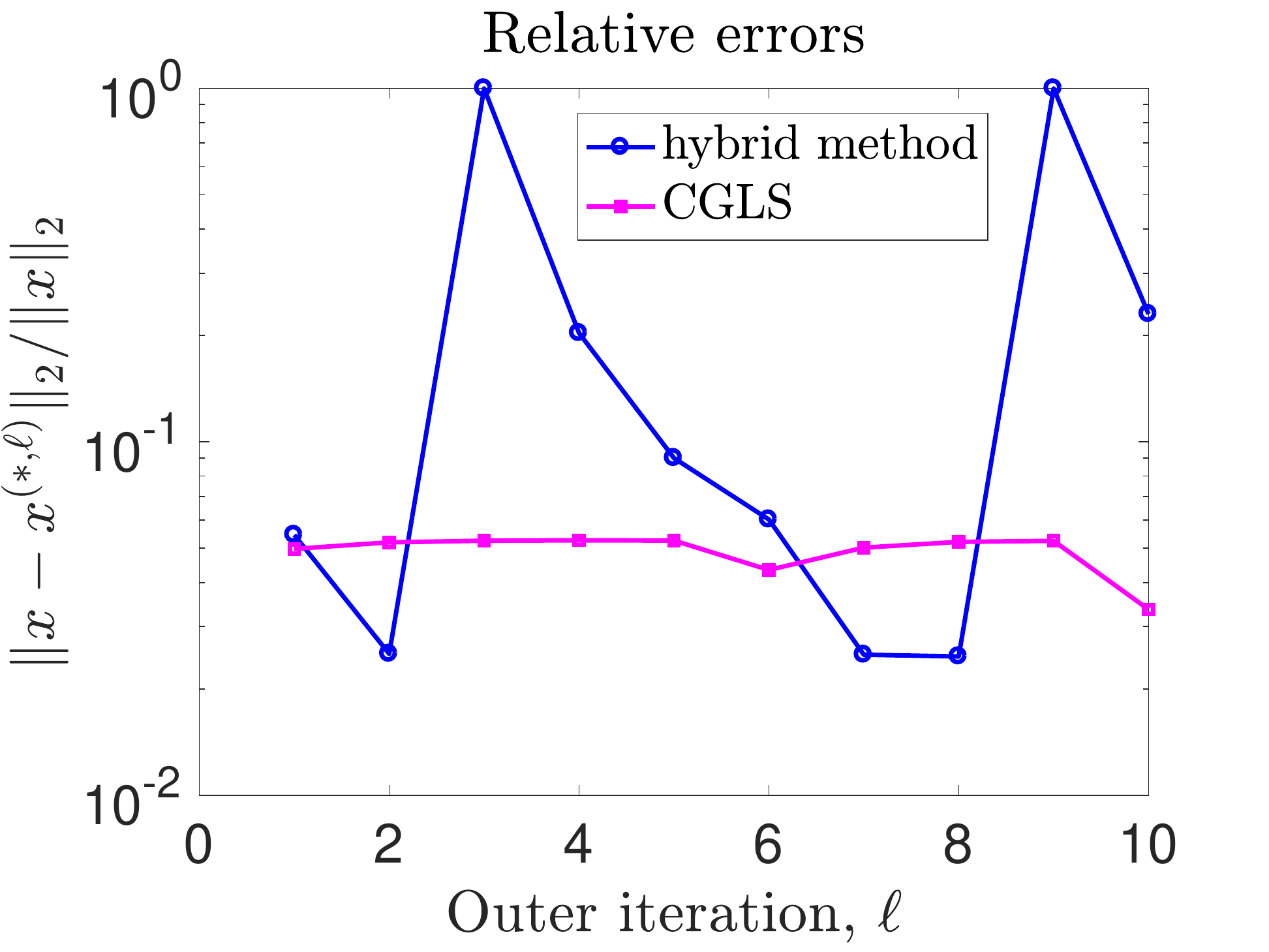}
\end{tabular}
\end{center}
\caption{\emph{Shaking blur} test problem. Relative errors versus number of (outer) iterations. Left frame: the new weights are used. Right frame: the IRN-TV weights are used.}
\label{fig:Blur1TV_comparisons}
\end{figure}

Next, we compare the new method to other inner-outer iterative methods for edge enhancement in imaging. Figure \ref{fig:Blur1TV} displays the relative error and regularization parameter values versus the number of outer iterations for the new method, and for a method that still employs the hybrid solver for general form Tikhonov regularization to handle the inner iterations while using the IRN-TV weights (\ref{eq:IRNTVweights}) at each outer iteration. For both methods, the ${\mathcal L}$-curve criterion is employed to adaptively choose the regularization parameter at each inner iteration. 
We can clearly see that, when the IRN-TV weights are used, the behavior of the relative errors is quite hectic and, although the quality of the solution computed at the second outer iteration is good, this trend is not maintained during the following outer iterations. This situation may improve if alternative stopping criteria for the inner iterations are performed; note also that the regularization parameter for IRN-TV almost immediately stabilizes around a value of the order of $10^{-3}$. 
%%%; indeed, for this test problem, setting the regularization parameter through the $\mathcal{L}$-curve results in a more accurate reconstructions
%% Similarly to the previous test problems, both methods perform well during the first (outer) iterations, with the IRN-TV weights being quite effective in reducing the error; however, the quality of the IRN-TV solutions eventually stagnates, while the new method keeps improving. 
%%%Similarly to the previous test problems, both methods perform well during the first (outer) iterations, with the IRN-TV weights being quite effective in reducing the error; however, the quality of the IRN-TV solutions eventually stagnates, while the new method keeps improving. 
%%The automatically selected regularization parameter for IRN-TV keeps oscillating between values of the order of $10^{-2}$ and values of the order of $10^{-4}$. 
% computes  the new method is amazing at the end.
Figure \ref{fig:Blur1TV_comparisons} assesses the influence of the inner iterative solver on the overall behavior of the method. Namely, we consider the inner-outer iterative schemes implemented with both the new and the IRN-TV weights, and with both the hybrid and the CGLS methods as inner solvers. The hybrid method chooses the regularization parameter adaptively according to the ${\mathcal L}$-curve criterion as the iterations proceed, and the value $\lambda_{\ast,\ell}$ selected when the $\ell$th inner iteration cycle terminates is taken to be the fixed regularization parameter to be set in advance of the $\ell$th CGLS iteration cycle. 
% requires a fixed value of the regularization parameter to be available in advance of each cycle of inner iterations, we first run the methods based on the hybrid solver for general form Tikhonov that adaptively chooses the regularization parameter at each inner iteration according to the discrepancy principle: in this way, a parameter 
% $\lambda_{\ast,\ell}$ is eventually set at the $\ell$th outer iteration. When running the methods based on CGLS, we take $\lambda_{\ast,\ell}$ as regularization parameter for the $\ell$th outer iteration.
Looking at both frames in Figure \ref{fig:Blur1TV_comparisons} we can see that CGLS performs quite well. When the new weights are adopted, the quality of the solution stagnates during the first outer iterations (this may be because the selected regularization parameter is quite small; see also Figure \ref{fig:Blur1Iterations}), then considerably improves, and eventually degenerates (the latter phenomenon can be avoided by experimenting more with different stopping criteria for the outer iterations). When the TV weights are adopted, no significant improvements in the quality of the solution are visible as the outer iterations proceed (again, this may be because of the almost constant value of the regularization parameter; see also Figure \ref{fig:Blur1TV}).

\begin{figure}[htbp]
\begin{center}
\begin{tabular}{ccc}
\hspace{-0.3cm}\includegraphics[width=5.5cm]{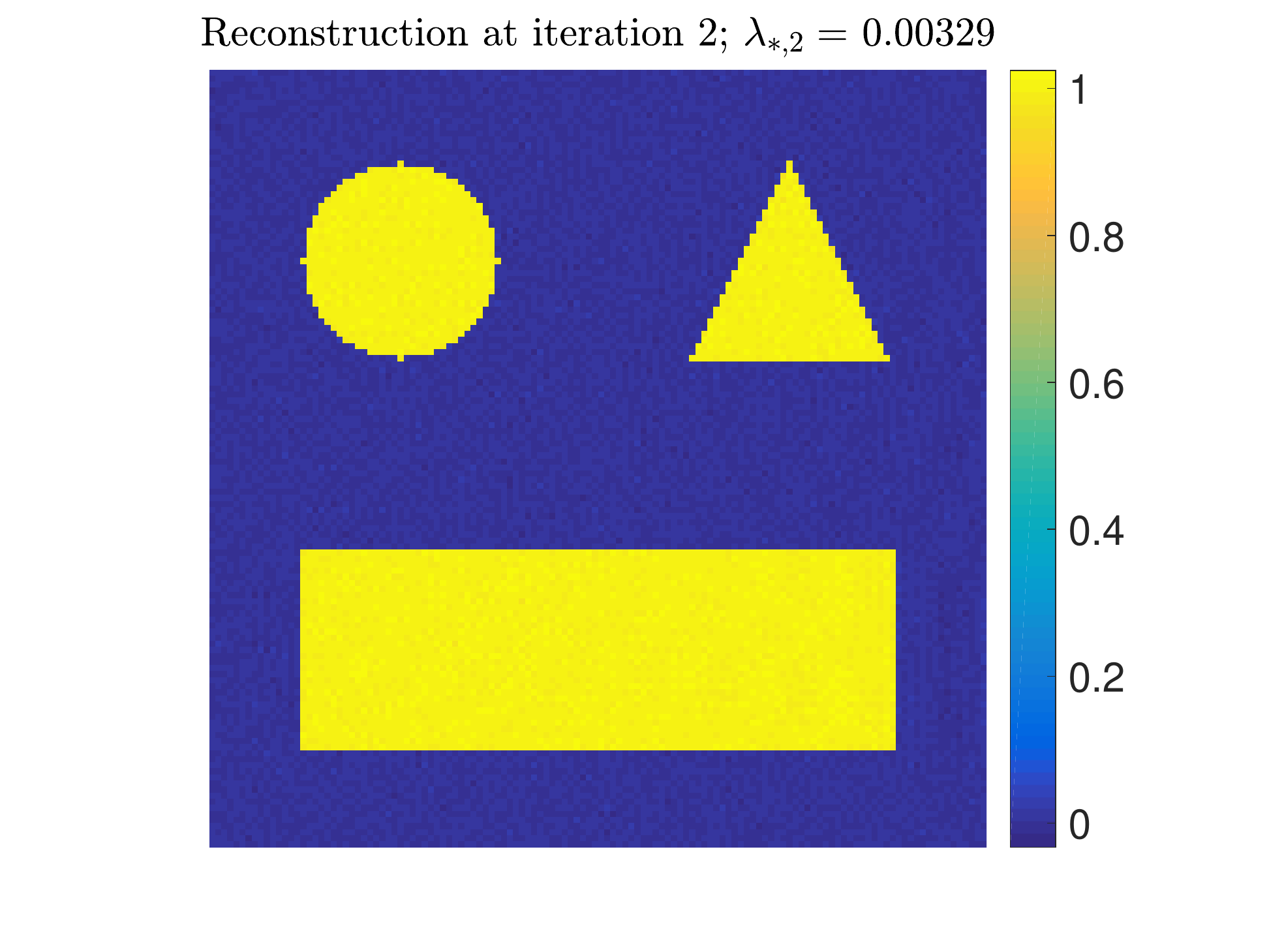} &
\hspace{-0.3cm}\includegraphics[width=5.5cm]{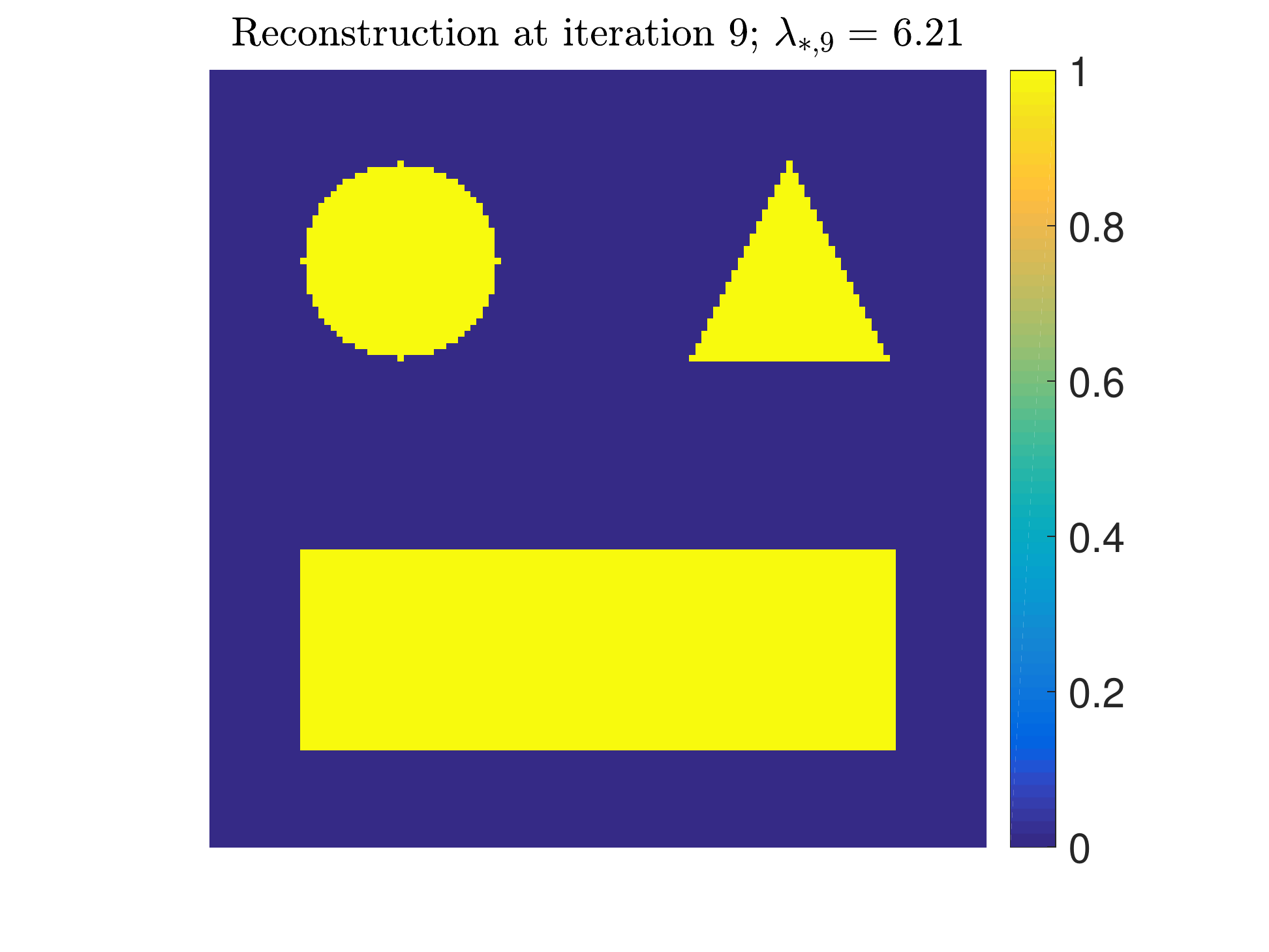} &
\hspace{-0.3cm}\includegraphics[width=5.5cm]{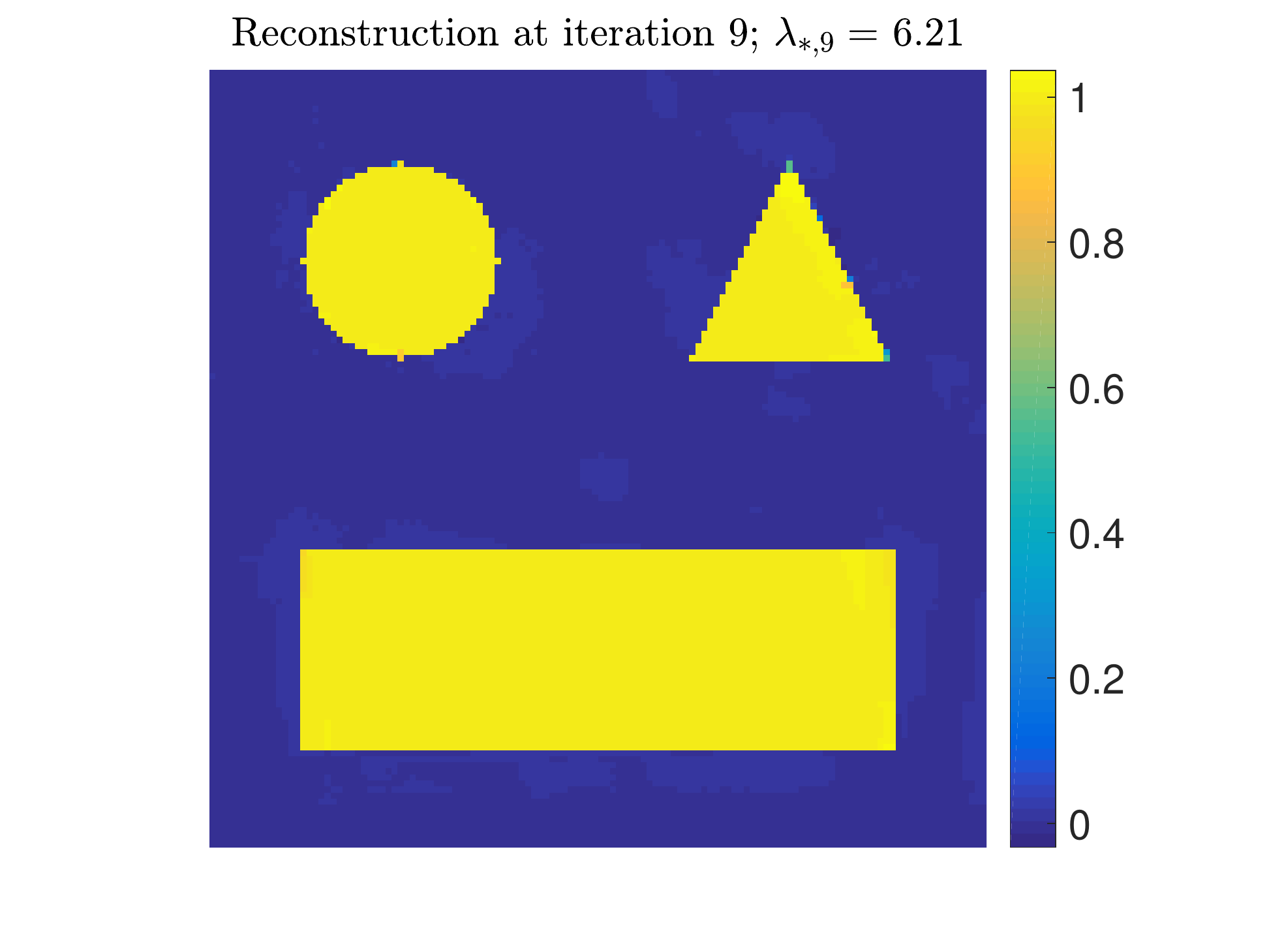}\\
\hspace{-0.3cm}\includegraphics[width=5.5cm]{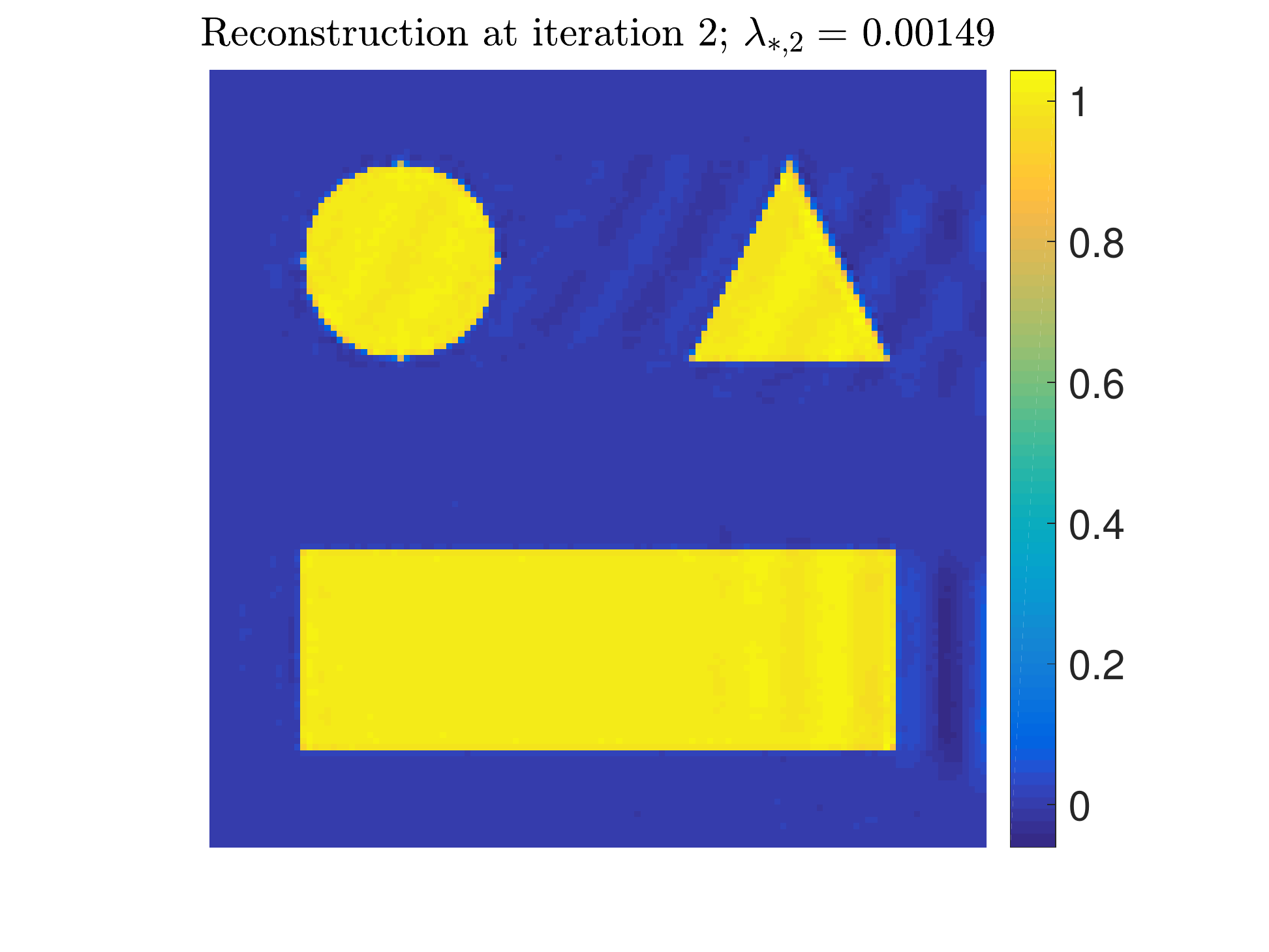} &
\hspace{-0.3cm}\includegraphics[width=5.5cm]{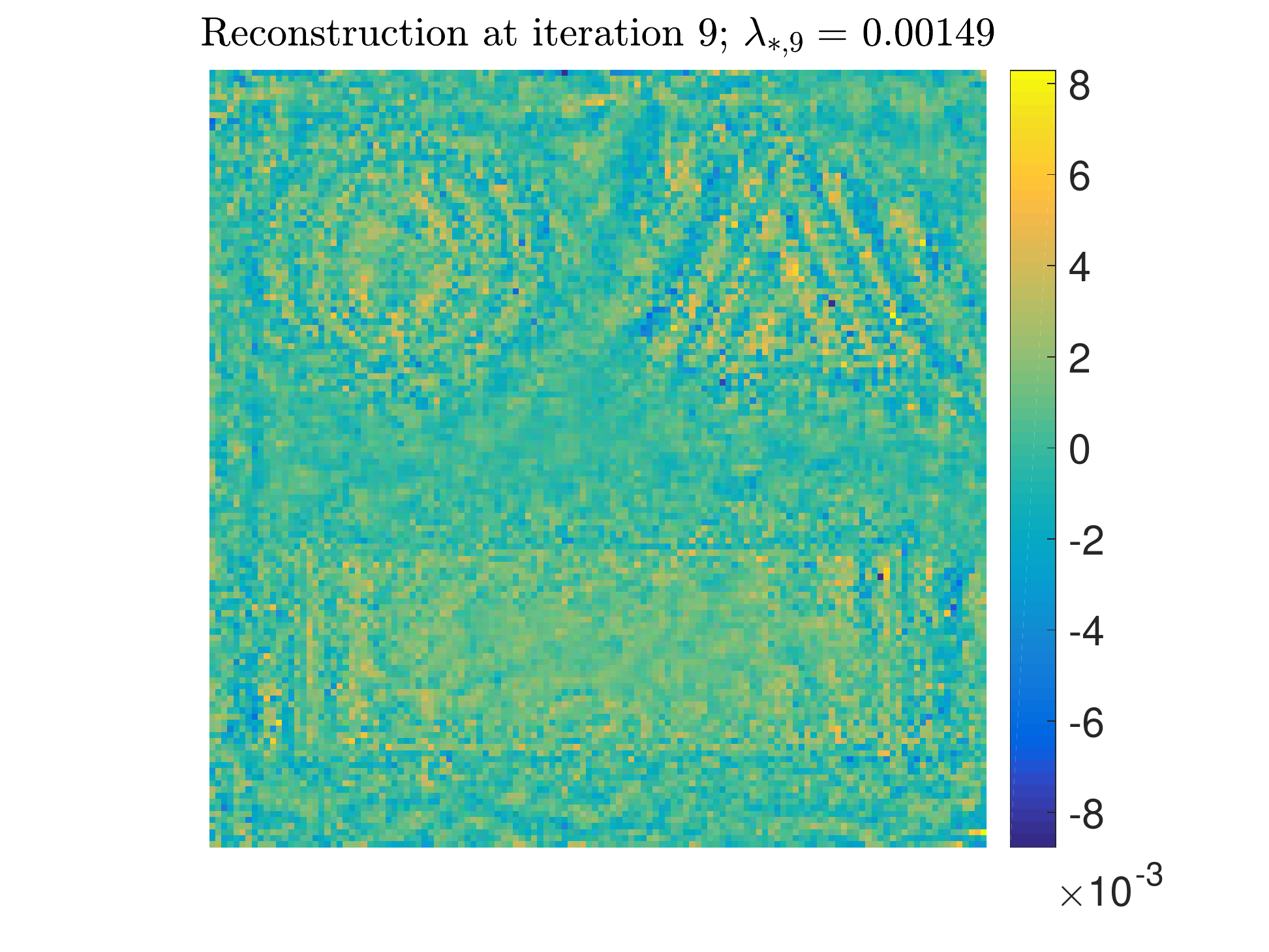} &
\hspace{-0.3cm}\includegraphics[width=5.5cm]{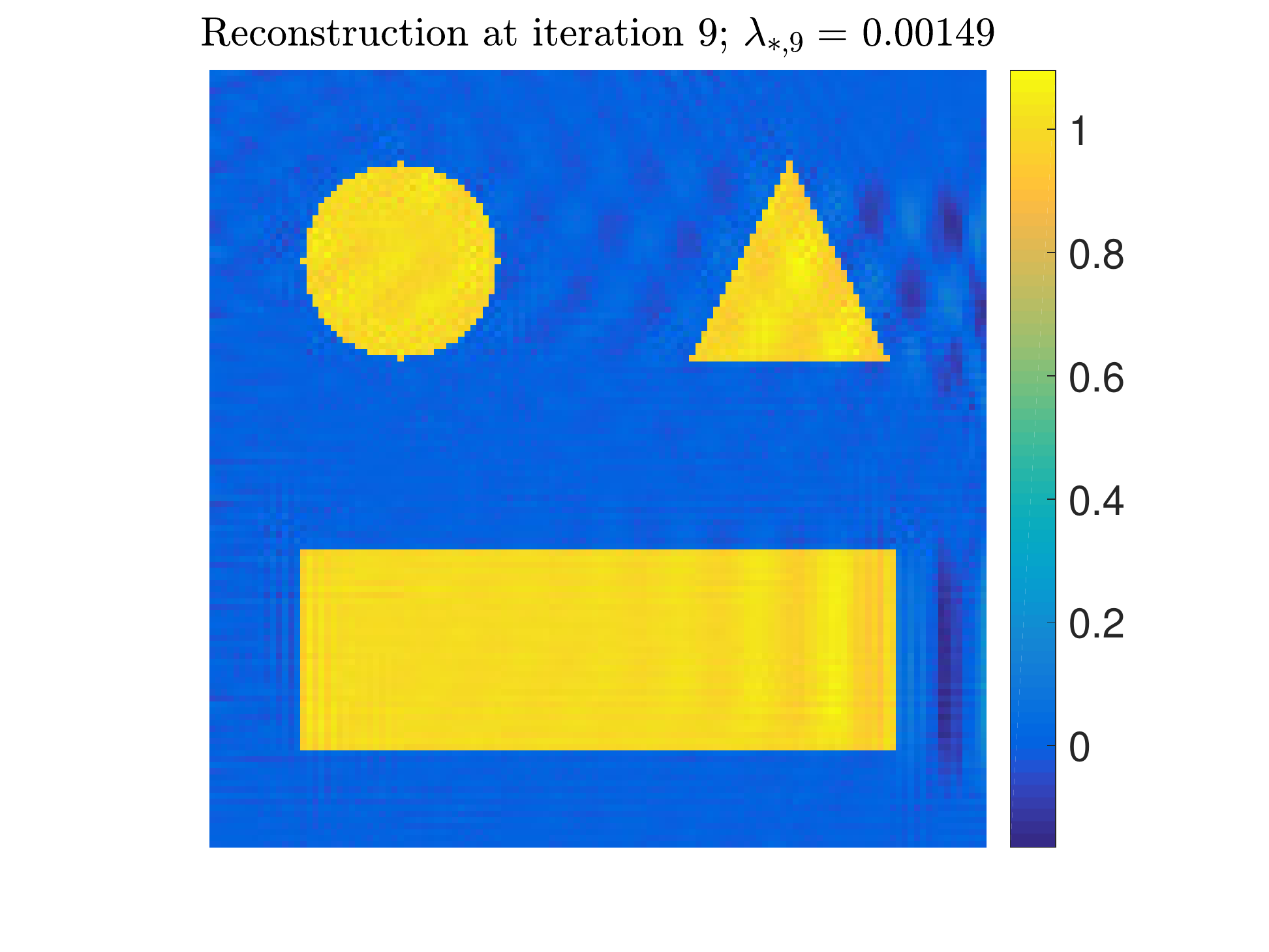}
\end{tabular}
\end{center}
\caption{\emph{Shaking blur} test problem. Upper row: reconstructions obtained using the new weights, at different outer iterations, and by different inner linear solvers. Lower row: reconstructions obtained using the IRN-TV weights, at different outer iterations, and by different inner linear solvers. The corresponding regularization parameters chosen
by the hybrid method, and used within CGLS, are displayed above each image.}
\label{fig:Blur1Solutions}
\end{figure}

\begin{figure}[htbp]
\begin{center}
\begin{tabular}{ccc}
\includegraphics[width=5cm]{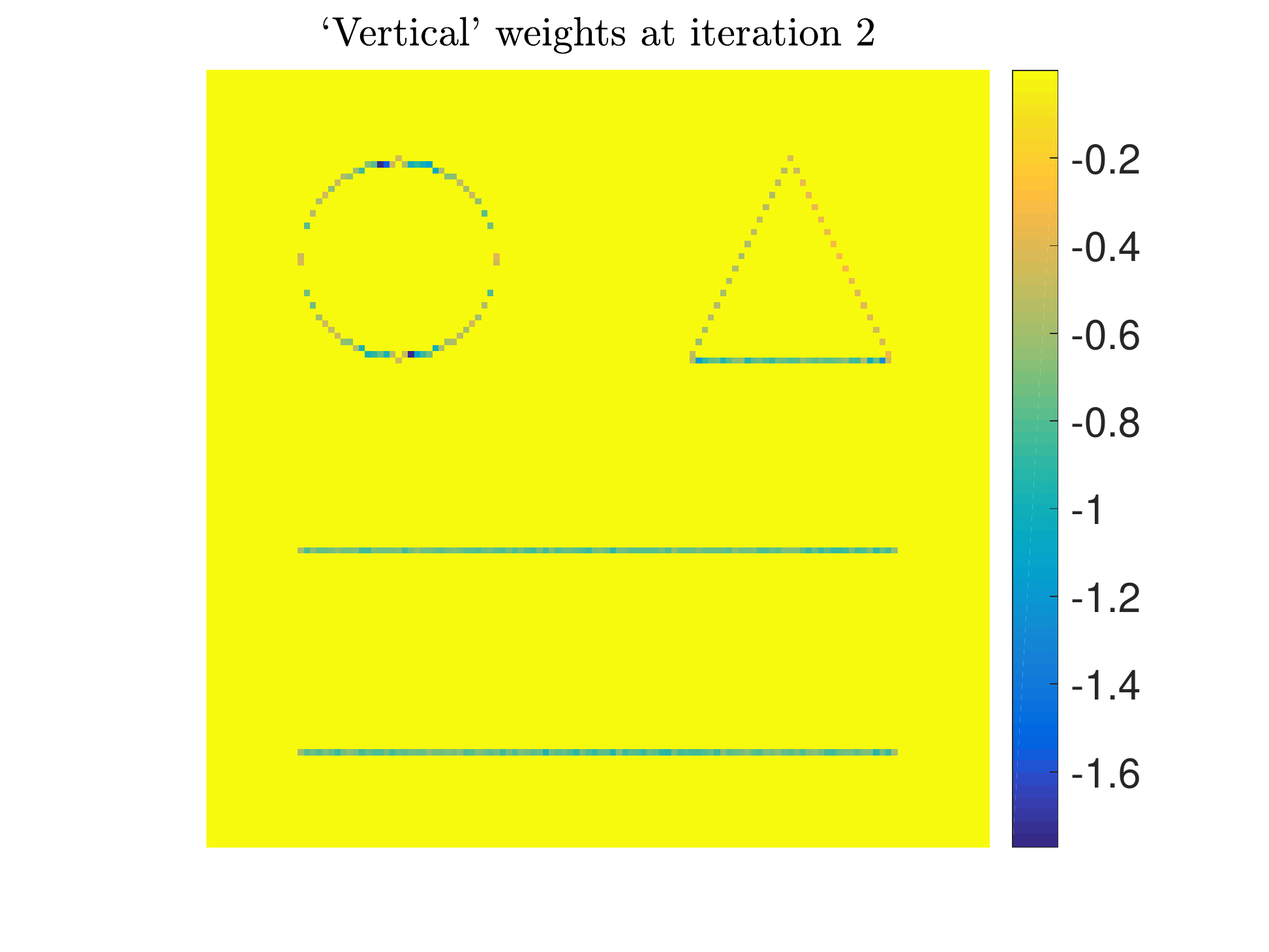} &
\includegraphics[width=5cm]{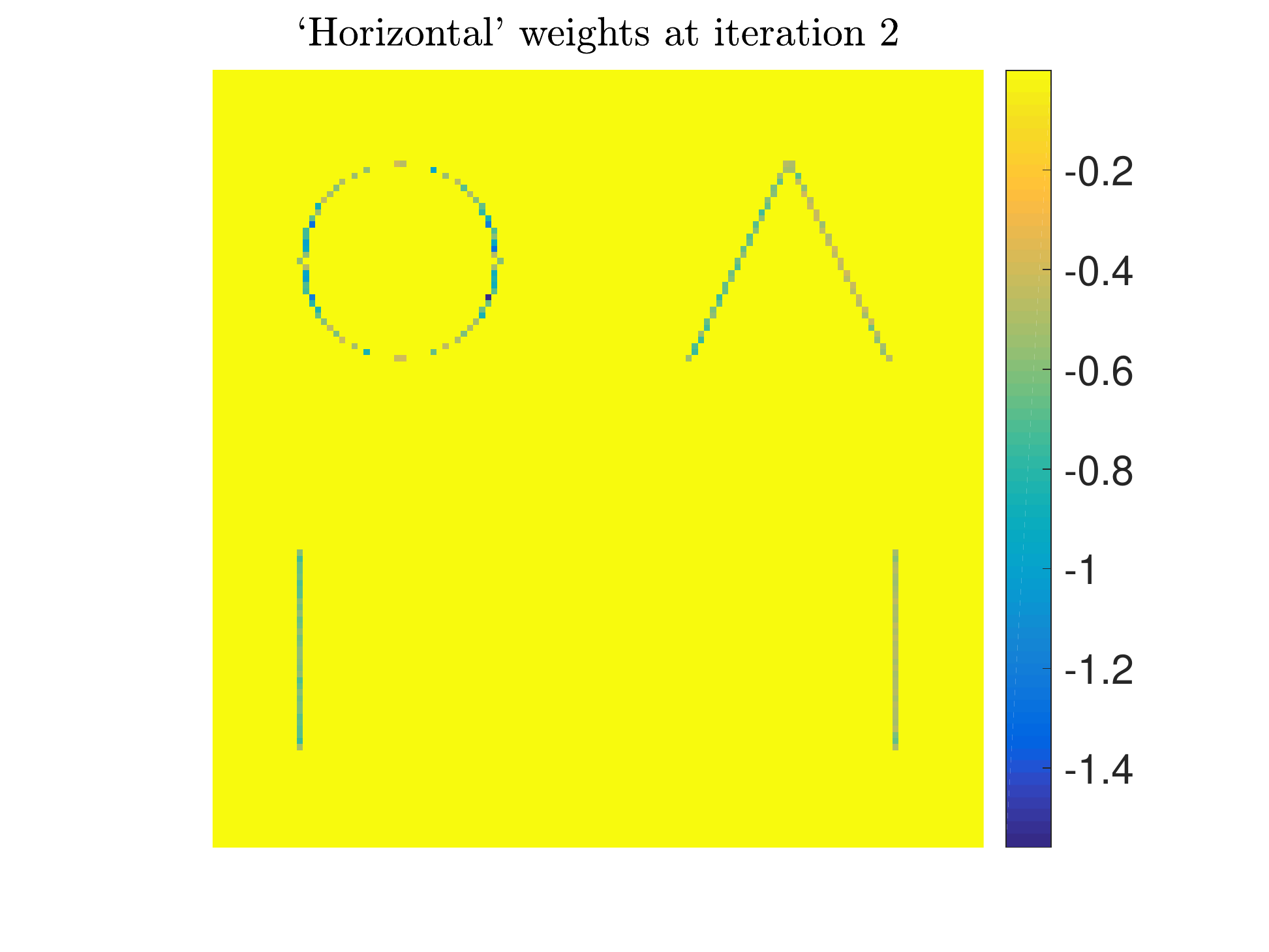} & 
\includegraphics[width=5cm]{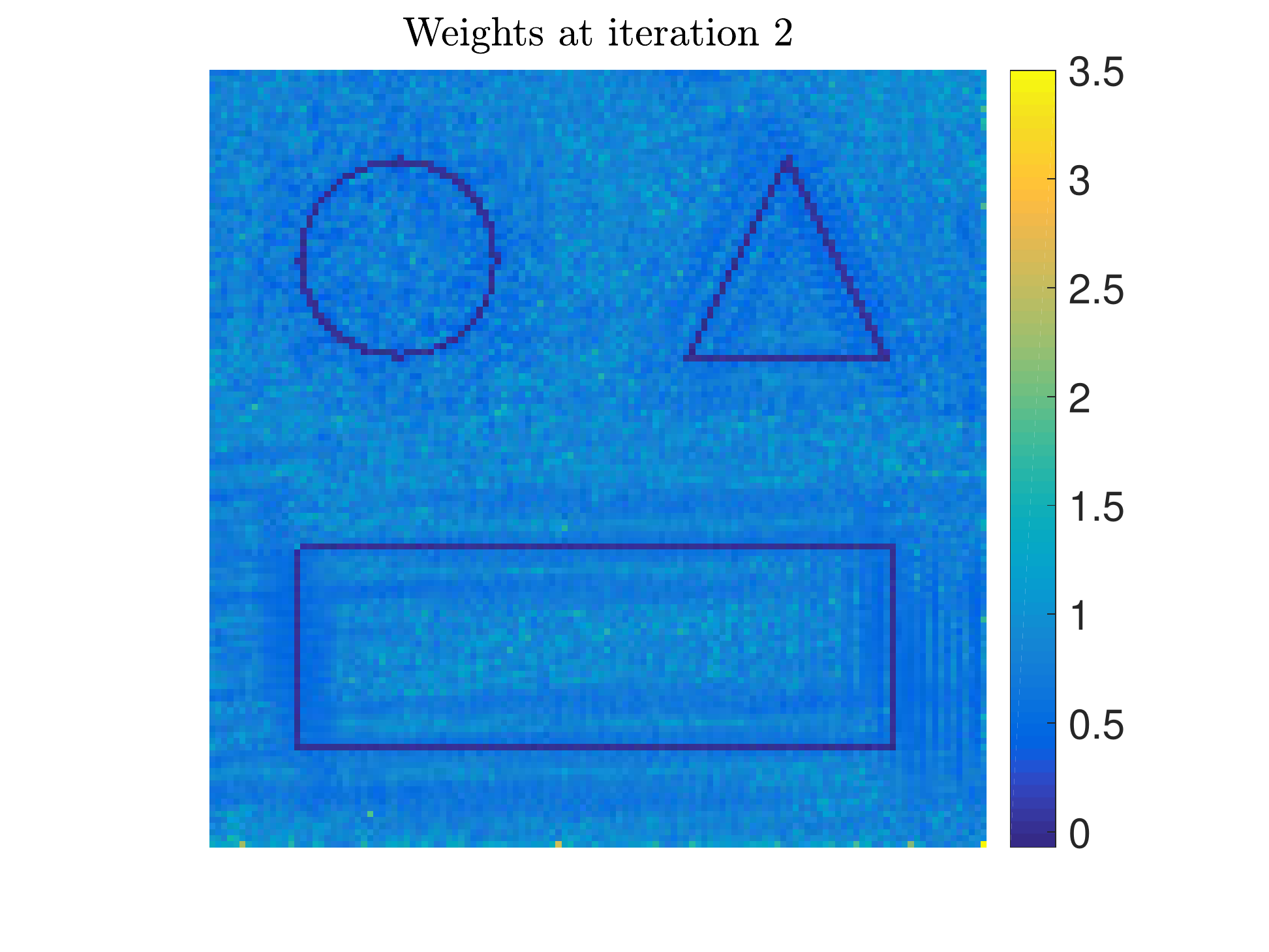}\\ 
\includegraphics[width=5cm]{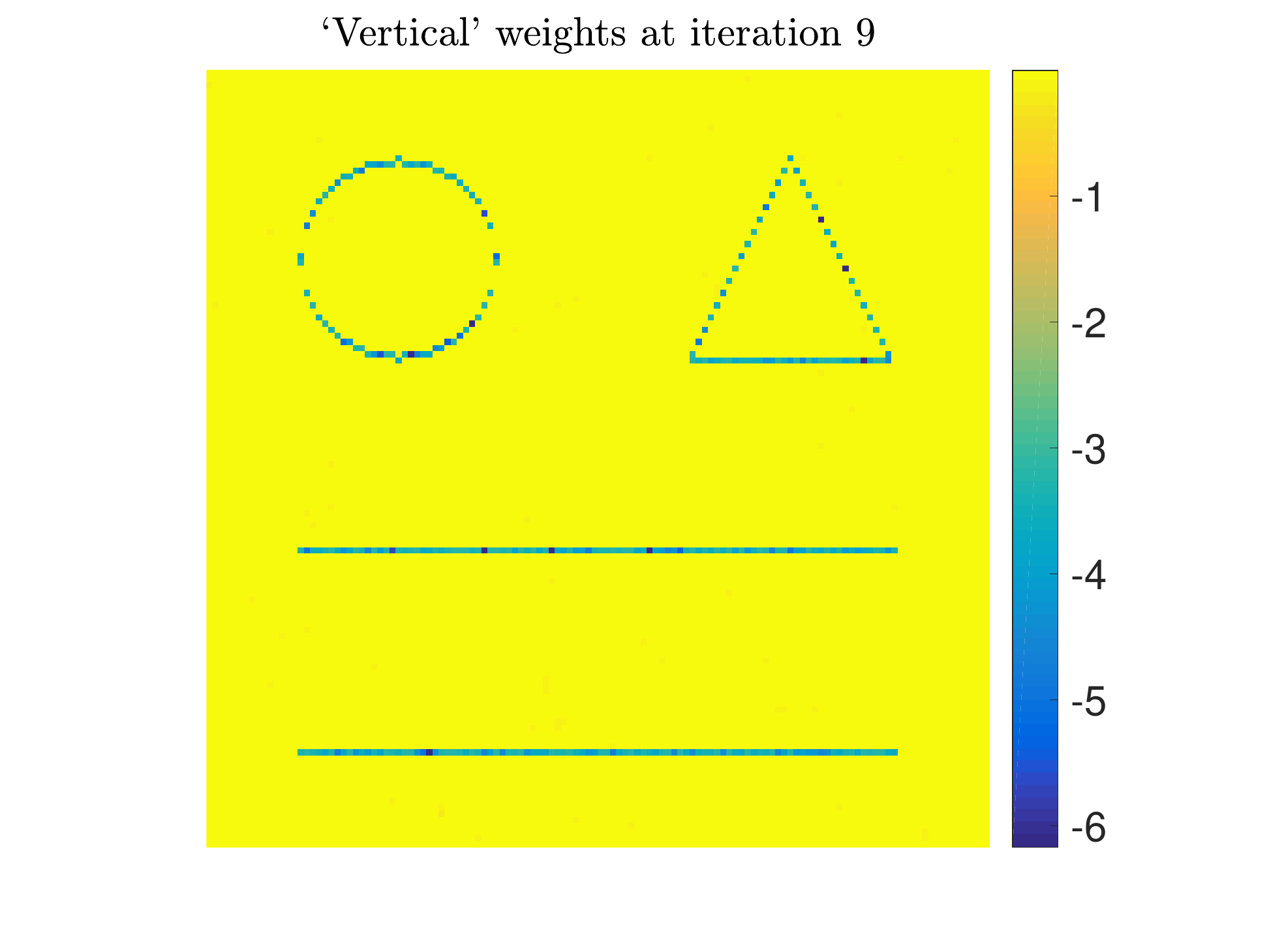} &
\includegraphics[width=5cm]{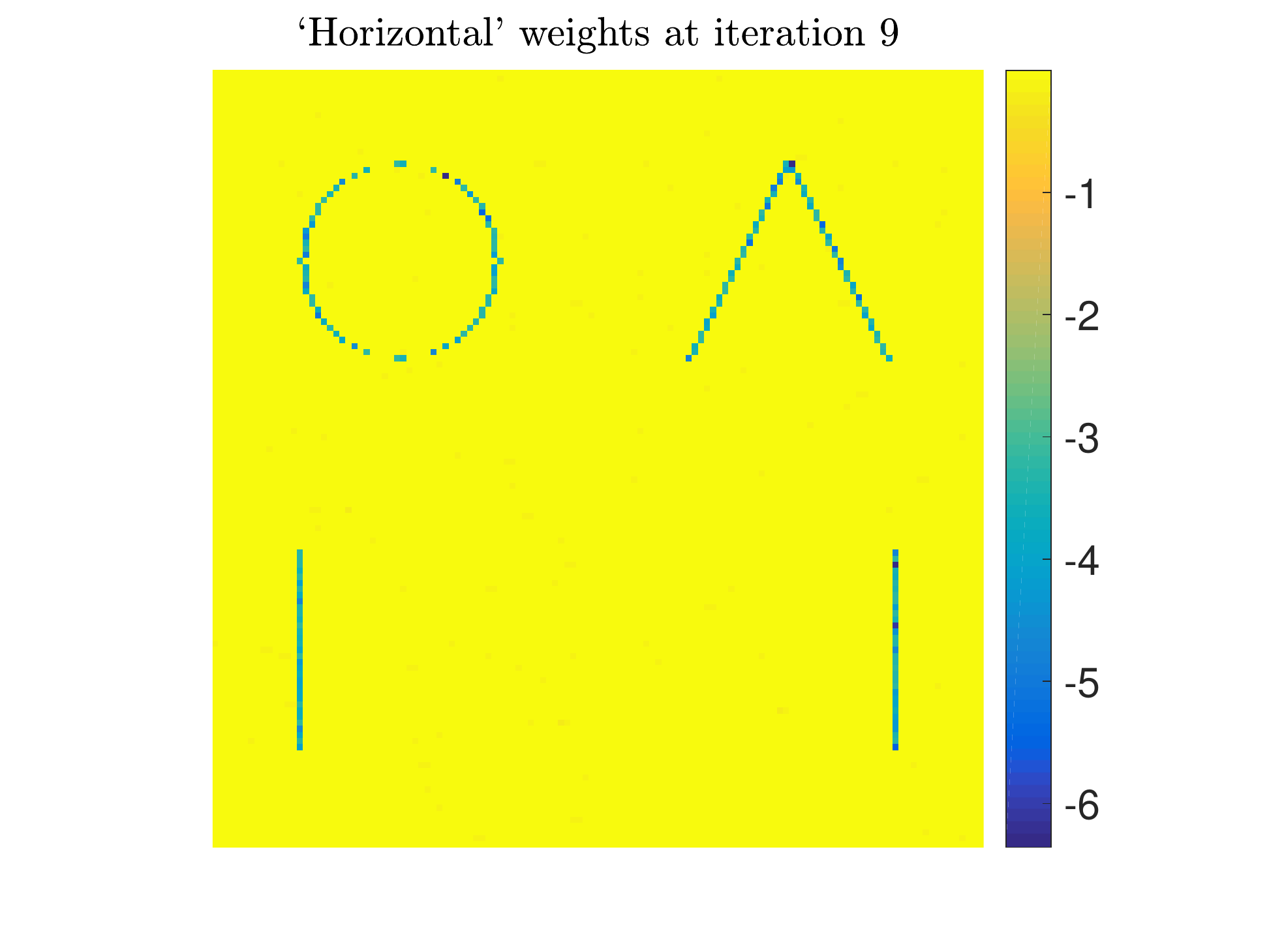} & 
\includegraphics[width=5cm]{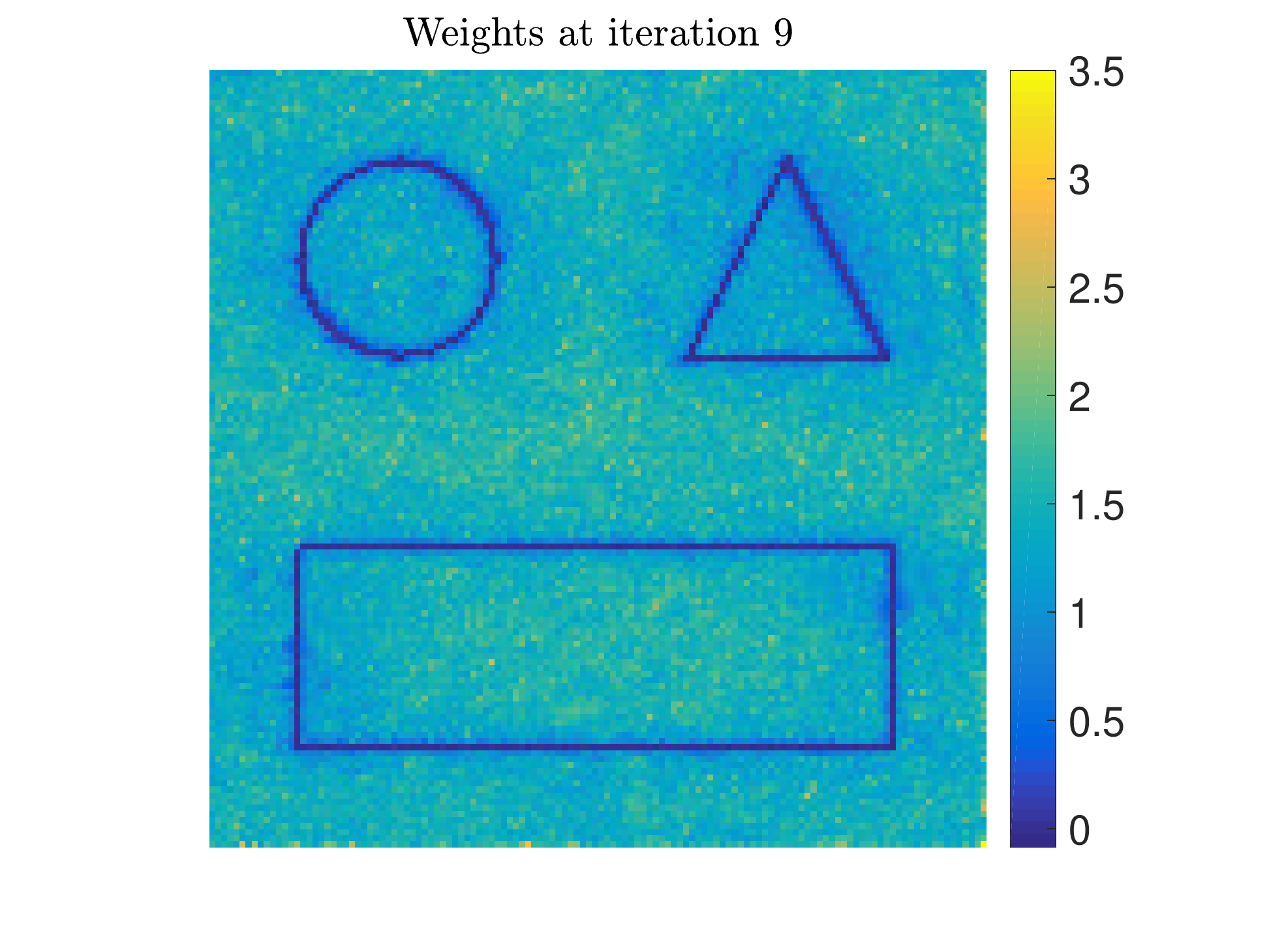}\\ \end{tabular}
\end{center}
\caption{\emph{Shaking blur} test problem. Left column: new weights, to be applied to the vertical derivatives. Middle column: new weights, to be applied to the horizontal derivatives. Right column: IRN-TV weights. The pixel values are displayed in logarithmic scale.}
\label{fig:Blur1TVweight}
\end{figure}

Finally, Figure \ref{fig:Blur1Solutions} displays some relevant reconstructions. We show the initial reconstructions $x^{(\ast,2)}$ obtained at the end of the second inner iteration cycle (i.e., as soon as the iterative reweighting of the regularization term is active), and the reconstructions $x^{(\ast,9)}$ obtained when the maximum number of outer iterations is performed. We consider the inner-outer iterative solvers that employ both the new and the IRN-TV weights, and both the hybrid method (with the ${\mathcal L}$-curve criterion) and CGLS. We can clearly see that, when the new weights together with the hybrid method are used, there is a good improvement in the reconstructions as the outer iterations proceed and, in particular, the final reconstruction is piecewise constant. when the new weights together with CGLS are used, the final reconstruction is overall good, but a few spurious constant patches show up in several locations. When considering the IRN-TV weights with both inner solvers, we cannot see any improvements as the outer iterations proceed (on the contrary, the reconstruction at the end of the outer iterations is much worse than the initial one), and none of the reconstructions properly looks piecewise constant. 
%
%\begin{figure}[htbp]
%\begin{center}
%\includegraphics[width=7cm]{FigsExampleRadon1b/LCurves1b} 
%\end{center}
%\caption{${\mathcal L}$-curves for each iteration, for the second test problem. As implied from
%the text in the plot, the top curve corresponds to the first outer iteration, $\ell = 1$, and the curves
%below this correspond sequentially to iterations $\ell = 2, 3, 4$.  The red circles
%denote corners of each ${\mathcal L}$-curve, which correspond to the chosen
%regularization parameter, $\lambda_{*,\ell}$ for the particular iteration.}
%\label{fig:Radon1bLCurves}
%\end{figure}
%
%Computed reconstructions for the first outer iteration (that is, $x^{(*,1)}$), and for the
%final outer iteration (that is, $x^{(*,4)}$) are shown 
%in Figure~\ref{fig:Radon1bSolutions}. As we can see from these plots, there is a 
%significant improvement in the reconstructions, and in particular 
%the edges at the final outer iteration are much sharper than in the initial outer iteration.
%
Figure \ref{fig:Blur1TVweight} displays the entries of the weight diagonal matrices at outer iterations $\ell=2$ and $\ell=9$, for both the new reweighting strategy (\ref{eq:GradientMap}) and the IRN-TV reweighting strategy. When considering the new reweighting strategy, we can clearly see that edges are properly detected in both the vertical and horizontal derivatives starting from the very first reweighting step: as explained in Section \ref{sec:algorithm}, this greatly contributes to the success of the new algorithm. 
%; the same is not true for the weights to be applied to the horizontal  derivatives, which seem to detect some spurious edges in the first reconstruction $x^{(\ast,1)}$, which are enhanced during the subsequent outer iterations. Despite this, spurious edges are hardly visible in the final reconstruction. 
Similarly to the previous test problems, the IRN-TV weights are not so effective in revealing the structure of the image; very small and oscillating weights are assigned to every location in the image.

\paragraph{Out-of-focus blur} 

\begin{figure}[htbp]
\begin{center}
\begin{tabular}{cc}
\footnotesize{True image} & \footnotesize{Corrupted image}\\
\includegraphics[height=4cm]{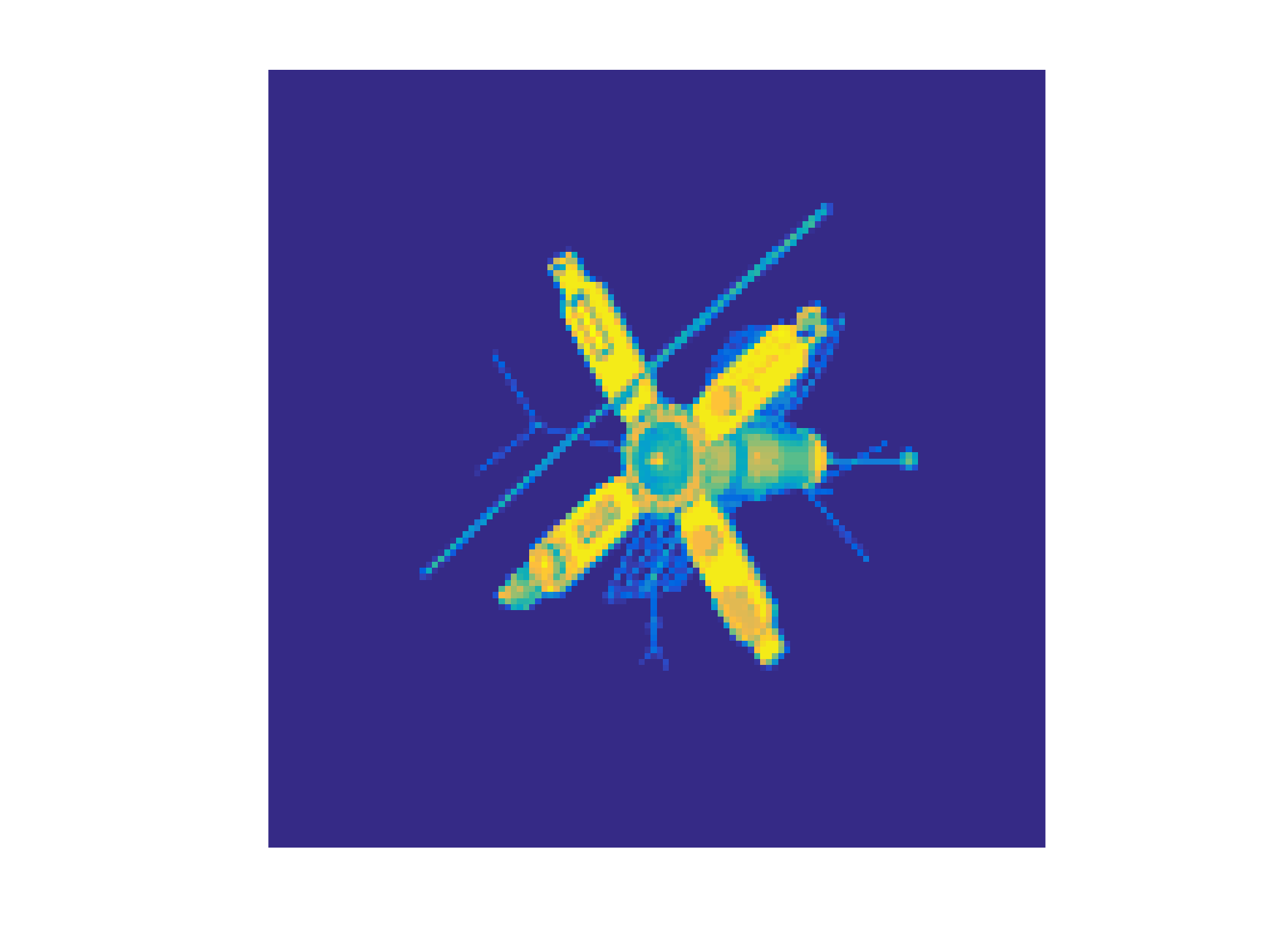} &
\includegraphics[height=4cm]{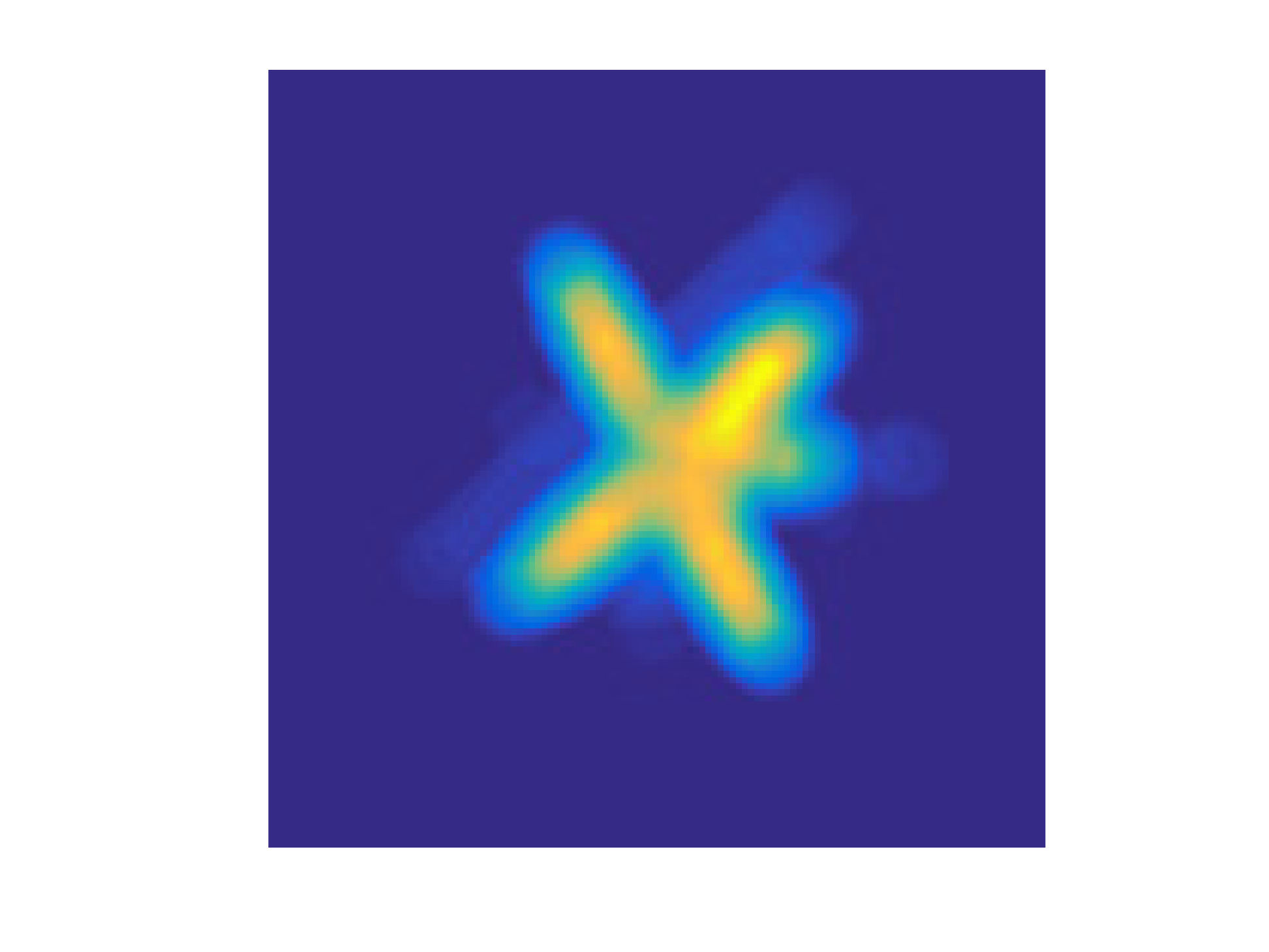}
\end{tabular}
\end{center}
\caption{\emph{Out-of-focus blur} test problem. Left frame: true sharp image $x$. Right frame: the measured data (corrupted image) $b$.}
\label{fig:Blur2Data}
\end{figure}

We use IR Tools to vary some of the options defining the previous image deblurring simulation. Firstly, we would like to change the sharp image to the well-known \emph{satellite} one, of size $128\times 128$ pixels. Secondly, we would like to consider an out-of-focus blur, still of medium intensity. Again, we generate this test problem within IR Tools, using the following MATLAB statements:
\begin{verbatim}
     ProblemOptions = PRset('trueImage', 'satellite', 'BlurLevel', medium); 
     [A, b_true, x_true, ProblemInfo] = PRblurdefocus(128, ProblemOptions);
     b = PRnoise(b_true, NoiseLevel);
\end{verbatim}
Figure~\ref{fig:Blur2Data} shows the true image, along with the
measured data $b$ (which are again corrupted by Gaussian white noise of level $10^{-3}$). 

\begin{figure}[htbp]
\begin{center}
\begin{tabular}{cc}
\includegraphics[width=6cm]{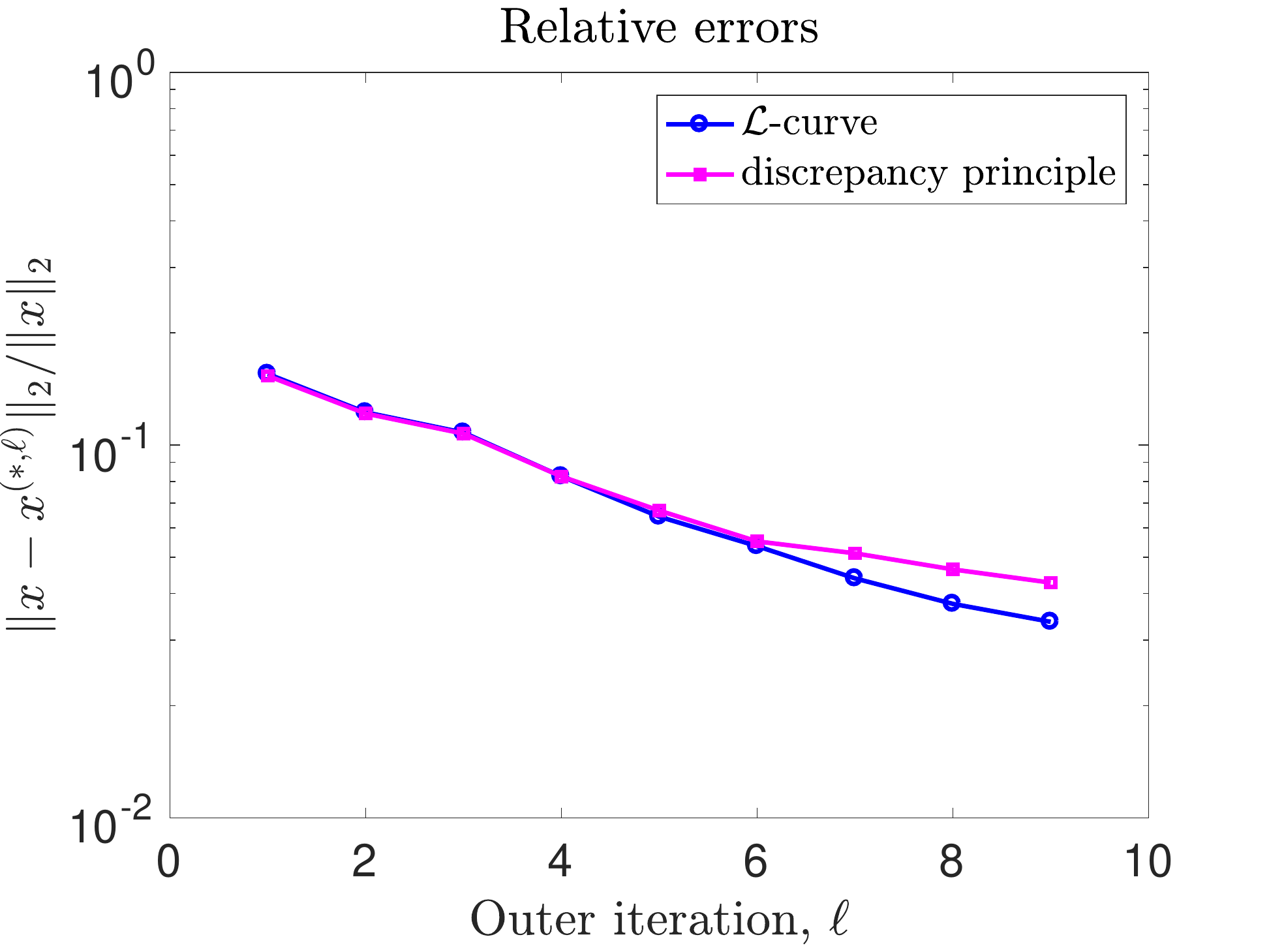} &
\includegraphics[width=6cm]{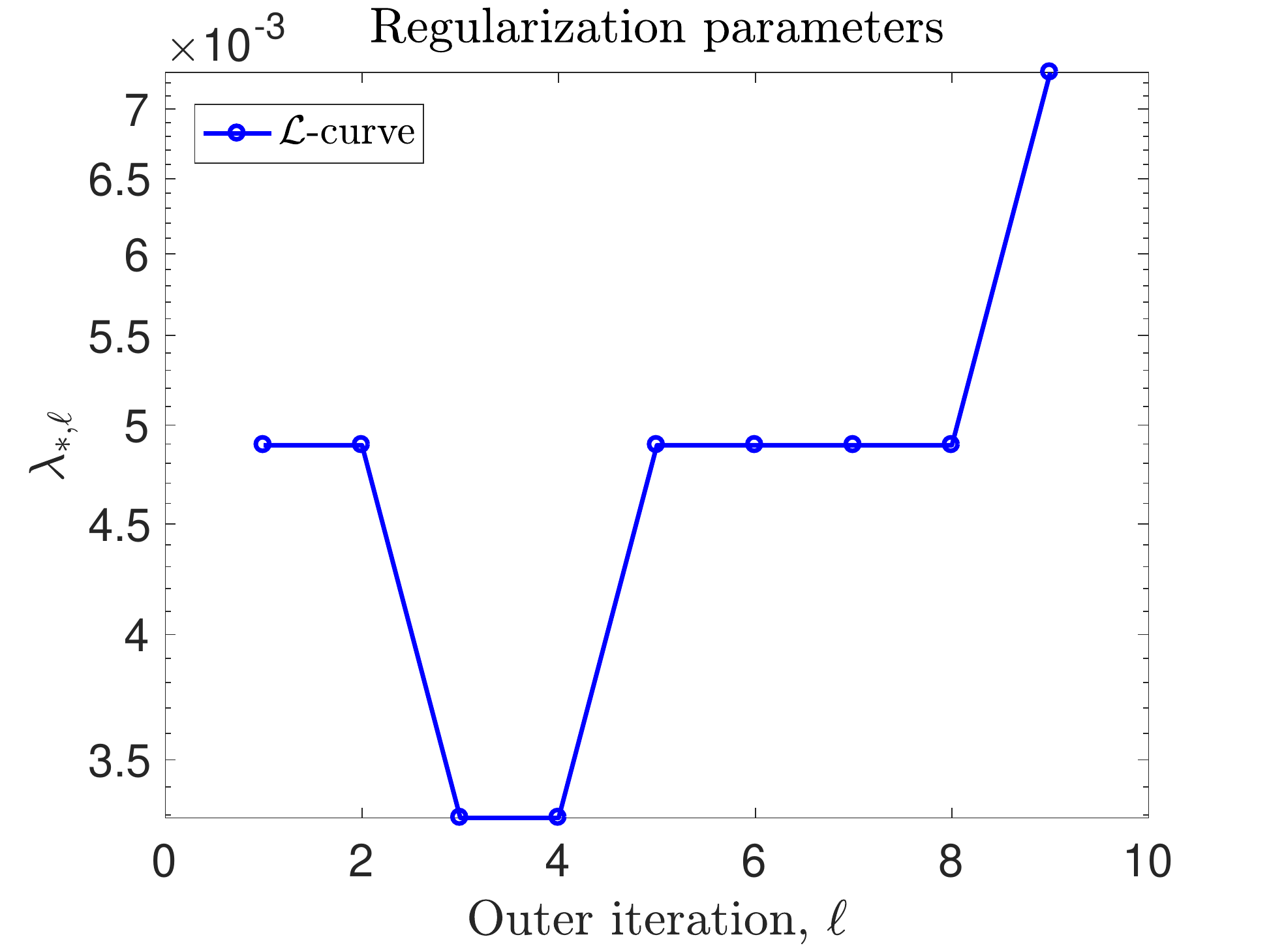}
\end{tabular}
\end{center}
\caption{\emph{Out-of-focus blur} test problem. Relative errors and regularization parameters values at each outer iterations $\ell$,until the stopping criterion is satisfied. Both the discrepancy principle and the ${\mathcal L}$-curve criterion are considered; the regularization parameter computed by the former is always numerically zero, and therefore it is not displayed in the rightmost frame.}
\label{fig:Blur2Iterations}
\end{figure}

As in the previous examples, we first run our algorithm with different parameter choice strategies within the inner hybrid scheme for generalized Tikhonov regularization: in order to fully assess how the reconstructions evolve when he outer iterations proceed, we perform 10 outer iterations, using both the discrepancy principle and the $\mathcal{L}$-curve criterion. Figure~\ref{fig:Blur2Iterations} shows a plot of the relative errors
and chosen regularization parameters at each outer iteration. We can clearly see that the reconstructions computed by both parameter selection strategies are of similar quality, and greatly improve as the outer iterations proceed. The regularization parameters selected by the $\mathcal{L}$-curve are of the order of $10^{-3}$ and slightly oscillate; as for the \emph{shaking blur} test problem, the regularization parameters computed by the discrepancy principle are always either zero or numerically zero (and therefore they are not displayed in the rightmost frame of Figure ~\ref{fig:Blur2Iterations}). 

The regularization parameter selected by the $\mathcal{L}$-curve criterion is always quite small and almost uniform across the outer iterations. Despite the different value of the regularization parameters selected by the two strategies, the quality of the corresponding solutions is almost identical. 

\begin{figure}[htbp]
\begin{center}
\begin{tabular}{cc}
\includegraphics[width=6cm]{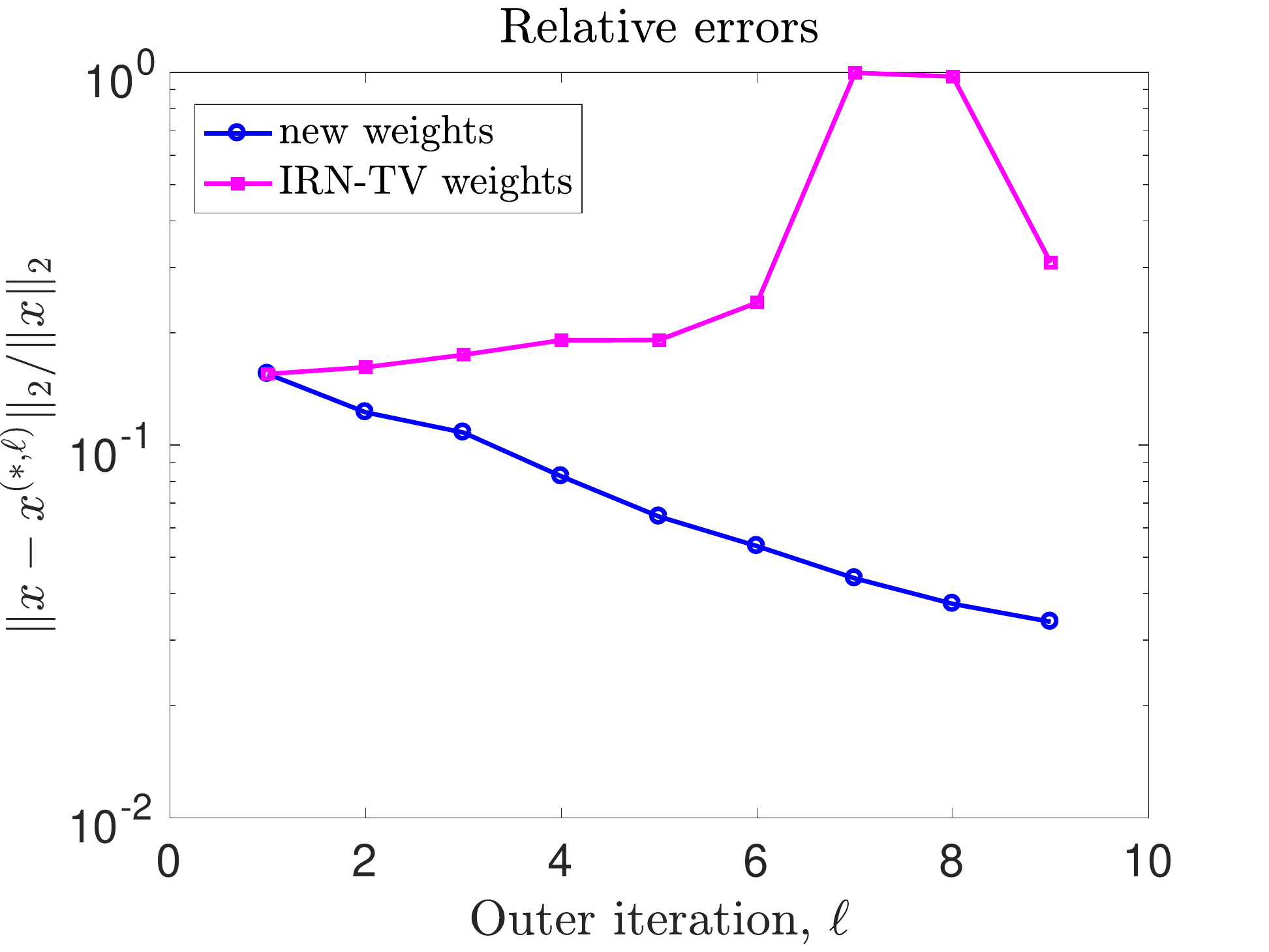} &
\includegraphics[width=6cm]{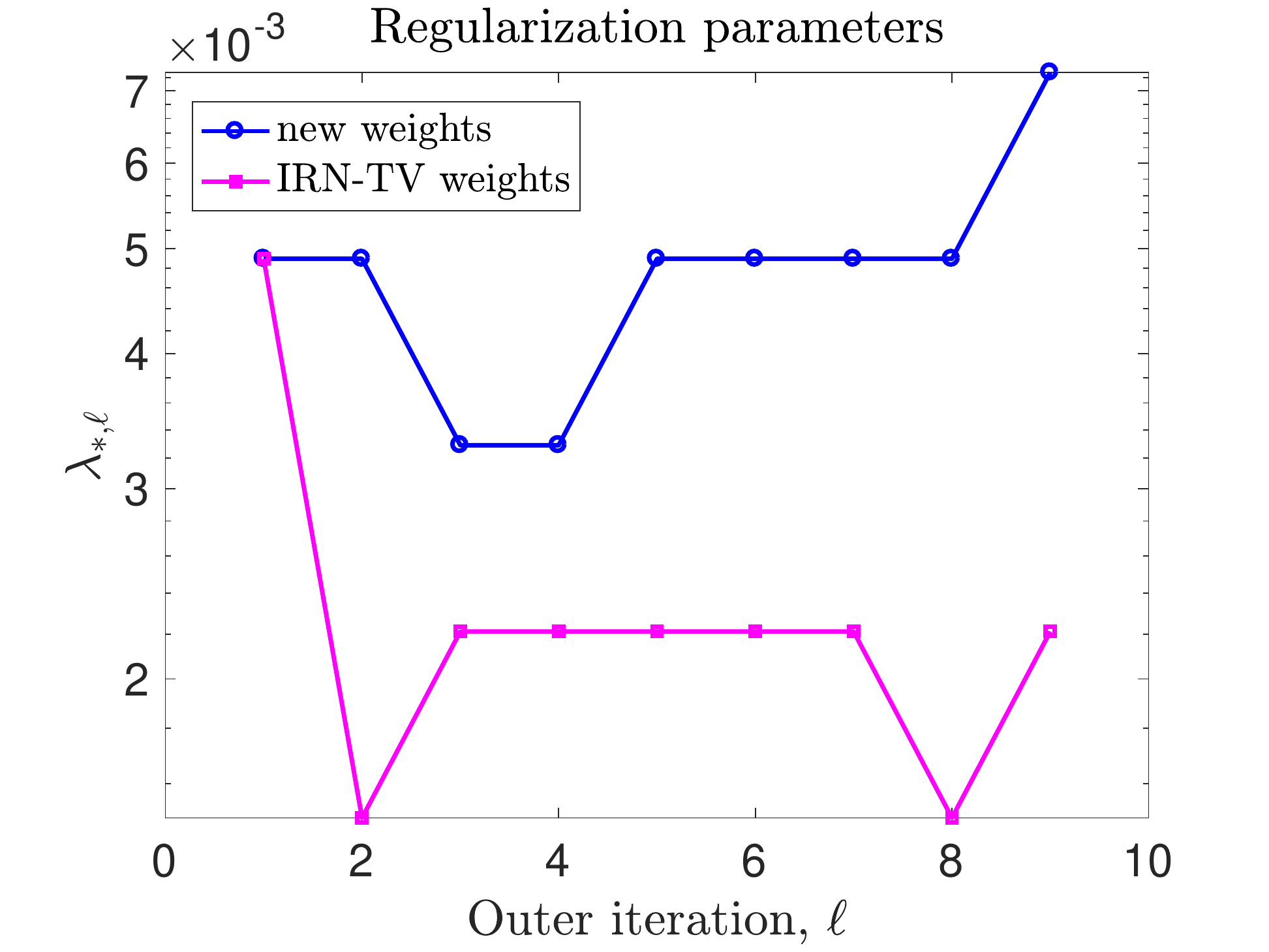}
\end{tabular}
\end{center}
\caption{\emph{Out-of-focus blur} test problem. Relative errors and regularization parameter versus number of (outer) iterations; both methods use the hybrid method for general form Tikhonov during the inner iterations, and adaptively select the regularization parameter according to the discrepancy principle.}
\label{fig:Blur2TV}
\end{figure}

Next, we compare the new method to other inner-outer iterative methods for edge enhancement in imaging. Figure \ref{fig:Blur2TV} displays the relative error and regularization parameter values versus the number of outer iterations for the new method, and for a method that still employs the hybrid solver for general form Tikhonov regularization to handle the inner iterations while using the IRN-TV weights (\ref{eq:IRNTVweights}) at each outer iteration. For both methods, the ${\mathcal L}$-curve criterion is employed to adaptively choose the regularization parameter at each inner iteration. We can clearly see that, when the IRN-TV weights are used, the behavior of the relative errors is quite hectic and the quality of the reconstructions keeps deteriorating as the outer iterations proceed (dramatically so during the last iterations). The regularization parameter for IRN-TV is always around a value of the order of $10^{-3}$ (as in the case of the new weights). 
%%%; indeed, for this test problem, setting the regularization parameter through the $\mathcal{L}$-curve results in a more accurate reconstructions
%% Similarly to the previous test problems, both methods perform well during the first (outer) iterations, with the IRN-TV weights being quite effective in reducing the error; however, the quality of the IRN-TV solutions eventually stagnates, while the new method keeps improving. 
%%%Similarly to the previous test problems, both methods perform well during the first (outer) iterations, with the IRN-TV weights being quite effective in reducing the error; however, the quality of the IRN-TV solutions eventually stagnates, while the new method keeps improving. 
%%The automatically selected regularization parameter for IRN-TV keeps oscillating between values of the order of $10^{-2}$ and values of the order of $10^{-4}$. 
% computes  the new method is amazing at the end.

\begin{figure}[htbp]
\begin{center}
\begin{tabular}{cc}
\includegraphics[width=6cm]{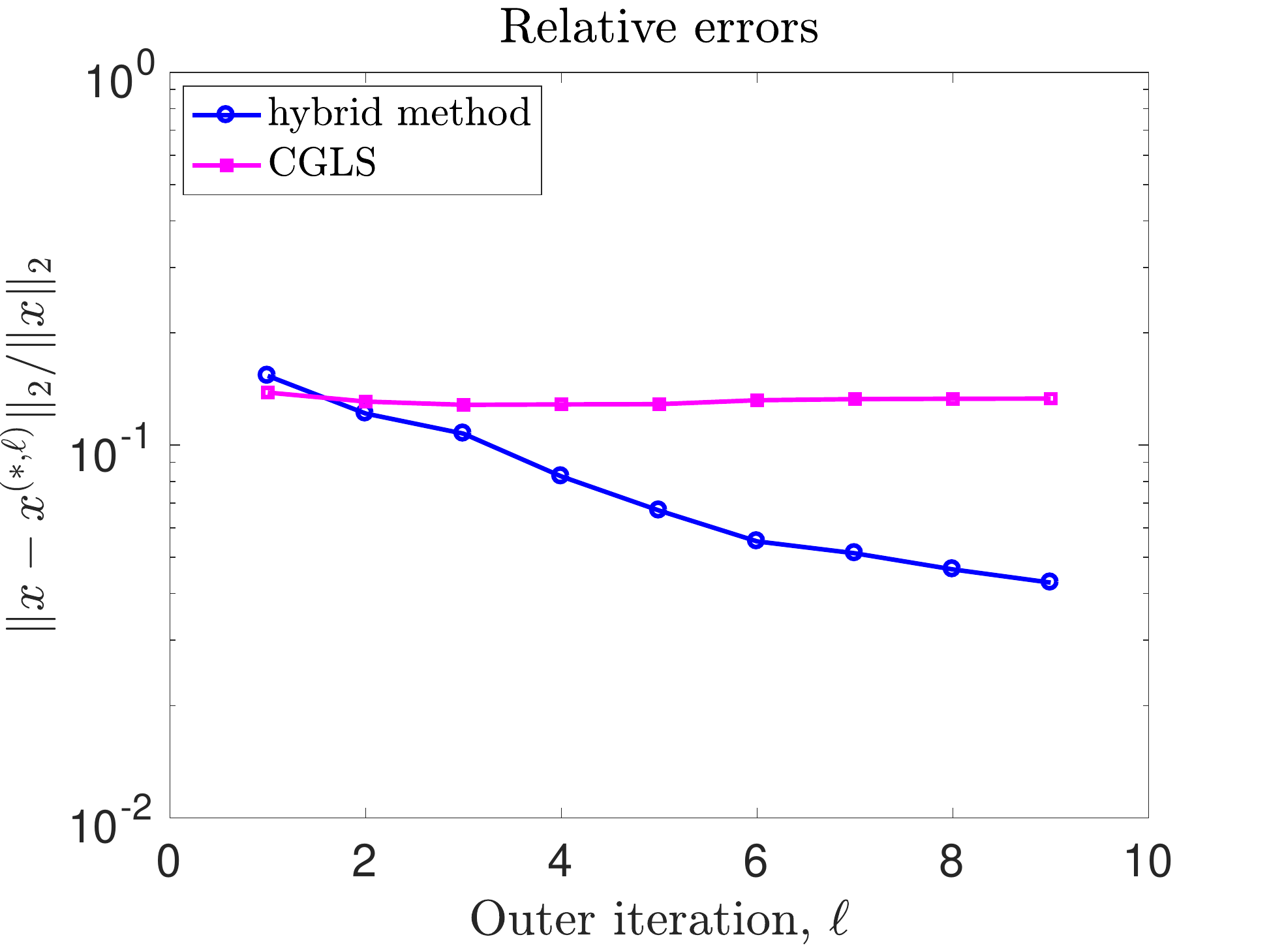} &
\includegraphics[width=6cm]{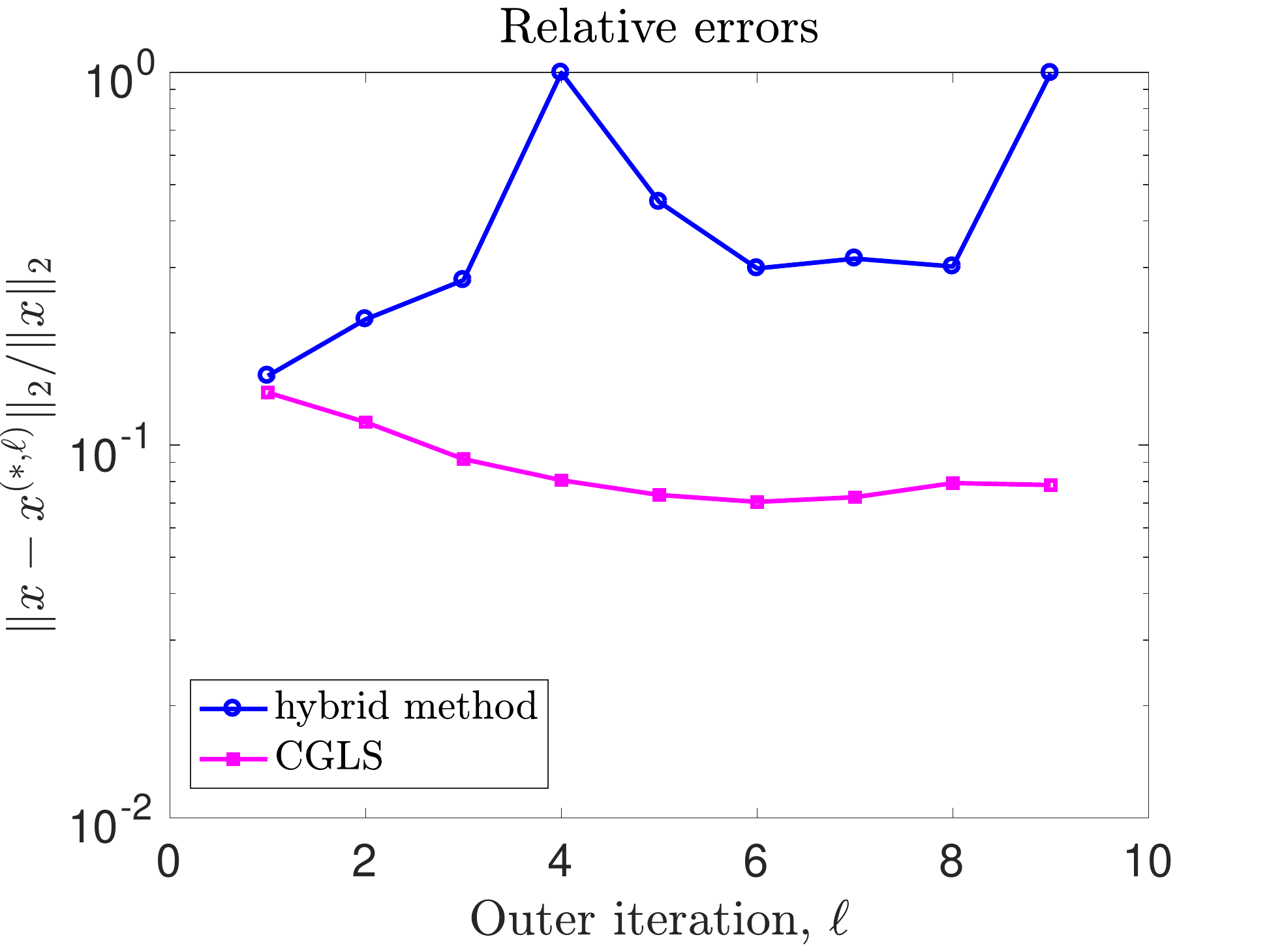}
\end{tabular}
\end{center}
\caption{\emph{Out-of-focus blur} test problem. Relative errors versus number of (outer) iterations. Left frame: the new weights are used. Right frame: the IRN-TV weights are used.}
\label{fig:Blur2TV_comparisons}
\end{figure}

Figure \ref{fig:Blur2TV_comparisons} assesses the influence of the inner iterative solver on the overall behavior of the method. Namely, we consider the inner-outer iterative schemes implemented with both the new and the IRN-TV weights, and with both the hybrid and the CGLS methods as inner solvers. The hybrid method chooses the regularization parameter adaptively according to the $\mathcal{L}$-curve criterion as the iterations proceed, and the value $\lambda_{\ast,\ell}$ selected when the $\ell$th inner iteration cycle terminates is taken to be the fixed regularization parameter to be set in advance of the $\ell$th CGLS iteration cycle. 
% requires a fixed value of the regularization parameter to be available in advance of each cycle of inner iterations, we first run the methods based on the hybrid solver for general form Tikhonov that adaptively chooses the regularization parameter at each inner iteration according to the discrepancy principle: in this way, a parameter 
% $\lambda_{\ast,\ell}$ is eventually set at the $\ell$th outer iteration. When running the methods based on CGLS, we take $\lambda_{\ast,\ell}$ as regularization parameter for the $\ell$th outer iteration.
Looking at the leftmost frame of Figure \ref{fig:Blur2TV_comparisons} we can see that CGLS stagnates, i.e., no sensible improvements in the quality of the reconstructions happen when the outer iterations proceed. Looking at the rightmost frame of Figure \ref{fig:Blur2TV_comparisons}, we can see that CGLS performs much better than the hybrid method when the IRN-TV weights are considered: although the quality of the reconstructions computed by the former are not dramatic, the latter degenerates. 

\begin{figure}[htbp]
\begin{center}
\begin{tabular}{ccc}
\hspace{-0.3cm}\includegraphics[width=5.5cm]{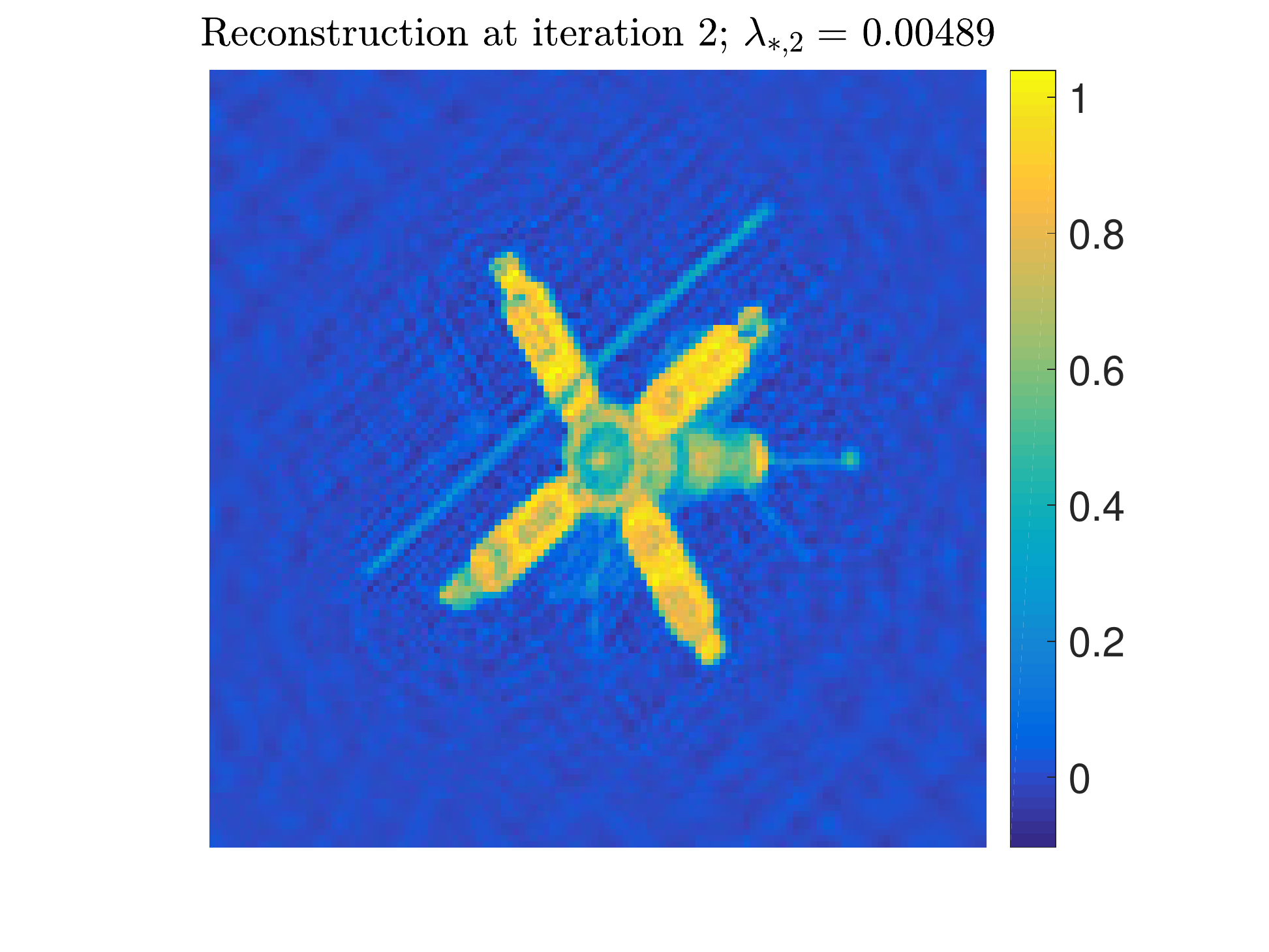} &
\hspace{-0.3cm}\includegraphics[width=5.5cm]{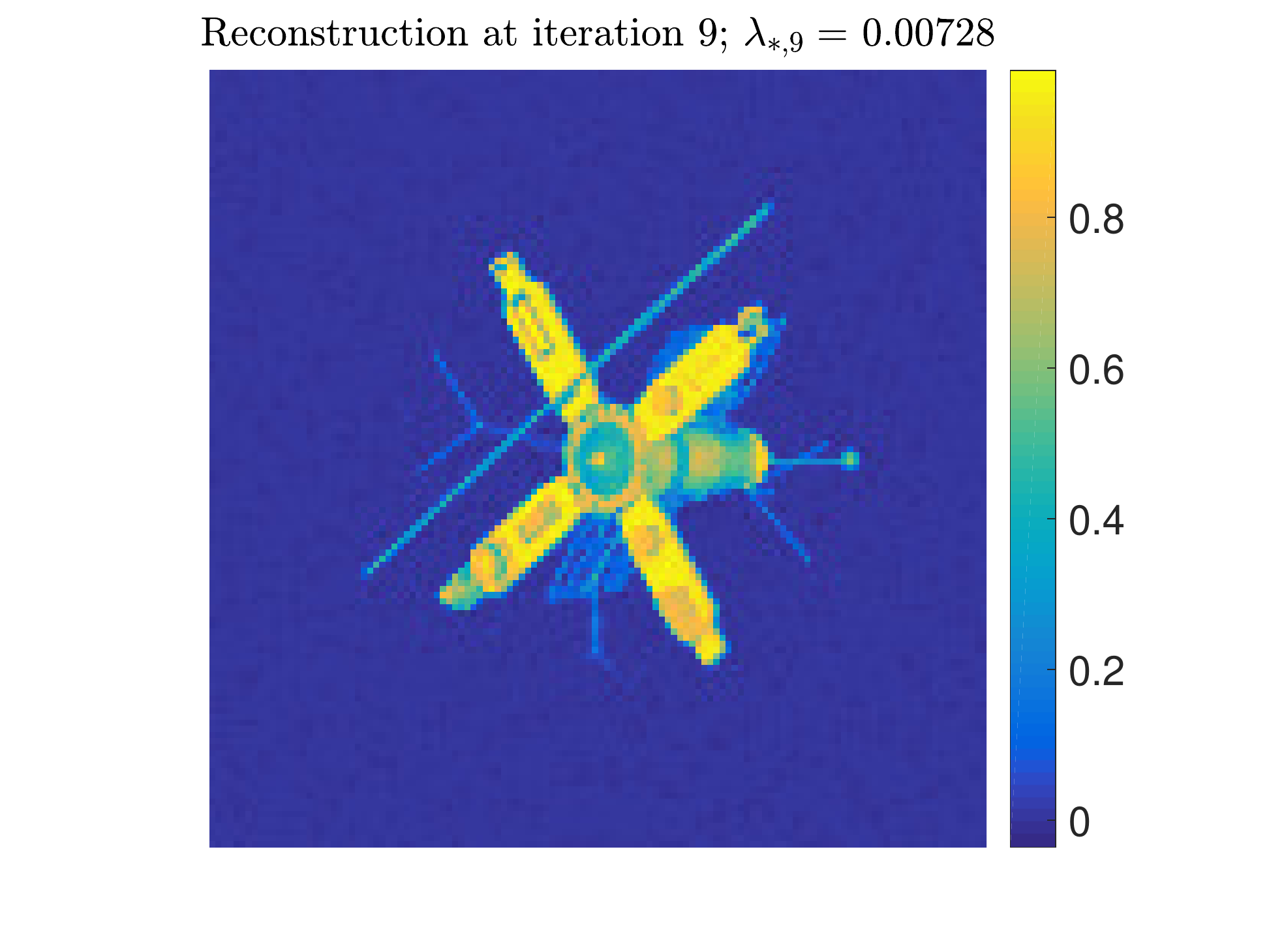} &
\hspace{-0.3cm}\includegraphics[width=5.5cm]{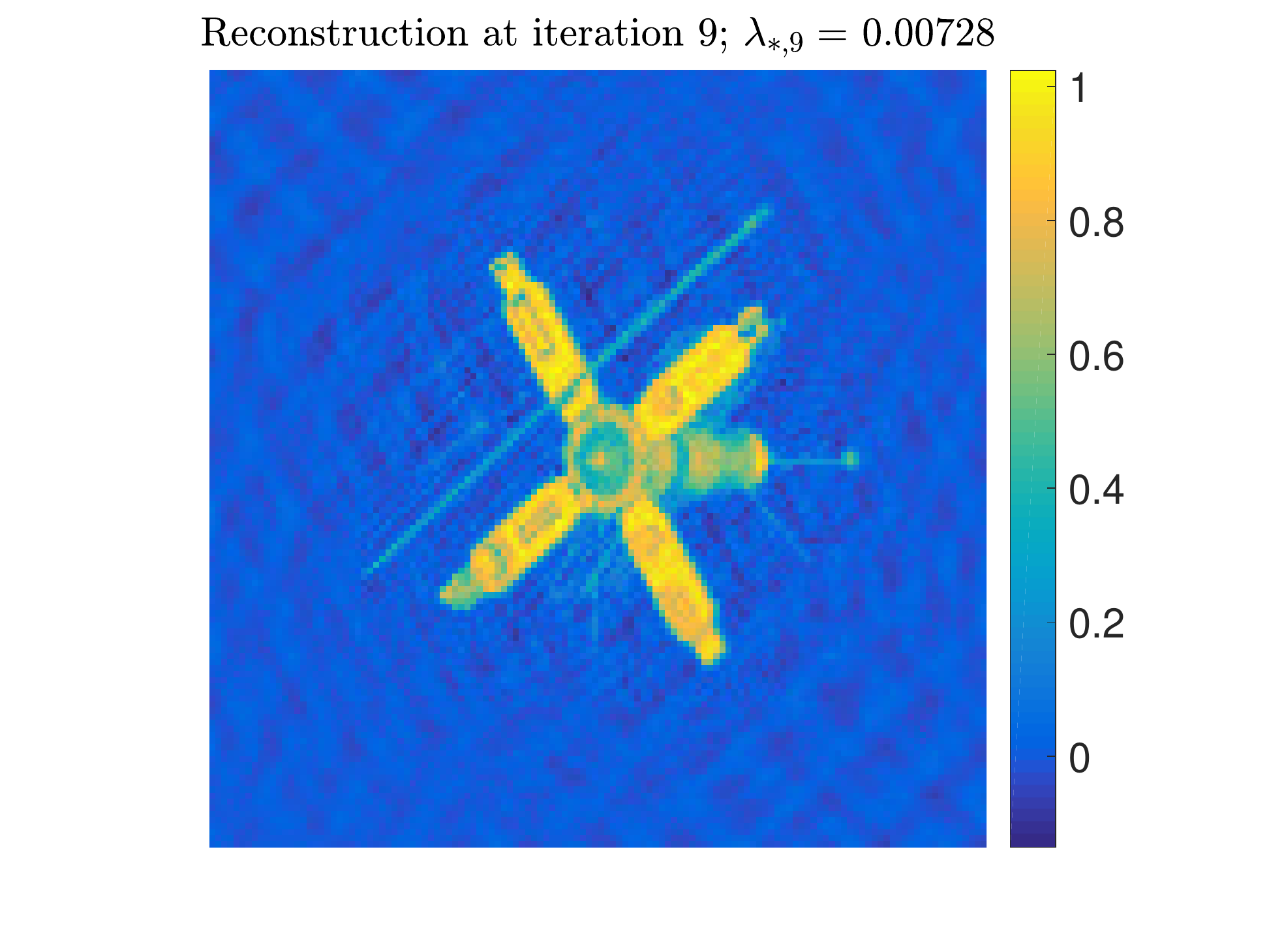}\\
\hspace{-0.3cm}\includegraphics[width=5.5cm]{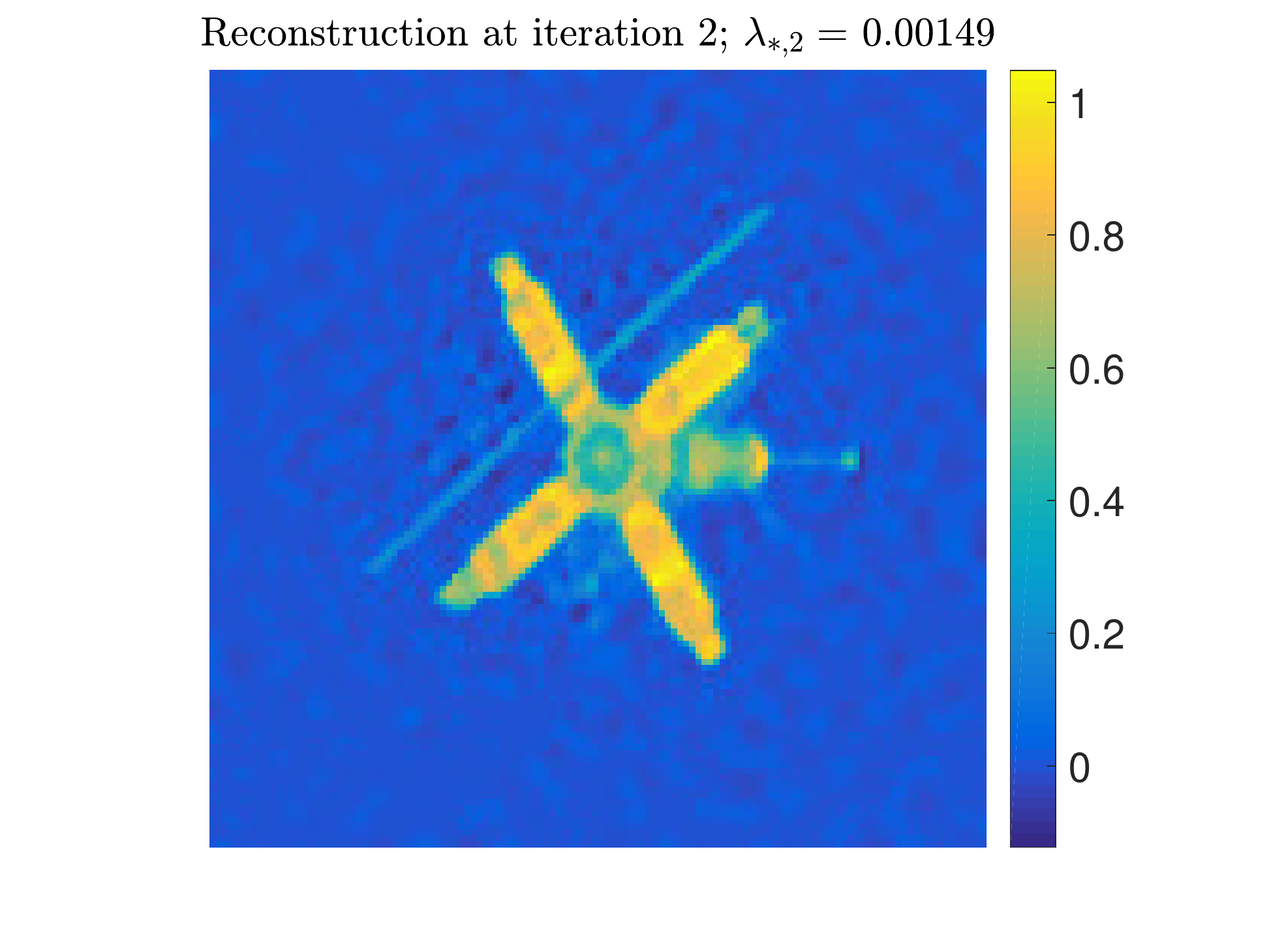} &
\hspace{-0.3cm}\includegraphics[width=5.5cm]{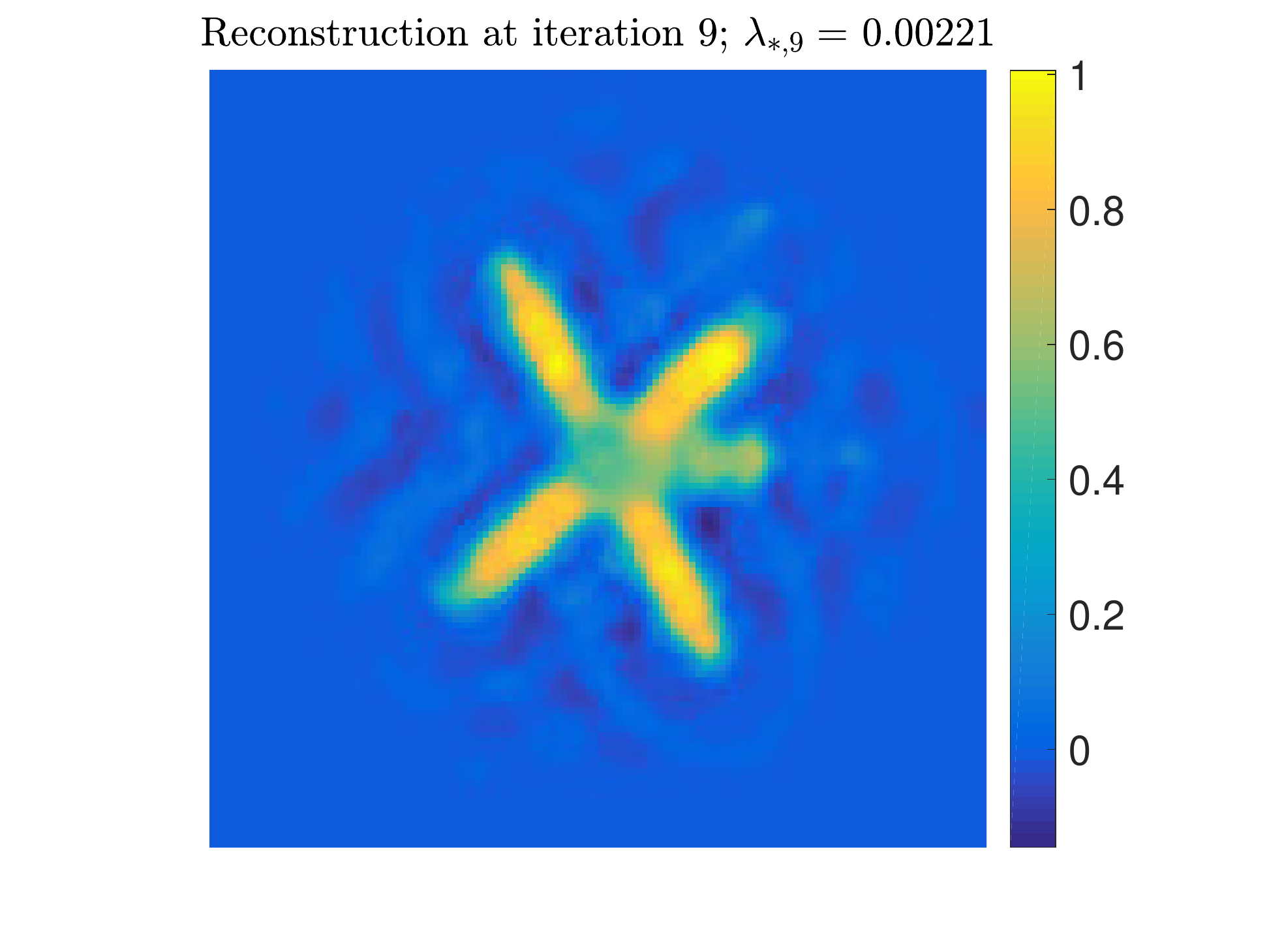} &
\hspace{-0.3cm}\includegraphics[width=5.5cm]{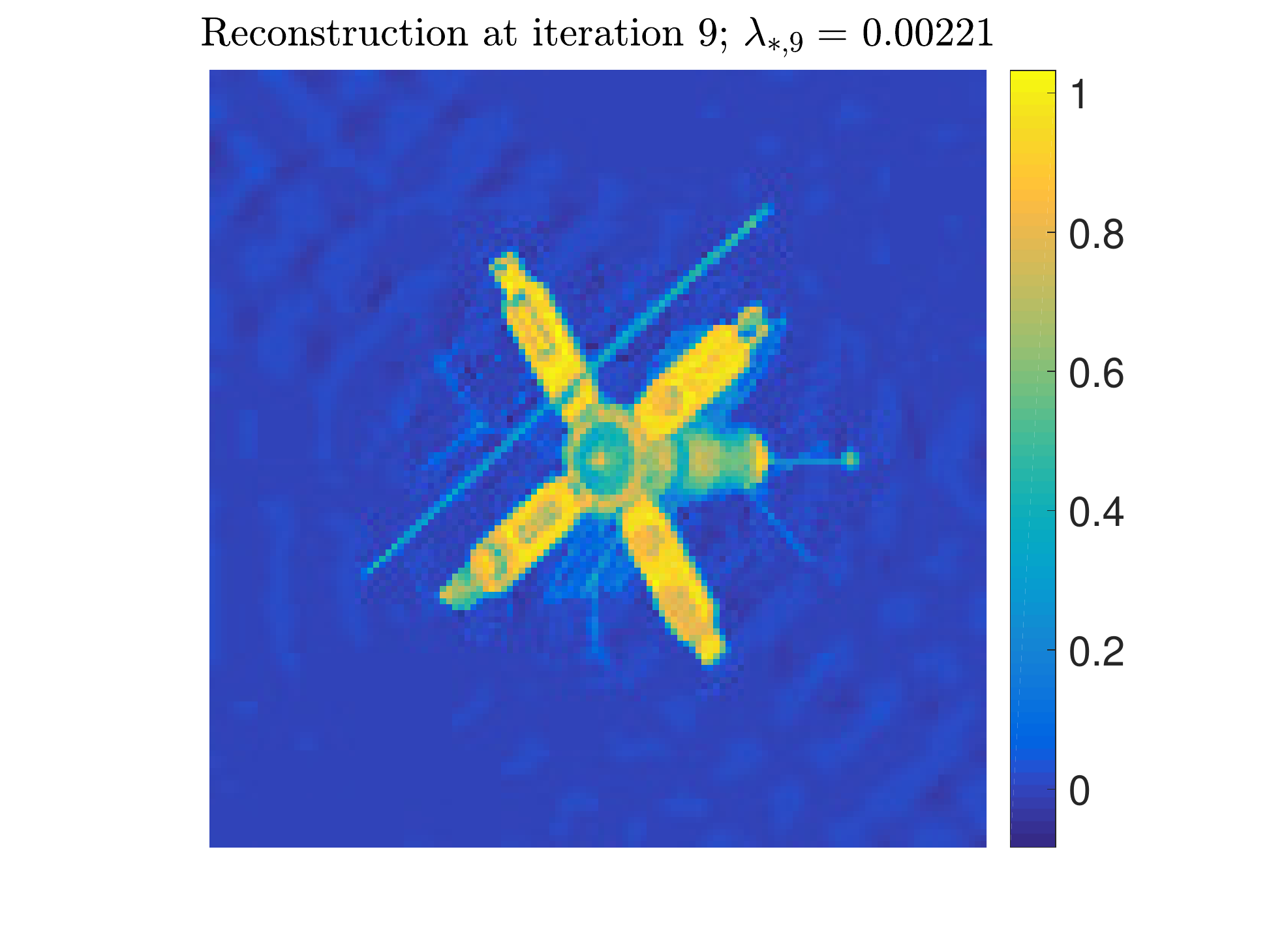}
\end{tabular}
\end{center}
\caption{\emph{Out-of-focus blur} test problem. Upper row: reconstructions obtained using the new weights, at different outer iterations, and by different inner linear solvers. Lower row: reconstructions obtained using the IRN-TV weights, at different outer iterations, and by different inner linear solvers. The corresponding regularization parameters chosen
by the hybrid method, and used within CGLS, are displayed above each image.}
\label{fig:Blur2Solutions}
\end{figure}

Finally, Figure \ref{fig:Blur2Solutions} displays some relevant reconstructions. We show the initial reconstructions $x^{(\ast,2)}$ obtained at the end of the second inner iteration cycle (i.e., as soon as the iterative reweighting of the regularization term is active), and the reconstructions $x^{(\ast,9)}$ obtained when the maximum number of outer iterations is performed. We consider the inner-outer iterative solvers that employ both the new and the IRN-TV weights, and both the hybrid method (with the ${\mathcal L}$-curve criterion) and CGLS. We can clearly see that, when the new weights together with the hybrid method are used, there is an excellent improvement in the reconstructions as the outer iterations proceed, although the last reconstruction is still slightly irregular (and one may consider performing additional outer iterations). When the new weights together with CGLS are used, the final reconstruction is almost identical to the one displayed in the leftmost top frame of Figure \ref{fig:Blur2Solutions} (as expected; see also Figure \ref{fig:Blur2TV_comparisons}). When considering the IRN-TV weights with the hybrid method, the solution ends up being oversmoothed. When considering the IRN-TV weights with CGLS, we have a slight improved final reconstruction (as expected; see again Figure \ref{fig:Blur2TV_comparisons}).  
%
%\begin{figure}[htbp]
%\begin{center}
%\includegraphics[width=7cm]{FigsExampleRadon1b/LCurves1b} 
%\end{center}
%\caption{${\mathcal L}$-curves for each iteration, for the second test problem. As implied from
%the text in the plot, the top curve corresponds to the first outer iteration, $\ell = 1$, and the curves
%below this correspond sequentially to iterations $\ell = 2, 3, 4$.  The red circles
%denote corners of each ${\mathcal L}$-curve, which correspond to the chosen
%regularization parameter, $\lambda_{*,\ell}$ for the particular iteration.}
%\label{fig:Radon1bLCurves}
%\end{figure}
%
%Computed reconstructions for the first outer iteration (that is, $x^{(*,1)}$), and for the
%final outer iteration (that is, $x^{(*,4)}$) are shown 
%in Figure~\ref{fig:Radon1bSolutions}. As we can see from these plots, there is a 
%significant improvement in the reconstructions, and in particular 
%the edges at the final outer iteration are much sharper than in the initial outer iteration.
%
\begin{figure}[htbp]
\begin{center}
\begin{tabular}{ccc}
\includegraphics[width=5cm]{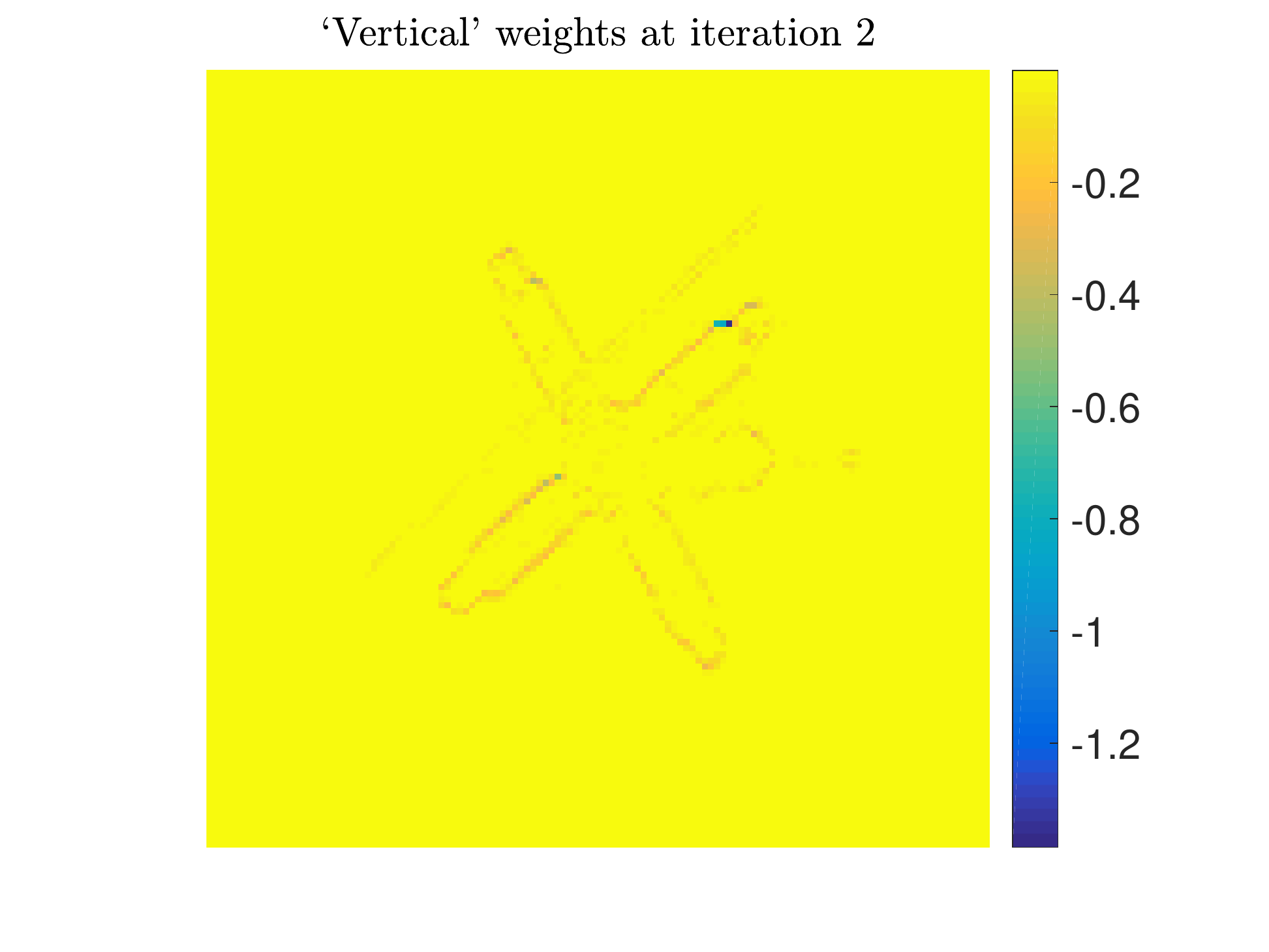} &
\includegraphics[width=5cm]{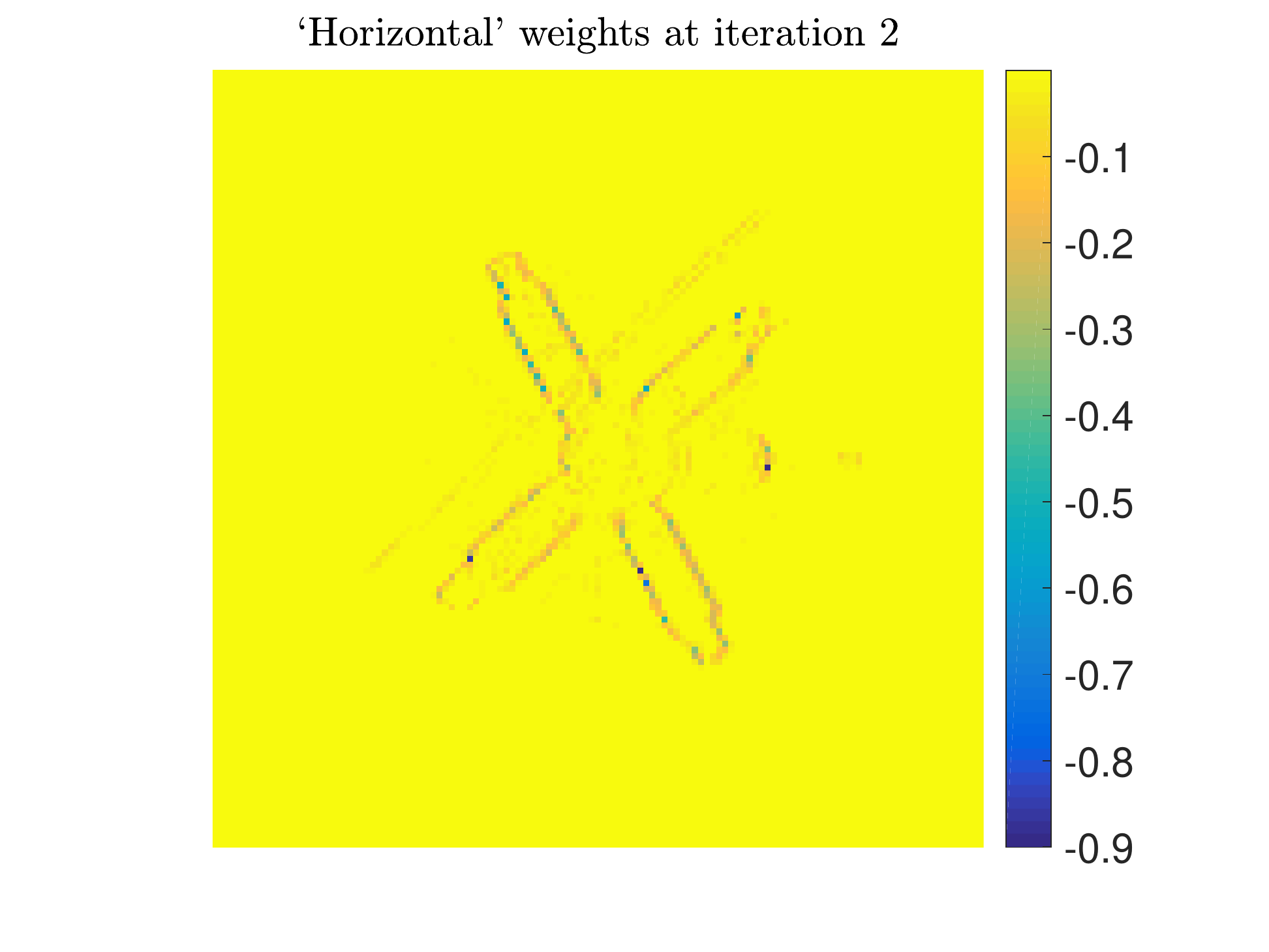} & 
\includegraphics[width=5cm]{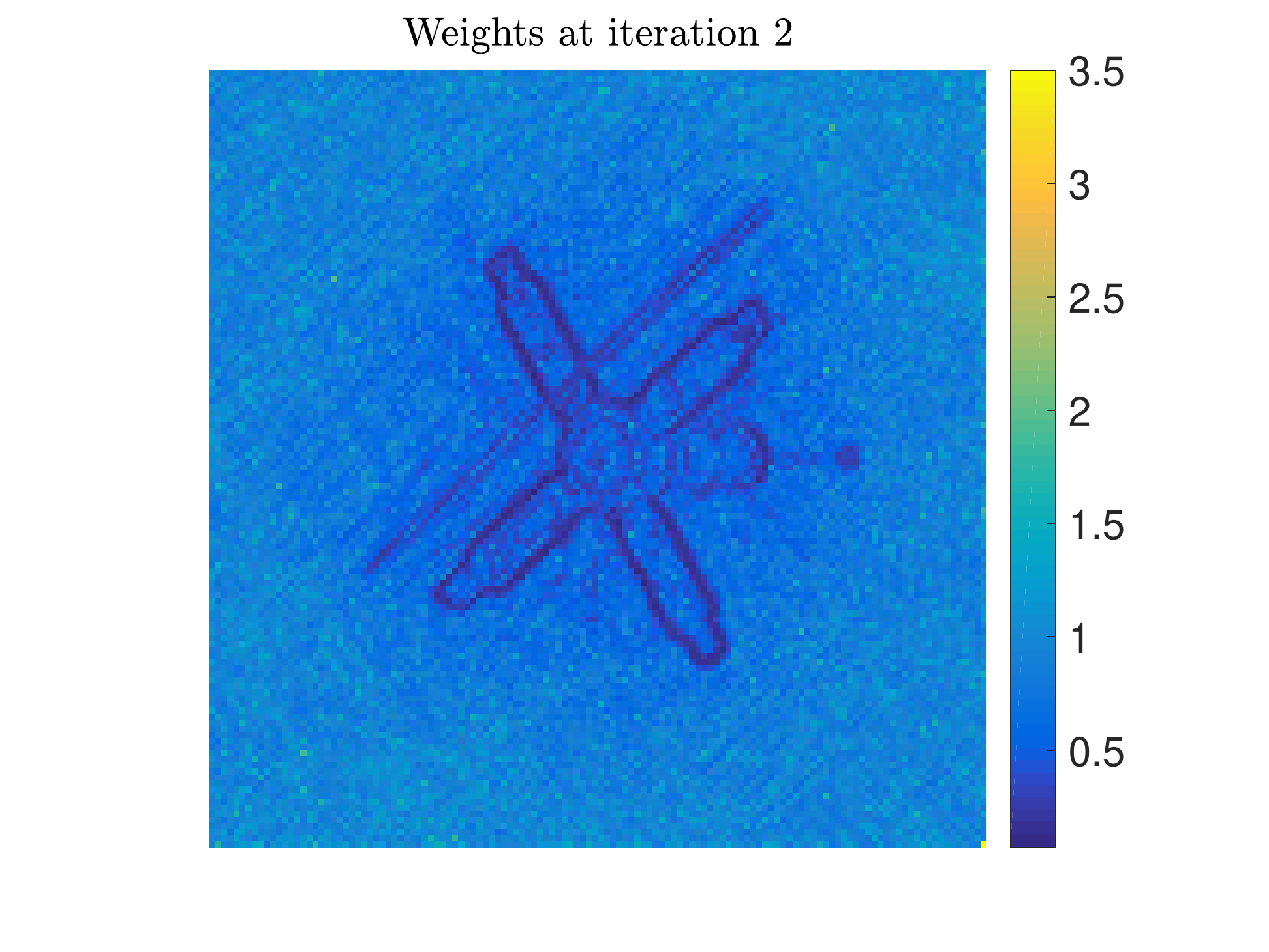}\\ 
\includegraphics[width=5cm]{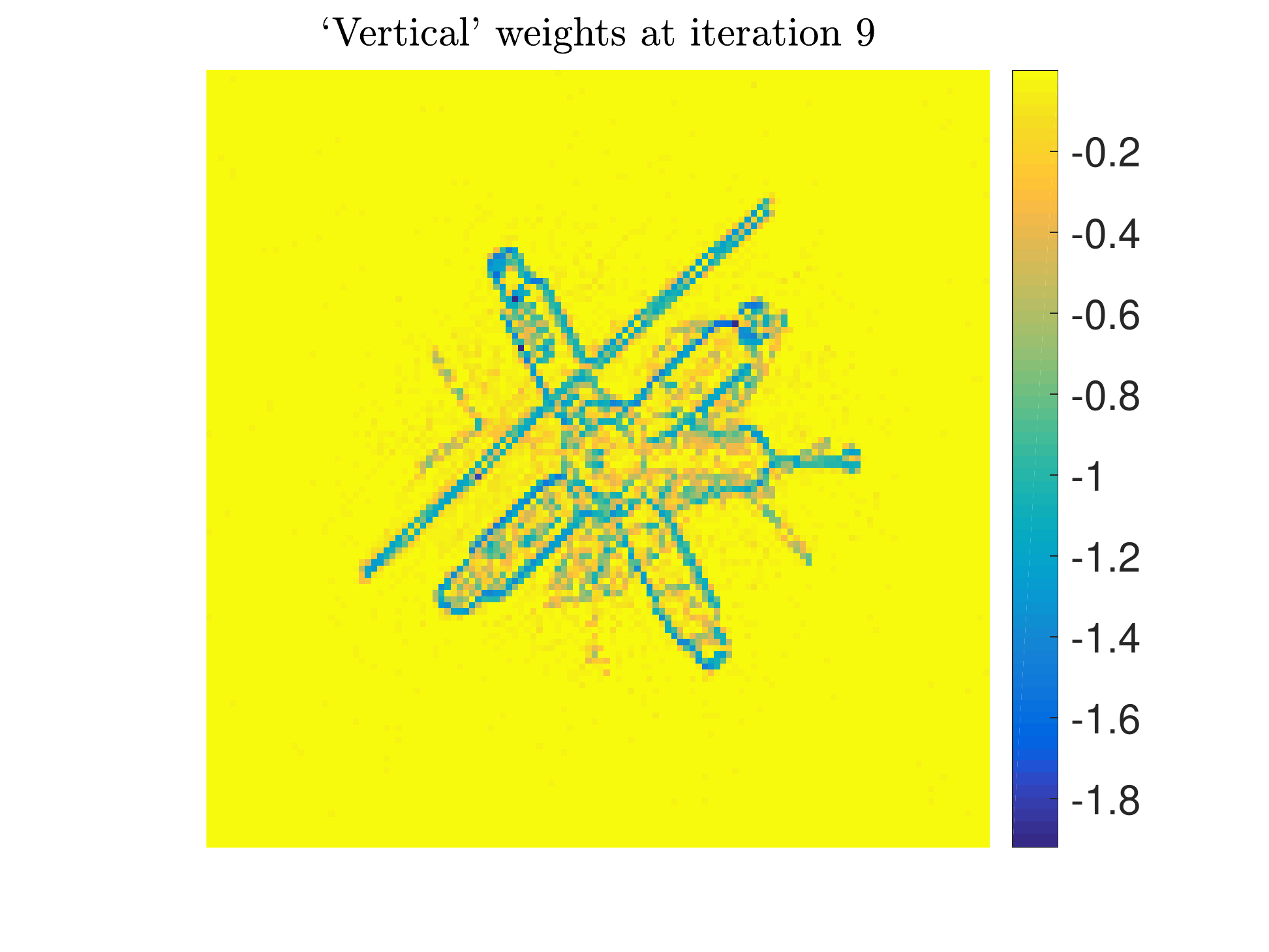} &
\includegraphics[width=5cm]{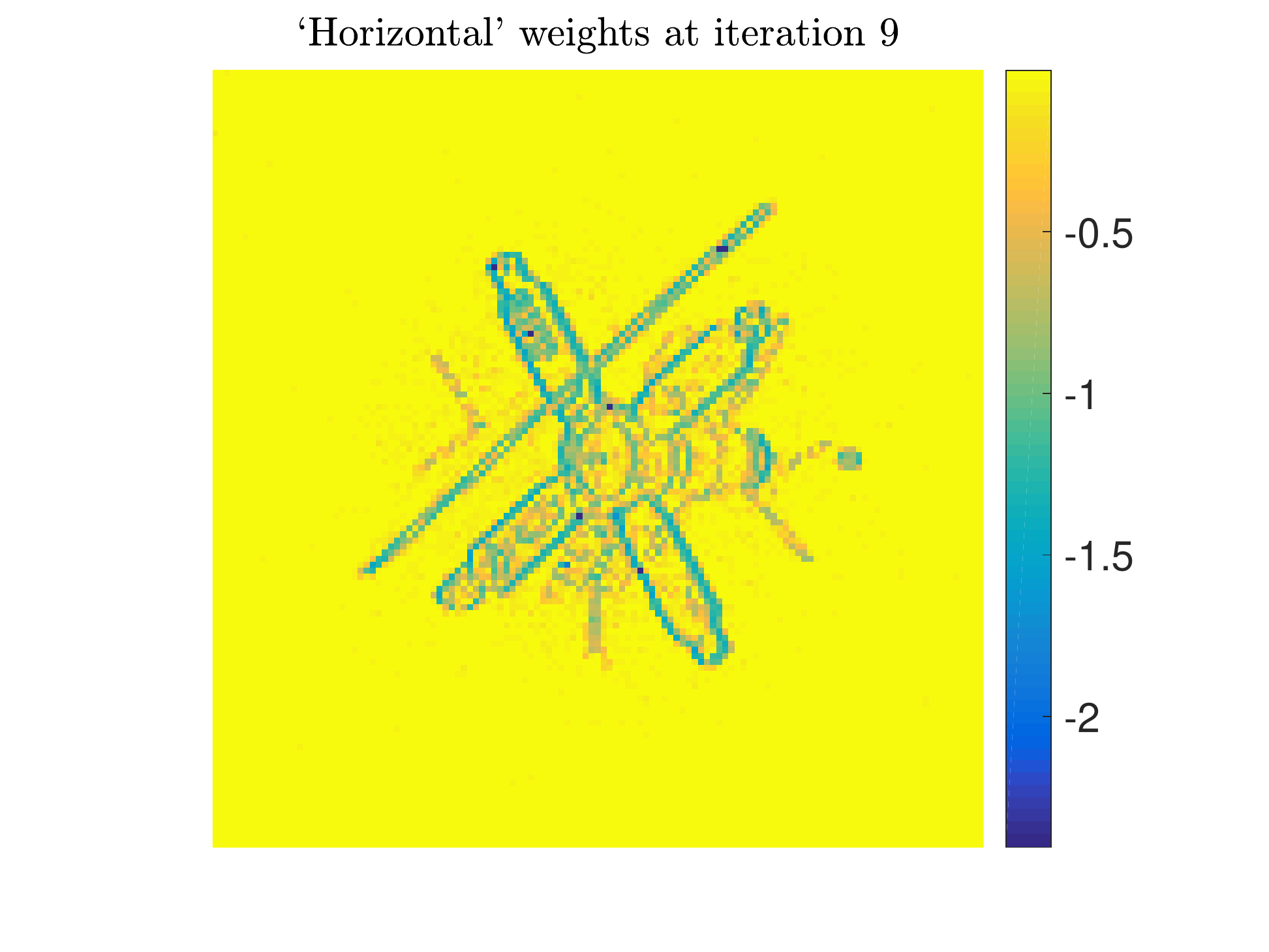} & 
\includegraphics[width=5cm]{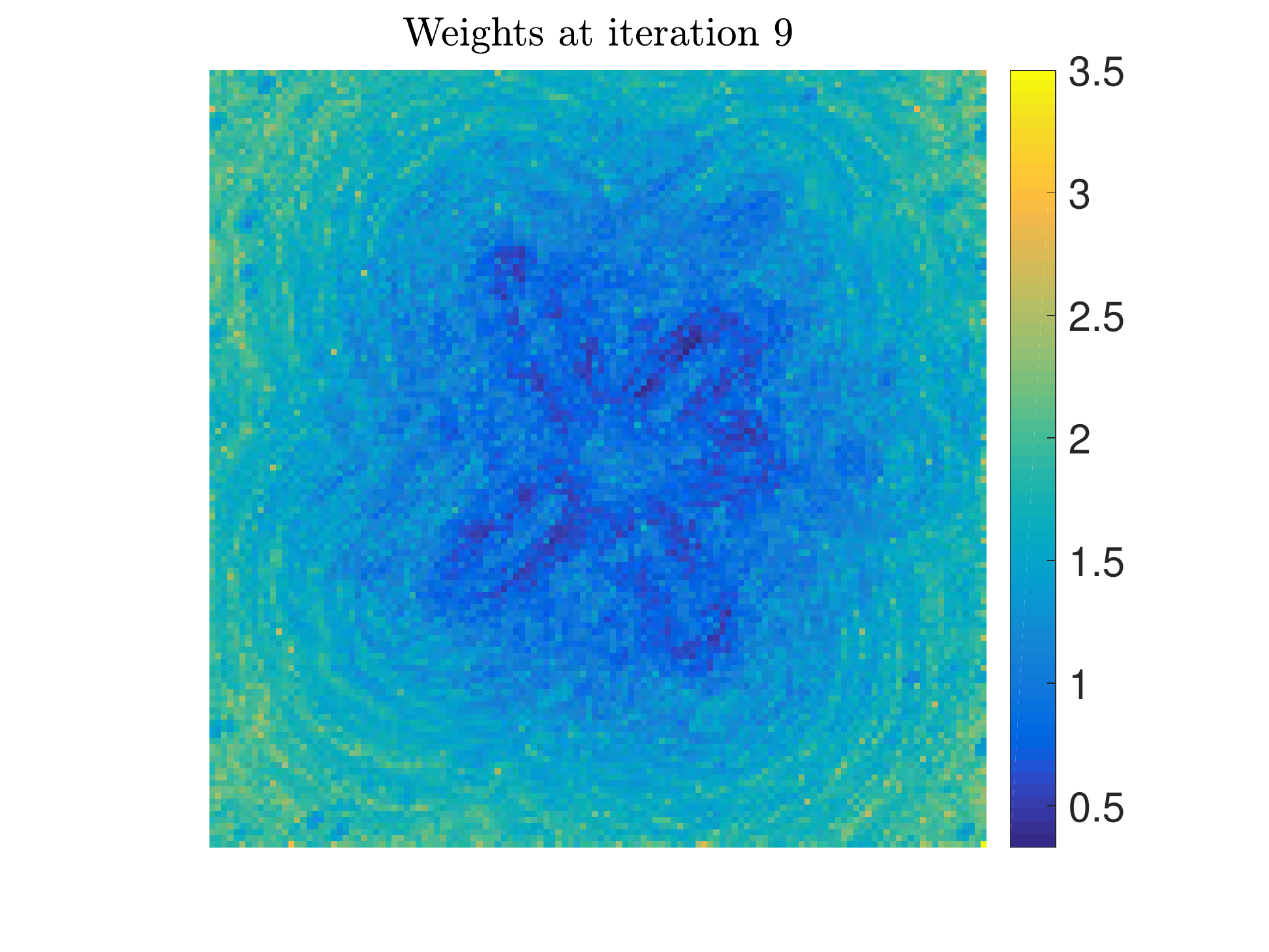}\\ \end{tabular}
\end{center}
\caption{\emph{Out-of-focus blur} test problem. Left column: new weights, to be applied to the vertical derivatives. Middle column: new weights, to be applied to the horizontal derivatives. Right column: IRN-TV weights. The pixel values are displayed in logarithmic scale.}
\label{fig:Blur2TVweight}
\end{figure} 
Figure \ref{fig:Blur2TVweight} displays the entries of the weight diagonal matrices at outer iterations $\ell=2$ and $\ell=9$, for both the new reweighting strategy (\ref{eq:GradientMap}) and the IRN-TV reweighting strategy. When considering the new reweighting strategy, we can clearly see that both the weights to be applied to the vertical and horizontal derivatives are appropriate and improve with  increasing outer iterations. 
%; the same is not true for the weights to be applied to the horizontal  derivatives, which seem to detect some spurious edges in the first reconstruction $x^{(\ast,1)}$, which are enhanced during the subsequent outer iterations. Despite this, spurious edges are hardly visible in the final reconstruction. 
More evidently than in the previous test problems, the IRN-TV weights are not so effective in revealing the structure of the image: they are oscillating but usually  very small, and the original structure of the image is totally lost in the last weights. 

\section{Conclusions and Future Work} \label{sec:conclusions}

In this paper we introduced a new inner-outer iterative algorithm for restoring and reconstructing images with enhanced edges, where 
% are enhanced. The highlights of the new method can be summarized as follows:
%\begin{enumerate}
%\item[(a)] A sequence of quadratic Tikhonov-regularized problems is solved (outer iterations): the regularization matrix for each problem is obtained by premultiplying the gradient by a weight matrix that encodes edge information disclosed within all the previous iterations, in such a way that edges are not penalized in the solution process.
%\item[(b)] A specific hybrid method is applied to solve each quadratic problem (inner iterations): this is a very efficient solver, where well-known strategies to set the regularization parameter can be succesfully incorporated. 
%\end{enumerate}
a sequence of quadratic Tikhonov-regularized problems is solved (outer iterations) and a specific hybrid method is applied to solve each quadratic problem in the sequence (inner iterations). The regularization matrix for each quadratic problem in the sequence is updated at each outer iteration and is obtained by premultiplying the gradient operator by a weight matrix that encodes edge information disclosed within all the previous outer iterations, in such a way that edges are not penalized in the solution process. The hybrid inner solver is very effective, and well-known strategies to set the regularization parameter can be successfully incorporated therein. The resulting strategy is innovative and very efficient when compared to available edge-recovery techniques, mainly because of the definition of the new weights and the incorporation of automatic regularization parameter choice techniques.

Future work will include the investigation of similar reweighting techniques, where the inner iterations are performed by a different hybrid method (for instance the one described in Section \ref{ssec:hybridstd}, with possible strategies to approximate the action of $M_A^\dagger$). Also, different expressions of the weighting matrices may be considered, in such a way that edge information is extracted by operators other than the gradient. Finally, similarly to \cite{chunggazzola2018,gazzolasabate2019,GaNa14}, possible hybrid variants that involve the use of flexible Krylov subspaces will be investigated, with the ideal goal of avoiding nested cycles of iterations. 

%\sg{Can this be reformulated within flexible Krylov subspaces?}

%\bibliographystyle{siam}

%\bibliography{edgebib}

\end{document}